\documentclass[12pt,a4paper,reqno]{amsart}
\pdfoutput=1
\usepackage{amsmath,amsfonts,amssymb,bm}   
\usepackage{amsthm}
\usepackage{mathrsfs}
\usepackage{graphicx,graphics}
\usepackage[latin2]{inputenc}
\usepackage{t1enc}
\usepackage{indentfirst}
\usepackage{pict2e}
\usepackage{epic}
\numberwithin{equation}{section}
\usepackage[margin=2.9cm]{geometry}
\usepackage{epstopdf} 
\usepackage[colorlinks,linkcolor=blue]{hyperref}
\usepackage[capitalise]{cleveref}

\crefformat{equation}{(#2#1#3)}
\crefrangeformat{equation}{(#3#1#4) to~(#5#2#6)}

\usepackage[normalem]{ulem}

\allowdisplaybreaks

\theoremstyle{plain}
\newtheorem{Thm}{Theorem}[section]
\newtheorem*{Thm*}{Theorem}
\newtheorem{Lem}[Thm]{Lemma}

\newtheorem{Prop}[Thm]{Proposition}
\newtheorem{Def}[Thm]{Definition}

\theoremstyle{definition}

\newtheorem{Rem}[Thm]{Remark}
\newtheorem{?}[Thm]{Problem}

\newcommand{\dv}{\mathrm{div}}
\newcommand{\sgnk}{\mathrm{sgn}(k)}
\newcommand{\sgnkone}{\mathrm{sgn}(k_1)}
\newcommand{\sgnktwo}{\mathrm{sgn}(k_2)}
\newcommand{\dt}[0]{\partial_t}
\newcommand{\vu}{\mathbf{u}}
\newcommand{\vB}{\mathbf{B}}
\newcommand{\vH}{\mathbf{H}}
\newcommand{\vE}{\mathbf{E}}

\newcommand{\R}{\mathbb{R}}
\newcommand{\e}{\varepsilon}
\newcommand{\sU}{\mathscr{U}}
\newcommand{\sB}{\mathscr{B}}
\newcommand{\sQ}{\mathscr{Q}}
\newcommand{\sH}{\mathscr{H}}
\newcommand{\sE}{\mathscr{E}}
\newcommand{\sL}{\mathscr{L}}

\newcommand{\sR}{\mathscr{R}}
\newcommand{\sA}{\mathscr{A}}
\newcommand{\sF}{\mathscr{F}}
\newcommand{\sG}{\mathscr{G}}
\newcommand{\bF}{\mathbb{F}}
\newcommand{\bG}{\mathbb{G}}
\newcommand{\fL}{\mathfrak{L}}

\newcommand{\norm}[2]{\left\lVert #1 \right\rVert_{#2}}
\newcommand{\abs}[2]{\left\lvert #1 \right\rvert^{#2}}

\begin{document}

%%%%%%%%%%%%%%%%%%%%%%%%%%%
% title,author,address, email
%%%%%%%%%%%%%%%%%%%%%%%%%%%

\title[Surface waves on plasma--vacuum interfaces]{Weakly nonlinear surface waves on the plasma--vacuum interface}

\author{Paolo Secchi}
\author{YUAN YUAN$ ^* $}
\address[P. Secchi]{DICATAM, Sezione di Matematica,
Universit\`a di Brescia, Via Valotti 9, 25133 Brescia, Italy}
\email{paolo.secchi@unibs.it}
\address[Y. Yuan]{South China Research Center for Applied Mathematics and Interdisciplinary Studies, South China Normal University, Guangzhou, Guangdong, China}
\email{yyuan2102@m.scnu.edu.cn}

\thanks{$ ^* $Corresponding author.}

%\subjclass[2010]{Primary: 05C??. Secondary: 05C??}

%\keywords{draft} 

\begin{abstract}
We consider the free boundary problem for a plasma--vacuum interface in ideal incompressible magnetohydrodynamics. Unlike the classical statement, where the vacuum magnetic field obeys the div-curl system of pre-Maxwell dynamics, we do not neglect the displacement current in the vacuum region and consider the Maxwell equations for electric and magnetic fields. 
Our aim is to construct highly oscillating surface wave solutions in weakly nonlinear regime to this plasma--vacuum interface problem.

Under a necessary and sufficient stability condition for a piecewise constant background state, we construct approximate solutions at any arbitrarily large order of accuracy to the free boundary problem in three space dimensions when the initial discontinuity displays high frequency oscillations. Moreover, such approximate surface waves have nontrivial residual non-oscillatory components.

%(As evidenced in earlier works, high frequency oscillations of the plasma--vacuum interface solution give rise to surface waves on either side of the interface. Such waves decay exponentially in the normal direction to the interface and, in the weakly nonlinear regime that we consider here, their leading amplitude is governed by a nonlocal Hamilton-Jacobi type equation, as for Rayleigh waves in elastodynamics and current-vortex sheets in MHD.)

\end{abstract}

\maketitle

%%%%%%%%%%%%%%%%%%%%%%%%%%%%%%%%%%%%%%%%%%%%%%%%%%%%%%%%%%%%%%%%
%%%%%Introduction #1
%%%%%%%%%%%%%%%%%%%%%%%%%%%%%%%%%%%%%%%%%%%%%%%%%%%%%%%%%%%%%%%%
\section{Introduction }\label{sec-intro}

In this paper, we are interested in the asymptotic analysis of the free boundary problem for a plasma--vacuum interface in ideal incompressible magnetohydrodynamics. 

plasma--vacuum interface problems are considered in the mathematical modeling of plasma confinement  in thermonuclear energy production (as in Tokamaks, Stellarators; see, e.g., \cite{BFKK}, \cite{Goed}). The plasma--vacuum interface appears as a typical phenomenon when the plasma is separated by a vacuum from outside wall, due to the effect of strong magnetic fields. There are also important applications in astrophysics, where the plasma--vacuum interface problem can be used for modeling the motion of a star or the solar corona when magnetic fields are taken into account. 

%

%

%\subsection{The plasma--vacuum interface problem}

%In this paper, we consider the weakly nonlinear surfaces waves on the free interface separating the incompressible plasma and vacuum in three-dimensional case.

 Assume that the plasma--vacuum interface is described by $\Gamma (t)=\{F(t,x)=0\}$, and that
$\Omega^\pm(t)=\{F(t,x)\gtrless 0\}$ are the space-time domains occupied by the plasma and the vacuum, respectively. The function $F$ is an unknown of the problem.

The motion of incompressible plasma in the region $\Omega^+(t)$ is governed by the incompressible magneto-hydrodynamics (MHD) equations
\begin{align}\label{MHD}
\left\{ 
\begin{aligned}
&\dt \vu + (\vu\cdot \nabla) \vu -(\vB \cdot \nabla ) \vB +\nabla q=0  &&\text{in~} \Omega^+(t),\\
&\dt \vB + (\vu\cdot \nabla) \vB -(\vB \cdot \nabla ) \vu =0&&\text{in~} \Omega^+(t),\\
&\dv \vu=0,\quad \dv \vB=0&&\text{in~} \Omega^+(t),
\end{aligned}
\right.
\end{align}
where $\vu=(u_1,u_2,u_3)$, $\vB=(B_1,B_2,B_3)$, $p$, $q =p+\frac{1}{2}|{\vB} |^2$ denote respectively the velocity, magnetic field, the pressure and the total pressure.

 In the classical description of \cite{BFKK} (see also \cite{Goed}) the plasma is described by the  equations of ideal compressible magnetohydrodynamics (MHD),
% \footnote{Here we do not write out explicitly the compressible MHD equations because in the sequel we are going to consider the incompressible MHD equations in the plasma region.}
whereas in the vacuum region $\Omega^-(t)$ one considers the so-called {\it pre-Maxwell dynamics}
\begin{equation}
\nabla \times \vH =0,\qquad {\rm div}\, {\vH}=0,\label{6}
\end{equation}
\begin{equation}
\nabla \times {\vE} =- \frac1{c}\partial_t{\vH},\qquad {\rm div}\, \vE=0,\label{6'}
\end{equation}
where $\vH=(H_1,H_2,H_3)$ , $\vE=(E_1,E_2,E_3)$ denote the magnetic and electric fields in vacuum and $c$ is the speed of light.
The equations \eqref{6}, \eqref{6'} are obtained from the Maxwell equations by neglecting the displacement current $(1/c)\,\partial_t\vE$.
From \eqref{6'} the electric field $E$ is a secondary variable that may be  computed from the magnetic field ${ \vH}$; thus it is enough to consider \eqref{6} for the magnetic field $\vH$.

The problem is completed by the boundary conditions at the free interface $\Gamma (t)$
\begin{equation}\label{7}
\frac{{\rm d}F }{{\rm d} t}=0,\quad
 [q]=0,\quad  \vB\cdot N=0,\quad
  \vH\cdot N=0 ,
 \end{equation}
where $\vB=(B_1,B_2,B_3)$ is the magnetic field in the plasma region, $[q]=q_{|\Gamma}-\frac12|\vH_{|\Gamma}|^2$ denotes the jump of the total pressure across the interface and $N=\nabla F$. The
first condition in \eqref{7} (where ${\rm d}/{\rm d}t= \partial_t +v\cdot\nabla$ denotes the material derivative) means that the interface moves with the velocity of plasma particles at the boundary.

Trakhinin \cite{Trakhinin2010} obtained a basic energy a priori estimate in Sobolev spaces for the linearized plasma--vacuum interface problem of \eqref{MHD}, \eqref{6}, \eqref{7} under the non-collinearity condition
\begin{equation}\label{non-coll}
\vB\times \vH\not=0 \qquad\mbox{on}\ \Gamma(t),
\end{equation}
satisfied on the interface for the unperturbed flow. Under the same non-collinearity condition \eqref{non-coll} satisfied at the initial time, the well-posedness of the nonlinear problem \eqref{MHD}, \eqref{6}, \eqref{7} was proved by Secchi and Trakhinin \cite{Secchi2013, Secchi2014} for compressible MHD equations in plasma region, and by Sun, Wang and Zhang \cite{Sun2019}  for incompressible MHD equations. It is interesting to observe that in \cite{Sun2019} the authors are able to avoid the loss of derivatives phenomenon, as in their previous work \cite{Sun2018} about the well-posedness of incompressible current-vortex sheets. See also \cite{Gu2019a, Hao2017, Morando2014}.

For the different model where the plasma is still described by the  equations \eqref{MHD},
but in the vacuum region $\Omega^-(t)$ it is assumed that $\vH\equiv0$,
Hao and Luo \cite{Hao2014} and Gu and Wang \cite{Gu2019b} studied the well-posedness for the nonlinear problem under the Rayleigh-Taylor sign condition $\partial q/\partial N_{|\Gamma}<0$. 
For the case in two space dimensions, Hao and Luo \cite{Hao2020} have recently proved the ill-posedness of the same problem when the Rayleigh-Taylor sign condition is violated at the initial time. 
%%%%

The linearized stability of the relativistic case was first addressed by Trakhinin in \cite{Trakhinin2012}, in the case of plasma expansion in vacuum.  Then the linearized stability for  the related non-relativistic problem was studied in the papers \cite{Catania2014, Catania2016, Mandrik2014} (see also \cite{Mandrik2018, Trakhinin2016}) by considering a model where in the vacuum region, instead of the pre-Maxwell dynamics, the displacement current is taken into account and the complete system of Maxwell equations for the electric and the magnetic fields is considered. The introduction of this model aims at investigating the influence of the electric field in vacuum on the well-posedness of the problem, since in the classical pre-Maxwell dynamics such an influence is hidden.

%%%%%%
%The plasma in the domain $\Omega^+(t)$ is assumed to be ideal and incompressible whereas in the vacuum region $\Omega^-(t)$
In this model the magnetic and electric fields $\vH$ and $\vE$
are governed by the Maxwell equations written in dimensionless form
\begin{align}\label{Maxwell}
\left\{ 
\begin{aligned}
&\nu\dt \vH + \nabla \times \vE=0&&\text{in~} \Omega^-(t),\\
&\nu\dt \vE - \nabla \times \vH=0&&\text{in~} \Omega^-(t), \\
&\dv \vH=0,\quad \dv \vE=0&&\text{in~} \Omega^-(t).
\end{aligned}
\right.
\end{align}
where 
$\nu\ll 1$ is a small parameter.
%We assume that the free interface $\Gamma (t)$ has the form of a graph and the domains $\Omega^{\pm}(t)$ occupied by the plasma and the vacuum are unbounded:
%\[
%\Gamma (t)=\{F(t,x)=x_1-\varphi (t,x')\},\quad \Omega^\pm(t)=\{\pm (x_1- \varphi(t,x'))>0,\ x'\in \mathbb{R}^2\},\quad x'=(x_2,x_3).
%\]
%For simplicity, we consider the plasma and the vacuum are occupied in $\{(x_1,x_2,x_3); x'=(x_1,x_2)\in \mathbb{T}^2,x_3\in \R\}$, 
%and the plasma and the vacuum are separated by a free interface $\Gamma(t):=\{(x_1,x_2,x_3)$$; x_3=\varphi(x_1,x_2) \}$ into two regions $\Omega^{\pm}(t):=\{(x_1,x_2,x_3);$$ x_3\gtrless\varphi(x_1,x_2) \}$, respectively. 
%On the interface $\Gamma(t)$ the variables $\vu,\vB,\vH,\vE$ satisfy 
%\begin{align}\label{interface}
%\begin{aligned}
%&\dt \varphi= \vu \cdot N,\quad [q]=0, \quad \vB\cdot N=\vH\cdot N=0,\\
%& \text{and}~ N\times \vE=\nu (\vu \cdot N)\vH \quad \text{on}~\Gamma(t),
%\end{aligned}
%\end{align} 
%where $N$ is the normal direction of free interface $N:=(-\partial_{x_1} \varphi, -\partial_{x_2} \varphi,1 )$ and $[q]=q|_{\Gamma}-\frac{1}{2}\abs{\vH}{2} |_{\Gamma}+\frac{1}{2}\abs{\vE}{2} |_{\Gamma}$.
On the interface $\Gamma(t)$ the variables $\vu,\vB,\vH,\vE$ satisfy 
\begin{align}\label{interface0}
\begin{aligned}
&\frac{{\rm d}F }{{\rm d} t}=0,\quad [q]=0, \quad \vB\cdot N=\vH\cdot N=0,\\
& \text{and}~ N\times \vE=\nu (\vu \cdot N)\vH \quad \text{on}~\Gamma(t),
\end{aligned}
\end{align} 
where now $[q]=q|_{\Gamma}-\frac{1}{2}\abs{\vH}{2} |_{\Gamma}+\frac{1}{2}\abs{\vE}{2} |_{\Gamma}$.

%%%%%%

%Note that for the last equation in \cref{interface}: 
%\begin{align*}
% - E_3\partial_{x_2} \varphi-E_2 &=\nu  H_1\dt \varphi,\\
% E_3\partial_{x_1} \varphi+E_1 &=\nu  H_2\dt \varphi,\\
% - E_2\partial_{x_1} \varphi+ E_1\partial_{x_2} \varphi&=\nu H_3\dt \varphi ,
%\end{align*}
%the third one can be derived by the first two equations and the slip condition of $\vH$ ($\vH\cdot N=0$).

The analysis of \cite{Trakhinin2012} suggests that ill-posedness occurs for problem \eqref{MHD}, \eqref{Maxwell}, \eqref{interface0} in the presence of a {\it sufficiently strong} vacuum electric field. On the other hand, it is shown in \cite{Catania2014, Mandrik2014} that a {\it sufficiently weak} vacuum electric field precludes ill-posedness and gives the well-posedness of the linearized problem, thus somehow justifying the practice of neglecting the displacement current in the classical pre-Maxwell formulation when the vacuum electric field is weak enough. Such smallness hypothesis is not required for linear well-posedness in the two dimensional case, see \cite{Catania2016}.

In \cite{Trakhinin2020b} Trakhinin investigates the linearized problem of \eqref{MHD}, \eqref{Maxwell}, \eqref{interface0} about a piecewise constant reference state and obtains a necessary and sufficient condition for the violent instability of a planar plasma--vacuum interface (the opposite of this condition is given in \eqref{H1}).
In particular, it is shown that as the unperturbed plasma and vacuum magnetic fields are collinear (i.e. when \eqref{non-coll} is violated), any nonzero unperturbed vacuum electric field makes the planar interface violently unstable. This shows the necessity of the corresponding non-collinearity condition \eqref{non-coll} for well-posedness and a crucial role of the vacuum electric field in the evolution of a plasma--vacuum interface.
In the recent paper \cite{Morando2020a}, the authors study the linearized constant coefficient problem (the same problem as in \cite{Trakhinin2020b}) and show that under the stability condition \eqref{H1} it admits an energy a priori estimate, showing the stability of the planar plasma--vacuum interface.

%%%%%%%%

In the present paper we are interested in the {\it qualitative behaviour} of exact solutions to \eqref{MHD}, \eqref{Maxwell}, \eqref{interface0} for {\it highly oscillating} initial data.
%%%
This problem has been first addressed by Al\`i and Hunter for the current-vortex sheets of two-dimensional incompressible flows in ref. \cite{Ali2003}, where they have shown that for some specific oscillation phase, and in the {\it weakly nonlinear regime}, the leading amplitude of the solution on either side of the current vortex sheet displays a {\it surface wave} structure: it oscillates with the same phase as the front and it is localized near the free surface with exponential decay in the normal direction to the free surface. The leading amplitude of the front is governed by a nonlocal Hamilton-Jacobi type equation, as for Rayleigh waves in elastodynamics, see \cite{Hamilton1995}. About amplitude equations and their Hamiltonian structure see, e.g., \cite{ali2002,benzoni2009a, benzoni2011, benzoni2009b, benzoni2012, hunter1989, hunter2006, hunter2011}. 
%%%
%(As evidenced in earlier works, high frequency oscillations of the plasma--vacuum interface solution give rise to surface waves on either side of the interface. Such waves decay exponentially in the normal direction to the interface and, in the weakly nonlinear regime that we consider here, their leading amplitude is governed by a nonlocal Hamilton-Jacobi type equation, as for Rayleigh waves in elastodynamics and current-vortex sheets in MHD.)
%%%
The occurrence of surface waves on the plasma--vacuum interface is studied in \cite{Secchi2015, secchi2016}. 

Such surface waves, exponentially decaying in the normal direction to the free boundary, are sometimes called {\it genuine} surface waves.
In other problems related to hydrodynamics, such as denotation waves in combustion \cite{majda1983} or supersonic compressible vortex sheets \cite{Artola1987, Coulombel2004a, Coulombel2008}, there occur non-decaying {\it radiative} or {\it leaky} surface waves that generate {\it bulk waves} propagating away from the boundary into the interior of the space domain; for a general description see \cite{benzoni2002, hunter2011}.
%%%

The aim of the present paper is to construct highly oscillating surface wave solutions in weakly nonlinear regime to the plasma--vacuum interface problem \eqref{MHD}, \eqref{Maxwell}, \eqref{interface0}, that are small perturbations of a reference solution. We follow the method of Pierre and Coulombel in \cite{Pierre2021}, that construct highly oscillating current-vortex sheets solutions to the MHD equations, by introducing a suitable WKB ansatz. As a matter of fact, the plasma--vacuum interface problem and the current-vortex sheet problem presents many similarities; in particular the evolution equation that governs the leading amplitude of the front is the same equation exhibited for current-vortex sheets and Rayleigh waves in elastodynamics, as shown in \cite{Secchi2015}.

Under the necessary and sufficient stability condition \eqref{H1} for a piecewise constant background state, we construct approximate solutions at any arbitrarily large order of accuracy to the free boundary problem in three space dimensions, when the initial discontinuity displays high frequency oscillations. 
Moreover, such approximate surface waves have nontrivial residual non-oscillatory components, which is referred to as \textit{rectification phenomenon} (\cite{Marcou2011}); see Section 1.1.5 for the definition.

Differently from \cite{Pierre2021}, in some parts the resolution of the plasma-vacuum problem requires new ideas, especially for the construction of slow means of correctors which leads to the investigation of the so-called {\it slow problem}, see \cref{sec-correctors}. While for the current-vortex sheets problem \cite{Pierre2021}, at this point the main tool is the coupled elliptic problem for the total pressures, here we are lead to investigate an overdetermined coupled elliptic-hyperbolic system for the plasma total pressure on one side and the extended total pressure in the vacuum region. For the resolution of this problem we also apply some results of \cite{Morando2020a,Trakhinin2020b}.

The plan of the paper is as follows: the present \cref{sec-intro} contains the formulation of the problem, the functional setting and the main result. In Section 2 we derive the WKB cascade, i.e. obtain the equations and boundary conditions for each order of the correctors. In Section 3 we solve the main issue of each step of the inductive processes, the so-called \emph{fast problem}. The leading amplitude and the correctors are constructed by induction in Section 4 and Section 5. In the final section, we discuss the rectification phenomenon, showing the necessity to introduce the residual components in the functional setting.

\subsection{Settings and initial data of the front}

We assume that the free interface $\Gamma(t)$ has the form of a graph 
\[
\Gamma(t):=\{(x_1,x_2,x_3);\,x'=(x_1,x_2)\in \mathbb{T}^2,\,  x_3=\varphi(x_1,x_2) \},
\]
where $\varphi$ denotes the unknown front function,
and assume that the plasma and the vacuum are contained in the unbounded region 
\[
\Omega:=\{(x_1,x_2,x_3);\,; x'\in \mathbb{T}^2,x_3\in \R \}.
\]
Denoting
\[
\Omega^{\pm}(t):=\{(x_1,x_2,x_3);\,; x'\in \mathbb{T}^2,\, x_3\gtrless\varphi(x_1,x_2) \}
\]
we assume that the plasma is contained in $\Omega^+(t)$ and is governed by equations \eqref{MHD}, while the vacuum region is $\Omega^{-}(t)$ where we consider the Maxwell equations \eqref{Maxwell}. At the free interface $\Gamma(t)$ we assume the boundary conditions \eqref{interface0}, which take the form
\begin{align}\label{interface}
\begin{aligned}
&\dt \varphi= \vu \cdot N,\quad [q]=0, \quad \vB\cdot N=\vH\cdot N=0,\\
& \text{and}~ N\times \vE=\nu (\vu \cdot N)\vH \quad \text{on}~\Gamma(t),
\end{aligned}
\end{align} 
where $N:=(-\partial_{x_1} \varphi, -\partial_{x_2} \varphi,1 )$ and $[q]=q|_{\Gamma}-\frac{1}{2}\abs{\vH}{2} |_{\Gamma}+\frac{1}{2}\abs{\vE}{2} |_{\Gamma}$.

Note that the last equation in \cref{interface} can be detailed as
\begin{align}\label{interface2}
 - E_3\partial_{x_2} \varphi-E_2 &=\nu  H_1\dt \varphi,\notag\\
 E_3\partial_{x_1} \varphi+E_1 &=\nu  H_2\dt \varphi,\\
 - E_2\partial_{x_1} \varphi+ E_1\partial_{x_2} \varphi&=\nu H_3\dt \varphi .\notag
\end{align}
In particular, the third equation can be derived by the first two ones and the slip condition $\vH\cdot N=0$.

\subsubsection{The reference plasma--vacuum interface}

The reference plasma--vacuum interface $\Gamma_0$, the plasma region $\Omega^+_0$ and the vacuum region $\Omega^-_0$ are respectively
\[\Gamma_0:=\{(x_1,x_2,x_3); x_3=0 \},\quad\Omega^{\pm}_0:=\{(x_1,x_2,x_3); x_3\gtrless0 \}.
\] 
We consider a piecewise constant solution to \cref{MHD}, 
\cref{Maxwell} and \cref{interface} as the reference plasma--vacuum interface:
\begin{align}
&\vu^0:=(u^0_1,u^0_2,0)^{\top},\quad  \vB^0:= (B^0_1,B^0_2,0)^{\top}, \notag\\
&\vH^0:=(H^0_1,H^0_2,0)^{\top},\quad  \vE^0:= (0,0,E^0_3)^{\top},\notag\\
&q^0:=\frac{1}{2}(H^0_1)^{2} +\frac{1}{2}(H^0_2)^{2}-\frac{1}{2}(E^0_3)^{2}, \label{reference-state} \\
&\varphi^0, \varphi^1:=0,\notag
\end{align}
where we use two functions $\varphi^0$ and $\varphi^1$ for the fixed reference front that could be rather thought of as $\varphi^0+\epsilon\varphi^1$, to be consistent with the notation below for the WKB ansatz.
Note that here $E^0_3$ may not be zero, which is different from the two dimensional case (\cite{Secchi2015}).

In the rest of this paper, we will use the notations \[U:=(\vu,\vB, q)^{\top},\quad V:=(\vH,\vE)^{\top},\]
\[U^0:=(\vu^0,\vB^0, q^0)^{\top},\quad  V^0:=(\vH^0,\vE^0)^{\top} .
\]
Moreover, we also assume for the reference plasma--vacuum interface \eqref{reference-state} the following stability criterion:
\begin{align}\label{H1} \tag{H1}
	\abs{\vE^0}{2}< \frac{\abs{\vB^0}{2}+\abs{\vH^0}{2}-\sqrt{(\abs{\vB^0}{2}+\abs{\vH^0}{2})^2-4\abs{\vB^0\times \vH^0}{2} } }{2}.
\end{align}
\cref{H1} was first proposed in \cite{Trakhinin2020b} by normal mode analysis, and proved to be the linear stability condition in \cite{Morando2020a}.
%The nonlinear surface wave solutions satisfy the decay conditions:
%\begin{align*}
%\lim\limits_{x_3\rightarrow+\infty}U =U^0, \quad
%\lim\limits_{x_3\rightarrow-\infty}V =V^0.
%\end{align*}

\subsubsection{The frequencies} \label{sec-frequencies}
In this paper, we use Einstein summation convention, and use $\alpha \in \{1,2,3 \}$ to represent the spatial coordinate and  $i,j \in \{1,2 \}$ to represent the tangential spatial coordinate. 
For tangential spatial variable $x'=(x_1,x_2)\in \mathbb{T}^2$, the corresponding frequency vector is defined by $\xi'=(\xi_1, \xi_2)\in \R^2$ with $\abs{\xi'}{} =1$.

It follows from \cite{Trakhinin2020b} that the stability condition \cref{H1} is equivalent to 
	\begin{align}\label{H1*} \tag{H1*}
		(E^0_3)^2< \min\limits_{|\xi|=1} \{ (\xi_j  B^0_j )^2+(\xi_j  H^0_j )^2\}.
	\end{align}
	
Given a (real) time frequency $\tau$ let us define the following parameters:
\begin{align}
	\begin{aligned} \label{simp-1}
%			\quad	&~\\
			~&  a^+=\xi_j u^0_j,  \quad b^+=\xi_j B^0_j, \quad c^+=\tau+a^+,\\
			&   a^-_1=\nu \tau H^0_1+\xi_2 E^0_3, \quad a^-_2=\nu \tau H^0_2-\xi_1 E^0_3 ,\quad b^-=\xi_j H^0_j. \quad%\\
%			&~
		\end{aligned}
\end{align}
In this way, $a^-_1, a^-_2, b^-$ have the relation $$\xi_j a^-_j =\nu\tau b^-.$$

Now we give {the assumptions on the frequencies} $\tau,\xi'$ that will be taken in the following main Theorem \ref{thm-main}.

\begin{enumerate}
%	\item 
%	It follows from \cite{Trakhinin2020b} that \cref{H1} is equivalent to 
%	\begin{align}\label{H1*} \tag{H1*}
%		(E^0_3)^2< \min\limits_{|\xi|=1} \{ (\xi_j  B^0_j )^2+(\xi_j  H^0_j )^2\}.
%	\end{align}
	
	\item Given a $2\pi$--periodic function $v=v(x',\theta)$ with respect to each of its arguments $(x',\theta)\in \mathbb{T}^3$, we shall require that functions of the form
\begin{equation*}
	\label{periodic}
	v_{\e}(x'):= v(x',\frac{\xi ' \cdot x'}{\e})
\end{equation*}
be $2\pi$--periodic with respect to $x'$. To do so we will choose as in \cite{Pierre2021}
\begin{align} \label{H2}\tag{H2}
		\xi'=\frac{1}{\sqrt{p^2+q^2} } (p,q)^{\top} ~\text{with~} (p,q)^\top \in \mathbb{Z}^2 \setminus \{0\}, \quad \text{and}\quad \e_{l}:=\frac{1}{l\sqrt{p^2+q^2} } \, .
	\end{align}
Then it follows from \eqref{H2} that $\frac{\xi'\cdot x'}{\e} \in \mathbb{T}^2 $ (for $\e=\e_l$), and the periodicity of $v$ w.r.t. $(x', \theta)$ implies the periodicity of $v_\e$ w.r.t. $x'$. In the following, the parameter $\e\in(0,1]$ stands for one element of the sequence $(\e_l)_{l\geq1}$ in \eqref{H2}. When writing $\e\to0$, it means that we consider $\e_l$ with $l\to+\infty$.
	
	\item We will assume
	\begin{equation}\label{H-key} \tag{H3}
		\sqrt{1-\nu^2 \tau^2 } ((c^+)^2 -(b^+)^2)+((a^-_1)^2+(a^-_2)^2-(b^-)^2)=0 \,.
	\end{equation}
	\vspace{0.2cm}
	It can easily be seen that \eqref{H-key} is equivalent to
	\begin{multline}\label{H-key2}
		\sqrt{1-\nu^2 \tau^2 } ( (\tau +\xi_j u^0_j)^2 -(\xi_j B^0_j)^2)+(\nu^2\tau^2 ((H^0_1)^2+(H^0_2)^2) \\
		-(\xi_j H^0_j)^2 + 2\nu \tau (\xi_2H^0_1-\xi_1 H^0_2)E^0_3 +(E^0_3)^2)=0\,.
	\end{multline}
	which is exactly the Lopatinskii determinant obtained in \cite{Trakhinin2020b,Morando2020a}.
		It has been proved in \cite{Trakhinin2020b} that under \cref{H1*}, given $\xi'\in \mathbb{S}^1$ and for $\nu$ sufficiently small there are only two different real roots $\tau=\tau(\xi)$ to \eqref{H-key} (or \cref{H-key2}), i.e. only the weak Lopatinskii condition is fulfilled.

	\item We will also assume
	\begin{equation} \label{H-a} \tag{H4}
		a^+\neq 0\,,\qquad \tau \neq 0\,.
	\end{equation}
	In fact, $a^+\neq 0\ $ can be easily achieved by excluding two directions on $\mathbb{S}^1$ for prescribed $u_j^0$   ($\xi'$ satisfying \cref{H2} is dense in $\mathbb{S}^1$),
	and $\tau\neq 0$ is fulfilled by choosing one of the two different roots of \cref{H-key}. 
	
\end{enumerate}

%\noindent\red{The existence of real roots to \cref{H-key} can also follow from the Implicit Function Theorem with \cref{H1*} and $a^+\neq 0$: when $\nu=0$, \cref{H-key} has two real roots $(\tau_0+a^+)^2=(b^+)^2+(b^-)^2-(E^0_3)^2$; if $\nu$ is sufficiently small, since the partial derivative w.r.t. $\tau$ at $(0,\tau_0$) is $2a^+\tau_0\neq 0$, then there exist two distinct real roots $\tau=\tau(\nu)$ to \cref{H-key} for any $\xi$ satisfying $a^+\neq 0$. }
%\todo{changed here-Apr: I give the details of the two proofs of \eqref{H-a}. I'm not sure whether to keep the Implicit Function Theorem part}
%\todo{Paolo, minor changes}

\noindent With the above assumptions, the frequencies satisfy: 
%\begin{itemize}
%\item[(i)]
$$(c^+)^2\neq (b^+)^2, \qquad c^+\neq 0\,.$$ 
Actually, if $(c^+)^2= (b^+)^2$, then the term of order $\mathcal{O}(\nu^2)$ as $\nu\to0$ in the left-hand side of \cref{H-key} is 
$$\nu^2\tau^2 ((H^0_1)^2+(H^0_2)^2)>0\,,$$
which could not vanish for $\nu\to0$; hence $(c^+)^2= (b^+)^2$ contradicts to \cref{H-key}. 
Moreover, if $c^+= 0$, then the term order $\mathcal{O}(1)$ in the left-hand side of \cref{H-key} is 
$$(E^0_3)^2-(\xi_j  B^0_j )^2-(\xi_j  H^0_j )^2<0\,$$
(according to \cref{H1*}); hence $c^+= 0$ also contradicts to \cref{H-key}. 

%\item[(ii)] 
%The frequencies also satisfy: 
%$$c^+\neq 0\,.$$ 
%Actually, if $c^+= 0$, then the order $\mathcal{O}(1)$ of \cref{H-key} is 
%$$(E^0_3)^2-(\xi_j  B^0_j )^2-(\xi_j  H^0_j )^2<0,$$
%which is again a contradiction.

%\end{itemize}
From now on, the reference plasma--vacuum interface \eqref{reference-state} and the frequencies $(\tau,\xi_1,\xi_2)\in\R^3$ satisfying \eqref{H1}, \eqref{H2}, \eqref{H-key},\eqref{H-a}, are {\bf fixed}.

\subsubsection{Initial front}
In the following we construct a solution to \eqref{MHD}, \eqref{Maxwell}, \eqref{interface} in the form of a small, highly oscillating perturbations of the reference constant state $(U_0,V_0, \varphi^0+ \e\varphi^1)$. The initial front is taken of the form
\begin{equation}
	\label{initial-front}
	\varphi_{\e}(t=0,x'):=\e^2 \varphi_0^2(x',\frac{\xi ' \cdot x'}{\e}),
\end{equation}
where $\varphi_0^2\in C^{\infty}(\mathbb{T}^3)$ is assumed to be of zero mean with respect to $\theta\in \mathbb{T}$. Recall that $\xi'$ and $\e=\e_l$ are as in \eqref{H2}, so that $\varphi_{\e}$ in \eqref{initial-front} is indeed $2\pi$-periodic w.r.t. $x'$.

%\begin{Rem}[just for the draft]
%	This method could not deal with the case that the initial front contains $\e^1 \varphi_0^1$.\\
%	While in \cite{Pierre2021}, the authors terms like $\e^m \varphi_0^m$ ($m\geq3)$ would not generate new difficulties but more complicated calculations. So we just assume in this way. 
%\end{Rem}

\subsubsection{WKB ansatz}
We look for an asymptotic expansion of the exact solution $(U_\e, V_\e, \varphi_{\e})$ to \eqref{MHD}, \eqref{Maxwell}, \eqref{interface} as a small, highly oscillating perturbation of the reference constant state. As in \cite{Pierre2021}, we aim at constructing an asymptotic expansion for $(U_\e, V_\e, \varphi_{\e})$ in the form:
\begin{equation}
	\label{WKB-ansatz}
\begin{aligned}
	U_\e(t,x) & \sim U^0+\sum\limits_{m\geq 1} \e^m U^m(t,x',y_3,Y_3,\theta),\\
	V_\e(t,x) & \sim V^0+\sum\limits_{m\geq 1} \e^m V^m(t,x',y_3,Y_3,\theta),\\
	\varphi_\e(t,x') & \sim \sum\limits_{m\geq 2} \e^m\varphi^m(t,x',\theta),
\end{aligned}
\end{equation}
where we have set
\begin{align*}
y_3:=x_3-\varphi_\e (t,x'),\quad Y_3:=\frac{x_3-\varphi_\e (t,x')}{\e},\quad \theta:=\frac{\tau t+\xi'\cdot x'}{\e}.
\end{align*}
By $\sim$ in \eqref{WKB-ansatz} we mean that the series has to be understood in the sense of asymptotic expansions in $\e$, see e.g. \cite{Rauch2012}.

Here the front $\varphi_\e$ is meant to oscillate with planar phase $\tau t+\xi'\cdot x'$ and a slow modulation in the variables $(t,x')$. The interior solution $(U_\e, V_\e)$ will display oscillations with the same planar phase $\tau t+\xi'\cdot x'$ and exponential decay with respect to the {\it fast normal variable} $(x_3-\varphi_\e (t,x'))/{\e}$, which describes the exponential decay of the surface wave. $\theta=(\tau t+\xi'\cdot x')/{\e}\in\mathbb{T}$ is the {\it fast tangential variable} and describes the oscillations. The {\it slow} variables are $(t,y',y_3)$, where $y'=x'\in\mathbb{T}^2$ are the tangential spatial variables and $y_3=x_3-\varphi_\e (t,x')\in\R^\pm$ is the normal variable. The oscillating plasma--vacuum interface $\Gamma_\e(t):=\{x_3=\varphi_\e (t,x')\}$ in the original space variables corresponds to the fixed interface $\{y_3=0\}$ in the straightened variables.

At the initial time we require for $\varphi_\e$:
\begin{equation}\label{initialphi}
\varphi^2(0,x',\theta)=\varphi^2_0(x',\theta)\quad \forall (x',\theta)\in\mathbb{T}^3,\quad\text{and}\quad 
\varphi^m(0,x',\theta)=0 \quad \forall m\geq3.
\end{equation}
Differently from $\varphi_\e$, the initial datum associated with $(U_\e, V_\e)$ {\it is not free}, as is well known in geometric optics because of {\it polarization}, see \cite{Rauch2012}. Actually, part of the profiles of $(U_\e, V_\e)$ will be determined for any time $t\in[0,T]$ by solving algebraic equations, so that part of the initial data for $(U_\e, V_\e)$ will be computed alongside the whole approximate solution. This restricts the choice of the initial data for $(U_\e, V_\e)$; nevertheless the choice of the initial profile $\varphi^2$ for the front is {\it free}, so this justifies \eqref{initialphi}.

The scaling in \eqref{WKB-ansatz} is analogous to the scaling of weakly nonlinear geometric optics, see e.g. \cite{Gues1993,hunter1989, Joly1993,Joly1995,Lescarret2007, Marcou2010, Pierre2021}.
Notice that the expansion \eqref{WKB-ansatz} of the front $\varphi_\e$ starts with the $O(\e^2)$ amplitude; this is because the gradient of $\varphi_\e$ is expected to have the same regularity of the trace of $(U_\e, V_\e)$ on $\Gamma_\e(t)$, as for current-vortex sheets, see \cite{Coulombel2012,Sun2018,Sun2019}.

In the following we shall refer to $(U^1,V^1,\varphi^{2})$ as the leading amplitude in the WKB expansion \eqref{WKB-ansatz}.

\subsubsection{Functional spaces}
The functional spaces used in this paper are defined similarly as in \cite{Pierre2021}, the only difference being that the slow variable $y_3\in \R^{\pm}$. These are the spaces for the profiles of the WKB ansatz \eqref{WKB-ansatz}.
\begin{Def}\label{def-space} We define the following functional spaces.
%	\begin{align*}
%		 \underline{S}^{\pm}:=& H^{\infty}([0,T]\times\mathbb{T}^2\times\R^{\pm}\times\mathbb{T}),\\
%		 {S}^{\pm}_{\star}:=& \{ u \in H^{\infty}([0,T]\times\mathbb{T}^2\times\R^{\pm}\times\R^{\pm}\times\mathbb{T}); ~\exists \delta >0, \forall \alpha \in \mathbb{N}^6,\exists C_\alpha >0, ~\text{s.t.}  \\
%		 &\quad\quad\quad \forall Y_3 > 0,
%		 \norm{\partial^{\alpha}u(\cdot,\cdot, Y_3,\cdot) }{L^{\infty}_{t,y,\theta}} \leq C_\alpha e^{- \delta |Y_3|}\},\\
%		 {S}^{\pm}:=&\underline{S}^{\pm} \oplus {S}^{\pm}_{\star}.
%	\end{align*}
%Functions in $\underline{S}^{\pm}$ depend on the slow variables $(t,y)$ and on the fast tangential variable $\theta$.
%
%Functions in ${S}^{\pm}_{\star}$ depend on both the slow variables  $(t,y)$ and the fast variables $(Y_3,\theta)$, and they decay exponentially as $Y_3\to\pm\infty$ as well as their derivatives, uniformly with respect to all other arguments.
\begin{align*}
	 \underline{S}^{+}:=&\{ U=(\vu,\vB,q)\, : \vu,\vB,{\nabla q} \in H^{\infty}([0,T]\times\mathbb{T}^2_{x'}\times\R^{+}_{y_3}\times\mathbb{T}_\theta) \},\\
		 {S}^{+}_{\star}:=&\{ U=(\vu,\vB,q)\, : \vu,\vB, {\nabla q} \in H^{\infty}([0,T]\times\mathbb{T}^2_{x'}\times\R^{+}_{y_3}\times\R^{+}_{Y_3}\times\mathbb{T}_\theta)),
	\\
	&   ~\exists \delta >0,\, \forall \alpha \in \mathbb{N}^6,\,\exists C_\alpha >0,\, \forall Y_3 > 0, 
	\,
	\norm{\partial^{\alpha}(\vu,\vB,{\nabla q})(\cdot,\cdot, Y_3,\cdot) }{L^{\infty}_{t,y,\theta}} \leq C_\alpha e^{- \delta Y_3}\},\\
	   \vspace*{3mm}
	 \underline{S}^{-}:=&\{ V=(\vE,\vH) : \vE,\vH \in H^{\infty}([0,T]\times\mathbb{T}^2_{x'}\times\R^{-}_{y_3}\times\mathbb{T}_\theta)) \},\\
{S}^{-}_{\star}:=&\{ V=(\vE,\vH) : \vE,\vH  \in H^{\infty}([0,T]\times\mathbb{T}^2_{x'}\times\R^{-}_{y_3}\times\R^{-}_{Y_3}\times\mathbb{T}_\theta)),
	\\
	&   ~\exists \delta >0,\, \forall \alpha \in \mathbb{N}^6,\,\exists C_\alpha >0,\, \forall Y_3 < 0, \,
	\norm{\partial^{\alpha}(\vE,\vH )(\cdot,\cdot, Y_3,\cdot) }{L^{\infty}_{t,y,\theta}} \leq C_\alpha e^{ \delta Y_3}\},\\
	   \vspace*{2mm}   %
    {S}^{\pm}:=&\underline{S}^{\pm} \oplus {S}^{\pm}_{\star}.
	\end{align*}

The space of profiles is ${S}^{\pm}:=\underline{S}^{\pm} \oplus {S}^{\pm}_{\star}$, where the sum is direct because  functions in $\underline{S}^{\pm} $ do not depend on $Y_3$ and functions in $ {S}^{\pm}_{\star}$ decay exponentially w.r.t. $Y_3$.
\end{Def}

%For simplicity,  we use the convections: $U=(\vu,\vB, q)\in S^+ $ is short for $ \vu,\vB,\nabla q \in S^+$, while $V=(\vE, \vH) \in S^-$ is short for $\vE, \vH \in S^-$. So  do $\underline{S}^\pm$ and $ S^\pm_{\star}$. 
%\todo{Here are the changes of def of the spaces}
%Generally the pressure $q$ cannot expect to decay to $0$ as $y_3$ or $Y_3$ tends to infinity.

The profiles of $(U^m,V^m)$, for $m\geq1$, will be sought in the functional space $S^+\times S^-$, while the profiles $\varphi^m,\, m\geq2$, will be sought in the function space $H^{\infty}([0,T]\times\mathbb{T}^2\times\mathbb{T})$.
Notice the particular definition for what concerns the total pressure, which appears in $\underline{S}^{+},{S}^{+}_{\star}$ through $\nabla q$, instead of $q$ itself. This is because $q$ will be found as the solution of an elliptic equation on the unbounded space domain $\mathbb{T}^2\times\R^+$. For each fixed $(t,Y_3,\theta)$, that appear in the problem as parameters, we will obtain $\nabla q\in H^{\infty}(\mathbb{T}^2\times\R^+)$. However, $q$ is not expected to decay sufficiently fast at infinity to get $ q\in L^2(\mathbb{T}^2\times\R^+)$.

The components on $\underline{S}^{\pm}$ and ${S}^{\pm}_{\star}$ are called the \textit{residual} components and the \textit{surface wave} component, respectively. It seems likely that the WKB cascade below can not be solved with profiles $(U^m, V^m)\in {S}^{+}_{\star}\times{S}^{-}_{\star}$ for all $m\geq1$, that is with only {surface wave} components. Namely, though the leading profile $(U^1,V^1)$ will belong to ${S}^{+}_{\star}\times{S}^{-}_{\star}$, it is likely that one corrector $(U^m,V^m)$ will have a nontrivial residual component, corresponding to a {\it rectification} phenomenon. We refer the reader to \cite{Marcou2011} for a rigorous proof of this phenomenon for a 2D model in elasticity.
Though we are mainly interested in the component on ${S}^{\pm}_{\star}$ of the leading amplitude, nevertheless determining the residual components is one major difficulty in the analysis below.

%\begin{align*}
%	%
%	\underline{S}^{+}:=&\{ U\in \R^7 ; \vu,\vB,\nabla q \in H^{\infty}([0,T]\times\mathbb{T}^2\times\R^{+}\times\mathbb{T}) \},\\
%	%
%	\underline{S}^{-}:=&\{ V\in \R^6 ; \vE,\vH \in H^{\infty}([0,T]\times\mathbb{T}^2\times\R^{-}\times\mathbb{T}) \},\\
%	%
%	{S}^{+}_{\star}:=&\{ U\in \R^7 ; \vu,\vB, \nabla q \in H^{\infty}([0,T]\times\mathbb{T}^2\times\R^{+}\times\R^{+}\times\mathbb{T}),
%	\\
%	&\quad   ~~\exists \delta >0, \forall \alpha \in \mathbb{N}^6,\exists C_\alpha >0, \forall Y_3 > 0, \\
%	&\quad\quad\quad
%	\norm{\partial^{\alpha}(\vu,\vB,\nabla q)(\cdot,\cdot, Y_3,\cdot) }{L^{\infty}_{t,y,\theta}} \leq C_\alpha e^{- \delta Y_3}\},\\
%	%
%	{S}^{-}_{\star}:=&\{ V\in \R^6 ; \vE,\vH  \in H^{\infty}([0,T]\times\mathbb{T}^2\times\R^{-}\times\R^{-}\times\mathbb{T}),
%	\\
%	&\quad   ~~\exists \delta >0, \forall \alpha \in \mathbb{N}^6,\exists C_\alpha >0, \forall Y_3 < 0, \\
%	&\quad\quad\quad
%	\norm{\partial^{\alpha}(\vE,\vH )(\cdot,\cdot, Y_3,\cdot) }{L^{\infty}_{t,y,\theta}} \leq C_\alpha e^{ \delta Y_3}\},\\
%	%
%	{S}^{\pm}:=&\underline{S}^{\pm} \oplus {S}^{\pm}_{\star}.
%\end{align*}

%The residue components are imposed because in \cite{Pierre2021}:
%``It seems likely that the WKB cascade below can not be solved with profiles $U^m, V^m\in {S}^{\pm}_{\star}$ for all $m\geq1$.''

\subsubsection{Notation for profiles}
We expand profiles $(U,V)\in S^+\times S^-$ into Fourier series in $\theta$ as follows
\begin{align*}
	U(t,y,Y_3,\theta) &=\sum\limits_{k\in \mathbb{Z}} \widehat{U} (t,y,Y_3, k) e^{2\pi i k \theta}=U_{\#}(t,y,Y_3,\theta)+ \widehat{U} (t,y,Y_3, 0),\\
	V(t,y,Y_3,\theta)&=\sum\limits_{k\in \mathbb{Z}} \widehat{V} (t,y,Y_3, k) e^{2\pi i k \theta}=V_{\#}(t,y,Y_3,\theta)+ \widehat{V} (t,y,Y_3, 0),
\end{align*}
where we use the subscript ``$\#$'' to represent the summation of nonzero modes with respect to $\theta$. According to the decomposition of $S^\pm$ we also write
\begin{align*}
U(t,y,Y_3,\theta) =\underbrace{\underline{U}(t,y,Y_3,\theta)}_{\in \underline{S}^{+}} + \underbrace{U_{\star}(t,y,Y_3,\theta)}_{\in {S}^{+}_{\star}}, \\ V(t,y,Y_3,\theta) =\underbrace{\underline{V}(t,y,Y_3,\theta)}_{\in \underline{S}^{-}} + \underbrace{V_{\star}(t,y,Y_3,\theta)}_{\in {S}^{-}_{\star}}.
\end{align*}
Consistently with this decomposition we split the zero Fourier modes as 
\begin{align*}
\widehat{U}(0)=\underline{\widehat{U}}(0)+\widehat{U}_{\star}(0), \quad \widehat{V}(0)=\underline{\widehat{V}}(0)+\widehat{V}_{\star}(0),
\end{align*}
with $(\underline{\widehat{U}}(0),\underline{\widehat{V}}(0))$ being referred to as the {\it slow mean} and $({\widehat{U}_{\star}}(0),{\widehat{V}_{\star}}(0))$ as the {\it fast mean}. For the plasma vector solution $U=(\vu,\vB, q)$ we define the projector
\begin{align*}
\Pi U := (u_1, u_2 ,0, B_1, B_2, 0,0)^{\top}
%\quad \Pi V := (0,0,0, 0, 0, 0)^{\top}
,
\end{align*}
consisting of the tangential components of the velocity $\vu$ and the magnetic field $\vB$. The remaining part $(I-\Pi) U=(0,0, u_3 ,0, 0, B_3, q)$ represents the non-characteristic components of $U$.

\subsection{The main result}

In this paper we show the existence of a sequence of profiles $(U^m,V^m,\varphi^{m+1})_{m\geq1}$ such that the asymptotic expansions \eqref{WKB-ansatz} satisfy \eqref{MHD}, \eqref{Maxwell}, \eqref{interface}, with accuracy $O(\e^\infty)$ in the sense of formal series, and the nontriviality of residual non-oscillatory components, $\underline{\widehat{U}}^3(0)$ or $\underline{\widehat{V}}^3(0)$.

\begin{Thm} \label{thm-main}
	Suppose that the reference plasma--vacuum interface defined by \cref{reference-state} satisfy the Assumption \cref{H1}, and the frequencies $\tau, \xi'$ satisfy Assumptions \cref{H2}, \cref{H-key} and \cref{H-a}. 
    Moreover, let $\varphi^2_0 \in  H^{\infty} (\mathbb{T}^2 \times \mathbb{T} )$ with zero mean with respect to $\theta$. \\
    Then there exists $\nu_*>0$ such that for $0<\nu<\nu_*$, there exists a time $T  >0$, only dependent on the norm $\norm{\varphi^2_0}{ H^4(\mathbb{T}^2 \times\mathbb{T} ) }  $, such that there exists a sequence of profiles $(U^m, V^m, \varphi^{m+1})_{m\geq1}$ in $S^+\times S^- \times H^{\infty} ([0,T]\times\mathbb{T}^2 \times\mathbb{T} ) $ satisfying the following properties:
    \begin{itemize}
    	\item for any $m\geq 2$, $\varphi^m|_{t=0}=\delta_{m2}\varphi^2_0$ (with $\delta_{m2}$ the Kronecker delta) and $\partial_t \widehat{\varphi}^m(0)|_{t=0}=0$.
    	\item $\underline{U}^1=\widehat{U}^1_\star(0)=\underline{U}^2=0$, $\underline{V}^1=\widehat{V}^1_\star(0)=\underline{V}^2=0$, while $\underline{U}^3=\underline{\widehat{U}}^3(0)\neq0$ or $\underline{V}^3=\underline{\widehat{V}}^3(0)\neq0$ for general $\varphi^2_0$. 
    	\item for any $m\geq2$, $\Pi\widehat{U}^{m}_\star(0) |_{t=0}=0$, 
    	\item for any integer $M\geq 1$, the functions
    	\begin{align*}
    		U_\e^{\textnormal{app},M}(t,x) & := U^0+\sum\limits_{m= 1}^M \e^m U^m(t,x',y_3,Y_3,\theta),\\
    		V_\e^{\textnormal{app},M}(t,x) & := V^0+\sum\limits_{m\geq 1}^M \e^m V^m(t,x',y_3,Y_3,\theta),\\
    		\varphi_\e^{\textnormal{app},M}(t,x') & := \sum\limits_{m\geq 2}^{M+1} \e^m\varphi^m(t,x',\theta),\\
    		\text{where we have set~~} y_3&:= x_3-\varphi_\e^{\textnormal{app},M} (t,x'),~ Y_3:=\frac{y_3}{\e},~ \theta:=\frac{\tau t+\xi'\cdot x'}{\e}, \text{~satisfy  }
		\end{align*}
   	\begin{align*}
       	&\left\{
       	\begin{aligned}
       		&\dt u_\e^{\textnormal{app},M} + (u_\e^{\textnormal{app},M}\cdot \nabla) u_\e^{\textnormal{app},M} -(B_\e^{\textnormal{app},M} \cdot \nabla ) B_\e^{\textnormal{app},M} +\nabla q_\e^{\textnormal{app},M}=R_\e^{1,+},
       		\\
       		&\dt B_\e^{\textnormal{app},M} + (u_\e^{\textnormal{app},M}\cdot \nabla) B_\e^{\textnormal{app},M} -(B_\e^{\textnormal{app},M} \cdot \nabla ) u_\e^{\textnormal{app},M} =R_\e^{2,+},
       		\\
       		&\dv u_\e^{\textnormal{app},M}=R_\e^{3,+},\qquad \dv B_\e^{\textnormal{app},M}=R_\e^{4,+}   
       		\qquad\qquad\qquad \text{in~} \Omega_\e^{\textnormal{app},M,+}(t),       		
       	 \end{aligned}
        \right.\\
        &\left\{
        \begin{aligned}
       	   &\nu\dt H_\e^{\textnormal{app},M} + \nabla \times E_\e^{\textnormal{app},M}=R_\e^{1,-},
       	   \\
       	   &\nu\dt E_\e^{\textnormal{app},M} - \nabla \times H_\e^{\textnormal{app},M}=R_\e^{2,-},
       	   \\
       	   &\dv H_\e^{\textnormal{app},M}=R_\e^{3,-},\qquad \dv E_\e^{\textnormal{app},M}=R_\e^{4,-}  
       	   \qquad\qquad\qquad\text{in~} \Omega_\e^{\textnormal{app},M,-}(t),       	   
         \end{aligned}
         \right.\\
        &\left\{
        \begin{aligned}
   	          &\dt \varphi_\e^{\textnormal{app},M}-u_\e^{\textnormal{app},M} \cdot N_\e^{\textnormal{app},M} = R_{b,\e}^{1},
   	          \\
   	          &q_\e^{\textnormal{app},M}-\frac{1}{2}\abs{H_\e^{\textnormal{app},M}}{2}+\frac{1}{2}\abs{E_\e^{\textnormal{app},M}}{2}=R_{b,\e}^{2}, \\
   	          &B_\e^{\textnormal{app},M}\cdot N_\e^{\textnormal{app},M}=R_{b,\e}^{3}, 
   	          \qquad H_\e^{\textnormal{app},M}\cdot N=R_{b,\e}^{4}, %&&
   	          \\
   	          & N_\e^{\textnormal{app},M} \times E_\e^{\textnormal{app},M}-\nu (u_\e^{\textnormal{app},M} \cdot N_\e^{\textnormal{app},M})H_\e^{\textnormal{app},M}=R_{b,\e}^{5} 
   	          \,\quad \text{on}~\Gamma_\e^{\textnormal{app},M}(t),
         \end{aligned}
        \right. 
       \end{align*}
       where we have denoted
           \begin{align*}
        \Omega_\e^{\textnormal{app},M,\pm}(t)&:=\{x\in \mathbb{T}^2\times \R; x_3\gtrless\varphi_\e^{\textnormal{app},M}(x') \},\\
        \Gamma_\e^{\textnormal{app},M}(t)&:=\{x\in \mathbb{T}^2\times \R; x_3=\varphi_\e^{\textnormal{app},M}(x') \},\\
        N_\e^{\textnormal{app},M}&:=(-\partial_{x_1} \varphi_\e^{\textnormal{app},M}, -\partial_{x_2} \varphi_\e^{\textnormal{app},M},1 )^{\top},
    	\end{align*}
and where the error terms satisfy the following bounds:
       \begin{align*}
        & \sup\limits_{ t\in [0,T], x\in \Omega_\e^{\textnormal{app},M,\pm}(t)} (\abs{R_\e^{1,\pm} }{} + \abs{R_\e^{2,\pm} }{}+\abs{R_\e^{3,\pm} }{}+\abs{R_\e^{4,\pm} }{})=\mathcal{O} (\e^M),\\
        & \sup\limits_{ t\in [0,T], x\in \Gamma_\e^{\textnormal{app},M}(t)} (\abs{ R_{b,\e}^{1} }{} + \abs{R_{b,\e}^{2} }{}+\abs{R_{b,\e}^{3} }{}+\abs{R_{b,\e}^{4} }{}+\abs{R_{b,\e}^{5} }{})=\mathcal{O} (\e^{M+1}).
       \end{align*}
    \end{itemize}
\end{Thm}

%%%%%%%%%%%%%%%%%%%%%%%%%%%%%%%%%%%%%%%%%%%%%%%%%%%%%%%%%%%%%%%%
%%%%%WKB cascade  #2
%%%%%%%%%%%%%%%%%%%%%%%%%%%%%%%%%%%%%%%%%%%%%%%%%%%%%%%%%%%%%%%%
%\newpage
\section{The WKB cascade} \label{sec-wkb}
In this section we derive the equations for the profiles $(U^m,V^m,\varphi^{m+1})_{m\geq0}$ in the WKB expansion \cref{WKB-ansatz}.
We first rewrite \cref{MHD} and \cref{Maxwell} in conservative form:
\begin{align}
&A^+_0 \dt U+\partial_{x_\alpha} f^+_{\alpha} (U)=0, \label{MHD-conserv}\\
&A^-_0 \dt V+ A^-_\alpha \partial_{x_\alpha} V=0, \label{Maxwell-conserv}
\end{align}
and for later use, we define the Jacobian matrices as 
\begin{align*}
A^+_\alpha := \text{d}  f^+_{\alpha} (U^0), \quad 
&\sA^\pm:=\tau A^\pm_0 +\xi_j A^\pm_j, \quad
\mathbb{A}_{\alpha}(\cdot,\cdot):=\text{d}^2  f^+_{\alpha} (U^0)(\cdot,\cdot),%\\
%A^-_\alpha := \text{d}  f^-_{\alpha} (U^0), \quad
% &\sA^-:=\tau A^-_0 +\xi_j A^-_j.
\end{align*}
where $f^+_{\alpha}, A^\pm_\alpha, \sA^\pm, \mathbb{A}_{\alpha}$ are given explicitly in \cref{appen-a}.
Note that $f_\alpha^+$ is a polynomial of degree at most 2, and thus $\mathbb{A}_{\alpha}$ does not depend on the state $U$, and $\text{d}^3  f^+_{\alpha}=0$.

Here we do not include in \cref{MHD-conserv} and \cref{Maxwell-conserv} the divergence-free constraints on the magnetic and electric fields, 
but keep the constraints as extra equations.
Actually, the divergence-free constraints on the magnetic and electric fields hold as long as they are satisfied at the initial data (see e.g. \cite{Sun2019}).

\subsection{Equations}
We substitute the WKB expansion \cref{WKB-ansatz} into \cref{MHD-conserv} and \cref{Maxwell-conserv}, 
collect the result in the form of a powers series with respect to the parameter $\e$; then we equate coefficients of all powers $\e^{m} (m\geq0)$ to zero, and obtain the so-called WKB cascade.
In the interior of the plasma and vacuum regions, $(U^m,V^m,\varphi^{m+1})$ satisfies the {\it fast system}
\begin{align}
&\forall m\geq 0,~ (t,y',y_3,Y_3, \theta) \in [0,T]\times\mathbb{T}^2\times\R^{\pm}\times\R^{\pm}\times\mathbb{T}, \notag
\\&\quad
\sL^+_f U^{m+1}=F^{m,+},\quad 
\sL^-_f V^{m+1}=F^{m,-},  \label{WKB-eq}
\end{align}
where the {\it fast} operators $\sL^{\pm}_f=\sL^{\pm}_f(\partial)$ are defined by
\begin{equation}
\label{def-sL-fast}
\sL^{\pm}_f =A^\pm_3 \partial_{Y_3}+\sA^\pm\partial_{\theta},
\end{equation}
and the source terms are
\begin{align}\label{Fm+}
F^{m,+} =& - L^+_s U^{m} 
+  \sum\limits_{\substack{l_1+l_2=m+2\\l_1\geq 1}} \partial_{\theta} \varphi^{l_2} \sA^+ \partial_{Y_3} U^{l_1}
+\sum\limits_{\substack{l_1+l_2=m+1\\l_1\geq 1}} (\dt \varphi^{l_2} A^+_0 + \partial_{y_j} \varphi^{l_2} A^+_j) \partial_{Y_3} U^{l_1}  \notag \\
&  \quad 
+  \sum\limits_{\substack{l_1+l_2=m+1\\l_1\geq 1}} \partial_{\theta} \varphi^{l_2} \sA^+ \partial_{y_3} U^{l_1}
+\sum\limits_{\substack{l_1+l_2=m\\l_1\geq 1}} (\dt \varphi^{l_2} A^+_0 + \partial_{y_j} \varphi^{l_2} A^+_j) \partial_{y_3} U^{l_1} \notag \\
&  \quad 
-\sum\limits_{\substack{l_1+l_2=m\\l_1,l_2\geq 1}} \mathbb{A}_{\alpha} (U^{l_1},\partial_{y_\alpha} U^{l_2})
-\sum\limits_{\substack{l_1+l_2=m+1\\l_1,l_2\geq 1}} \xi_j \mathbb{A}_{j} (U^{l_1},\partial_{\theta} U^{l_2})  
-\sum\limits_{\substack{l_1+l_2=m+1\\l_1,l_2\geq 1}} \mathbb{A}_{3} (U^{l_1},\partial_{Y_3} U^{l_2}) \notag  \\
&  \quad 
+\sum\limits_{\substack{l_1+l_2+l_3=m+2\\l_1,l_2\geq 1}} \partial_{\theta} \varphi^{l_3} \xi_j \mathbb{A}_{j}  (U^{l_1},\partial_{Y_3} U^{l_2})
+\sum\limits_{\substack{l_1+l_2+l_3=m+1\\l_1,l_2\geq 1}} \partial_{y_j} \varphi^{l_3}  \mathbb{A}_{j}  (U^{l_1},\partial_{Y_3} U^{l_2})  \notag \\
&  \quad 
+\sum\limits_{\substack{l_1+l_2+l_3=m+1\\l_1,l_2\geq 1}} \partial_{\theta} \varphi^{l_3} \xi_j \mathbb{A}_{j}  (U^{l_1},\partial_{y_3} U^{l_2})
+\sum\limits_{\substack{l_1+l_2+l_3=m\\l_1,l_2\geq 1}} \partial_{y_j} \varphi^{l_3}  \mathbb{A}_{j}  (U^{l_1},\partial_{y_3} U^{l_2})   ,
\end{align}
\begin{align}\label{Fm-}
	F^{m,-} =& - L^-_s V^{m} 
	+  \sum\limits_{\substack{l_1+l_2=m+2\\l_1\geq 1}} \partial_{\theta} \varphi^{l_2} \sA^- \partial_{Y_3} V^{l_1}
	+\sum\limits_{\substack{l_1+l_2=m+1\\l_1\geq 1}} (\dt \varphi^{l_2} A^-_0 + \partial_{y_j} \varphi^{l_2} A^-_j) \partial_{Y_3} V^{l_1}  \notag \\
	&  \quad 
	+  \sum\limits_{\substack{l_1+l_2=m+1\\l_1\geq 1}} \partial_{\theta} \varphi^{l_2} \sA^- \partial_{y_3} V^{l_1}
	+\sum\limits_{\substack{l_1+l_2=m\\l_1\geq 1}} (\dt \varphi^{l_2} A^-_0 + \partial_{y_j} \varphi^{l_2} A^-_j) \partial_{y_3} V^{l_1}.   
\end{align}
In \eqref{Fm+}, \eqref{Fm-}, the {\it slow} operators $L^{\pm}_s=L^{\pm}_s(\partial)$ are defined by
\begin{equation}
	\label{def-sL-slow}
		L^{\pm}_s =A^\pm_0 \dt+A^\pm_{\alpha}\partial_{y_{\alpha}}.
\end{equation}
When $m=0,1$, $F^{m,\pm}$ can be simplified as 
\begin{align}
&
\label{F01+}
\left\{
\begin{aligned}
& F^{0,+}=0,\\
& F^{1,+}=- L^+_s U^{1} + \partial_{\theta} \varphi^2 \sA^+ \partial_{Y_3} U^{1}
-\frac{1}{2} (\xi_j  \partial_{\theta} \mathbb{A}_{j} (U^{1}, U^{1}) + \partial_{Y_3} \mathbb{A}_{3} (U^{1}, U^{1})),
\end{aligned}
\right.\\
&
\label{F01-}
\left\{
\begin{aligned}
& F^{0,-}=0,\\
& F^{1,-}=- L^-_s V^{1} + \partial_{\theta} \varphi^2 \sA^- \partial_{Y_3} V^{1}.
\end{aligned}
\right.
\end{align}

Substituting the WKB expansion \cref{WKB-ansatz} into the divergence-free condition for $\vu$ yields
\begin{align}\label{F8m}
\partial_{Y_3}u^{m+1}_3+\xi_j\partial_{\theta} u^{m+1}_j
&=F^{m,+}_7:=-\nabla \cdot u^{m} + \sum\limits_{\substack{l_1+l_2=m+2\\l_1\geq 1}} \partial_{\theta} \varphi^{l_2} \xi_j \partial_{Y_3} u_{j}^{l_1} 
+ \sum\limits_{\substack{l_1+l_2=m+1\\l_1\geq 1}} \partial_{y_j} \varphi^{l_2}  \partial_{Y_3} u_{j}^{l_1}  \notag  \\
&\quad + \sum\limits_{\substack{l_1+l_2=m+1\\l_1\geq 1}} \partial_{\theta} \varphi^{l_2} \xi_j \partial_{y_3} u_{j}^{l_1} 
+ \sum\limits_{\substack{l_1+l_2=m\\l_1\geq 1}} \partial_{y_j} \varphi^{l_2}  \partial_{y_3} u_{j}^{l_1}.
\end{align}
Similar equations are derived from the  divergence-free conditions for $\vB, \vH, \vE$, and read:
\begin{equation}
	\label{WKB-div} 
	\begin{aligned}
	\partial_{Y_3}B^{m+1}_3+\xi_j\partial_{\theta} B^{m+1}_j
	&=F^{m,+}_8, \\
	\partial_{Y_3}H^{m+1}_3+\xi_j\partial_{\theta} H^{m+1}_j
	&=F^{m,-}_7,
	\\
	\partial_{Y_3}E^{m+1}_3+\xi_j\partial_{\theta} E^{m+1}_j
	&=F^{m,-}_8, 
	\end{aligned}
\end{equation}
where the source terms $F^{m,+}_8, F^{m,-}_7$ $F^{m,-}_8$, are defined analogously to $F^{m,+}_7$ in \eqref{F8m}. Equations \eqref{F8m}, \eqref{WKB-div} are the {\it fast divergences} constraints.
When $m=0,1$, $F_7^{m,+}$ can be simplified as 
\begin{align}
\label{F701}
\left\{
\begin{aligned}
& F_7^{0,+}=0,\\
& F_7^{1,+}=-\nabla \cdot u^1 + \xi_j \partial_{\theta} \varphi^2  \partial_{Y_3} u_j^1.
\end{aligned}
\right.
\end{align}
Similarly, $F_7^{0,-}=F_8^{0,\pm}=0,$ while $F_7^{1,-},  F_8^{1,\pm}$ have expressions similar to $F_7^{1,+}$, with $\vu$ respectively substituted by $ \vH, \vB,\vE$.

\subsection{Jump conditions}
We plug \cref{WKB-ansatz} in the jump condition \cref{interface} and equate coefficients of $\e^{m} (m\geq0)$ to zero, 
and thus obtain the jump conditions of $(U^m,V^m,\varphi^{m+1})$:
\begin{align*}
%equation 1
u^{m+1}_3-c^+\partial_{\theta} \varphi^{m+2}=&\dt \varphi^{m+1} +u^0_j \partial_{y_j} \varphi^{m+1}\\
& +\sum\limits_{\substack{l_1+l_2=m+2\\l_1\geq 1}} \xi_j   u^{l_1}_j \partial_{\theta}\varphi^{l_2}
+\sum\limits_{\substack{l_1+l_2=m+1\\l_1\geq 1}}  u^{l_1}_j \partial_{y_j}\varphi^{l_2},\\
%equation 2
B^{m+1}_3-b^+\partial_{\theta} \varphi^{m+2}=&B^0_j \partial_{y_j} \varphi^{m+1}+ \sum\limits_{\substack{l_1+l_2=m+2\\l_1\geq 1}}  \partial_{\theta}\varphi^{l_2} \xi_j   B^{l_1}_j
+\sum\limits_{\substack{l_1+l_2=m+1\\l_1\geq 1}}  \partial_{y_j}\varphi^{l_2} B^{l_1}_j,\\
%equation 3
H^{m+1}_3-b^-\partial_{\theta} \varphi^{m+2}=&H^0_j \partial_{y_j} \varphi^{m+1}+ \sum\limits_{\substack{l_1+l_2=m+2\\l_1\geq 1}} \partial_{\theta}\varphi^{l_2} \xi_j   H^{l_1}_j
+\sum\limits_{\substack{l_1+l_2=m+1\\l_1\geq 1}}  \partial_{y_j}\varphi^{l_2} H^{l_1}_j,
\end{align*}
\begin{align}
%equation 4
E^{m+1}_2+\xi_2 E^0_3\partial_{\theta}\varphi^{m+2} +\nu\tau H^0_1 \partial_{\theta}\varphi^{m+2} 
=&-E^0_3\partial_{y_2}\varphi^{m+1} -\nu H^0_1 \dt \varphi^{m+1}\notag\\
-\sum\limits_{\substack{l_1+l_2=m+2\\l_1\geq 1}}  \partial_{\theta}\varphi^{l_2}\xi_2   E^{l_1}_3
&-\sum\limits_{\substack{l_1+l_2=m+1\\l_1\geq 1}}    \partial_{y_2}\varphi^{l_2} E^{l_1}_3\notag\\
-\sum\limits_{\substack{l_1+l_2=m+2\\l_1\geq 1}} \partial_{\theta}\varphi^{l_2} \nu \tau  H^{l_1}_1
&-\sum\limits_{\substack{l_1+l_2=m+1\\l_1\geq 1}} \dt\varphi^{l_2}\nu  H^{l_1}_1 ,\notag\\
%equation 5
-E^{m+1}_1-\xi_1 E^0_3\partial_{\theta}\varphi^{m+2} +\nu\tau H^0_2 \partial_{\theta}\varphi^{m+2} 
=&E^0_3\partial_{y_1}\varphi^{m+1} -\nu H^0_2 \dt \varphi^{m+1}\notag\\
+\sum\limits_{\substack{l_1+l_2=m+2\\l_1\geq 1}} \partial_{\theta}\varphi^{l_2} \xi_1   E^{l_1}_3
&+\sum\limits_{\substack{l_1+l_2=m+1\\l_1\geq 1}}  \partial_{y_1}\varphi^{l_2}   E^{l_1}_3\notag\\
-\sum\limits_{\substack{l_1+l_2=m+2\\l_1\geq 1}} \partial_{\theta}\varphi^{l_2} \nu \tau  H^{l_2}_2
&-\sum\limits_{\substack{l_1+l_2=m+1\\l_1\geq 1}}  \dt\varphi^{l_2} \nu H^{l_1}_2,\notag\\
%equation 6
q^{m+1}-H^0_jH^{m+1}_j+E^0_3E^{m+1}_3 =
\frac{1}{2} \sum\limits_{\substack{l_1+l_2=m+1\\l_1,l_2\geq 1}}& H^{l_1}_{\alpha}H^{l_2}_{\alpha}-\frac{1}{2} \sum\limits_{\substack{l_1+l_2=m+1\\l_1,l_2\geq 1}} E^{l_1}_{\alpha}E^{l_2}_{\alpha}.\label{WKB-jump-q}
\end{align}
For simplicity, the above equations can be rewritten in the form:
\begin{align}\label{WKB-jump}
	\forall m\geq 0, \qquad
\mathbb{B}^+ U^{m+1} |_{y_3=Y_3=0}+\mathbb{B}^- V^{m+1} |_{y_3=Y_3=0}+\underline{b} \partial_{\theta} \varphi^{m+2} =G^m.
\end{align}
where the matrices $\mathbb{B}^\pm, \underline{b}$ and the source term $G^m$ are
\begin{align*}
	&\mathbb{B}^+ =\begin{pmatrix}
		0 & 0 &1 &0 &0 &0 & 0\\
		0 & 0 &0 &0 &0 &1 & 0\\
		0 & 0 &0 &0 &0 &0 & 0\\
		0 & 0 &0 &0 &0 &0 & 0\\
		0 & 0 &0 &0 &0 &0 & 0\\
		0 & 0 &0 &0 &0 &0 & 1
	\end{pmatrix}, \quad 
	\mathbb{B}^- =\begin{pmatrix}
		0 & 0 &0 &0 &0 &0 \\
		0 & 0 &0 &0 &0 &0 \\
		0 & 0 &1 &0 &0 &0 \\
		0 & 0 &0 &0 &1 &0 \\
		0 & 0 &0 &-1 &0 &0 \\
		-H^0_1 & -H^0_2 &0 &0 &0 &E^0_3 
	\end{pmatrix},\\
	& \underline{b}=-(c^+,b^+,b^-,-a^-_1, -a^-_2, 0)^\top,
\end{align*}
\begin{subequations}\label{Gm}
	\begin{align}
	G^m_1=&
	\dt \varphi^{m+1} +u^0_j \partial_{y_j} \varphi^{m+1}
	+\sum\limits_{\substack{l_1+l_2=m+2\\l_1\geq 1}} \partial_{\theta}\varphi^{l_2} \xi_j   u^{l_1}_j |_{y_3=Y_3=0}
	+\sum\limits_{\substack{l_1+l_2=m+1\\l_1\geq 1}} \partial_{y_j}\varphi^{l_2} u^{l_1}_j |_{y_3=Y_3=0},\\
	G^m_2=&
	B^0_j \partial_{y_j} \varphi^{m+1}+ \sum\limits_{\substack{l_1+l_2=m+2\\l_1\geq 1}} \partial_{\theta}\varphi^{l_2} \xi_j   B^{l_1}_j |_{y_3=Y_3=0} 
	+\sum\limits_{\substack{l_1+l_2=m+1\\l_1\geq 1}} \partial_{y_j}\varphi^{l_2} B^{l_1}_j |_{y_3=Y_3=0} ,\\
	G^m_3=&H^0_j \partial_{y_j} \varphi^{m+1}+ \sum\limits_{\substack{l_1+l_2=m+2\\l_1\geq 1}}\partial_{\theta}\varphi^{l_2} \xi_j   H^{l_1}_j |_{y_3=Y_3=0} 
	+\sum\limits_{\substack{l_1+l_2=m+1\\l_1\geq 1}} \partial_{y_j}\varphi^{l_2} H^{l_1}_j |_{y_3=Y_3=0},&
	\\
	G^m_4=&
	-E^0_3\partial_{y_2}\varphi^{m+1} -\nu H^0_1 \dt \varphi^{m+1}-\sum\limits_{\substack{l_1+l_2=m+2\\l_1\geq 1}}\partial_{\theta}\varphi^{l_2} \xi_2   E^{l_1}_3 |_{y_3=Y_3=0} 
	-\sum\limits_{\substack{l_1+l_2=m+1\\l_1\geq 1}} \partial_{y_2}\varphi^{l_2}     E^{l_1}_3 |_{y_3=Y_3=0}  
	\notag\\
	&-\sum\limits_{\substack{l_1+l_2=m+2\\l_1\geq 1}}\partial_{\theta}\varphi^{l_2} \nu \tau  H^{l_1}_1 |_{y_3=Y_3=0} 
	-\sum\limits_{\substack{l_1+l_2=m+1\\l_1\geq 1}}\dt\varphi^{l_2} \nu  H^{l_1}_1 |_{y_3=Y_3=0} ,\\
	G^m_5=&
	E^0_3\partial_{y_1}\varphi^{m+1} -\nu H^0_2 \dt \varphi^{m+1}
	+\sum\limits_{\substack{l_1+l_2=m+2\\l_1\geq 1}}\partial_{\theta}\varphi^{l_2} \xi_1   E^{l_1}_3 |_{y_3=Y_3=0} 
	+\sum\limits_{\substack{l_1+l_2=m+1\\l_1\geq 1}} \partial_{y_1}\varphi^{l_2}    E^{l_1}_3 |_{y_3=Y_3=0} 
	\notag \\
	&-\sum\limits_{\substack{l_1+l_2=m+2\\l_1\geq 1}} \partial_{\theta}\varphi^{l_2} \nu \tau  H^{l_2}_2 |_{y_3=Y_3=0}
	-\sum\limits_{\substack{l_1+l_2=m+1\\l_1\geq 1}} \dt\varphi^{l_2} \nu H^{l_1}_2 |_{y_3=Y_3=0} ,\\
	G^m_6=& (\frac{1}{2} \sum\limits_{\substack{l_1+l_2=m+1\\l_1,l_2\geq 1}} H^{l_1}_{\alpha}H^{l_2}_{\alpha}-\frac{1}{2} \sum\limits_{\substack{l_1+l_2=m+1\\l_1,l_2\geq 1}} E^{l_1}_{\alpha}E^{l_2}_{\alpha}) |_{y_3=Y_3=0}.\label{Gm6}
	\end{align}
\end{subequations}

\vspace{0.2cm}
In the rest of the paper, we will determine $(U^m,V^m,\varphi^{m+1})_{m\geq1}$ by the fast system \cref{WKB-eq}, together with the fast divergence constrains \cref{WKB-div} and the jump condition \cref{WKB-jump}.
In agreement with \cref{initial-front}, the inital data for the front profiles $\varphi^{m+1}$ are 
\begin{equation}
	\label{WKB-front}
	\forall (x',\theta) \in \mathbb{T}^3, ~
	\varphi^2(0,y',\theta)=\varphi^2_0(y',\theta), \quad 
	\text{and}~
	\forall m\geq3,~
	\varphi^m(0,y',\theta)=0.
\end{equation}
Recalling the setting of $\varphi_0^2$, at the initial time $\varphi^m$ has zero mean with respect to $\theta$ for all $m\geq2$.

%%%%%%%%%%%%%%%%%%%%%%%%%%%%%%%%%%%%%%%%%%%%%%%%%%%%%%%%%%%%%%%%
%%%%%Fast problem  #3
%%%%%%%%%%%%%%%%%%%%%%%%%%%%%%%%%%%%%%%%%%%%%%%%%%%%%%%%%%%%%%%%
%\newpage
\vspace{0.5cm}
\section{The fast problem}\label{sec-fast}
The construction of a solution to the WKB cascade \cref{WKB-eq}, \cref{WKB-div}, \cref{WKB-jump} is done by induction. At each step of the inductive process 
the main issue is to solve the following so-called \textit{fast problem}:
\begin{align}\label{f}
\left\{
\begin{aligned}
&\sL^+_f U =F^+, &y_3,Y_3>0,\\
&\partial_{Y_3}B_3+\xi_j \partial_{\theta} B_j =F^+_8,&y_3,Y_3>0,\\
&\sL^-_f V =F^-,&y_3,Y_3<0,\\
&\partial_{Y_3}H_3+\xi_j \partial_{\theta}H_j =F^-_7,&y_3,Y_3<0,\\
&\partial_{Y_3}E_3+\xi_j \partial_{\theta}E_j =F^-_8,&y_3,Y_3<0,\\
&\mathbb{B}^+ U  |_{y_3=Y_3=0}+\mathbb{B}^- V  |_{y_3=Y_3=0}+\underline{b} \partial_{\theta} \varphi =G.
\end{aligned}
\right.
\end{align}
where $\sL^{\pm}_f$ is defined in \cref{def-sL-fast}, the matrices $\mathbb{B}^\pm$ and vector $\underline{b}$ can be found in \cref{appen-a}. 

Notice that the slow variables $(t,y')$ and $y_3$ act as parameters in \cref{f}, and that the boundary conditions on $\Gamma_0$ requires $y_3=Y_3=0$.
An important feature of the fast problem \eqref{f} (and some other problems in the subsequent sections) is that, when  ignoring the front profile $\varphi$ or just regarding it as a given term, the equations and the boundary conditions in \cref{f} for the plasma and those for vacuum can be solved independently. Since the equations for the plasma are almost the same as those in \cite{Pierre2021}, most of times we can directly adopt the results and conclusions in \cite{Pierre2021} for the plasma part.

In this section, we prove the solvability of the fast problem \cref{f} for an arbitrary final time $T>0$.
\begin{Thm}\label{thm-fast}
	Suppose that Assumptions \cref{H1,H2,H-key,H-a} hold, and let $F^{\pm}, F^-_7, F^{\pm}_8$$\in S^{\pm}$ and $G\in H^{\infty} ([0,T]\times\mathbb{T}^2\times\mathbb{T})$. Then the following facts hold.
	
	(i) \cref{f} has a solution $(U, V, \varphi)$ in $S^+\times S^- \times H^{\infty} ([0,T]\times\mathbb{T}^2 \times\mathbb{T} ) $ if and only if the following compatibility conditions are satisfied:
	\begin{subequations}\label{f-comp}
		\begin{align} 	
		&
		\left\{\begin{aligned}
		&F^+_6 |_{y_3=Y_3=0}=-b^+ \partial_{\theta}G_1+c^+\partial_{\theta}G_2,\\
		&F^-_3 |_{y_3=Y_3=0}=\nu\tau \partial_{\theta}G_3+\xi_j\partial_{\theta}G_{3+j},
		\end{aligned}\right. \label{comp-a}
		\\[1em]
		&
		\underline{\widehat{F}}^{\pm}(t,y,0)=0, ~\underline{\widehat{F}}_7^-(t,y,0)= \underline{\widehat{F}}_8^{\pm}(t,y,0)=0,
		\label{comp-b}\\[1em]
		&
		\left\{
		\begin{aligned}
			&u_j^0 \widehat{F}^+_{7,\star}(0) -B_j^0 \widehat{F}^+_{8,\star}(0)=\widehat{F}^+_{j,\star} (0),\quad j=1,2,\\
			&B_j^0 \widehat{F}^+_{7,\star}(0) -u_j^0 \widehat{F}^+_{8,\star}(0)=\widehat{F}^+_{3+j,\star} (0),\quad j=1,2,\\
			%	&\widehat{F}^-_{3,\star} (0)=0.
		\end{aligned}
		\right.\label{comp-c}\\[1em]
		&
		\left\{
		\begin{aligned}
		&\partial_{Y_3} F_6^+ + \xi_j \partial_{\theta} F^+_{3+j}-\tau \partial_{\theta} F^+_8=0,\\
		&\partial_{Y_3} F_3^- + \xi_j \partial_{\theta} F^-_{j}-\nu\tau \partial_{\theta} F^-_7=0,\\
		&\partial_{Y_3} F_6^- + \xi_j \partial_{\theta} F^-_{3+j}-\nu\tau \partial_{\theta} F^-_8=0,
		\end{aligned}
		\right.
		\label{comp-d}\\[1em]
		&
		\begin{aligned}
		&\quad \int_{0}^{+\infty}e^{-|k|Y_3}\sL^+(k) \cdot \widehat{F}^+(t,y',0,Y_3,k) dY_3 + l^+_j \widehat{G}_j(t,y',k)\\
		& -\int_{-\infty}^{0} e^{|k|\sqrt{1-\nu^2\tau^2}Y_3}\sL^-(k) \cdot \widehat{F}^-(t,y',0,Y_3,k)  dY_3 -\frac{-\nu\tau a^-_j +\xi_j b^-}{\nu\sqrt{1-\nu^2\tau^2}}  \widehat{G}_{3+j}(t,y',k) \\
		& -i\sgnk \tau \widehat{G}_6(t,y',k) +\frac{E^0_3}{|k|\nu} \widehat{F}^-_6(t,y',0,0,k)=0, \quad \text{for all~} k\in \mathbb{Z}\setminus \{0\}.
		\end{aligned}\label{comp-e}
		\end{align}
	\end{subequations}
	Here we  use  ``$\cdot$'' to denote the inner product of vectors of complex numbers : $X\cdot Y = \sum_{\alpha} \overline{X}_\alpha Y_\alpha$. $\sL^{\pm}(k)$ is defined in \cref{vector-L} and $l^+_j$, $j=1,2$ are defined in \cref{l+}. 
	
	(ii) If the compatibility conditions \eqref{f-comp} are satisfied, then \cref{f} has a solution of the form $(\mathbb{U},\mathbb{V},0)$ with $\Pi\, \widehat{\mathbb{U}}(0)=0$. %\todo{need condition for $\Pi \widehat{\mathbb{U}}(0)$}
	
	(iii) Any solution to \cref{f} is in the form
	\begin{align}
		U(t,y,Y_3,\theta)=&\mathbb{U}(t,y,Y_3,\theta)+\underline{U}_h(t,y)+(\widehat{u}_{1,\star}(0),\widehat{u}_{2,\star}(0),0,\widehat{B}_{1,\star}(0),\widehat{B}_{2,\star}(0),0,0)^{\top} (t,y,Y_3)  \notag \\
		&\quad + \sum\limits_{k\in \mathbb{Z} \setminus \{0\} }\gamma^+(t,y,k) e^{-|k|Y_3+ik\theta} \sR^+(k), \label{decomp-1} \\
		%line 
		V(t,y,Y_3,\theta)=&\mathbb{V}(t,y,Y_3,\theta)+\underline{V}_h(t,y)+\sum\limits_{k\in \mathbb{Z} \setminus \{0\}}\gamma^-_j(t,y,k) e^{|k|\sqrt{1-\nu^2\tau^2}Y_3+ik\theta} \sR^-_j(k), \notag
	\end{align}
	where $\underline{U}_h \in \underline{S}^{+}$ and $\underline{V}_h \in \underline{S}^{-}$ are independent of $\theta$ and satisfy
	the boundary conditions on $\Gamma_0$:
	\begin{align*}
	&\underline{u}_{h,3} |_{y_3=0}=\underline{B}_{h,3} |_{y_3=0}=\underline{H}_{h,3} |_{y_3=0}=\underline{E}_{h,1} |_{y_3=0}=\underline{E}_{h,2} |_{y_3=0}=0,\\
	& \underline{q}_h |_{y_3=0}-H^0_j \underline{H}_{h,j} |_{y_3=0}+E^0_3 \underline{E}_{h,3} |_{y_3=0}=0.
	\end{align*}
	Moreover,
	$\widehat{u}_{1,\star}(0),\widehat{u}_{2,\star}(0), \widehat{B}_{1,\star}(0),\widehat{B}_{2,\star}(0) \in S^+_{\star}$, the vectors $\sR^+(k),\sR^-_1(k), \sR^-_2(k)$ are explicitely given in \cref{vector-R}, and the coefficients $\gamma^+, \gamma^-_1,\gamma^-_2$ satisfy, for all $(t,y',k)\in [0,T]\times \mathbb{T}^2\times (\mathbb{Z}\setminus \{0\})$,
	\begin{equation}
		\label{def-gamma}
		\begin{aligned}
			\quad\quad \gamma^+(t,y',0,k)=|k|\widehat{\varphi}(t,y',k),
			\qquad\gamma^-_j(t,y',0,k)=\widetilde{\gamma}_j|k|\widehat{\varphi}(t,y',k),
		\end{aligned}
	\end{equation}
where $\widetilde{\gamma}_j, j=1,2,$  are two constants satisfying the relation
 	\begin{equation}
 	\label{def-gamma2} a^-_j\widetilde{\gamma}_j=\frac{(a_1^-)^2+(a_2^-)^2-(b^-)^2}{1-\nu^2\tau^2 }.
 		\end{equation}
 	The difference $(U-\mathbb{U},V-\mathbb{V},\varphi)$ is a solution to \cref{f-h}, that is problem \cref{f} with zero source terms.
 	
\end{Thm}

\begin{Rem}\label{rem-decomp}
	The solution to \cref{f} can also be decomposed in the following way: 
	\begin{align}
	U(t,y,Y_3,\theta)=&\mathbb{U}_{\#}(t,y,Y_3,\theta)+(0,0,\widehat{\sU}_{3,\star}(0),0,0,\widehat{\sB}_{3,\star}(0),\widehat{\sQ}_{\star}(0))^{\top}(t,y,Y_3) \notag \\
	&  +(\widehat{u}_{1,\star}(0),\widehat{u}_{2,\star}(0),0,\widehat{B}_{1,\star}(0),\widehat{B}_{2,\star}(0),0,0)^{\top} (t,y,Y_3) \notag \\
	& 
	+\underline{\widehat{U}}(0)(t,y)+ \sum\limits_{k\in \mathbb{Z} \setminus \{0\} }\gamma^+(t,y,k) e^{-|k|Y_3+ik\theta} \sR^+(k), \label{decomp-2} \\
	%line 
	V(t,y,Y_3,\theta)=&\mathbb{V}_{\#}(t,y,Y_3,\theta)+\widehat{\mathbb{V}}_{\star}(0) (t,y,Y_3)\notag \\
	& 
	+\underline{\widehat{V}}(0)(t,y)+\sum\limits_{k\in \mathbb{Z} \setminus \{0\}}\gamma^-_j(t,y,k) e^{|k|\sqrt{1-\nu^2\tau^2}Y_3+ik\theta} \sR^-_j(k). \notag
	\end{align}
	In \cref{decomp-2} $\mathbb{U}_{\#}$ and $ (0,0,\widehat{\sU}_{3,\star}(0),0,0,\widehat{\sB}_{3,\star}(0),\widehat{\sQ}_{\star}(0))$ (and similarly $\mathbb{V}_{\#}$ and $\widehat{\mathbb{V}}_{\star}(0)$, see \cref{def-0-Ustar} and \cref{def-Usharp}) originate from $\mathbb{U}$ (or $\mathbb{V}$), and they are determined by the source terms.
	$\underline{\widehat{U}}(0)$ (or $\underline{\widehat{V}}(0)$ ) comes from $\mathbb{U}$ and $\underline{U}_h$ (or $\mathbb{V}$ and $\underline{V}_h$), and they just have constraints at the boundary.
	$\gamma^+$ and $\gamma^-_j$ are determined in terms of $\varphi$, and we {have one degree of freedom} in the choice of the coefficients $\widetilde{\gamma}_j$ because they are linked by \cref{def-gamma2}.  $\Pi\, \widehat{U}_\star(0)$ can be freely chosen. 
\end{Rem}

As in \cite{Pierre2021}, the compatibility conditions \eqref{f-comp} do not involve the mean $\widehat{G}(0)$ of the boundary source term. The mean $\widehat{G}(0)$ will appear in the induction process when determining the mean of the front profiles $\widehat{\varphi}^m$, $m\geq2$.

\vspace{0.3cm}
\subsection{The homogeneous case}
We first consider the homogeneous case of the fast problem \cref{f}, i.e.
\begin{align}\label{f-h}
\left\{
\begin{aligned}
&\sL^+_f U =0, &y_3,Y_3>0,\\
&\partial_{Y_3}B_3+\xi_j \partial_{\theta} B_j =0,&y_3,Y_3>0,\\
&\sL^-_f V =0,&y_3,Y_3<0,\\
&\partial_{Y_3}H_3+\xi_j \partial_{\theta}H_j =0,&y_3,Y_3<0,\\
&\partial_{Y_3}E_3+\xi_j \partial_{\theta}E_j =0,&y_3,Y_3<0,\\
&\mathbb{B}^+ U  |_{y_3=Y_3=0}+\mathbb{B}^- V  |_{y_3=Y_3=0}+\underline{b} \partial_{\theta} \varphi =0.
\end{aligned}
\right.
\end{align}
%We will discuss the residual components and surface wave components of the solutions separately. 
The solutions to \cref{f-h} are characterized as follows.
\begin{Prop}\label{prop-h}
	Given $\varphi\in H^{\infty}([0,T]\times\mathbb{T}^2\times\mathbb{T})$, the solutions $(U,V,\varphi) \in S^+ \times S^-\times H^{\infty}([0,T]\times\mathbb{T}^2\times\mathbb{T})$ to \cref{f-h} are 
	in the form:
	\begin{align*}
	U(t,y,Y_3,\theta)=&\underline{U}(t,y)+(\widehat{u}_{1,\star}(0),\widehat{u}_{2,\star}(0),0,\widehat{B}_{1,\star}(0),\widehat{B}_{2,\star}(0),0,0)^{\top} (t,y,Y_3)\\
	&\quad + \sum\limits_{k\in \mathbb{Z} \setminus \{0\} }\gamma^+(t,y,k) e^{-|k|Y_3+ik\theta} \sR^+(k),\\
	%line 
	V(t,y,Y_3,\theta)=&\underline{V}(t,y)+\sum\limits_{k\in \mathbb{Z} \setminus \{0\}}\gamma^-_j(t,y,k) e^{|k|\sqrt{1-\nu^2\tau^2}Y_3+ik\theta} \sR^-_j(k),
	\end{align*}
	where $\underline{U} \in \underline{S}^{+}$ and $\underline{V} \in \underline{S}^{-}$ are independent of $\theta$ and satisfy
	the boundary conditions on $\Gamma_0$:
	\begin{align*}
	&\underline{u}_3 |_{y_3=0}=\underline{B}_3 |_{y_3=0}=\underline{H}_3 |_{y_3=0}=\underline{E}_1 |_{y_3=0}=\underline{E}_2 |_{y_3=0}=0,\\
	& \underline{q} |_{y_3=0}-H^0_j \underline{H}_j |_{y_3=0}+E^0_3 \underline{E}_3 |_{y_3=0}=0.
	\end{align*}
	Moreover,
	$\widehat{u}_{1,\star},\widehat{u}_{2,\star}, \widehat{B}_{1,\star},\widehat{B}_{2,\star} \in S^+_{\star}$, 
	the vectors $\sR^+(k),\sR^-_1(k), \sR^-_2(k)$ are explicitly given in \cref{vector-R}, and the coefficients $\gamma^+, \gamma^-_1,\gamma^-_2$ satisfy, for all $(t,y',k)\in [0,T]\times \mathbb{T}^2\times (\mathbb{Z}\setminus \{0\})$,
	\begin{equation*}
		\begin{aligned}
			&\quad\quad \gamma^+(t,y',0,k)=|k|\widehat{\varphi}(t,y',k),
			\qquad\gamma^-_j(t,y',0,k)=\widetilde{\gamma}_j|k|\widehat{\varphi}(t,y',k), 
		\end{aligned}
	\end{equation*}
where $\widetilde{\gamma}_j$ are two constants satisfying the relation
\begin{equation*} a^-_j\widetilde{\gamma}_j=\frac{(a_1^-)^2+(a_2^-)^2-(b^-)^2}{1-\nu^2\tau^2 }.
		\end{equation*}
\end{Prop}

\begin{proof}[Proof of \cref{prop-h}]
	%The proof is basically algebraic manipulations. It suffices to
	We prove that all the solutions to \cref{f-h} are in the form given in the proposition, since the other side is just a direct verification. We solve the linear PDE \cref{f-h} by decompositing the solution into three parts: the residual components, the surface waves components with zero Fourier modes and the surface waves components with nonzero Fourier modes.
	
	\textbullet~ \textit{The residual components.}
	When $Y_3\rightarrow \pm \infty$, it yields that
	\begin{align*}
		&\sA^+ \partial_{\theta} \underline{U}=0,\quad \xi_j \partial_{\theta} \underline{B}_j =0, &y_3>0,\\
		&\sA^- \partial_{\theta} \underline{V}=0,\quad \xi_j \partial_{\theta} \underline{H}_j =0,\quad \xi_j \partial_{\theta} \underline{E}_j =0, &y_3<0.
	\end{align*}
	where $\sA^+$ are invertible when $\tau\neq0 $, and $\sA^-$ are invertible when $\nu\ll1$ and $\tau\neq0 $.
	Therefore, both $\underline{U}$ and $\underline{V}$ are independent of $\theta$. 
	The constraints of $\underline{U}$ and $\underline{V}$ only come from the boundary conditions, and are determined after the analysis of the surface wave components.

	\textbullet~\textit{The surface wave components with zero Fourier mode.} 
	Projecting the system \cref{f-inh} onto the surface wave components and taking the average with respect to $\theta$,  we have 
	\begin{align*}
		&A^+_3 \partial_{Y_3} \widehat{U}_{\star}(0)=0,\quad \partial_{Y_3} \widehat{B}_{3,\star}(0) =0, &y_3>0,\\
		&A^-_3 \partial_{Y_3} \widehat{V}_{\star}(0)=0,\quad \partial_{Y_3} \widehat{H}_{3,\star}(0) =0,\quad \partial_{Y_3} \widehat{E}_{3,\star}(0) =0, &y_3<0.
	\end{align*}
	Since the surface wave components decay exponentially as $Y_3\rightarrow \pm \infty$, we have that
	\begin{align*}
		& \widehat{u}_{3,\star}(0)=\widehat{B}_{3,\star}(0)=\widehat{q}_{\star}(0)\equiv 0, \quad \widehat{V}_{\star}(0)\equiv0.
	\end{align*}
	Therefore, $ \widehat{U}_{\star}(0)$ is in the form:
	\begin{align*}
		(\widehat{u}_{1,\star}(0),\widehat{u}_{2,\star}(0),0,\widehat{B}_{1,\star}(0),\widehat{B}_{2,\star}(0),0,0)^{\top} (t,y,Y_3),
	\end{align*}
	while $\widehat{V}_{\star}(0)$ is equal to $0$.
	
	With the above expressions of $ \widehat{U}_{\star}(0)$ and $ \widehat{V}_{\star}(0)$, it follows from the Fourier transform of the boundary conditions \cref{f-h} for $k=0$ that 
	\begin{align*}
		&\underline{u}_3 |_{y_3=0}=\underline{B}_3 |_{y_3=0}=\underline{H}_3 |_{y_3=0}=\underline{E}_1 |_{y_3=0}=\underline{E}_2 |_{y_3=0}=0,\\
		& \underline{q} |_{y_3=0}-H^0_j \underline{H}_j |_{y_3=0}+E^0_3 \underline{E}_3 |_{y_3=0}=0.
	\end{align*}
	
	\textbullet~\textit{The surface wave components with non-zero Fourier mode.} 
	Since the residual component does not depend on $\theta$,
	the Fourier transform of the solutions to \cref{f-h} with respect to $\theta$ for $k\neq 0$ yields
	\begin{align}\label{f-h-n0}
		\left\{
		\begin{aligned}
			&A_3^+ \partial_{Y_3} \widehat{U}_{\star}(k)+ik\sA^+ \widehat{U}_{\star}(k) =0, &y_3,Y_3>0,\\
			&\partial_{Y_3}\widehat{B}_{3,\star}(k)+ik\xi_j  \widehat{B}_{j,\star}(k) =0,&y_3,Y_3>0,\\
			&A_3^- \partial_{Y_3} \widehat{V}_{\star}(k)+ik\sA^- \widehat{V}_{\star}(k) =0, &y_3,Y_3<0,\\
			&\partial_{Y_3}\widehat{H}_{3,\star}(k)+ik\xi_j\widehat{H}_{j,\star} (k)=0,&y_3,Y_3<0,\\
			&\partial_{Y_3}\widehat{E}_{3,\star}(k)+ik\xi_j\widehat{E}_{j,\star}(k)=0,&y_3,Y_3<0,\\
			&\mathbb{B}^+ \widehat{U}_{\star}(k)  |_{y_3=Y_3=0}+\mathbb{B}^- \widehat{V}_{\star}(k)  |_{y_3=Y_3=0}+ik\underline{b}\widehat{\varphi}(k)=0.
		\end{aligned}
		\right.
	\end{align}

For simplicity, in this part %we omit ``(k)'' in $\widehat{U}_\star(k), \widehat{V}_\star(k)$ and use the notation $$q^+:=\widehat{q}_{\star},~ q^-:=H^0_j \widehat{H}_{j,\star}-E^0_3 \widehat{E}_{3,\star}.$$
we use the notation $$q^+:=\widehat{q}_{\star}(k),~ q^-:=H^0_j \widehat{H}_{j,\star}(k)-E^0_3 \widehat{E}_{3,\star}(k).$$
Using the divengence-free constraints, one can obtain the equation of $q^+$ and $q^-$ as in \cite{Pierre2021}:
\begin{equation}
	\label{f-h-q}
\left\{
 \begin{aligned}
	&	\partial_{Y_3}^2 q^+ -k^2 q^+ =0, && y_3, Y_3>0,\\
	&	\partial_{Y_3}^2 q^- -k^2(1-\nu^2 \tau^2) q^- =0, && y_3, Y_3<0,\\
	& q^+|_{y_3=Y_3=0}=q^-|_{y_3=Y_3=0}, &&\\
	& \partial_{Y_3}q^+|_{y_3=Y_3=0}= k^2 ((c^+)^2 -(b^+)^2) \widehat{\varphi},&&\\
	& \partial_{Y_3}q^+|_{y_3=Y_3=0}= k^2 ((a^-_1)^2+(a^-_2)^2-(b^-)^2)\widehat{\varphi}.&&
\end{aligned}
\right.
\end{equation}
Therefore, the exponential decaying profiles of $q^\pm$ in $Y_3$ are in the following form: 
\begin{align*}
	&q^+(t,y,Y_3)=((b^+)^2 -(c^+)^2) \gamma^+(t,y,k) e^{-|k|Y_3}, \\
	&q^-(t,y,Y_3)=((b^-)^2 -(a^-_1)^2-(a^-_2)^2)  \gamma^-(t,y,k) e^{|k|\sqrt{1-\nu^2\tau^2}Y_3},
\end{align*}
for any $(t,y,Y_3)\in[0,T]\times \mathbb{T}^2 \times \mathbb{R}^\pm \times \mathbb{R}^\pm$.
Here $\gamma^\pm (t,y,k)$ are coefficients to be determined such that the boundary conditions in \cref{f-h-q} are satisfied:
$$
\left\{
\begin{aligned}
&~~ \quad	((b^+)^2 -(c^+)^2) \gamma^+(t,y',0,k) =((b^-)^2 -(a^-_1)^2-(a^-_2)^2) \gamma^-(t,y',0,k),\\
&	-k((b^+)^2 -(c^+)^2) \gamma^+(t,y',0,k)= k^2 ((c^+)^2 -(b^+)^2) \widehat{\varphi}(t,y',k),\\
&	k\sqrt{1-\nu^2 \tau^2}((b^-)^2-(a^-_1)^2-(a^-_2)^2) \gamma^-(t,y',0,k) \\
&\qquad\qquad\qquad= k^2 ((a^-_1)^2+(a^-_2)^2-(b^-)^2)\widehat{\varphi}(t,y',k).
\end{aligned}
\right.
$$
Therefore, the boundary conditions are satisfied if and only if 
\begin{equation}\label{Lopatinskii}
	\Delta(\tau,\xi)=\sqrt{1-\nu^2 \tau^2}((c^+)^2-(b^+)^2) + ((a^-_1)^2+(a^-_2)^2-(b^-)^2)=0.
\end{equation}
The left-hand side is exactly the Lopatinskii determinant (\cite{Morando2020a, Trakhinin2020b})
and the equality to zero is fulfilled by our assumption \cref{H-key},
As mentioned in the introduction, under \cref{H1*} there are only simple real roots $\tau=\tau(\xi)$ to $\Delta(\tau,\xi)=0$.
We refer to \cite{Pierre2021} and the references therein for more discussions about connections between these linear surface waves and the normal mode analysis.

Once the total pressures are settled down, the exponential decaying profiles of $\widehat{U}_\star(k), \widehat{V}_\star(k)$ in $Y_3$ can be obtained in the form: 

\begin{align*}
	&\widehat{U}_\star(t,y,Y_3, k)=\gamma^+ (t,y,k)e^{-|k|Y_3} \sR^+(k), \\
	&\widehat{V}_{\star}(t,y,Y_3, k)=\gamma^-_j (t,y,k)e^{|k|\sqrt{1 -\nu^2\tau^2}Y_3} \sR^-_j(k)
\end{align*}
for any $(t,y,Y_3)\in[0,T]\times \mathbb{T}^2 \times \mathbb{R}^\pm \times \mathbb{R}^\pm$, where $\sR^+(k),\sR^-_j(k)$ satisfy
$$
(-A_3^+ +i \sgnk \sA^+ ) \sR^+=0,\qquad
(\sqrt{1-\nu^2\tau^2}A_3^- +i \sgnk \sA^- ) \sR^-=0,
$$
and are given explicity in \cref{vector-R}.
The boundary conditions of \cref{f-h-q} determine 
$$
\gamma^+ (t,y',0,k)=|k|\widehat{\varphi}(t,y',k), \qquad 
\gamma^-_j (t,y',0,k)=\widetilde{\gamma}_j|k|\widehat{\varphi}(t,y',k), 
$$
where $\widetilde{\gamma}_j$ are two constants satisfying the relation $a^-_j\widetilde{\gamma}_j=\frac{(a_1^-)^2+(a_2^-)^2-(b^-)^2}{1-\nu^2\tau^2 }$.
Therefore, the proof of \cref{prop-h} is completed.
\end{proof}

\subsection{The inhomogeneous case}
Since the fast problem \cref{f} is linear, 
after subtracting a solution to the homogeneous problem \cref{f-h} as in \cref{prop-h} ($\varphi$ in \cref{f-h} should be the same as in \cref{f}), 
it suffices to show that \cref{f} admits one solution of the form $(\mathbb{U},\mathbb{V},0)$ with $\mathbb{U} \in S^+,\mathbb{V} \in S^-$, i.e.
\begin{align}\label{f-inh}
\left\{
\begin{aligned}
&\sL^+_f \mathbb{U} =F^+, &y_3,Y_3>0,\\
&\partial_{Y_3}\sB_3+\xi_j \partial_{\theta} \sB_j =F^+_8,&y_3,Y_3>0,\\
&\sL^-_f \mathbb{V} =F^-,&y_3,Y_3<0,\\
&\partial_{Y_3}\sH_3+\xi_j \partial_{\theta}\sH_j =F^-_7,&y_3,Y_3<0,\\
&\partial_{Y_3}\sE_3+\xi_j \partial_{\theta}\sE_j =F^-_8,&y_3,Y_3<0,\\
&\mathbb{B}^+ \mathbb{U}  |_{y_3=Y_3=0}+\mathbb{B}^- \mathbb{V}  |_{y_3=Y_3=0}=G.
\end{aligned}
\right.
\end{align}
Here we denote $\mathbb{U}=\{\sU_1,\sU_2,\sU_3,\sB_1,\sB_2,\sB_3,\sQ\}^\top$,
$\mathbb{V}=\{\sH_1,\sH_2,\sH_3,\sE_1,\sE_2,\sE_3\}^\top$.
We will prove that \eqref{f-comp} are the necessary and sufficient conditions on the source terms for the existence of a solution to \cref{f-inh}.
 
\subsubsection{Necessary conditions}
Suppose that $(\mathbb{U},\mathbb{V},0)$ is a solution to \cref{f-inh}.
By the boundary condition, the values of $\sU_3, \sB_3, \sH_3, \sE_1,\sE_2$ on the boundary can be expressed by $G$. Substituting these expressions into the 6th equation in the plasma side and 3rd equation in the vacuum side yields \eqref{comp-a}.

 Letting $Y_3\rightarrow \pm \infty$ in the system \cref{f-inh} and considering the average with respect to $\theta$, we easily obtain \eqref{comp-b}.

% Projecting the system \cref{f-inh} onto the surface wave components and taking the average with respect to $\theta$,  we get 
%	\begin{align*}
%	%line 1
%	&\partial_{Y_3} \widehat{\sU}_{3,\star}(0)=\widehat{F}^+_{7,\star}(0),
%	\quad\partial_{Y_3} \widehat{\sB}_{3,\star}(0)=\widehat{F}^+_{8,\star}(0),
%	\quad\partial_{Y_3} \widehat{\sQ}_{\star}(0)=\widehat{F}^+_{3,\star}(0),\\
%	%\line 2
%	&\partial_{Y_3} \widehat{\sH}_{1,\star}(0)=-\widehat{F}^-_{5,\star}(0),
%	\quad\partial_{Y_3} \widehat{\sH}_{2,\star}(0)=\widehat{F}^-_{4,\star}(0),
%	\quad\partial_{Y_3} \widehat{\sH}_{3,\star}(0)=\widehat{F}^-_{7,\star}(0),\\
%	%\line 3
%	&\partial_{Y_3} \widehat{\sE}_{1,\star}(0)=\widehat{F}^-_{2,\star}(0),
%	\quad\partial_{Y_3} \widehat{\sE}_{2,\star}(0)=-\widehat{F}^-_{1,\star}(0),
%	\quad\partial_{Y_3} \widehat{\sE}_{3,\star}(0)=\widehat{F}^-_{8,\star}(0).
%	\end{align*}
Projecting the equations on the plasma side of system \cref{f-inh} onto the surface wave components and taking the average with respect to $\theta$, we obtain
\begin{align}\label{surf}
		\quad A_3^+\partial_{Y_3} \widehat{\mathbb{U}}_{\star}(0)=\widehat{F}^+_{\star}(0),
	\quad\partial_{Y_3} \widehat{\sB}_{3,\star}(0)=\widehat{F}^+_{8,\star}(0).
\end{align}
In particular, from the 7th equation of the vector equality in \cref{surf} we get
\begin{align}\label{nd}
	\partial_{Y_3} \widehat{\sU}_{3,\star}(0)=\widehat{F}^+_{7,\star}(0),
	\quad\partial_{Y_3} \widehat{\sB}_{3,\star}(0)=\widehat{F}^+_{8,\star}(0),
\end{align}
and from the equations for the tangential components we have
\begin{equation}
	\left\{
\begin{aligned}
	&u_j^0 \partial_{Y_3} \widehat{\sU}_{3,\star}(0) -B_j^0 \partial_{Y_3} \widehat{\sB}_{3,\star}(0)=\widehat{F}^+_{j,\star} (0),\quad j=1,2,\\
	&B_j^0 \partial_{Y_3} \widehat{\sU}_{3,\star}(0) -u_j^0 \partial_{Y_3} \widehat{\sB}_{3,\star}(0)=\widehat{F}^+_{3+j,\star} (0),\quad j=1,2,\\
	%	&\widehat{F}^-_{3,\star} (0)=0.
\end{aligned}
\right.\label{comp-tang}
\end{equation}
Substituting \cref{nd} into \cref{comp-tang} yields \eqref{comp-c}.
	Notice that the 6th equation of the vector equality in \cref{surf} necessarily gives $\widehat{F}^+_{6,\star} (0)=0$. However, as in \cite{Pierre2021}, we don't include this equation among the compatibility conditions \eqref{f-comp} because it can be derived from the next compatibility condition \eqref{comp-d}. Similarly, in the vacuum side the equality $A_3^-\partial_{Y_3} \widehat{\mathbb{V}}_{\star}(0)=\widehat{F}^-_{\star}(0)$ necessarily yields  $\widehat{F}^-_{3,\star} (0)=\widehat{F}^-_{6,\star} (0)=0$ that we do not include in \eqref{f-comp} because it can also be derived from \eqref{comp-d}.
	
Applying $\partial_{Y_3}$ to the 3rd equation in vacuum side of \cref{f-inh}, $\xi_1\partial_\theta$ to the 1st equation,  $\xi_2\partial_\theta$ to the 2nd equation, and summing up yields 
	$$
	\nu\tau \partial_{\theta} \left( \partial_{Y_3} \sH_3 + \xi_j \partial_{\theta} \sH_j \right) =\partial_{Y_3} F^-_3+\xi_j \partial_{\theta}F^-_j.
	$$
	Similar formulas hold for $\sB, \sE$. Substituting the divergence constraints in \cref{f-inh}, these relations give the compatibility condition \eqref{comp-d}. 
	We refer to \cite{Pierre2021} for more discussions on this compatibility condition.

We now derive the last compatibility condition \cref{comp-e}. Considering the nonzero Fourier modes of \cref{f-inh},  with the exclusion of the three divergence-free conditions, letting
	$y_3=0$ and regarding $t,y',k$ as parameters, we have the following ODE equations with respect to $Y_3$:
	\begin{align}\label{f-comp-ode}
	\left\{
	\begin{aligned}
	&A_3^+ \partial_{Y_3} \widehat{\mathbb{U}}(k) |_{y_3=0}+ik\sA^+ \widehat{\mathbb{U}}(k) |_{y_3=0} =\widehat{F}^+(k) |_{y_3=0}, &Y_3>0,\\
	&A_3^- \partial_{Y_3} \widehat{\mathbb{V}}(k) |_{y_3=0}+ik\sA^- \widehat{\mathbb{V}}(k) |_{y_3=0} =\widehat{F}^-(k) |_{y_3=0}, &Y_3<0,\\
	&\mathbb{B}^+ \widehat{\mathbb{U}}(k)  |_{y_3=Y_3=0}
	+\mathbb{B}^- \widehat{\mathbb{V}}(k)  |_{y_3=Y_3=0}
	=\widehat{G}(k).
	\end{aligned}
	\right.
	\end{align}
Let us introduce the operators $L^{\pm}_k:=A^{\pm} \partial_{Y_3} +ik\sA^{\pm}$ and their adjoint operators $(L^{\pm}_k)^*=-(A^{\pm})^\top \partial_{Y_3} +ik(\sA^{\pm})^\top$. For any smooth test function $W^{\pm}$ with exponential decay at $Y_3=\pm \infty$, from \eqref{f-comp-ode} it follows that:
	\begin{align}
	%line 1
	&\int_{0}^{+\infty} W^+ \cdot (L^+_k \widehat{\mathbb{U}}(k) |_{y_3=0} ) dY_3-\int_{-\infty}^{0} W^- \cdot (L^-_k \widehat{\mathbb{V}}(k) |_{y_3=0} ) dY_3 \notag \\
	%	 line 2
	%	 &\quad = \int_{0}^{+\infty} (L^+_k)*W^+ \cdot  \widehat{\mathbb{U}}(k) |_{y_3=0} dY_3
	%	 -\int_{-\infty}^{0} (L^-_k)*W^- \cdot \widehat{\mathbb{V}}(k) |_{y_3=0}  dY_3 
	%	 - W^+(0) \cdot A^+_3 \widehat{\mathbb{U}}(k) |_{y_3=Y_3=0}
	%	 - W^-(0) \cdot A^-_3 \widehat{\mathbb{V}}(k) |_{y_3=Y_3=0}\\
	%line 2,3,4
	& = \uwave{\int_{0}^{+\infty} (L^+_k)^*W^+ \cdot  \widehat{\mathbb{U}}(k) |_{y_3=0} dY_3
	-\int_{-\infty}^{0} (L^-_k)^*W^- \cdot \widehat{\mathbb{V}}(k) |_{y_3=0}  dY_3 }
	\notag\\
	& \uwave{+\left( \nu\tau \overline{W^-_5(0)}-a^-_1 \overline{W^+_3(0)}\right)\frac{1}{\nu\tau} \widehat{\sH}_1(k) |_{y_3=Y_3=0}}  \notag\\
	& \uwave{+\left( -\nu\tau \overline{W^-_4(0)}-a^-_2\overline{W^+_3(0)}\right)\frac{1}{\nu\tau} \widehat{\sH}_2(k) |_{y_3=Y_3=0}} \label{eq-test}\\
	%line 5,6,7
	& - \left(u^0_j \overline{W^+_j(0)} +B^0_j \overline{W^+_{3+j}(0)} + \overline{W^+_7(0)} \right) \widehat{G}_1(k)\notag\\
	& + \left(B^0_j \overline{W^+_j(0)} +u^0_j \overline{W^+_{3+j}(0)} \right) \widehat{G}_2(k)\notag\\
	& + \overline{W^-_j(0)}\widehat{G}_{3+j}(k)  - \overline{W^+_3(0)}\widehat{G}_6(k) + \frac{E_3^0}{ik\nu\tau}\overline{W^+_3(0)} \widehat{F}_6^- (k) |_{y_3=Y_3=0}	 \,.
	\end{align}
	Here we  use  ``$\cdot$'' to denote the inner product of vectors of complex numbers : $X\cdot Y = \sum_{\alpha} \overline{X}_\alpha Y_\alpha$. 
	Moreover, we also use the 6th equations of \eqref{f-comp-ode} in the vacuum side:  
	$$ ik\xi_2 \widehat{\sH}_1(k)-ik\xi_1 \widehat{\sH}_2(k)+ik\nu\tau\widehat{\sE}_3(k) =  \widehat{F}^-_6(k) \quad y_3,Y_3\leq 0.$$
	To choose test functions $W^\pm$ such that underwave terms in \cref{eq-test} vanish, we consider the {\it dual} problem
	\begin{align}\label{test}
	\left\{
	\begin{aligned}
	&(L^{\pm}_k)^*W^{\pm} =0, \quad \pm Y_3 >0,\\
	&  \nu\tau W^-_5(0)-a^-_1W^+_3(0) =0\\
	&-\nu\tau W^-_4(0) -a^-_2 W^+_3(0)=0.
	\end{aligned}
	\right.
	\end{align}
    After calculations similar to those in the proof of Proposition \ref{prop-h} we find 
	\begin{align}
	W^+=e^{-|k| Y_3} \sL^+(k), W^-=e^{|k|\sqrt{1-\nu^2\tau^2} Y_3} \sL^-(k),
	\end{align}
	where the explicit definitions of $\sL^\pm(k)$ are given in \cref{vector-L}. 
	Then \eqref{comp-e} can be derived by defining \begin{align}\label{l+}
	l^+_1:=2(b^+)^2-\tau c^+,\quad l^+_2:= -(a^+ +c^+)b^+.
	\end{align}

\subsubsection{Solving \cref{f-inh} on zero Fourier modes}
To solve \cref{f-inh} we firstly consider its zero Fourier modes, namely we construct functions $\widehat{\mathbb{U}}(0) \in S^+,\widehat{\mathbb{V}}(0) \in S^-$, independent of $\theta$, and that satisfy 
\begin{align}\label{f-inh-0}
	\left\{
	\begin{aligned}
		&A_3^+ \partial_{Y_3} \widehat{\mathbb{U}}(0) =\widehat{F}^+(0),\quad
		\partial_{Y_3}\widehat{\sB}_3(0) =\widehat{F}^+_8(0),&y_3,Y_3>0,\\
		&A_3^- \partial_{Y_3} \widehat{\mathbb{V}}(0) =\widehat{F}^-(0), &y_3,Y_3<0,\\
		&\partial_{Y_3}\widehat{\sH}_3(0) (k)=\widehat{F}^-_7(0),\quad
		\partial_{Y_3}\widehat{\sE}_3(0)=\widehat{F}^-_8(0),\quad&y_3,Y_3<0,\\
		&\mathbb{B}^+ \widehat{\mathbb{U}}(0)  |_{y_3=Y_3=0}+\mathbb{B}^- \widehat{\mathbb{V}}(0)  |_{y_3=Y_3=0}=\widehat{G}(0).
	\end{aligned}
	\right.
\end{align}
Using \cref{comp-b} and integrating in $Y_3$ with zero limit at $Y_3=\pm\infty$, from \cref{f-inh-0} we obtain
\begin{equation}\label{def-0-Ustar}
\begin{aligned}
%line 1+2
&\widehat{\sU}_{3,\star}(t,y,Y_3,0)=-\int_{Y_3}^{+\infty} \widehat{F}^+_{7,\star}(t,y,Y,0) dY, 
&& \widehat{\sB}_{3,\star}(t,y,Y_3,0)=-\int_{Y_3}^{+\infty} \widehat{F}^+_{8,\star}(t,y,Y,0) dY,\\
&\widehat{\sQ}_{\star}(t,y,Y_3,0)=-\int_{Y_3}^{+\infty} \widehat{F}^+_{3,\star}(t,y,Y,0) dY,\\
%line 3+4
&\widehat{\sH}_{1,\star}(t,y,Y_3,0)=-\int_{-\infty}^{Y_3} \widehat{F}^-_{5,\star}(t,y,Y,0) dY, 
&& \widehat{\sH}_{2,\star}(t,y,Y_3,0)=\int_{-\infty}^{Y_3} \widehat{F}^-_{4,\star}(t,y,Y,0) dY,\\
& \widehat{\sH}_{3,\star}(t,y,Y_3,0)=\int_{-\infty}^{Y_3} \widehat{F}^-_{7,\star}(t,y,Y,0) dY,\\
%line 5+6
&\widehat{\sE}_{1,\star}(t,y,Y_3,0)=\int_{-\infty}^{Y_3} \widehat{F}^-_{2,\star}(t,y,Y,0) dY, 
&& \widehat{\sE}_{2,\star}(t,y,Y_3,0)=-\int_{-\infty}^{Y_3} \widehat{F}^-_{1,\star}(t,y,Y,0) dY,\\
& \widehat{\sE}_{3,\star}(t,y,Y_3,0)=\int_{-\infty}^{Y_3} \widehat{F}^-_{8,\star}(t,y,Y,0) dY.
\end{aligned}
\end{equation}
We remark that the above functions do not depend on $\theta$.
By \eqref{comp-c} and \eqref{comp-d} (which implies $\widehat{F}^\pm_{6,\star} (0)=\widehat{F}^-_{3,\star} (0)=0$), it is easy to verify that 
\begin{align*}
&\widehat{\mathbb{U}}_{\star}(0)=(0,0,\widehat{\sU}_{3,\star}(0),0,0,\widehat{\sB}_{3,\star}(0),\widehat{\sQ}_{\star}(0))^\top, \\
&
\widehat{\mathbb{V}}_{\star}(0)=(\widehat{\sH}_{1,\star}(0),\widehat{\sH}_{2,\star}(0),\widehat{\sH}_{3,\star}(0),\widehat{\sE}_{1,\star}(0),\widehat{\sE}_{2,\star}(0),\widehat{\sE}_{3,\star}(0))^\top
\end{align*}
satisfies \cref{f-inh-0} projected on $S^{\pm}_{\star}$. It remains to add suitable slow functions in order to satisfy the boundary conditions in \cref{f-inh-0}. Indeed, adding any function independent of $Y_3$ will not modify the fullfillment of the fast equations in \cref{f-inh-0}.
The residual components $\underline{\widehat{\mathbb{U}}}(0), \underline{\widehat{\mathbb{V}}}(0)$ are therefore defined as:
\begin{align*}
&\underline{\widehat{\sU}}_3 (t,y',0,0)=\underline{\widehat{G}}_1(t,y',0)-\widehat{\sU}_{3,\star}(t,y',0,0,0),\\
&\underline{\widehat{\sB}}_3 (t,y',0,0)=\underline{\widehat{G}}_2(t,y',0)-\widehat{\sB}_{3,\star}(t,y',0,0,0),\\
&\underline{\widehat{\sH}}_3 (t,y',0,0)=\underline{\widehat{G}}_3(t,y',0)-\widehat{\sH}_{3,\star}(t,y',0,0,0),\\
&\underline{\widehat{\sE}}_1 (t,y',0,0)=-\underline{\widehat{G}}_5(t,y',0)-\widehat{\sE}_{1,\star}(t,y',0,0,0),\\
&\underline{\widehat{\sE}}_2 (t,y',0,0)=\underline{\widehat{G}}_4(t,y',0)-\widehat{\sE}_{2,\star}(t,y',0,0,0),\\
&(\underline{\widehat{\sQ}}-H^0_j\underline{\widehat{\sH}}_j+E^0_3\underline{\widehat{\sE}}_3) (t,y',0,0)=\underline{\widehat{G}}_6(t,y',0)-(\widehat{\sQ}_{\star}-H^0_j\widehat{\sH}_{j,\star}+E^0_3\widehat{\sE}_{3,\star})(t,y',0,0,0).
\end{align*}
This completes the construction of the particular solution $\widehat{\mathbb{U}}(0) \in S^+,\widehat{\mathbb{V}}(0) \in S^-$; the above construction shows that it is always possible to choose $\widehat{\mathbb{U}}(0)$ such that $\Pi\,\widehat{\mathbb{U}}(0)=0$ in the inhomogeneous case.

\subsubsection{Solving \cref{f-inh} on non-zero Fourier modes}
Now we consider the non-zero Fourier modes of \cref{f-inh}.
The $k$-Fourier coefficient of \cref{f-inh} for $k\neq 0$ is 
\begin{align}\label{f-inh-n0}
\left\{
\begin{aligned}
&A_3^+ \partial_{Y_3} \widehat{\mathbb{U}}(k)+ik\sA^+ \widehat{\mathbb{U}}(k) =\widehat{F}^+(k), &y_3,Y_3>0,\\
&\partial_{Y_3}\widehat{\sB}_3(k)+ik\xi_j  \widehat{\sB}_j(k) =\widehat{F}^+_8(k),&y_3,Y_3>0,\\
&A_3^- \partial_{Y_3} \widehat{\mathbb{V}}(k)+ik\sA^- \widehat{\mathbb{V}}(k) =\widehat{F}^-(k), &y_3,Y_3<0,\\
&\partial_{Y_3}\widehat{\sH}_3(k)+ik\xi_j\widehat{\sH}_j (k)=\widehat{F}^-_7(k),&y_3,Y_3<0,\\
&\partial_{Y_3}\widehat{\sE}_3(k)+ik\xi_j\widehat{\sE}_j(k)=\widehat{F}^-_8(k),&y_3,Y_3<0,\\
&\mathbb{B}^+ \widehat{\mathbb{U}}(k)  |_{y_3=Y_3=0}+\mathbb{B}^- \widehat{\mathbb{V}}(k)  |_{y_3=Y_3=0}=\widehat{G}(k).
\end{aligned}
\right.
\end{align}
We look for a solution to the above system under the compatibility conditions \eqref{comp-a}, \eqref{comp-d} and \eqref{comp-e} (in fact \eqref{comp-b} and \eqref{comp-c} are just conditions on zero Fourier modes).

Firstly, by taking the limit as $Y_3\rightarrow\pm\infty$, we obtain that the residual components are necessarily determined by
$$
\underline{\widehat{\mathbb{U}}}(k):=\frac{1}{ik} (\sA^+)^{-1} \underline{\widehat{F}}^+(k),\quad \underline{\widehat{\mathbb{V}}}(k):=\frac{1}{ik} (\sA^-)^{-1} \underline{\widehat{F}}^-(k).
$$
The three divergence-free conditions can be verified by the compatibility condition \eqref{comp-d} as in \cite{Pierre2021}.
Given the above functions $(\underline{\widehat{\mathbb{U}}},\underline{\widehat{\mathbb{V}}})$, we define
$$
\underline{\widehat{G}}(k):=\mathbb{B}^+\underline{\widehat{\mathbb{U}}}(k) |_{y_3=0}+\mathbb{B}^-\underline{\widehat{\mathbb{V}}}(k) |_{y_3=0}, \qquad \widehat{G}_{\star}(k):= \widehat{G}(k)-\underline{\widehat{G}}(k).
$$ 
Then $(\underline{\widehat{\mathbb{U}}},\underline{\widehat{\mathbb{V}}})$ solves a problem which is the limit of \cref{f-inh-n0} as $Y_3\rightarrow\pm\infty$,
%, and the compatibility conditions \eqref{comp-a-n0*} 
%and \eqref{comp-e-n0*}, 
with the source terms  $\underline{\widehat{F}}^{\pm}(k),\underline{\widehat{F}}^{-}_7(k),\underline{\widehat{F}}^{\pm}_8(k),\underline{\widehat{G}}(k)$. Such problem is the residual component of \cref{f-inh-n0}, except for the different boundary source term.

\vspace{0.3cm}

Therefore, by subtracting the residual components from \cref{f-inh-n0}, we look for surface wave components that are solutions to
\begin{align}\label{f-inh-n0*}
\left\{
\begin{aligned}
&A_3^+ \partial_{Y_3} \widehat{\mathbb{U}}_{\star}(k)+ik\sA^+ \widehat{\mathbb{U}}_{\star}(k) =\widehat{F}^+_{\star}(k), &y_3,Y_3>0,\\
&\partial_{Y_3}\widehat{\sB}_{3,\star}(k)+ik\xi_j  \widehat{\sB}_{j,\star}(k) =\widehat{F}^+_{8,\star}(k),&y_3,Y_3>0,\\
&A_3^- \partial_{Y_3} \widehat{\mathbb{V}}_{\star}(k)+ik\sA^- \widehat{\mathbb{V}}_{\star}(k) =\widehat{F}^-_{\star}(k), &y_3,Y_3<0,\\
&\partial_{Y_3}\widehat{\sH}_{3,\star}(k)+ik\xi_j\widehat{\sH}_{j,\star} (k)=\widehat{F}^-_{7,\star}(k),&y_3,Y_3<0,\\
&\partial_{Y_3}\widehat{\sE}_{3,\star}(k)+ik\xi_j\widehat{\sE}_{j,\star}(k)=\widehat{F}^-_{8,\star}(k),&y_3,Y_3<0,\\
&\mathbb{B}^+ \widehat{\mathbb{U}}_{\star}(k)  |_{y_3=Y_3=0}+\mathbb{B}^- \widehat{\mathbb{V}}_{\star}(k)  |_{y_3=Y_3=0}=\widehat{G}_{\star}(k).
\end{aligned}
\right.
\end{align}
By linearity from \eqref{comp-a}, \eqref{comp-e}, and projection on the surface waves components of  \eqref{comp-d}, the source terms of \cref{f-inh-n0*} satisfy the following compatibility conditions: 
\begin{align} \label{comp-a-n0*}
\left\{
\begin{aligned}
&\widehat{F}^+_{6,\star} |_{y_3=Y_3=0}=-ikb^+ \widehat{G}_{1,\star}+ikc^+\widehat{G}_{2,\star},\\
&\widehat{F}^-_{3,\star} |_{y_3=Y_3=0}=ik\nu\tau \widehat{G}_{3,\star}+ik\xi_j\widehat{G}_{3+j,\star}.
\end{aligned}
\right.
\end{align}
\begin{align}
\label{comp-d-n0*}
\left\{
\begin{aligned}
&\partial_{Y_3} \widehat{F}_{6,\star}^+ + ik\xi_j  \widehat{F}^+_{3+j,\star}-ik\tau \widehat{F}^+_{8,\star}=0,\\
&\partial_{Y_3} \widehat{F}_{3,\star}^- + ik\xi_j  \widehat{F}^-_{j,\star}-ik\nu\tau  \widehat{F}^-_{7,\star}=0,\\
&\partial_{Y_3} \widehat{F}_{6,\star}^- + ik\xi_j  \widehat{F}^-_{3+j,\star}-ik\nu\tau  \widehat{F}^-_{8,\star}=0.
\end{aligned}
\right.
\end{align}
\begin{align}
\label{comp-e-n0*}
\begin{aligned}
&\quad \int_{0}^{+\infty} e^{-|k|Y_3} \sL^+(k) \cdot \widehat{F}^+_{\star}(t,y',0,Y_3,k) dY_3
+ l^+_j \widehat{G}_{j,\star}(t,y',k)
 \\
&-\int_{-\infty}^{0} e^{|k|\sqrt{1-\nu^2\tau^2}Y_3} \sL^-(k) \cdot \widehat{F}_{\star}^-(t,y',0,Y_3,k)  dY_3 
-\frac{-\nu\tau a^-_j +\xi_j b^-}{\nu\sqrt{1-\nu^2\tau^2}}  \widehat{G}_{3+j,\star}(t,y',k)\\
&  -i\sgnk \tau \widehat{G}_{6,\star}(t,y',k) +\frac{E^0_3}{|k|\nu} \widehat{F}^-_{6,\star}(t,y',0,0,k)=0, \quad \text{for~} k\in \mathbb{Z}\setminus\{0\}.
\end{aligned}
\end{align}
In the following we will solve \cref{f-inh-n0*} in the plasma and the vacuum separately: we will consider the  two equations in $\{y_3, Y_3>0\}$ and the first two boundary conditions for the plasma, and consider the three equations in $\{y_3, Y_3<0\}$ and the 3rd - 5th boundary condition for the vacuum.  
Eventually, the 6th boundary condition (on the pressures) will be verified by the compatibility conditions.  
Solving \cref{f-inh-n0*} in the plasma region is the same as in \cite{Pierre2021}, and thus we omit it here.  
Solving \cref{f-inh-n0*} in the vacuum region consists of two steps: constructing the solutions to some of the equations and boundary conditions, and then verifying the remaining equations by the compatibility conditions.

\vspace{0.3cm}

\noindent \textit{1) Solving \cref{f-inh-n0*} in vacuum.}

In the following we will omit ``$(k)$'' in $\widehat{\mathbb{V}}_{\star}(k)$.
It follows from \cref{f-inh-n0*} that 
\begin{align}
\partial_{Y_3}^2 \widehat{\sE}_{1,\star}=
&-ik\nu\tau \partial_{Y_3} \widehat{\sH}_{2,\star} + ik\xi_1 \partial_{Y_3} \widehat{\sE}_{3,\star} + \partial_{Y_3} \widehat{F}_{2,\star}^- \notag  \\
=&-ik\nu \tau (-ik\nu\tau  \widehat{\sE}_{1,\star} + ik\xi_2  \widehat{\sH}_{3,\star} +  \widehat{F}_{4,\star}^- )  \notag\\
&\qquad +ik\xi_1 (-ik\xi_1  \widehat{\sE}_{1,\star} - ik\xi_2  \widehat{\sE}_{2,\star} +  \widehat{F}_{8,\star}^-)+\partial_{Y_3} \widehat{F}_{2,\star}^-   \notag\\
=&k^2(-\nu^2\tau^2+\xi_1^2) \widehat{\sE}_{1,\star} -ik\xi_2 (ik\nu\tau \widehat{\sH}_{3,\star} +ik \xi_1 \widehat{\sE}_{2,\star})  \notag\\
&\qquad +\partial_{Y_3}\widehat{F}_{2,\star}^-+ ik(-\nu\tau \widehat{F}_{4,\star}^- +\xi_1 \widehat{F}_{8,\star}^- ) \notag\\
=&k^2(1-\nu^2\tau^2) \widehat{\sE}_{1,\star} +\partial_{Y_3}\widehat{F}_{2,\star}^-+ ik(-\nu\tau \widehat{F}_{4,\star}^- +\xi_1 \widehat{F}_{8,\star}^- -\xi_2\widehat{F}_{3,\star}^-).\label{ode-deri}
\end{align}
Therefore, we can conclude that $\widehat{\sE}_{1,\star}^-$ satisfies a second ODE with respect to $Y_3$. With similar calculations for the other variables in the vacuum we obtain:
\begin{align}\label{ode-V*}
\partial_{Y_3}^2 \widehat{\mathbb{V}}_{\star}-k^2(1-\nu^2\tau^2)\widehat{\mathbb{V}}_{\star}
=\underbrace{\partial_{Y_3} \begin{pmatrix}
	-\widehat{F}^-_{5,\star}\\
	\widehat{F}^-_{4,\star}\\
	\widehat{F}^-_{7,\star}\\
	\widehat{F}^-_{2,\star}\\
	-\widehat{F}^-_{1,\star}\\
	\widehat{F}^-_{8,\star}
	\end{pmatrix}
	+ik\begin{pmatrix}
	-\nu \tau \widehat{F}^-_{1,\star}+\xi_1 \widehat{F}^-_{7,\star} +\xi_2 \widehat{F}^-_{6,\star}\\
	-\nu \tau \widehat{F}^-_{2,\star}+\xi_2 \widehat{F}^-_{7,\star} -\xi_1 \widehat{F}^-_{6,\star}\\
	-\nu \tau \widehat{F}^-_{3,\star}+\xi_1 \widehat{F}^-_{5,\star} -\xi_2 \widehat{F}^-_{4,\star}\\
	-\nu \tau \widehat{F}^-_{4,\star}+\xi_1 \widehat{F}^-_{8,\star} -\xi_2 \widehat{F}^-_{3,\star}\\
	-\nu \tau \widehat{F}^-_{5,\star}+\xi_2 \widehat{F}^-_{8,\star} +\xi_1 \widehat{F}^-_{3,\star}\\
	-\nu \tau \widehat{F}^-_{6,\star}-\xi_1 \widehat{F}^-_{2,\star} +\xi_2 \widehat{F}^-_{1,\star}
	\end{pmatrix} }_{=:\mathscr{F}^-}.
\end{align}
First, we construct $\widehat{\sE}_{1,\star}, \widehat{\sE}_{2,\star}$ by {the 4th and 5th ODEs in \cref{ode-V*} and the 4th and 5th boundary conditions in \cref{f-inh-n0*}:}
\begin{align}
&\widehat{\sE}_{j,\star}:=\kappa^-_{3+j} e^{|k|\sqrt{1-\nu^2\tau^2}Y_3}-\frac{1}{2|k|\sqrt{1-\nu^2\tau^2}} \Big( \int_{-\infty}^{Y_3} e^{-|k|\sqrt{1-\nu^2\tau^2}(Y_3-Y)} \mathscr{F}^-_{3+j}(Y) dY  \notag\\
&\quad\quad +\int_{Y_3}^{0} e^{-|k|\sqrt{1-\nu^2\tau^2}(Y-Y_3)}  \mathscr{F}^-_{3+j}(Y) dY \Big), \label{def-E*}\\
%kappa
&\kappa^-_{3+j} |_{y_3=0}=\frac{1}{2|k|\sqrt{1-\nu^2\tau^2}}  \int_{-\infty}^{0} e^{|k|\sqrt{1-\nu^2\tau^2}Y} \mathscr{F}^-_{3+j} |_{y_3=0}(Y) dY+ \widehat{\sE}_{j,\star} |_{y_3=Y_3=0}\,, \notag\\
%G
&\widehat{\sE}_{1,\star} |_{y_3=Y_3=0}= -\widehat{G}_{5,\star}, \quad  \widehat{\sE}_{2,\star} |_{y_3=Y_3=0}= \widehat{G}_{4,\star}\,, \notag
\end{align}
where $\kappa^-_4$ and $\kappa^-_5$ are constants w.r.t $Y_3$ (functions of $(t,y,k)$,  each one having a constraint at $y_3=0$). 
 By \cref{def-E*} and {the 3rd equation of \cref{f-inh-n0*} in the vacuum}, one can obtain $\widehat{\sH}_{3,\star}$.
It follows from {the 1st, 2nd and 6th equations of \cref{f-inh-n0*} in the vacuum} that 
\begin{align*}
\begin{pmatrix}
ik\nu\tau &0  & ik\xi_2 \\
0 &ik\nu\tau & -ik\xi_1 \\
ik\xi_2 & -ik \xi_1 & ik \nu\tau
\end{pmatrix} 
\begin{pmatrix}
\widehat{\sH}_{1,\star}\\
\widehat{\sH}_{2,\star}\\
\widehat{\sE}_{3,\star}
\end{pmatrix}=
\begin{pmatrix}
\widehat{F}^-_{1,\star}+\partial_{Y_3} \widehat{\sE}^-_{2,\star}\\
\widehat{F}^-_{2,\star}-\partial_{Y_3} \widehat{\sE}^-_{1,\star}\\
\widehat{F}^-_{6,\star}
\end{pmatrix}.
\end{align*}
Since the above matrix is invertible, $\widehat{\sH}_{1,\star}, \widehat{\sH}_{2,\star}, \widehat{\sE}_{3,\star}$ are determined. It remains to verify that the solutions constructed above satisfy the other equations and boundary conditions of \cref{f-inh-n0*} in the vacuum part.

\vspace{0.5cm}
{The 3rd boundary condition} can be proved by \eqref{comp-a-n0*} directly, while
{the 7th equation of \cref{f-inh-n0*} in the vacuum} can be easily verified by the 1st, 2nd, 3rd ones and \eqref{comp-d-n0*}.

By substituting $\widehat{\sH}_{3,\star}$, obtained from the 3rd equation in the vacuum part,  {the 4th equation of \cref{f-inh-n0*} in the vacuum} reduces to
\begin{align}\label{4th}
\partial_{Y_3} \widehat{\sH}_{2,\star}
+\frac{ik}{\nu\tau} (\nu^2\tau^2-\xi_2^2) \widehat{\sE}_{1,\star}
+\frac{ik}{\nu\tau} \xi_1\xi_2 \widehat{\sE}_{2,\star} =\widehat{F}^-_{4,\star} +\frac{\xi_2}{\nu\tau} \widehat{F}^-_{3,\star}.
\end{align}
Analogously, by substituting $\widehat{\sE}_{3,\star}$ obtained from the 8th equation of \cref{f-inh-n0*} in the vacuum part, the 1st and 2nd ones reduce to
\begin{align*}
-\partial_{Y_3} \widehat{\sE}_{2,\star}
+\frac{ik}{\nu\tau} (\nu^2\tau^2-\xi_2^2) \widehat{\sH}_{1,\star}
+\frac{ik}{\nu\tau} \xi_1\xi_2 \widehat{\sH}_{2,\star} =\widehat{F}^-_{1,\star} -\frac{\xi_2}{\nu\tau} \widehat{F}^-_{6,\star},\\
\partial_{Y_3} \widehat{\sE}_{1,\star}
+\frac{ik}{\nu\tau} \xi_1\xi_2 \widehat{\sH}_{1,\star}
+\frac{ik}{\nu\tau} (\nu^2\tau^2-\xi_1^2) \widehat{\sH}_{2,\star} =\widehat{F}^-_{2,\star} +\frac{\xi_1}{\nu\tau} \widehat{F}^-_{6,\star}.
\end{align*}
Eliminating $\widehat{\sH}_{1,\star}$ gives
\begin{align*}
&(\frac{ik}{\nu\tau})^2 ((\nu^2\tau^2-\xi_2^2)(\nu^2\tau^2-\xi_1^2)-\xi_1^2\xi_2^2)  \widehat{\sH}_{2,\star}
+\frac{ik}{\nu\tau} (\nu^2\tau^2-\xi_2^2) \partial_{Y_3} \widehat{\sE}_{1,\star}
+\frac{ik}{\nu\tau} \xi_1\xi_2 \partial_{Y_3} \widehat{\sE}_{2,\star} \\
&\qquad\qquad
=\frac{ik}{\nu\tau} \Big((\nu^2\tau^2-\xi_2^2)\widehat{F}^-_{2,\star}-\xi_1\xi_2\widehat{F}^-_{1,\star} +\nu\tau \xi_1 \widehat{F}^-_{6,\star}\Big)
.
\end{align*}
Since $\widehat{\sE}_{1,\star}$ and $\widehat{\sE}_{2,\star}$ defined in \cref{def-E*} satisfy the ODEs, 
\begin{align*}
&k^2(1-\nu^2 \tau^2 ) \Big(\partial_{Y_3} \widehat{\sH}_{2,\star}
+\frac{ik}{\nu\tau} (\nu^2\tau^2-\xi_2^2) \widehat{\sE}_{1,\star}
+\frac{ik}{\nu\tau} \xi_1\xi_2 \widehat{\sE}_{2,\star} \Big)\\
&= \frac{ik}{\nu\tau} \Big((\nu^2\tau^2-\xi_2^2)\partial_{Y_3}\widehat{F}^-_{2,\star}-\xi_1\xi_2\partial_{Y_3}\widehat{F}^-_{1,\star} +\nu\tau \xi_1 \partial_{Y_3}\widehat{F}^-_{6,\star}\Big)\\
& \quad -\frac{ik}{\nu\tau}(\nu^2\tau^2-\xi_2^2) \Big( \partial_{Y_3} \widehat{F}^-_{2,\star}+ik(-\nu\tau\widehat{F}^-_{4,\star}+\xi_1\widehat{F}^-_{8,\star} - \xi_2 \widehat{F}^-_{3,\star})\Big)\\
& \quad -\frac{ik}{\nu\tau}\xi_1\xi_2 \Big(- \partial_{Y_3} \widehat{F}^-_{1,\star}+ik(-\nu\tau\widehat{F}^-_{5,\star}+\xi_2\widehat{F}^-_{8,\star} + \xi_1 \widehat{F}^-_{3,\star})\Big)\\
&= \frac{ik}{\nu\tau}\Big\{
 \nu\tau \xi_1 \partial_{Y_3} \widehat{F}^-_{6,\star} + ik\Big( 
 \nu\tau (\nu^2\tau^2-\xi_2^2)\widehat{F}^-_{4,\star}
 +\nu\tau \xi_1\xi_2\widehat{F}^-_{5,\star}\\
&\qquad
  -\nu^2\tau^2\xi_1\widehat{F}^-_{8,\star} +(\nu^2\tau^2-1)\xi_2 \widehat{F}^-_{3,\star}\Big)
\Big\}\\
&=k^2(1-\nu^2 \tau^2 ) (\widehat{F}^-_{4,\star} +\frac{\xi_2}{\nu\tau}\widehat{F}^-_{3,\star}),
\end{align*}
where in the last step \eqref{comp-d-n0*} has been applied. 
So {the 4th equation of \cref{f-inh-n0*} in the vacuum} is proved, and the proof of the 5th one is similar to this one.
{The last equation in the vacuum} is proved by the compatibility for the divergences \eqref{comp-d-n0*},  just like the 7th one.

Now we have found the solution to \cref{f-inh-n0*} in the vacuum. With the same solution to \cref{f-inh-n0*} in the plasma part of \cite{Pierre2021}, it remains to verify that they satisfy the 6th boundary condition.

\vspace{0.5cm}

\noindent \textit{2) Verifying the 6th boundary condition.}

Since the above constructed $\widehat{\mathbb{V}}_{\star}=(\widehat{\sH}_{1,\star},\widehat{\sH}_{2,\star},\widehat{\sH}_{3,\star},\widehat{\sE}_{1,\star},\widehat{\sE}_{2,\star},\widehat{\sE}_{3,\star})$ solves \cref{f-inh-n0*}, then it also solves \cref{ode-V*}.
Setting $\widetilde{\sQ}_{\star}=H^0_j \widehat{\sH}_{j,\star}-E^0_3 \widehat{\sE}_{3,\star}$, by \cref{f-inh-n0*}, \eqref{comp-d-n0*} and \cref{ode-V*}, we have
\begin{subequations}
	\label{sQ+}
	\begin{align}
	%equation
	&\partial_{Y_3}^2 \widehat{\sQ}_{\star}-k^2 \widehat{\sQ}_{\star}  = \partial_{Y_3} \widehat{F}^+_{3,\star}-ik(a^++c^+)\widehat{F}^+_{7,\star}+2ikb^+\widehat{F}^+_{8,\star} +ik\xi_j \widehat{F}^+_{j,\star},  \notag \\
	&\quad\quad =\underbrace{\partial_{Y_3} \widehat{F}^+_{3,\star}+\frac{2b^+}{\tau}\partial_{Y_3} \widehat{F}^+_{6,\star}-ik(a^++c^+)\widehat{F}^+_{7,\star} +\frac{2ikb^+}{\tau}\xi_j\widehat{F}^+_{3+j,\star} +ik\xi_j \widehat{F}^+_{j,\star} }_{=:\mathscr{F}^+}   
	 \\
	%boundary condition
	&\partial_{Y_3} \widehat{\sQ}_{\star} |_{y_3=Y_3=0}= \widehat{F}^+_{3,\star} |_{y_3=Y_3=0}-ikc^+\widehat{G}_{1,\star}+ikb^+\widehat{G}_{2,\star},
	\end{align}
\end{subequations}
\begin{subequations}
	\label{sQ-}
	\begin{align}
	% equation 
	&\partial_{Y_3}^2 \widetilde{\sQ}_{\star}-k^2(1-\nu^2\tau^2) \widetilde{\sQ}_{\star}  = -H^0_1 \partial_{Y_3} \widehat{F}^-_{5,\star} + H^0_2 \partial_{Y_3} \widehat{F}^-_{4,\star}-E^0_3\partial_{Y_3} \widehat{F}^-_{8,\star} - ika^-_j  \widehat{F}^-_{j,\star} \notag  \\
	&\qquad\qquad  +ik b^- \widehat{F}^-_{7,\star} +ik(\xi_2 H^0_1-\xi_1 H^0_2+\nu\tau E^0_3) \widehat{F}^-_{6,\star} \notag \\
	&\quad\quad   
	\left.
	\begin{aligned}
	&=\partial_{Y_3} ( -\frac{a^-_1}{\nu\tau} \widehat{F}^-_{5,\star} +\frac{a^-_2}{\nu\tau}\widehat{F}^-_{4,\star} -\frac{E^0_3}{ik\nu\tau}\partial_{Y_3}\widehat{F}^-_{6,\star}  + \frac{b^-}{\nu\tau}\widehat{F}^-_{3,\star} )+\frac{ik}{\nu\tau} (-\nu\tau a^-_j+\xi_j b^-)\widehat{F}^-_{j,\star}  \\
	&\quad   +ik(\xi_2 H^0_1-\xi_1 H^0_2+\nu\tau E^0_3) \widehat{F}^-_{6,\star}	\end{aligned}
	\right\}
	=:\mathscr{F}^-_7,    \\
	%boundary condition
	& \partial_{Y_3} \widetilde{\sQ}_{\star} |_{y_3=Y_3=0}=
	(-H^0_1  \widehat{F}^-_{5,\star}   + H^0_2 \widehat{F}^-_{4,\star} -E^0_3 \widehat{F}^-_{8,\star}) |_{y_3=Y_3=0}  + ika^-_j \widehat{G}_{3+j,\star}+ikb^- \widehat{G}_{3,\star} \notag \\
	&\quad\quad=(-\frac{a^-_1}{\nu\tau}\widehat{F}^-_{5,\star} +\frac{a^-_2}{\nu\tau}\widehat{F}^-_{4,\star} -\frac{E^0_3}{ik\nu\tau}\partial_{Y_3}\widehat{F}^-_{6,\star}  ) |_{y_3=Y_3=0}+ik(a^-_j \widehat{G}_{3+j,\star}+b^- \widehat{G}_{3,\star}).
	\end{align}
\end{subequations}
The solutions $\widehat{\sQ}_{\star}$ and $\widetilde{\sQ}_{\star}$ are in the form:
\begin{align}\label{sQ-sol}
\begin{aligned}
\widehat{\sQ}_{\star}&=\kappa^+ e^{-|k|Y_3}-\frac{1}{2|k|} \left( \int_{0}^{Y_3} e^{-|k|(Y_3-Y)} \mathscr{F}^+(Y) dY +\int_{Y_3}^{+\infty} e^{-|k|(Y-Y_3)}  \mathscr{F}^+(Y) dY \right),\\
\widetilde{\sQ}_{\star}&=\kappa^-_7 e^{|k|\sqrt{1-\nu^2\tau^2}Y_3}-\frac{1}{2|k|\sqrt{1-\nu^2\tau^2}} \Big( \int_{-\infty}^{Y_3} e^{-|k|\sqrt{1-\nu^2\tau^2}(Y_3-Y)} \mathscr{F}^-_7(Y) dY \\
&\quad \quad +\int_{Y_3}^{0} e^{-|k|\sqrt{1-\nu^2\tau^2}(Y-Y_3)}  \mathscr{F}^-_7(Y) dY \Big),
\end{aligned}
\end{align}
where
\begin{align}\label{kappa-sol}
\begin{aligned}
\kappa^+ |_{y_3=0}=&-\frac{1}{2|k|} \int_{0}^{+\infty} e^{-|k|Y} \mathscr{F}^+ |_{y_3=0}(Y)  dY 
   -\frac{1}{|k|}\widehat{F}^+_{3,\star} |_{y_3=Y_3=0}+ i\sgnk (c^+\widehat{G}_{1,\star}- b^+\widehat{G}_{2,\star})\\
\kappa^-_7 |_{y_3=0}=&-\frac{1}{2|k|\sqrt{1-\nu^2\tau^2}}  \int_{-\infty}^{0} e^{|k|\sqrt{1-\nu^2\tau^2}Y} \mathscr{F}^-_7 |_{y_3=0}(Y) dY+ \frac{1}{|k| \sqrt{1-\nu^2\tau^2}} \Big( (-\frac{a^-_1}{\nu\tau}\widehat{F}^-_{5,\star} \\
                   &\quad+\frac{a^-_2}{\nu\tau}\widehat{F}^-_{4,\star} -\frac{E^0_3}{ik\nu\tau}\partial_{Y_3}\widehat{F}^-_{6,\star} ) |_{y_3=Y_3=0} + ik(a^-_j \widehat{G}_{3+j,\star}+b^- \widehat{G}_{3,\star})  \Big).
\end{aligned}
\end{align}
Note that here $\kappa^{+}$ and $\kappa_7^-$ are functions of $(t,y,k)$,  $\kappa^{+}$ only has constraints on $y_3=0$, 
while $\kappa_7^-$ is actually determined by $\widehat{\sH}_{1,\star}, \widehat{\sH}_{2,\star}$ and $\widehat{\sE}_{3,\star}$, or just $\kappa_4^-$ and $\kappa_5^-$ (especially in the interior of the vacuum) by the definition of $\widetilde{\sQ}_{\star}$.
However, to verify the 6th boundary condition, we will just use the value of $\kappa_7^-$ on $y_3=0$, which is determined by $\mathscr{F}^-_7$ on $y_3=0$ and $\widehat{G}_{3,\star},\widehat{G}_{4,\star},\widehat{G}_{5,\star}$, which have been verified above.
%Namely, $\kappa_7^-|_{y_3=0}$ is independent of the constructions of $\kappa_4^-$ and $\kappa_5^-$.

With the explicit expressions of $\widehat{\sQ}_{\star}$ and $\widetilde{\sQ}_{\star}$, one has
\begin{align*}
&(\widehat{\sQ}_{\star} -\widetilde{\sQ}_{\star} ) |_{y_3=Y_3=0} =\Big(-\frac{1}{|k|} \int_{0}^{+\infty} e^{-|k|Y} \mathscr{F}^+ |_{y_3=0}(Y)  dY 
	-\frac{1}{|k|}\widehat{F}^+_{3,\star} |_{y_3=Y_3=0}\\
&\quad \quad + i\sgnk( c^+\widehat{G}_{1,\star}-  b^+\widehat{G}_{2,\star})\Big)
%2nd integral
+\Big( \frac{1}{|k|\sqrt{1-\nu^2\tau^2}}  \int_{-\infty}^{0} e^{|k|\sqrt{1-\nu^2\tau^2}Y} \mathscr{F}^-_7 |_{y_3=0}(Y) dY\\
&\quad\quad - \frac{1}{|k| \sqrt{1-\nu^2\tau^2}} ( -\frac{a^-_1}{\nu\tau}\widehat{F}^-_{5,\star}  +\frac{a^-_2}{\nu\tau}\widehat{F}^-_{4,\star} -\frac{E^0_3}{ik\nu\tau}\partial_{Y_3}\widehat{F}^-_{6,\star}  ) |_{y_3=Y_3=0} \\
&\quad\quad   -\frac{i\sgnk}{\sqrt{1-\nu^2\tau^2}} (a^-_j \widehat{G}_{3+j,\star} +b^- \widehat{G}_{3,\star} )\Big) =: (I_1)+(I_2),
\end{align*}
where $\mathscr{F}$ is defined in \eqref{sQ-}.
As shown in \cite{Pierre2021}, applying integration by parts gives
\begin{align*}
I_1=&-\frac{i\sgnk}{\tau} \left( \int_{0}^{+\infty} e^{-|k|Y_3} \sL^+(k) \cdot \widehat{F}^+_{\star} |_{y_3=0}(Y) dY+ l^+_j \widehat{G}_{j,\star} \right)\\
I_2=& -\frac{i\sgnk }{\nu\tau\sqrt{1-\nu^2\tau^2}}  \Big\{-
\int_{-\infty}^{0} e^{|k|\sqrt{1-\nu^2\tau^2}Y} 
\Big(-i\sgnk a^-_1  \sqrt{1-\nu^2\tau^2}\widehat{F}^-_{5,\star} \\
&\quad 
+ i\sgnk a^-_2  \sqrt{1-\nu^2\tau^2} \widehat{F}^-_{4,\star} 
-\frac{\sqrt{1-\nu^2\tau^2}}{|k|}E_3^0 \partial_{Y_3} \widehat{F}^-_{6,\star} +i\sgnk \sqrt{1-\nu^2\tau^2} b^- \widehat{F}^-_{3,\star} \\
&\quad 
+  (-\nu\tau a^-_j+\xi_j b^-) \widehat{F}^-_{j,\star} + \nu\tau(\xi_2 H^0_1-\xi_1 H^0_2+\nu\tau E^0_3) \widehat{F}^-_{6,\star}\Big)\Big|_{y_3=0}(Y) dY\\
&\quad 
-\frac{b^-}{ik} \widehat{F}^-_{3,\star} |_{y_3=Y_3=0} + \nu\tau(a^-_j \widehat{G}_{3+j,\star}+b^- \widehat{G}_{3,\star})\Big\}\\
=&
-\frac{i\sgnk}{\tau} \Big( -\int_{-\infty}^{0} e^{|k|\sqrt{1-\nu^2\tau^2}Y} \sL^-(k) \cdot \widehat{F}^-_{\star} |_{y_3=0}(Y) dY \\
&\quad +\frac{E^0_3}{|k|\nu} \widehat{F}^-_{6,\star} |_{y_3=Y_3=0}(Y) 
-\frac{1}{\nu\sqrt{1-\nu^2\tau^2}}(-\nu\tau a^-_j +\xi_j b^-)  \widehat{G}_{3+j,\star}  \Big)
\end{align*}
Here in the last step, we  use  ``$\cdot$'' to denote the inner product of vectors of complex numbers: $X\cdot Y = \sum_{\alpha} \overline{X}_\alpha Y_\alpha$, integration by parts to substitute $\partial_{Y_3} \widehat{F}_{6,\star}^-$, and \eqref{comp-a-n0*}  to substitute  $\widehat{F}_{3,\star}^-$.
By \eqref{comp-e-n0*}, one has
\begin{align*}
(\widehat{\sQ}_{\star} -\widetilde{\sQ}_{\star} ) |_{y_3=Y_3=0} = -\frac{i\sgnk}{\tau} \cdot i\sgnk\tau \widehat{G}_{6,\star}=\widehat{G}_{6,\star}.
\end{align*}
Therefore, the 6th boundary condition is verified. Once $\widehat{\sQ}$ is constructed, the other unknowns in the plasma part are found exactly as in \cite{Pierre2021}, so we omit the details.

\vspace{0.5cm}

With the above two steps we have completed the construction of a solution to the fast problem \cref{f-inh-n0*}. It remains to verify that the sums
$$\sum\limits_{k\in \mathbb{Z} \setminus \{0\} } \widehat{\mathbb{U} }_{\star}   e^{ik\theta},~\qquad~ 
\sum\limits_{k\in \mathbb{Z} \setminus \{0\} }\widehat{\mathbb{V} }_{\star}  e^{ik\theta}$$
belong to $S^{\pm}_{\star}$. For this part we refer the readers to \cite[Lemma~4.2]{Marcou2010}, where this question is addressed in full details.
In conclusion,
\begin{equation}\label{def-Usharp}
\mathbb{U}_{\#}:= \sum\limits_{k\in \mathbb{Z} \setminus \{0\} } (\underline{\widehat{\mathbb{U} } }+\widehat{\mathbb{U} }_{\star}  )(k)\, e^{ik\theta},
\quad 
\mathbb{V}_{\#}:= \sum\limits_{k\in \mathbb{Z} \setminus \{0\} } (\underline{\widehat{\mathbb{V} } }+\widehat{\mathbb{V} }_{\star}  )(k)\, e^{ik\theta}
\end{equation}
belong to $S^{\pm}$ and are the solution to \cref{f-inh-n0},that is the solution to \cref{f-inh} for non-zero Fourier modes. The proof of \cref{thm-fast} is complete.

%%%%%%%%%%%%%%%%%%%%%%%%%%%%%%%%%%%%%%%%%%%%%%%%%%%%%%%%%%%%%%%%
%%%%%Solving the WKB cascade: the leading amplitude  #4
%%%%%%%%%%%%%%%%%%%%%%%%%%%%%%%%%%%%%%%%%%%%%%%%%%%%%%%%%%%%%%%%
\vspace{0.5cm}
\section{Solving the WKB cascade: the leading amplitude} \label{sec-leading}

In this section we begin to solve the WKB cascade \cref{WKB-eq}, \cref{WKB-div} and \cref{WKB-jump}, that is we construct the leading amplitude $(U^1,V^1,\varphi^{2})$ in the WKB ansatz \cref{WKB-ansatz}. To do that, we first need to identify the degrees of freedom at our disposal, and then, to complete the construction of the solution, we determine all functions by imposing (some of) the necessary solvability conditions \eqref{f-comp} for the fast problem to be solved by the first corrector $(U^2,V^2,\varphi^{3})$.
%A brief description of the inductive construction, for instances when constructing $(U^m,V^m,\varphi^{m+1})$, is to identify the freedoms in solving the fast problem for 
% , and then to determine the rest freedoms such that the necessary solvability conditions \eqref{f-comp} on the fast problem for $(U^{m+1},V^{m+1},\varphi^{m+2})$ hold.
This procedure is standard in geometric optics, see \cite{ Rauch2012}, with the main feature here, as in \cite{Pierre2021},  that some of the conditions \eqref{f-comp} must {\it come for free} in the WKB cascade because we have more constraints than degrees of freedom
at our disposal.

Noting that in \cref{F01+}, \cref{F01-}, \cref{F701} and \eqref{Gm} the source terms and boundary terms for $m=0$ are exactly zero, the leading amplitude $(U^1, V^1, \varphi^2)$ satisfy the homogeneous fast problem
\begin{equation}\label{f-leading}
	\left\{
	 \begin{aligned}
	 	&\sL^+_f U^1 =0, &y_3,Y_3>0,\\
	 	&\partial_{Y_3}B_3^1+\xi_j \partial_{\theta} B_j^1 =0,&y_3,Y_3>0,\\
	 	&\sL^-_f V^1 =0,&y_3,Y_3<0,\\
	 	&\partial_{Y_3}H_3^1+\xi_j \partial_{\theta}H_j^1 =0,&y_3,Y_3<0,\\
	 	&\partial_{Y_3}E_3^1+\xi_j \partial_{\theta}E_j^1 =0,&y_3,Y_3<0,\\
	 	&\mathbb{B}^+ U^1  |_{y_3=Y_3=0}+\mathbb{B}^- V^1  |_{y_3=Y_3=0}+\underline{b} \partial_{\theta} \varphi^2 =0.
	 \end{aligned}
	\right.
\end{equation}
It follows from \cref{prop-h} that the leading amplitude $(U^1, V^1, \varphi^2)$ is in the form
\begin{align}
	U^1(t,y,Y_3,\theta)=&\underline{U}^1(t,y)+(\widehat{u}_{1,\star}^1,\widehat{u}_{2,\star}^1,0,\widehat{B}_{1,\star}^1,\widehat{B}_{2,\star}^1,0,0)^{\top} (t,y,Y_3) \notag\\
	&\quad + \sum\limits_{k\in \mathbb{Z} \setminus {0} }\gamma^{1,+}(t,y,k) e^{-|k|Y_3+ik\theta} \sR^+(k), \label{sol-leading-general}\\
	%line 
	V^1(t,y,Y_3,\theta)=&\underline{V}^1(t,y)+\sum\limits_{k\in \mathbb{Z} \setminus \{0\}}\gamma^{1,-}_j(t,y,k) e^{|k|\sqrt{1-\nu^2\tau^2}Y_3+ik\theta} \sR^-_j(k), \notag
\end{align}
where $\underline{U}^1 \in \underline{S}^{+}$ and $\underline{V}^1 \in \underline{S}^{-}$ are independent of $\theta$ and satisfy the boundary conditions:
\begin{align*}
	&\underline{u}^1_3 |_{y_3=0}=\underline{B}^1_3 |_{y_3=0}=\underline{H}^1_3 |_{y_3=0}=\underline{E}^1_1 |_{y_3=0}=\underline{E}^1_2 |_{y_3=0}=0,\\
	& \underline{q}^1 |_{y_3=0}-H^0_j \underline{H}^1_j |_{y_3=0}+E^0_3 \underline{E}^1_3 |_{y_3=0}=0.
\end{align*}
Moreover,
$\widehat{u}_{1,\star}^1(0),\widehat{u}_{2,\star}^1(0), \widehat{B}_{1,\star}^1(0),\widehat{B}_{2,\star}^1(0) \in S^+_{\star}$, and can be freely chosen. 
$\sR^+(k),\sR^-_1(k), \sR^-_2(k)$ are explicitly given in \cref{vector-R}, the coefficients $\gamma^{1,+}, \gamma^{1,-}_1,\gamma^{1,-}_2$ satisfy that for all $(t,y',k)\in [0,T]\times \mathbb{T}^2\times (\mathbb{Z}\setminus \{0\})$,
\begin{equation}
	\label{def-gamma1}
	\gamma^{1,+}(t,y',0,k)=|k|\widehat{\varphi}^2(t,y',k),
	\quad\gamma^{1,-}_j(t,y',0,k)=\widetilde{\gamma}^1_j|k|\widehat{\varphi}^2(t,y',k), 
\end{equation}
where $a^-_j\widetilde{\gamma}^1_j =\frac{(a_1^-)^2+(a_2^-)^2-(b^-)^2}{1-\nu^2\tau^2}$.

All of $\underline{U}^1, \underline{V}^1$, $\Pi \widehat{U}_\star^1(0)$, and $\varphi^2$ are the degrees of freedom of the fast problem for $(U^1,V^1,\varphi^2)$,
and they will be determined by the necessary solvability conditions \eqref{f-comp} on the fast problem for $(U^{2},V^{2},\varphi^{3})$ in the subsequent parts.
It will be proved that $\underline{U}^1, \underline{V}^1$ and $\Pi \,\widehat{U}_\star^1(0)$ can be chosen equal to $0$, 
and that $\varphi^2$ is determined by 
a nonlocal Hamilton-Jacobi equation. The final time $T>0$ of local existence will be fixed in this step.

\subsection{The slow mean and fast mean of the leading profile}
By \cref{WKB-eq}, \cref{WKB-div} and \cref{WKB-jump}, the first corrector $(U^{2},V^{2},\varphi^{3})$ must satisfy
\begin{equation}\label{f-corrector1}
	\left\{
	\begin{aligned}
		&\sL^+_f U^2 =F^{1,+}, &y_3,Y_3>0,\\
		&\partial_{Y_3}B_3^2+\xi_j \partial_{\theta} B_j^2	 =F^{1,+}_8,&y_3,Y_3>0,\\
		&\sL^-_f V^2 =F^{1,-},&y_3,Y_3<0,\\
		&\partial_{Y_3}H_3^2+\xi_j \partial_{\theta}H_j^2 =F^{1,-}_7,&y_3,Y_3<0,\\
		&\partial_{Y_3}E_3^2+\xi_j \partial_{\theta}E_j^2 =F^{1,-}_8,&y_3,Y_3<0,\\
		&\mathbb{B}^+ U^2  |_{y_3=Y_3=0}+\mathbb{B}^- V^2  |_{y_3=Y_3=0}+\underline{b} \partial_{\theta} \varphi^3 =G^1,
	\end{aligned}
	\right.
\end{equation}
where $F^{1,\pm}, F^{1,-}_7, F^{1,\pm}_8$ are given in \cref{F01+}, \cref{F01-} and \cref{F701}; for $G^1$ see definition \eqref{Gm}.
The necessary solvability condition \eqref{comp-b} of \cref{f-corrector1} is 
\begin{equation}
	\label{comp-b1}
		\underline{\widehat{F}}^{1,\pm}(t,y,0)=0, ~\underline{\widehat{F}}_7^{1,-}(t,y,0)= \underline{\widehat{F}}_8^{1, \pm}(t,y,0)=0.
\end{equation}
From \cref{F01+}, \cref{F01-}, \cref{F701}, after taking the limit as $Y_3\to\pm\infty$ and the mean in $\theta$, it follows that $\underline{U}^1$ and $\underline{V}^1$ should satisfy the \textit{slow} problem
\begin{equation}
	\label{s-corrector1}
	\left\{
	\begin{aligned}
		&L_s^+ \underline{\widehat{U}}^{1}(0)=0, \quad  \nabla\cdot \underline{\widehat{B}}^{1}(0) = 0,  & y\in \Omega_0^+, \\
		&L_s^- \underline{\widehat{V}}^{1}(0)=0, \quad \nabla\cdot \underline{\widehat{H}}^{1}(0) = 0, \quad
		\nabla\cdot \underline{\widehat{E}}^{1}(0) = 0, & y\in \Omega_0^-,
	\end{aligned}
	\right.
\end{equation}
with boundary conditions on $\Gamma_0$  
\begin{equation}
	\label{bc-corrector1}
	\left\{
	\begin{aligned}
	&\underline{u}^1_3 |_{y_3=0}=\underline{B}^1_3 |_{y_3=0}=\underline{H}^1_3 |_{y_3=0}=\underline{E}^1_1 |_{y_3=0}=\underline{E}^1_2 |_{y_3=0}=0\,,\\
	& \underline{q}^1 |_{y_3=0}-H^0_j \underline{H}^1_j |_{y_3=0}+E^0_3 \underline{E}^1_3 |_{y_3=0}=0\,.
	\end{aligned}
\right.
\end{equation}
Problem \cref{s-corrector1}, \cref{bc-corrector1} is simply the linearized plasma--vacuum problem with zero source terms; the choice of the initial data is free, up to the obvious compatibility conditions for the divergences and the boundary conditions at the initial time. One can find that the trivial solution is exactly the unique one to \cref{s-corrector1} with zero initial data (the easiest possible choice)
\begin{equation*}
	(\underline{\widehat{u}}^{1}, \underline{\widehat{B}}^1,\underline{\widehat{E}}^1,\underline{\widehat{H}}^1 ) |_{t=0}\equiv0\,. 
\end{equation*}
This choice gives zero slow means of the leading profile \cref{sol-leading-general}: $\underline{\widehat{U}}^{1}(0)=\underline{U}^1(t,y)=0, \underline{\widehat{V}}^{1}(0)=\underline{V}^1(t,y)=0$.
%In \cref{sec-correctors} we will consider more complicated slow problems.

The fast mean of the leading profile is determined by imposing the solvability condition \eqref{comp-c} on problem \cref{f-corrector1}. As in \cite{Pierre2021}, we prescribe zero initial values for the fast mean of the velocity and magnetic fields in the plasma region; this choice implies $\widehat{u}_{1,\star}^1,\widehat{u}_{2,\star}^1,\widehat{B}_{1,\star}^1,\widehat{B}_{2,\star}^1\equiv0$. Since this is just a repetition of the argument
in \cite{Pierre2021} we omit the details here. 
Therefore, $\underline{U}^1, \underline{V}^1$ and $\Pi \widehat{U}_\star^1(0)$ can be chosen as $0$, and the leading profile $(U^1, V^1, \varphi^2)$ can be simplified as 
\begin{equation}\label{sol-leading0}
	\begin{aligned}
		U^1(t,y,Y_3,\theta)=& \sum\limits_{k\in \mathbb{Z} \setminus \{0\} }\gamma^{1,+}(t,y,k) e^{-|k|Y_3+ik\theta} \sR^+(k), \\
		%line 
		V^1(t,y,Y_3,\theta)=&\sum\limits_{k\in \mathbb{Z} \setminus \{0\}}\gamma^{1,-}_j(t,y,k) e^{|k|\sqrt{1-\nu^2\tau^2}Y_3+ik\theta} \sR^-_j(k).
	\end{aligned}
\end{equation} 
It remains to determine $\varphi^2$ in $\gamma^{1,+}(t,y,k)$ and $\gamma^{1,-}_j(t,y,k)$, $j=1,2$ in \cref{def-gamma1} by the solvability condition \eqref{comp-e}.

\subsection{The nonlocal Hamilton-Jacobi equation for the leading front}

\subsubsection{Derivation of the equation of the leading front}
We focus again on the inhomogeneous fast problem \cref{f-corrector1} for the first corrector $(U^{2},V^{2},\varphi^{3})$. From \cref{thm-fast} we know that a necessary solvability condition for \cref{f-corrector1} is the orthogonality condition \eqref{comp-e} for all non-zero Fourier modes. 
In the present case it reads:
\begin{equation}
	\begin{aligned}
		&\quad \int_{0}^{+\infty}e^{-|k|Y_3}\sL^+(k) \cdot \widehat{F}^{1,+} (k)|_{y_3=0}  dY_3 + l^+_j \widehat{G}^1_j(k)\\
		& -\int_{-\infty}^{0} e^{|k|\sqrt{1-\nu^2\tau^2}Y_3}\sL^-(k) \cdot \widehat{F}^{1,-} (k)|_{y_3=0}   dY_3 -\frac{-\nu\tau a^-_j +\xi_j b^-}{\nu\sqrt{1-\nu^2\tau^2}}  \widehat{G}^1_{3+j}(k) \\
		& -i\sgnk \tau \widehat{G}^1_6(k) +\frac{E^0_3}{|k|\nu} \widehat{F}^{1,-}_6 (k)|_{y_3=Y_3=0} =0, \quad \text{for all~} k\in \mathbb{Z}\setminus \{0\}.
	\end{aligned}\label{comp-e1}
\end{equation}
%Now we will rewrite the above equation in terms of $\widehat{\varphi}^2$, i.e. replacing $\widehat{F}^{1,\pm} (k), \widehat{G}^1(k)$ by $\widehat{U}^1(k), \widehat{V}^1(k), \widehat{\varphi}^2(k)$, and finally by the only unknown $\widehat{\varphi}^2(k)$.

We will see that the boundary source term $\widehat{G}^1(k)$ can be explicitely written in terms of the leading front $\widehat{\varphi}^2$. Moreover, the source terms $\widehat{F}^{1,\pm} (k)|_{y_3=0}$ will be written in terms of the traces $\widehat{U}^1(k)|_{y_3=0}, \widehat{V}^1(k)|_{y_3=0}$, which in turn can also be written in terms of $\widehat{\varphi}^2$ by \cref{sol-leading0}. Thus all \cref{comp-e1} will be expanded in terms of $\widehat{\varphi}^2$, and this will lead to the nonlocal Hamilton-Jacobi equation for the leading front.

\vspace{0.5cm}

The terms for the plasma side
are exactly the same as those in \cite{Pierre2021}, and thus here we omit the calculations and just list the result:
\begin{align}
	&\int_{0}^{+\infty}e^{-|k|Y_3}\sL^+(k) \cdot \widehat{F}^{1,+} (k)|_{y_3=0}  dY_3 + l^+_j \widehat{G}^1_j(k)\notag\\
	&=-2\tau c^+ \dt \widehat{\varphi}^2 (k) -2\tau (c^+ u^0_j -b^+ B^0_j) \partial_{y_j} \widehat{\varphi}^2 (k) \label{front-1}  \\
	&\quad -2i\tau \Big( (c^+)^2 -(b^+)^2   \Big) \sgnk \sum\limits_{k_1+k_2 = k}  \frac{|k||k_1||k_2|}{|k|+|k_1|+|k_2|}\widehat{\varphi}^2 (k_1) \widehat{\varphi}^2(k_2)=0\,.
	\notag
\end{align}

\textbullet~
\textit{The boundary terms  for the vacuum of \eqref{comp-e1}.} Recall in \cref{sol-leading0} that 
$$\widehat{V}^1 (k) |_{y_3=Y_3=0} = \widetilde{\gamma}^1_j |k| \widehat{\varphi}^2 (k)  \sR^-_j(k),
\quad 
\partial_{Y_3}\widehat{V}^1 (k) |_{y_3=Y_3=0} = \widetilde{\gamma}^1_j |k|^2 \sqrt{1-\nu^2\tau^2} \widehat{\varphi}^2 (k)  \sR^-_j(k). $$
Therefore, by \eqref{Gm} and \eqref{vector-R}, one has that 
\begin{align*}
\widehat{G}^1_4(k)=& -E^0_3 \partial_{y_2} \widehat{\varphi}^2(k)-\nu H^0_1 \partial_t \widehat{\varphi}^2 (k)  +   i\sqrt{1-\nu^2 \tau ^2 } \big\{ -   \xi_2 (-\xi_2 \widetilde{\gamma}^1_1 +\xi_1 \widetilde{\gamma}^1_2) \\
&\quad -    \nu^2 \tau^2 \widetilde{\gamma}^1_1 \big\} \sum\limits_{k_1+k_2 = k}  \frac{k_1|k_2|+k_2|k_1|}{2 }\widehat{\varphi}^2 (k_1) \widehat{\varphi}^2(k_2),
\\%new equation
\widehat{G}^1_5(k)=& E^0_3 \partial_{y_1} \widehat{\varphi}^2(k)-\nu H^0_2 \partial_t \widehat{\varphi}^2 (k) +  i\sqrt{1-\nu^2 \tau ^2 }   \big\{  \xi_1 (-\xi_2 \widetilde{\gamma}^1_1 +\xi_1 \widetilde{\gamma}^1_2) \\
&\quad -  \nu^2 \tau^2 \widetilde{\gamma}^1_2 \big\} \sum\limits_{k_1+k_2 = k}  \frac{k_1|k_2|+k_2|k_1|}{2 }\widehat{\varphi}^2 (k_1) \widehat{\varphi}^2(k_2),
\\%new equation
\widehat{G}^1_6(k)=& \frac{1}{2} \sum\limits_{k_1+k_2 = k}  \big(\widehat{H}^1_\alpha(k_1)\widehat{H}^1_\alpha(k_2)-\widehat{E}^1_\alpha(k_1)\widehat{E}^1_\alpha(k_2) \big)\\
=&\frac{1}{2} \Big\{ (1-\nu^2\tau^2)\nu^2\tau^2 ((\widetilde{\gamma}^1_1)^2+(\widetilde{\gamma}^1_2)^2 ) - \sgnkone\sgnktwo \nu^2\tau^2 ( \xi_1 \widetilde{\gamma}^1_1 + \xi_2 \widetilde{\gamma}^1_2)^2\\
& \quad +\sgnkone\sgnktwo  ( \xi_1 \xi_2  \widetilde{\gamma}^1_1 + (\nu^2\tau^2-\xi_1^2) \widetilde{\gamma}^1_2)^2
 \\
&\quad +\sgnkone\sgnktwo  (  (\nu^2\tau^2-\xi_2^2) \widetilde{\gamma}^1_1+\xi_1 \xi_2  \widetilde{\gamma}^1_2)^2\\
&\quad - (1-\nu^2\tau^2) (-\xi_2 \widetilde{\gamma}^1_1 +\xi_1 \widetilde{\gamma}^1_2 )^2
  \Big\} |k_1||k_2| \widehat{\varphi}^2 (k_1) \widehat{\varphi}^2(k_2)\\
  =&  (1-\nu^2\tau^2)  \big\{\nu^2\tau^2 ((\widetilde{\gamma}^1_1)^2+(\widetilde{\gamma}^1_2)^2 )- (-\xi_2 \widetilde{\gamma}^1_1 +\xi_1 \widetilde{\gamma}^1_2 )^2
  \big\}\\
  & \quad \cdot \sum\limits_{k_1+k_2 = k} \frac{1-\sgnkone\sgnktwo}{2} |k_1||k_2| \widehat{\varphi}^2 (k_1) \widehat{\varphi}^2(k_2),
\end{align*}
\begin{align*}
\widehat{F}^{1,-}_6 (k)|_{y_3=Y_3=0} =& - (\nu \dt \widehat{E}^1_3 -\partial_{y_1} \widehat{H}^1_2 +\partial_{y_2} \widehat{H}^1_1   ) (k)|_{y_3=Y_3=0}  \\
&\quad +\sum\limits_{k_1+k_2=k} ik_2 \widehat{\varphi}^2(k_2) \partial_{Y_3} (\nu\tau \widehat{E}_3^1 -\xi_1 \widehat{H}_2^1 +\xi_2 \widehat{H}_1^1 ) |_{y_3=Y_3=0}(k_1)\\
=&-\sqrt{1-\nu^2 \tau^2} \nu (-\xi_2 \widetilde{\gamma}^1_1 +\xi_1 \widetilde{\gamma}^1_2 ) |k| \dt \widehat{\varphi}^2(k)\\
&\quad +\sqrt{1-\nu^2 \tau^2} \nu\tau|k| (\widetilde{\gamma}^1_2 \partial_{y_1} \widehat{\varphi}^2(k) - \widetilde{\gamma}^1_1 \partial_{y_2} \widehat{\varphi}^2(k)).
\end{align*}
We observe, also using the relation $\xi_j a^-_j =\nu\tau b^-$, that the coefficient of the nonlinear terms of $-\frac{-\nu\tau a^-_j +\xi_j b^-}{\nu\sqrt{1-\nu^2\tau^2}}  \widehat{G}^1_{3+j}(k)$ can be simplified as 
\begin{align*}
&\frac{i}{\nu} \Big\{
(-\nu\tau a^-_1 +\xi_1 b^-) \big((\nu^2\tau^2 -\xi_2^2) \widetilde{\gamma}^1_1 + \xi_1\xi_2 \widetilde{\gamma}^1_2 \big) \\
&\quad + (-\nu\tau a^-_2 +\xi_2 b^-) \big(\xi_1\xi_2 \widetilde{\gamma}^1_1 +  (\nu^2\tau^2 -\xi_1^2)\widetilde{\gamma}^1_2 \big)
\Big\} \\
=&\frac{i}{\nu} \Big\{
(-\nu^3\tau^3a^-_1 +\nu\tau\xi_2^2 a^-_1 -\nu\tau \xi_1\xi_2  a^-_2   +\nu^2\tau^2 \xi_1 b^- )\widetilde{\gamma}^1_1 \\
&\quad + 
(-\nu^3\tau^3a^-_2 +\nu\tau\xi_1^2 a^-_2 -\nu\tau \xi_1\xi_2  a^-_1   +\nu^2\tau^2 \xi_2 b^- )\widetilde{\gamma}^1_2
\Big\}\\
=&i\tau (1-\nu^2\tau^2)  (a^-_j \widetilde{\gamma}^1_j) \,.
\end{align*}
Therefore, the boundary terms of \eqref{comp-e1} in vacuum part are
\begin{align}
&-\frac{-\nu\tau a^-_j +\xi_j b^-}{\nu\sqrt{1-\nu^2\tau^2}}  \widehat{G}^1_{3+j}(k)  -i\sgnk \tau \widehat{G}^1_6(k) +\frac{E^0_3}{|k|\nu} \widehat{F}^{1,-}_6 (k)|_{y_3=Y_3=0} 
\notag\\% new line
=& \Big(
\frac{-\nu\tau a^-_j +\xi_j b^-}{\sqrt{1-\nu^2\tau^2}} H^0_j - \sqrt{1-\nu^2\tau^2} (-\xi_2 \widetilde{\gamma}^1_1 +\xi_1 \widetilde{\gamma}^1_2 ) E^0_3  
\Big) \dt \widehat{\varphi}^2(k)
\notag\\% new line
&\quad +\Big(-
\frac{-\nu\tau a^-_2 +\xi_2 b^-}{\nu\sqrt{1-\nu^2\tau^2}} E^0_3+\tau  \sqrt{1-\nu^2\tau^2}\widetilde{\gamma}^1_2 E^0_3  
\Big) \partial_{y_1} \widehat{\varphi}^2(k)
\notag\\% new line
&\quad +\Big(
\frac{-\nu\tau a^-_1 +\xi_1 b^-}{\nu\sqrt{1-\nu^2\tau^2}} E^0_3- \tau\sqrt{1-\nu^2\tau^2} \widetilde{\gamma}^1_1 E^0_3  
\Big) \partial_{y_2} \widehat{\varphi}^2(k)
\notag\\% new line
&\quad +i\tau (1-\nu^2\tau^2)  (a^-_j \widetilde{\gamma}^1_j)
\sum\limits_{k_1+k_2 = k}  \frac{k_1|k_2|+k_2|k_1|}{2 }\widehat{\varphi}^2 (k_1) \widehat{\varphi}^2(k_2)
\notag\\% new line
& \quad -i \tau (1-\nu^2\tau^2)  \Big( 
\nu^2\tau^2 ((\widetilde{\gamma}^1_1)^2+(\widetilde{\gamma}^1_2)^2 )- (-\xi_2 \widetilde{\gamma}^1_1 +\xi_1 \widetilde{\gamma}^1_2 )^2
\Big)
\notag\\% new line
& \quad \cdot \sgnk\sum\limits_{k_1+k_2 = k} \frac{1-\sgnkone\sgnktwo}{2} |k_1||k_2| \widehat{\varphi}^2 (k_1) \widehat{\varphi}^2(k_2)\,. \label{front-3}
\end{align}

\textbullet~
\textit{The normal slow derivatives for the vacuum in \eqref{comp-e1}. } 
Now we decompose $\widehat{F}^{1,-} (k)|_{y_3=0} $ in the integration \cref{comp-e1} in several terms.
Since $\sL^-(k) \cdot A_3^- \sR^-_{j}(k)=0,~j=1,2$ (noting that  ``$\cdot$'' between two complex vectors denotes the inner product, i.e. $X\cdot Y = \sum_{\alpha} \overline{X}_\alpha Y_\alpha$), it follows from \cref{sol-leading0} that 
$$\sL^-(k) \cdot A_3^- \widehat{V}^{1,-}(Y_3,k)=0.$$
Therefore, differentiating the above equality with respect to $y_3$, taking trace on $y_3=0$ and integrating over $Y_3\in(-\infty,0)$ yields
\begin{align}
	\int_{-\infty}^{0}e^{|k|\sqrt{1-\nu^2\tau^2}Y_3}\sL^-(k) \cdot A^-_3 \partial_{y_3} \widehat{V}^1 (k)|_{y_3=0}  dY_3=0. 
	\label{front-2}
\end{align}

\textbullet~
\textit{The other linear terms for the vacuum in \eqref{comp-e1}. }
\begin{align}
& \int_{-\infty}^{0}e^{|k|\sqrt{1-\nu^2\tau^2}Y_3}\sL^-(k) \cdot (A^-_0 \dt+A^-_j \partial_{y_j}) \widehat{V}^1 (k)|_{y_3=0}  dY_3  \notag\\
&\quad  = \frac{\widetilde{\gamma}^1_i}{2\sqrt{1-\nu^2 \tau^2}} \Big\{
\big(\sL^-(k) \cdot A_0 \sR^-_i (k)\big) \dt \widehat{\varphi}^2 (k)+
\big(\sL^-(k) \cdot A_j \sR^-_i (k)\big) \partial_{y_j} \widehat{\varphi}^2 (k)
\Big\}   \notag\\
&\quad= \frac{1}{\sqrt{1-\nu^2 \tau^2}} \Big\{ (-a^-_j+2\nu\tau \xi_j b^-)  \widetilde{\gamma}^1_j\Big\}\dt \widehat{\varphi}^2 (k)  \notag \\
&\quad\quad +\frac{1}{\nu\sqrt{1-\nu^2 \tau^2}} \Big\{ - (\nu^2\tau^2-\xi_2^2 ) b^-\widetilde{\gamma}^1_1 + ( \nu\tau \xi_1 a^-_2 -\nu\tau \xi_2 a^-_1 -\xi_1\xi_2 b^- )\widetilde{\gamma}^1_2 \Big\}\partial_{y_1} \widehat{\varphi}^2 (k) \notag \\
&\quad\quad +\frac{1}{\nu\sqrt{1-\nu^2 \tau^2}} \Big\{( \nu\tau \xi_2 a^-_1 -\nu\tau \xi_1 a^-_2 -\xi_1\xi_2 b^- )  \widetilde{\gamma}^1_1 - (\nu^2\tau^2-\xi_1^2 ) b^- \widetilde{\gamma}^1_2 \Big\}\partial_{y_2} \widehat{\varphi}^2 (k).
\label{front-4}
\end{align}

\textbullet~
\textit{The quadratic terms for the vacuum in \eqref{comp-e1}. }
\begin{align}
&-\int_{-\infty}^{0} e^{|k|\sqrt{1-\nu^2\tau^2}Y_3}\sL^-(k) \cdot (\widehat{ \partial_{\theta} \varphi^2  \sA^- \partial_{Y_3} V^1  }) (k)|_{y_3=0}  dY_3 \notag\\
&\quad = -\int_{-\infty}^{0} e^{|k|\sqrt{1-\nu^2\tau^2}Y_3}\sL^-(k)  \notag\\ 
&\quad \quad \cdot \sum\limits_{k_1+k_2 = k} ik_2\widehat{\varphi}^2(k_2) 
\widetilde{\gamma}^1_j |k_1|^2 \sqrt{1-\nu^2\tau^2} e^{|k_1|\sqrt{1-\nu^2\tau^2}Y_3}\widehat{\varphi}^2 (k_1)  \sA^- \sR^-_j(k_1) dY_3  \notag\\
&\quad = -\sum\limits_{k_1+k_2 = k} \frac{ik_1^2k_2}{|k|+|k_1|} \widetilde{\gamma}^1_j \big( \sL^-(k) \cdot \sA^- \sR^-_j(k_1)
\big) \widehat{\varphi}^2(k_1)\widehat{\varphi}^2(k_2)  \notag\\
&\quad = -i\tau (1-\nu^2\tau^2) (a^-_j \widetilde{\gamma}^1_j ) \sgnk \sum\limits_{k_1+k_2 = k} \Big(
\frac{\sgnk -\sgnkone}{2}\frac{k_1^2k_2}{|k|+|k_1|}   \notag\\
&\quad \quad +\frac{\sgnk- \sgnktwo}{2}\frac{k_1k_2^2}{|k|+|k_2|}\Big) \widehat{\varphi}^2(k_1)\widehat{\varphi}^2(k_2). 
\label{front-5}
\end{align}

Before adding up all the integrations, we first deal with the quadratic terms appearing in \cref{front-3} and \cref{front-5}.
Let us define 
\begin{equation}\label{def-Lambda}
\begin{aligned}
\Lambda(k_1,k_2)=& \sgnk  \frac{k_1|k_2|+k_2|k_1|}{2 } + \frac{1-\sgnkone\sgnktwo}{2} |k_1||k_2|\\
&\quad -\Big(
\frac{\sgnk -\sgnkone}{2}\frac{k_1^2k_2}{|k|+|k_1|}    +\frac{\sgnk- \sgnktwo}{2}\frac{k_1k_2^2}{|k|+|k_2|}\Big),
\end{aligned}
\end{equation}
where $k=k_1+k_2$. The kernel $\Lambda$ satisfies the following properties:
\begin{align}
&\Lambda(k_1,k_2)=\Lambda(k_2,k_1)            &\text{(symmetry)}, \notag\\
&\Lambda(k_1,k_2)=\overline{\Lambda(-k_1,-k_2)}            &\text{(reality)}, \label{prop-Lambda}\\
&\Lambda(\alpha k_1,\alpha k_2)=\alpha^2\Lambda(k_1,k_2)   \quad \forall \alpha>0        &\text{(homogeneity)}. \notag
\end{align}
With these properties, to simplify $\Lambda$ it suffices to consider the two particular cases:
\begin{equation*}
\Lambda(k_1,k_2)=
\begin{cases}
k_1k_2 &\text{if~}k_1,k_2>0,\\
-k_1k_2-k_2^2 &\text{if~}k_1+k_2>0,k_2<0.
\end{cases}
\end{equation*}
Therefore, one has 
\begin{equation}\label{def-Lambda2}
\Lambda(k_1,k_2)=\frac{2|k_1+k_2||k_1||k_2|}{|k_1+k_2|+|k_1|+|k_2|}.
\end{equation}
As a result, to simplify the summation of the quadratic terms of the equation for the $\widehat{\varphi}^2(k)$, we take $(\widetilde{\gamma}^1_1, \widetilde{\gamma}^1_2)=(\iota_1, \iota_2)$ satisfying
\begin{equation}
\label{eq-gamma-special}
-\Big(
\nu^2\tau^2((\widetilde{\gamma}^1_1)^2+(\widetilde{\gamma}^1_2)^2)-(-\xi_2\widetilde{\gamma}^1_1+\xi_1\widetilde{\gamma}^1_2)^2
\Big)
=a^-_j\widetilde{\gamma}^1_j
=\frac{(a^-_1)^2+(a^-_2)^2-(b^-)^2}{1-\nu^2\tau^2},
\end{equation}
i.e. 
\begin{equation}\label{def-gamma-special}
\iota_j=\frac{\nu\tau a^-_j -\xi_j b^-}{\nu\tau(1-\nu^2\tau^2)}, \quad j=1,2.
\end{equation}
(It is worth noting that \cref{def-gamma-special} is the only solution to \cref{eq-gamma-special}.)

In the case of \cref{def-gamma-special}, summing up \cref{front-1,front-2,front-3,front-4,front-5} and dividing by $-2\tau$ gives
\begin{align}
&(c^+ + d_0) \dt \widehat{\varphi}^2(k)  +(c^+ u^0_j -b^+ B^0_j+d_j) \partial_{y_j} \widehat{\varphi}^2(k)\notag\\
&\qquad+ i \Big( (c^+)^2 -(b^+)^2 - (a_1^-)^2-(a_2^-)^2+(b^-)^2 \Big)\sgnk  \label{eq-HJ}\\
&\qquad \cdot\sum\limits_{k_1+k_2 = k}  \frac{|k||k_1||k_2|}{|k|+|k_1|+|k_2|}\widehat{\varphi}^2 (k_1) \widehat{\varphi}^2(k_2)=0\quad \text{for all~} k\in \mathbb{Z}\setminus \{0\}, \notag
\end{align}
where the constants $d_0\in \mathcal{O}(\nu), d_1, d_2 $ are defined as follows:
\begin{align*}
d_0=&\frac{\nu}{ 2(1-\nu^2\tau^2)^{3/2}  } \{\nu\tau (2-\nu^2\tau^2)((H^0_1)^2+(H^0_2)^2) +2(\xi_2H^0_1-\xi_1H^0_2)E^0_3\\
&\quad+\nu\tau (E^0_3)^2-\nu\tau (\xi_j H^0_j)^2
 \}, \\
d_1=&-\frac{1}{2(1-\nu^2\tau^2)^{3/2}}\{  \xi_1 (1-\nu^2\tau^2 +\xi_2^2)(H^0_1)^2 - 2  \xi_2 (\nu^2\tau^2-\xi_2^2)H^0_1H^0_2 + \xi_1 (\nu^2\tau^2-\xi_2^2)(H^0_2)^2 \\
&\quad+2\nu\tau \xi_1\xi_2 H^0_1E^0_3-2\nu\tau (\nu^2\tau^2 -\xi^2_2) H^0_2E^0_3-(1-2\nu^2\tau^2)\xi_1(E^0_3)^2
   \},\\
d_2=&-\frac{1}{2(1-\nu^2\tau^2)^{3/2}}\{ \xi_2 (\nu^2\tau^2-\xi_1^2)  (H^0_1)^2 - 2  \xi_1 (\nu^2\tau^2-\xi_1^2)H^0_1H^0_2 + \xi_2(1-\nu^2\tau^2 +\xi_1^2) (H^0_2)^2 \\
 &\quad+2\nu\tau (\nu^2\tau^2 -\xi^2_1) H^0_1E^0_3-2\nu\tau \xi_1\xi_2 H^0_2E^0_3-(1-2\nu^2\tau^2)\xi_2(E^0_3)^2
 \}.
\end{align*}

\subsubsection{Solving the leading amplitude equation}
Since $c^+\neq 0$ and $\nu \ll 1$, the coefficient of $\dt \widehat{\varphi}^2(k) $ in \cref{eq-HJ} is different from zero. 
Moreover, since $\tau$ is real, the coefficients of $\partial_{y_j} \widehat{\varphi}^2(k) $ are real as well.
Therefore, the first line of \cref{eq-HJ} is just a transport operator with respect to the tangential variables $y_1, y_2$. 

The kernel $\Lambda$ is exactly the kernel proposed by Hamilton \textit{et al.} \cite{Hamilton1995} to study the nonlinear Rayleigh waves (see also\cite{hunter2006}).
Such type of equation \cref{eq-HJ} also appears in (weakly) nonlinear surfaces waves on vortex sheets in magnetohydrodynamics \cite{Ali2003, Pierre2021}, and in nonlinear surfaces waves on the plasma--vacuum interface \cite{Secchi2013}.
Its local well-posedness easily follows by adapting the proof of Hunter \cite{hunter2006} to our case:
%Therefore, it follows from \cite[Theorem~4.1]{Pierre2021} and the references therein that there exists a unique local-in-time solution to \cref{eq-HJ}.

Denote $H^s_\#,~s\geq0 $ as the Sobolev space of all $g(y',\theta) \in H^s (\mathbb{T}^3 )$ with zero mean with respect to $\theta$.
\begin{Thm}[\cite{hunter2006}]
	\label{thm-HJ}
	Let $s>2$ and $g\in H^s_\#$. Then there exist $T>0$ and a unique local solution to \eqref{eq-HJ} on $[0,T]$
	$$\varphi^2 \in C([0,T];H^s_\#)\cap C^1([0,T];H^{s-1}_\#)$$
	 with $ \varphi^2 |_{t=0}=g$, where  $T=\frac{1}{K\norm{g}{H^4_\#}}$ for a suitable constant $K>0$.
	 
	 Moreover, if $g\in H^\infty_\#$, then the unique solution $\varphi^2 \in C^{\infty}([0,T];H^\infty_\#)$.
\end{Thm}
In \cite[Theorem~4.2]{hunter2006}, the constant $K$ depends on $s$, but it is shown to be independent of $s$ in \cite{Wheeler2018}. The global existence for \eqref{eq-HJ} is an open question. 
We refer the interested readers to \cite[Theorem~2.2]{Pierre2021} and its discussion for more details.

\subsubsection{Construction of the leading front profile}
We decompose the leading front profile $\varphi^2$ as $$\varphi^2=\varphi^2_\# +\widehat{\varphi}^2(0),$$ where $\varphi^2_\#$ has zero mean with respect to $\theta$. Observing that \eqref{eq-HJ} only involves the nonzero modes, 
we can obtain the oscillating part $\varphi^2_\#$ as the solution to \cref{eq-HJ} with initial condition
$$\varphi^2_\# |_{t=0} = \varphi_0^2 \,$$ given  by \cref{thm-HJ}.
Notice that  $\varphi_0^2 \in H^\infty_\#$ in the setting of initial front \cref{initial-front}. The mean
$\widehat{\varphi}^2(0)$ will be determined in the next step of the induction process.

\subsection{Summary on the construction of the leading profile}

Before we introduce the induction process for the correctors, let us make a short summary of our construction of the leading profile.
\begin{enumerate}
	\item Due to \cref{WKB-eq}, \cref{WKB-div} and \cref{WKB-jump}, the leading profile  $(U^1, V^1, \varphi^2)$ must satisfy the fast problem \cref{f-leading}, which gives the general decomposition \cref{sol-leading-general}.
	
	\item By the compatibility conditions \eqref{comp-b} and \eqref{comp-c} on the fast problem for $(U^2, V^2, \varphi^3)$,  $\underline{U}^1, \underline{V}^1$ and $\Pi \,\widehat{U}_\star^1(0)$ can be chosen as $0$.
	
	\item By the compatibility condition \eqref{comp-e} on the fast problem for $(U^2, V^2, \varphi^3)$,  $\widehat{\varphi}^2(k)$ ($k\neq 0$) can be determined by the nonlocal Hamilton-Jacobi equation \eqref{eq-HJ} when taking $(\widetilde{\gamma}^1_1, \widetilde{\gamma}^1_2)=(\iota_1, \iota_2)$ in \cref{def-gamma-special}.
	
\end{enumerate}
With the above constructions, the leading profile  $(U^1, V^1, \varphi^2)$ can be defined as 
\begin{equation}\label{sol-leading1}
	\begin{aligned}
		U^1(t,y,Y_3,\theta)=& \sum\limits_{k\in \mathbb{Z} \setminus {0} }|k|\widehat{\varphi}^2(t,y',k) \chi(y_3) e^{-|k|Y_3+ik\theta} \sR^+(k), \\
		%line 
		V^1(t,y,Y_3,\theta)=&\sum\limits_{k\in \mathbb{Z} \setminus \{0\}}
		|k| \widehat{\varphi}^2(t,y',k) \chi(y_3) e^{|k|\sqrt{1-\nu^2\tau^2}Y_3+ik\theta} \sR^-(k),
	\end{aligned}
\end{equation} 
where $\sR^-(k)=\iota_j\sR^-_j(k) (j=1,2)$ is defined in \cref{vector-R3}, 
%\todo{changed here-Apr: use new $\sR^-(k)$ here} 
and $\chi \in H^\infty (\mathbb{R})$ is a cut-off function such that $\chi=1$ near $y_3=0$. The final time $T>0$ of local existence is determined by \cref{thm-HJ}.

Up to now the compatibility conditions \eqref{comp-a} and \eqref{comp-d} of
the fast problem for $(U^2, V^2, \varphi^3)$ have not been satisfied yet. %i.e.
%\begin{align*}
%	&
%	\left\{\begin{aligned}
%		&F^{1,+}_6 |_{y_3=Y_3=0}=-b^+ \partial_{\theta}G_1^1+c^+\partial_{\theta}G_2^1,\\
%		&F^{1,-}_3 |_{y_3=Y_3=0}=\nu\tau \partial_{\theta}G_3^1+\xi_j\partial_{\theta}G_{3+j}^1,
%	\end{aligned}\right.\\
%    &
%    \left\{
%    \begin{aligned}
%    	&\partial_{Y_3} F_6^{1,+} + \xi_j \partial_{\theta} F^{1,+}_{3+j}-\tau \partial_{\theta} F^{1,+}_8=0,\\
%    	&\partial_{Y_3} F_3^{1,-} + \xi_j \partial_{\theta} F^{1,-}_{j}-\nu\tau \partial_{\theta} F^{1,-}_7=0,\\
%    	&\partial_{Y_3} F_6^{1,-} + \xi_j \partial_{\theta} F^{1,-}_{3+j}-\nu\tau \partial_{\theta} F^{1,-}_8=0.
%    \end{aligned}
%    \right.
%\end{align*}
Actually, we will see that they automatically hold true in the induction process, see the proofs of \cref{lem-comp-am} and \cref{lem-comp-dm}.

Even though all the compatibility conditions \eqref{f-comp} on the fast problem \cref{f-corrector1} for $(U^2, V^2, \varphi^2)$ are satisfied, the mean $\widehat{\varphi}^2(0)$ has not been fixed yet (in \cref{eq-HJ} only $\widehat{\varphi}^2(k)$, $k\neq 0$ are concerned). 
$\widehat{\varphi}^2(0)$ will be determined in the next steps of the induction process.

%%%%%%%%%%%%%%%%%%%%%%%%%%%%%%%%%%%%%%%%%%%%%%%%%%%%%%%%%%%%%%%%
%%%%%Solving the WKB cascade: the correctors  #5
%%%%%%%%%%%%%%%%%%%%%%%%%%%%%%%%%%%%%%%%%%%%%%%%%%%%%%%%%%%%%%%%
\vspace{0.5cm}
\section{Solving the WKB cascade: the correctors} 
\label{sec-correctors}
This section is devoted to the induction process for the WKB cascade, 
i.e. the iterative construction of the corrector $(U^{m+1}, V^{m+1}, \varphi^{m+2})$, satisfying a collection of properties H(m+1), for given profiles $(U^{1}, V^{1}, \varphi^{2})$, $(U^{2}, V^{2}, \varphi^{3})$, $\cdots$, $(U^{m}, V^{m}, \varphi^{m+1})$ that satisfy the induction assumption H(m).

The induction assumption H(m) is as follows: there exists a positive time $T$ (independent of $m$) and profiles $(U^1,V^1,\varphi^2)$, $(U^2,V^2,\varphi^3)$,$\cdots$, $(U^{m-1},V^{m-1},\varphi^{m})$, $(U^m,V^m,\varphi^{m+1}_{\#})$ such that the following five properties are satisfied:
\begin{align} 
&(U^1,V^1,\varphi^2),\cdots, (U^{m-1},V^{m-1},\varphi^{m}),(U^m,V^m,\varphi^{m+1}_{\#}) \notag\\
& \quad 
\in S^+ \times S^- \times H^{\infty}([0,T]\times\mathbb{T}^2\times\mathbb{T}), 
\tag{H(m)-1}\label{H(m)-1}
\end{align}
\begin{align}
& \left\{
\begin{aligned}
&\sL^+_f U^{\mu} =F^{\mu-1,+}, &y\in\Omega_0^+,Y_3>0,\\
&\partial_{Y_3}B_3^{\mu}+\xi_j \partial_{\theta} B_j^{\mu} =F^{\mu-1,+}_8,&y\in\Omega_0^+,Y_3>0,\\
&\sL^-_f V^{\mu} =F^{\mu-1,-},&y\in\Omega_0^-,Y_3<0,\\
&\partial_{Y_3}H_3^{\mu}+\xi_j \partial_{\theta}H_j^{\mu} =F^{\mu-1,-}_7,&y\in\Omega_0^-,Y_3<0,\\
&\partial_{Y_3}E_3^{\mu}+\xi_j \partial_{\theta}E_j^{\mu} =F^{\mu-1,-}_8,&y\in\Omega_0^-,Y_3<0,\\
&\mathbb{B}^+ U^{\mu}  |_{y_3=Y_3=0}+\mathbb{B}^- V^{\mu}  |_{y_3=Y_3=0}+\underline{b} \partial_{\theta} \varphi^{\mu+1} =G^{\mu-1},
\end{aligned}
\right.
\label{fm}
\tag{H(m)-2}\\
&\text{for~} \mu=1,\cdots,m, \notag
\end{align}
\begin{equation}
\underline{\widehat{F}}^{m,\pm}(t,y,0)=0, \underline{\widehat{F}}_7^{m,-}(t,y,0)= \underline{\widehat{F}}_8^{m,\pm}(t,y,0)=0.
\label{comp-bm}
\tag{H(m)-3}
\end{equation}
\begin{equation}
\left\{
\begin{aligned}
&u_j^0 \widehat{F}^{m,+}_{7,\star}(0) -B_j^0 \widehat{F}^{m,+}_{8,\star}(0)=\widehat{F}^{m,+}_{j,\star} (0),\quad j=1,2,\\
&B_j^0 \widehat{F}^{m,+}_{7,\star}(0) -u_j^0 \widehat{F}^{m,+}_{8,\star}(0)=\widehat{F}^{m,+}_{3+j,\star} (0),\quad j=1,2,
%	&\widehat{F}^-_{3,\star} (0)=0.
\end{aligned}
\right.
\label{comp-cm}
\tag{H(m)-4}
\end{equation}
\begin{equation}
\begin{aligned}
&\quad \int_{0}^{+\infty}e^{-|k|Y_3}\sL^+(k) \cdot \widehat{F}^{m,+}(t,y',0,Y_3,k) dY_3 + l^+_j \widehat{G}_j^{m}(t,y',k)\\
& -\int_{-\infty}^{0} e^{|k|\sqrt{1-\nu^2\tau^2}Y_3}\sL^-(k) \cdot \widehat{F}^{m,-}(t,y',0,Y_3,k)  dY_3 -\frac{-\nu\tau a^-_j +\xi_j b^-}{\nu\sqrt{1-\nu^2\tau^2}}  \widehat{G}_{3+j}^{m}(t,y',k) \\
& -i\sgnk \tau \widehat{G}_6^{m}(t,y',k) +\frac{E^0_3}{|k|\nu} \widehat{F}^{m,-}_6(t,y',0,0,k)=0, \quad \text{for all~} k\in \mathbb{Z}\setminus \{0\},
\end{aligned}\label{comp-em}
\tag{H(m)-5}
\end{equation}
where $S^\pm$ are defined in \cref{def-space}, the source terms $F^{\nu-1,\pm}, F^{\nu-1,-}_7, F^{\nu-1,\pm}_8$ and the boundary terms $G^{\nu-1}$ are given in \cref{Fm+}, \cref{Fm-}, \cref{F8m} and \eqref{Gm}, respectively.

Note that the time $T$ in the induction assumption H(m) is actually determined by the solvability of the leading front to \cref{eq-HJ}. Thus it is independent of m.
This is because except the equation for the leading front, all other partial differential equations in our constructions are linear, for which the lifespan is not an issue.
Moreover, in H(m) the mean of $\varphi^{m+1}$ with respect to $\theta$ is not assigned, because it does not enter in the definition of the source terms $F^{m,\pm}, F^{m,-}_7, F^{m,\pm}_8, G^m$.
The mean will be determined in order to fulfill (H(m+1)-3).

\subsection{The initial step of the induction}\label{sec-initial step}
We first verify H(1) by the results obtained in the last section for the leading profiles.

From \cref{thm-HJ}, the leading front $\varphi^2_\# \in H^\infty([0,T]\times\mathbb{T}^2 \times\mathbb{T} )$ for some time $T>0$.
Our constructions \cref{sol-leading1} ensures that $(U^1, V^1) \in S^+ \times S^-$, and the leading front satisfy \cref{f-leading}, which show that (H(1)-1) and (H(1)-2) hold true.
During our construction of the leading profiles, condition \cref{comp-b1}, the solvability condition \eqref{comp-c} on problem \cref{f-corrector1}, and \cref{comp-e1} have been enforced, and they correspond to (H(1)-3), (H(1)-4) and (H(1)-5), respectively. 
Thus, we have shown that there exists a positive time $T>$ such that H(1) is satisfied.

\vspace{0.3cm}

In the subsequent part of this section, given H(m) for some time $T>0$, we will show that H(m+1) is satisfied for the same time $T$.
We will construct $(U^{m+1}, V^{m+1}, \varphi^{m+2}_\sharp)$ by its fast problem, i.e. $\mu=m+1$ in (H(m+1)-2): 
\begin{align}
	& \left\{
	\begin{aligned}
		&\sL^+_f U^{m+1} =F^{m,+}, &y\in\Omega_0^+,Y_3>0,\\
		&\partial_{Y_3}B_3^{m+1}+\xi_j \partial_{\theta} B_j^{m+1} =F^{m,+}_8,&y\in\Omega_0^+,Y_3>0,\\
		&\sL^-_f V^{m+1} =F^{m,-},&y\in\Omega_0^-,Y_3<0,\\
		&\partial_{Y_3}H_3^{m+1}+\xi_j \partial_{\theta}H_j^{m+1} =F^{m,-}_7,&y\in\Omega_0^-,Y_3<0,\\
		&\partial_{Y_3}E_3^{m+1}+\xi_j \partial_{\theta}E_j^{m+1} =F^{m,-}_8,&y\in\Omega_0^-,Y_3<0,\\
		&\mathbb{B}^+ U^{m+1}  |_{y_3=Y_3=0}+\mathbb{B}^- V^{m+1}  |_{y_3=Y_3=0}+\underline{b} \partial_{\theta} \varphi^{m+2} =G^{m}.
	\end{aligned}
	\right.
	\label{fm+1},
\end{align}
We need to verify the solvability conditions \eqref{f-comp} of \cref{fm+1}.
In fact \eqref{comp-bm}, \eqref{comp-cm} and \eqref{comp-em} are exactly the solvability conditions \eqref{comp-b}, \eqref{comp-c} and \eqref{comp-e} of \eqref{fm+1},
while \eqref{comp-a} and \eqref{comp-d} of \eqref{fm+1} will be automatically fulfilled by the assumption H(m) in our construction, see \cref{lem-comp-am} and \cref{lem-comp-dm}.

\subsection{Preliminaries}

The following lemmata show that H(m) yields the compatibility conditions \eqref{comp-a}, \eqref{comp-d} of problem \eqref{fm+1}.

\begin{Lem}\label{lem-comp-am}
Under the induction assumptions H(m), it holds 
\begin{equation}
\left\{\begin{aligned}
&F^{m,+}_6 |_{y_3=Y_3=0}=-b^+ \partial_{\theta}G_1^{m}+c^+\partial_{\theta}G_2^{m},\\
&F^{m,-}_3 |_{y_3=Y_3=0}=\nu\tau \partial_{\theta}G_3^{m}+\xi_j\partial_{\theta}G_{3+j}^{m}.
\end{aligned}\right. \label{comp-am}
\end{equation}
\end{Lem}
The proof of \cref{lem-comp-am} is postponed to \cref{appen-proofs}.

\begin{Lem}\label{lem-comp-dm}
	Under the induction assumptions H(m), it holds 
	\begin{equation}
	\left\{
	\begin{aligned}
	&\partial_{Y_3} F_6^{m,+} + \xi_j \partial_{\theta} F_{3+j}^{m,+}-\tau \partial_{\theta} F_8^{m,+}=0,\\
	&\partial_{Y_3} F_3^{m,-} + \xi_j \partial_{\theta} F_{j}^{m,-}-\nu\tau \partial_{\theta} F_7^{m,-}=0,\\
	&\partial_{Y_3} F_6^{m,-} + \xi_j \partial_{\theta} F_{3+j}^{m,-}-\nu\tau \partial_{\theta} F_8^{m,-}=0.
	\end{aligned}
	\right.
	\label{comp-dm}
	\end{equation}
\end{Lem}
The proof of \cref{lem-comp-dm} is postponed to \cref{appen-proofs}.

\vspace{0.3cm}

%%%%%%%%%%%% slow problem %%%%%%%%%%%%%%%%

Given the five compatibility conditions \eqref{comp-bm}, \eqref{comp-cm}, \eqref{comp-em}, \eqref{comp-am} and \eqref{comp-dm}, by \cref{thm-fast} and according to \cref{rem-decomp}, there is a solution to the fast problem \cref{fm+1} in the form
\begin{align}
U^{m+1}(t,y,Y_3,\theta)=&\mathbb{U}_{\#}^{m+1}(t,y,Y_3,\theta)+(0,0,\widehat{\sU}_{3,\star}^{m+1}(0),0,0,\widehat{\sB}_{3,\star}^{m+1}(0),\widehat{\sQ}_{\star}^{m+1}(0))^{\top}(t,y,Y_3) \notag \\
&  +(\widehat{u}_{1,\star}^{m+1}(0),\widehat{u}_{2,\star}^{m+1}(0),0,\widehat{B}_{1,\star}^{m+1}(0),\widehat{B}_{2,\star}^{m+1}(0),0,0)^{\top} (t,y,Y_3)  \notag\\
& 
+\underline{\widehat{U}}^{m+1}(0)(t,y)+ \sum\limits_{k\in \mathbb{Z} \setminus {0} }|k|\widehat{\varphi}^{m+2}(t,y',k) \chi(y_3) e^{-|k|Y_3+ik\theta} \sR^+(k), \notag \\
%line 
V^{m+1}(t,y,Y_3,\theta)=&\mathbb{V}_{\#}^{m+1}(t,y,Y_3,\theta)+\widehat{\mathbb{V}}_{\star}^{m+1}(0)(t,y,Y_3)
+\underline{\widehat{V}}^{m+1}(0)(t,y)\notag \\
& 
+\sum\limits_{k\in \mathbb{Z} \setminus \{0\}}\widetilde{\gamma}^{m+1}_j |k| \widehat{\varphi}^{m+2}(t,y',k) \chi(y_3) e^{|k|\sqrt{1-\nu^2\tau^2}Y_3+ik\theta} \sR^-_j(k). \label{decomp-fm+1}
\end{align}
The functions $\mathbb{U}_{\#}^{m+1}$, $ (0,0,\widehat{\sU}_{3,\star}^{m+1}(0),0,0,\widehat{\sB}_{3,\star}^{m+1}(0),\widehat{\sQ}_{\star}^{m+1}(0))$, $\mathbb{V}_{\#}^{m+1}$ and $\widehat{\mathbb{V}}_{\star}^{m+1}(0)$ are determined by the source terms $F^{m,\pm}, G^{m}$ (see \cref{def-0-Ustar} and \cref{def-Usharp}).
%$\underline{\widehat{U}}^{m+1}$ (or $\underline{\widehat{V}}^{m+1}$ ) just have constrains on the boundary.
%$\gamma^{m+1,+}$ and $\gamma^{m+1,-}_j$ are determined by $\varphi_{\#}^{m+2}$, but have one freedom from the coefficients $\widetilde{\gamma}_j^{m+1}$.  
Moreover, $\sR^+(k)$, $\sR^-_1(k)$, $\sR^-_2(k)$ are explicitly given in \cref{vector-R},
$\chi \in H^\infty (\mathbb{R})$ is an arbitrary cut-off function such that $\chi=1$ near $y_3=0$,
and the coefficients $\widetilde{\gamma}^{m+1}_j$ ($j=1,2$) satisfy $a^-_j\widetilde{\gamma}^{m+1}_j =\frac{(a_1^-)^2+(a_2^-)^2-(b^-)^2}{1-\nu^2\tau^2}$.

To complete the construction of \cref{decomp-fm+1}, it remains to 
 determine the fast means of the tangential components $\widehat{u}_{1,\star}^{m+1}(0)$, $\widehat{u}_{2,\star}^{m+1}(0)$, $ \widehat{B}_{1,\star}^{m+1}(0)$, $\widehat{B}_{2,\star}^{m+1}(0)$, the slow means $\underline{\widehat{U}}^{m+1}(0)$, $\underline{\widehat{V}}^{m+1}(0)$  
and $\varphi_{\#}^{m+2}$. To complete H(m+1) we also need to determine $\widehat{\varphi}^{m+1}(0)$.
All these functions will be found by enforcing H((m+1)-3), H((m+1)-4), H((m+1)-5) as follows.

\subsection{The slow means}\label{sec-slow means}
In this part, we will construct  $\underline{\widehat{U}}^{m+1}(0)$, $\underline{\widehat{V}}^{m+1}(0)$ and
$\widehat{\varphi}^{m+1}(0)$ by enforcing H((m+1)-3).

\subsubsection{Deriving the equations}
Recalling the general definitions \cref{Fm+}, \cref{Fm-}, \cref{F8m} of ${F}^{m,\pm}, {F}^{m,+}_8$, we obtain that the residual components of ${F}^{m+1,\pm}$ are given by
\begin{align*}
\underline{F}^{m+1,+} =& - L^+_s \underline{U}^{m+1} 
+  \sum\limits_{\substack{l_1+l_2=m+2\\l_1\geq 1}} \partial_{\theta} \varphi^{l_2} \sA^+ \partial_{y_3} \underline{U}^{l_1}
+\sum\limits_{\substack{l_1+l_2=m+1\\l_1\geq 1}} (\dt \varphi^{l_2} A^+_0 + \partial_{y_j} \varphi^{l_2} A^+_j) \partial_{y_3} \underline{U}^{l_1} \notag \\
&  \quad 
-\sum\limits_{\substack{l_1+l_2=m+1\\l_1,l_2\geq 1}} \mathbb{A}_{\alpha} (\underline{U}^{l_1},\partial_{y_\alpha} \underline{U}^{l_2})
-\sum\limits_{\substack{l_1+l_2=m+2\\l_1,l_2\geq 1}} \xi_j \mathbb{A}_{j} (\underline{U}^{l_1},\partial_{\theta} \underline{U}^{l_2})  
 \notag  \\
&  \quad 
+\sum\limits_{\substack{l_1+l_2+l_3=m+2\\l_1,l_2\geq 1}} \partial_{\theta} \varphi^{l_3} \xi_j \mathbb{A}_{j}  (\underline{U}^{l_1},\partial_{y_3} \underline{U}^{l_2})
+\sum\limits_{\substack{l_1+l_2+l_3=m+1\\l_1,l_2\geq 1}} \partial_{y_j} \varphi^{l_3}  \mathbb{A}_{j}  (\underline{U}^{l_1},\partial_{y_3} \underline{U}^{l_2}),  \notag	\\	 
\underline{F}^{m+1,-}
=& - L^-_s \underline{V}^{m+1} 
+  \sum\limits_{\substack{l_1+l_2=m+2\\l_1\geq 1}} \partial_{\theta} \varphi^{l_2} \sA^- \partial_{y_3} \underline{V}^{l_1}
+\sum\limits_{\substack{l_1+l_2=m+1\\l_1\geq 1}} (\dt \varphi^{l_2} A^-_0 + \partial_{y_j} \varphi^{l_2} A^-_j) \partial_{y_3} \underline{V}^{l_1}. \notag	
\end{align*}
Therefore, to ensure H((m+1)-3), that is
$\underline{\widehat{F}}^{m+1,\pm}(t,y,0)=0$, $\underline{\widehat{F}}_7^{m+1,-}(t,y,0)=0$ and $\underline{\widehat{F}}_8^{m+1,\pm}(t,y,0)=0$,
we impose the \textit{slow problem}:
\begin{equation}
\label{s}
\left\{
\begin{aligned}
&L_s^+ \underline{\widehat{U}}^{m+1}(0)=\bF^{m,+},  \nabla\cdot \underline{\widehat{B}}^{m+1}(0) = \bF_8^{m,+},  & y\in \Omega_0^+, \\
&L_s^- \underline{\widehat{V}}^{m+1}(0)=\bF^{m,-},  \nabla\cdot \underline{\widehat{H}}^{m+1}(0) = \bF_7^{m,-},
 \nabla\cdot \underline{\widehat{E}}^{m+1}(0) = \bF_8^{m,-}, & y\in \Omega_0^-,
\end{aligned}
\right.
\end{equation}
where
\begin{equation}
\label{s-source}
\begin{aligned}
\bF^{m,+}:=
&\mathbf{c_0} \Big\{ \sum\limits_{\substack{l_1+l_2=m+2\\l_1\geq 1}} \partial_{\theta} \varphi^{l_2} \sA^+ \partial_{y_3} \underline{U}^{l_1}
+\sum\limits_{\substack{l_1+l_2=m+1\\l_1\geq 1}} (\dt \varphi^{l_2} A^+_0 + \partial_{y_j} \varphi^{l_2} A^+_j) \partial_{y_3} \underline{U}^{l_1} \notag \\
&   
-\sum\limits_{\substack{l_1+l_2=m+1\\l_1,l_2\geq 1}} \mathbb{A}_{\alpha} (\underline{U}^{l_1},\partial_{y_\alpha} \underline{U}^{l_2})
-\sum\limits_{\substack{l_1+l_2=m+2\\l_1,l_2\geq 1}} \xi_j \mathbb{A}_{j} (\underline{U}^{l_1},\partial_{\theta} \underline{U}^{l_2})  
\notag  \\
&   
+\sum\limits_{\substack{l_1+l_2+l_3=m+2\\l_1,l_2\geq 1}} \partial_{\theta} \varphi^{l_3} \xi_j \mathbb{A}_{j}  (\underline{U}^{l_1},\partial_{y_3} \underline{U}^{l_2})
+\sum\limits_{\substack{l_1+l_2+l_3=m+1\\l_1,l_2\geq 1}} \partial_{y_j} \varphi^{l_3}  \mathbb{A}_{j}  (\underline{U}^{l_1},\partial_{y_3} \underline{U}^{l_2})
\Big\}, \\
\bF_8^{m,+}:=
&\mathbf{c_0} \Big\{ \sum\limits_{\substack{l_1+l_2=m+2\\l_1\geq 1}} \partial_{\theta} \varphi^{l_2} \xi_j \partial_{y_3} \underline{B}_j^{l_1}
+\sum\limits_{\substack{l_1+l_2=m+1\\l_1\geq 1}} \partial_{y_j} \varphi^{l_2}  \partial_{y_3} \underline{B}_j^{l_1}
\Big\}, \\
\bF^{m,-}:=
&\mathbf{c_0} \Big\{ \sum\limits_{\substack{l_1+l_2=m+2\\l_1\geq 1}} \partial_{\theta} \varphi^{l_2} \sA^- \partial_{y_3} \underline{V}^{l_1}
+\sum\limits_{\substack{l_1+l_2=m+1\\l_1\geq 1}} (\dt \varphi^{l_2} A^-_0 + \partial_{y_j} \varphi^{l_2} A^-_j) \partial_{y_3} \underline{V}^{l_1}
\Big\}, \\
\bF_7^{m,-}:=
&\mathbf{c_0} \Big\{ \sum\limits_{\substack{l_1+l_2=m+2\\l_1\geq 1}} \partial_{\theta} \varphi^{l_2} \xi_j \partial_{y_3} \underline{H}_j^{l_1}
+\sum\limits_{\substack{l_1+l_2=m+1\\l_1\geq 1}} \partial_{y_j} \varphi^{l_2}  \partial_{y_3} \underline{H}_j^{l_1}
\Big\}, \\
\bF_8^{m,-}:=
&\mathbf{c_0} \Big\{ \sum\limits_{\substack{l_1+l_2=m+2\\l_1\geq 1}} \partial_{\theta} \varphi^{l_2} \xi_j \partial_{y_3} \underline{E}_j^{l_1}
+\sum\limits_{\substack{l_1+l_2=m+1\\l_1\geq 1}} \partial_{y_j} \varphi^{l_2}  \partial_{y_3} \underline{E}_j^{l_1}
\Big\}.
\end{aligned}
\end{equation}
Here $\mathbf{c_0}(\cdot)$ represents the mean of "$\cdot$" with respect to $\theta$.

From the fast problem \eqref{fm+1} we see that $\underline{\widehat{U}}^{m+1}(0)$ and  $\underline{\widehat{V}}^{m+1}(0)$ have to satisfy the boundary condition:
$$
\mathbb{B}^+ {\widehat{U}}^{m+1} (0) |_{y_3=Y_3=0}+\mathbb{B}^- {\widehat{V}}^{m+1}(0)  |_{y_3=Y_3=0} ={\widehat{G}}^{m}(0),
$$
which equivalently reads:
\begin{equation}
	\label{s-boundary}
	\left\{
	\begin{aligned}
		&
		\underline{\widehat{u}}_3^{m+1}(0) |_{y_3=0}= \widehat{G}_1^{m}(0)-\widehat{u}_{3,\star}^{m+1}(0) |_{y_3=Y_3=0} =: \bG^{m}_1,  \\
		&
		\underline{\widehat{B}}_3^{m+1}(0) |_{y_3=0}= \widehat{G}_2^{m}(0)-\widehat{B}_{3,\star}^{m+1}(0) |_{y_3=Y_3=0} =: \bG^{m}_2, \\
		&
		\underline{\widehat{H}}_3^{m+1}(0) |_{y_3=0}= \widehat{G}_3^{m}(0)-\widehat{H}_{3,\star}^{m+1}(0) |_{y_3=Y_3=0} =: \bG^{m}_3, \\
		&
		\underline{\widehat{E}}_2^{m+1}(0) |_{y_3=0}= \widehat{G}_4^{m}(0)-\widehat{E}_{2,\star}^{m+1}(0) |_{y_3=Y_3=0}  =: \bG^{m}_4, \\
		&
		-\underline{\widehat{E}}_1^{m+1}(0) |_{y_3=0}= \widehat{G}_5^{m}(0)+\widehat{E}_{1,\star}^{m+1}(0) |_{y_3=Y_3=0} =: \bG^{m}_5, \\
		&
		\{\underline{\widehat{q}}^{m+1}(0)-H_j^0\underline{\widehat{H}}_j^{m+1}(0) + E_3^0 \underline{\widehat{E}}_3^{m+1}(0) \}  |_{y_3=0} \\
		&
		\qquad\qquad= \widehat{G}_6^{m}(0) -\{\widehat{q}_{\star}^{m+1}(0)-H_j^0\widehat{H}_{j,\star}^{m+1}(0) + E_3^0 \widehat{E}_{3,\star}^{m+1}(0) \}  |_{y_3=Y_3=0}
		=: \bG^{m}_6.
	\end{aligned}
	\right.
\end{equation}
Then $\underline{\widehat{U}}^{m+1}(0)$, $\underline{\widehat{V}}^{m+1}(0)$ will be determined as the solutions to \cref{s}, \cref{s-boundary} supplemented with initial conditions (to be chosen):
\begin{equation}
	\label{s-initial}
	(\underline{\widehat{u}}^{m+1}(0), \underline{\widehat{B}}^{m+1}(0),\underline{\widehat{H}}^{m+1}(0),\underline{\widehat{E}}^{m+1}(0)) |_{t=0}=(\mathfrak{u}_0^{m+1}, \mathfrak{B}_0^{m+1}, \mathfrak{H}_0^{m+1}, \mathfrak{E}_0^{m+1}). 
\end{equation}

\subsubsection{The linearized plasma--vacuum interface problem}
Let us consider the linearized plasma--vacuum interface problem \cref{s}, \cref{s-boundary}, \cref{s-initial} for the unknowns $U=(u,B,q), V=(H,E)$, where for the sake of shortness $U$ stands for $\underline{\widehat{U}}^{m+1}(0)$ and $V$ stands for $\underline{\widehat{V}}^{m+1}(0)$. For shortness we also omit the index $m$ in the source terms $\bF^{m,\pm}$, $\bG^m_i$. Written in more explicit form the problem reads:
\begin{align}
	&\label{s*}
	\left\{
	\begin{aligned}
		&(\dt+u^0_j\partial_{y_j})u_\alpha+u^0_\alpha\nabla\cdot u-B^0_j\partial_{y_j}B_\alpha -B^0_\alpha\nabla\cdot B +\partial_{x_\alpha}q =\bF^{+}_\alpha,&&\quad\alpha=1,2,3, \\ &(\dt+u^0_j\partial_{y_j})B_\alpha-u^0_\alpha\nabla\cdot B-B^0_j\partial_{y_j}u_\alpha +B^0_\alpha\nabla\cdot u  =\bF^{+}_{3+\alpha},&&\quad\alpha=1,2,3, \\ &\nabla\cdot u = \bF_7^{+}, \quad \nabla\cdot B = \bF_8^{+},  && \quad y\in \Omega_0^+,\\
		&(\nu\dt H + \nabla \times E)_\alpha=\bF^{-}_\alpha,&&\quad\alpha=1,2,3, \\
		&(\nu\dt E - \nabla \times H)_\alpha=\bF^{-}_{3+\alpha},&&\quad\alpha=1,2,3, \\
		& \nabla\cdot H = \bF_7^{-},    
		\quad \nabla\cdot E = \bF_8^{-}, && \quad y\in \Omega_0^-,
	\end{aligned}
	\right.\\
	&\label{s-boundary*}
	\left\{
	\begin{aligned}
		&	u_3|_{y_3=0}=  \bG_1\,,  	B_3|_{y_3=0}=  \bG_2\,,    \\
		&		H_3|_{y_3=0}=  \bG_3\,,E_2|_{y_3=0}=  \bG_4\,,  	-E_1|_{y_3=0}=  \bG_5,  \\
		&   (q-H_j^0H_j + E_3^0 E_3)|_{y_3=0} = \bG_6,
	\end{aligned}
	\right.\\
	&\label{s-initial*}\quad
	(u,B,E,H) |_{t=0}=(u_0,B_0,E_0,H_0). 
\end{align}
The source terms $\bF^{\pm}, \dots, \bF_8^{\pm}$ in the differential equations are completely determined by previous steps in the induction process, while the boundary terms $\bG_1,\dots,\bG_5$ also depend on the unknown slow mean $\widehat{\varphi}^{m+1}(0)$.
The source terms in the equations, the boundary source terms and the initial data are related each other by several compatibility conditions that will be detailed in the following. 
From the analysis of the problem it turns out that there is one more boundary condition than needed for the resolution of \cref{s*}, \cref{s-boundary*}. 
Without the last boundary condition of \cref{s-boundary*}, the equations in the plasma or vacuum region can be solved individually.
This additional condition will be used for the determination of $\widehat{\varphi}^{m+1}(0)$.

\vspace{0.3cm}

\subsubsection{Determining the slow mean of the front profile}
Let us first consider the total pressure $q$ %$q=\underline{\widehat{q}}^{m+1}(0)$ 
in plasma. By applying the divergence operator to the equation for the velocity in \cref{s*} and substituting the divergence constraints for $u,B$, it follows that $q$ solves
\begin{align}
	\label{s-q}
	&\left\{ 
	\begin{aligned}
		-\Delta q&=  -\partial_{y_{\alpha}} \bF^+_{\alpha}+(\dt+2u^0_j\partial_{y_j}) \bF^+_7-2B^0_j \partial_{y_j} \bF^+_8,  \qquad y\in \Omega_0^+, \\
		\partial_{y_3}q |_{y_3=0} &=  \bF^+_{3}  |_{y_3=0}  -(\dt+u^0_j\partial_{y_j}) \bG_1+B^0_j \partial_{y_j} \bG_2,   
	\end{aligned}
	\right.
\end{align}
where the Neumann condition follows from the equation for $\alpha=3$ restricted at the boundary.
Let us define the vacuum total pressure 
\begin{equation}\label{def-tilde-q}
	\widetilde{q}:= H_j^0 H_j -E_3^0 E_3.
\end{equation}
%(recall $H_j=\underline{\widehat{H}}_j^{m+1}(0),  E_3= \underline{\widehat{E}}_3^{m+1}(0)$).
By classical manipulation of the Maxwell equations in \cref{s*}, it is easily seen that $H_j, E_3$ solve the wave equation with suitable source terms. It follows from \cref{def-tilde-q} that 
\begin{align}\label{tildeq}
	\begin{aligned}
		(\nu^2 \dt^2-\Delta) \widetilde{q} &= H^0_1 \partial_{y_3} \bF^{-}_5-H^0_2 \partial_{y_3} \bF^{-}_4 +E^0_3 \partial_{y_3} \bF^{-}_8 \\
		& \qquad +(-E^0_3 \nu\dt+ H^0_2\partial_{y_1}-H^0_1 \partial_{y_2}  ) \bF^{-}_6 -H^0_j \partial_{y_j} \bF^{-}_7 \\
		&\qquad + (H^0_1 \nu \dt +E^0_3 \partial_{y_2} ) \bF^{-}_1
		+(H^0_2 \nu \dt -E^0_3 \partial_{y_1} ) \bF^{-}_2   .
	\end{aligned}
\end{align}
Let us compute the normal derivative of $\widetilde{q}$ at $\Gamma_0$.
From the evolution equation for $H$ and the divergence constraint for $E$ in \cref{s*vac} it follows that
\begin{align*}
	&\partial_{y_3}H_1=-\bF_5^- +\nu\dt E_2+\partial_{y_1}H_3,\quad
	\partial_{y_3}H_2=\bF_4^- -\nu\dt E_1+\partial_{y_2}H_3,
	\\
	&\partial_{y_3}E_3=\bF_8^- -\partial_{y_1} E_1-\partial_{y_2}E_2.
\end{align*}
Restricting the above equations at the boundary and substituting the boundary values in \cref{s-boundary*} gives
\begin{align}\label{s-tildeq-N}
	\begin{aligned}
		\partial_{y_3}\widetilde{q}|_{y_3=0}&=(H_2^0 \bF_4^- - H_1^0 \bF_5^-- E_3^0 \bF_8^-)|_{y_3=0} 
		+H_j^0 \partial_{y_j} \bG_3
		\\ 
		&+(H_1^0 \nu\dt+ E_3^0 \partial_{y_2})\bG_4 + (H_2^0 \nu\dt- E_3^0 \partial_{y_1})\bG_5.
	\end{aligned}
\end{align}
Summarizing, from \cref{s-boundary*}, \cref{tildeq}, \cref{s-tildeq-N} the vacuum total pressure $\widetilde{q}$ solves the problem
\begin{align}\label{s-tildeq}
	&\left\{
	\begin{aligned}
		(\nu^2 \dt^2-\Delta) \widetilde{q} &= H^0_1 \partial_{y_3} \bF^{-}_5-H^0_2 \partial_{y_3} \bF^{-}_4 +E^0_3 \partial_{y_3} \bF^{-}_8 \\
		& \qquad +(-E^0_3 \nu\dt+ H^0_2\partial_{y_1}-H^0_1 \partial_{y_2}  ) \bF^{-}_6 -H^0_j \partial_{y_j} \bF^{-}_7 \\
		&\qquad + (H^0_1 \nu \dt +E^0_3 \partial_{y_2} ) \bF^{-}_1
		+(H^0_2 \nu \dt -E^0_3 \partial_{y_1} ) \bF^{-}_2,   \\
		\partial_{y_3}\widetilde{q}|_{y_3=0}&=(H_2^0 \bF_4^- - H_1^0 \bF_5^-- E_3^0 \bF_8^-)|_{y_3=0} 
		+H_j^0 \partial_{y_j} \bG_3
		\\ 
		&+(H_1^0 \nu\dt+ E_3^0 \partial_{y_2})\bG_4 + (H_2^0 \nu\dt- E_3^0 \partial_{y_1})\bG_5\\
		   (q-\widetilde{q})|_{y_3=0} &= \bG_6.
	\end{aligned}
	\right.
\end{align}
Notice that the problems \cref{s-q} and \cref{s-tildeq} are coupled through the last boundary condition in \cref{s-tildeq}.

As in \cite{Pierre2021} we solve the elliptic-hyperbolic problem \cref{s-q}, \cref{s-tildeq} by Fourier expansion with respect to $y'\in \mathbb{T}^2$, i.e., 
$$q =\sum\limits_{j'\in\mathbb{Z}^2} \textbf{c}_{j'}(q)(y_3) e^{ij'\cdot y'},\quad \widetilde{q} =\sum\limits_{j'\in\mathbb{Z}^2} \textbf{c}_{j'}(\widetilde{q})(y_3) e^{ij'\cdot y'}, $$
where we use  $\textbf{c}_{j'}$ to represent the $j'$ Fourier coefficient with respect to $y'$. 
For simplicity, we denote by $\sF^{m,\pm}, \sG^{m,\pm}$ the source terms and the boundary terms of \cref{s-q}, \cref{s-tildeq}, respectively. 

The calculations can be divided into two parts: the nonzero Fourier modes and the zero Fourier mode.

\vline

$\blacklozenge$ \underline{The analysis of nonzero Fourier modes}. 

For $j'\neq 0$, we look for the solution such that $\textbf{c}_{j'}(q)(y_3)\rightarrow0$ as $y_3\rightarrow +\infty$. By taking the Fourier expansion in $y'$ of \cref{s-q} and both the Fourier transform in $y'$ and the Laplace transform ($\fL$) with respect to $t$ of \cref{s-tildeq}, with Laplace dual variable $z=\gamma+ i\delta$, we obtain  
\begin{align*}
	\textbf{c}_{j'}(q)(y_3)= & -\frac{1}{|j'|} \textbf{c}_{j'}(\sG^{m,+}) e^{-|j'| y_3} 
	+\int_{0}^{\infty} \frac{\cosh (|j'|y_3)}{|j'|} e^{-|j'|z}  \textbf{c}_{j'}(\sF^{m,+})(z) dz \\
	&\quad -\int_{0}^{y_3} \frac{\sinh (|j'|(y_3-z))}{|j'|}  \textbf{c}_{j'}(\sF^{m,+})(z) dz,\\
	\fL (\textbf{c}_{j'}(\widetilde{q}))(y_3)=& \frac{1}{\sqrt{\nu^2 z^2 +|j'|^2}} \fL\textbf{c}_{j'}(\sG^{m,-})  e^{(\sqrt{\nu^2 z^2 +|j'|^2}y_3)} \\
	&\quad 
	+\int_{-\infty}^{0} \frac{ \cosh (\sqrt{\nu^2 z^2 +|j'|^2}y_3)    }{\sqrt{\nu^2 z^2 +|j'|^2} } e^{ \sqrt{\nu^2 z^2 +|j'|^2} \tau  }\\
	&\quad\quad\quad\quad \cdot
	\Big( \fL\textbf{c}_{j'}(\sF^{m,-} ) +\nu^2 z  \textbf{c}_{j'}(\widetilde{q})|_{t=0} +\nu^2\partial_t \textbf{c}_{j'}(\widetilde{q})|_{t=0}  \Big)d \tau\\
	&\quad + \int_{y_3}^{0} \frac{ \sinh (\sqrt{\nu^2 z^2 +|j'|^2}(y_3-\tau))    }{\sqrt{\nu^2 z^2 +|j'|^2} } \\
	&\quad\quad\quad\quad \cdot
	\Big(\fL \textbf{c}_{j'}(\sF^{m,-} )+\nu^2 z  \textbf{c}_{j'}(\widetilde{q})|_{t=0} +\nu^2\partial_t \textbf{c}_{j'}(\widetilde{q})|_{t=0}   \Big) d \tau.
\end{align*}
Therefore, the pressure condition $q^+-q^-=\sG_6^m:=\bG_6$ on the boundary gives
\begin{align}
	&\fL(\textbf{c}_{j'}(q-\widetilde{q}))|_{y_3=0}=	-\frac{\fL\textbf{c}_{j'}(\sG^{m,+}) }{|j'|} 
	+\int_{0}^{\infty} \frac{e^{-|j'|\tau}}{|j'|}   \fL\textbf{c}_{j'}(\sF^{m,+})  (\tau) d\tau \notag\\
	&\qquad\qquad
	-\frac{\fL\textbf{c}_{j'}(\sG^{m,-}) }{\sqrt{\nu^2 z^2 +|j'|^2}}   
	-\int_{-\infty}^{0} \frac{   e^{ \sqrt{\nu^2 z^2 +|j'|^2} \tau  }  }{\sqrt{\nu^2 z^2 +|j'|^2} }
	\Big( \fL\textbf{c}_{j'}(\sF^{m,-} ) \notag\\
	&\qquad\qquad +\nu^2 z  \textbf{c}_{j'}(\widetilde{q})|_{t=0} +\nu^2\partial_t \textbf{c}_{j'}(\widetilde{q})|_{t=0}  \Big)(\tau) d \tau \label{s-phi}\\
	&\qquad\qquad\qquad\qquad
	=\fL\textbf{c}_{j'}(\sG_6^m) . \notag
\end{align}
The initial values $\widetilde{q}|_{t=0}$ and $\partial_t\widetilde{q}|_{t=0}$ are determined in terms of $(\underline{\widehat{H}}^{m+1}(0),\underline{\widehat{E}}^{m+1}(0))|_{t=0}
$ and will be fixed in \cref{s-initial-values}.
After collecting the front profile $\widehat{\varphi}^{m+1}(0)$ from $\fL\textbf{c}_{j'}(\sG^{m,\pm})$, one can rewrite \cref{s-phi} in the equivalent form
\begin{align}
	&\Delta(z, j') \fL\textbf{c}_{j'} (\widehat{\varphi}^{m+1}(0))= \fL\textbf{c}_{j'} (\mathfrak{A}^m)  \label{s-phi*}
\end{align}
where 
\begin{align}
	\Delta(z, j'):= &\sqrt{\nu^2 z^2 +|j'|^2}	((z+iu_j^0 j_j')^2 -(iB_j^0 j_j')^2) \notag\\
	&\quad + |j'| ((H_1^0\nu z +iE_3^0j_2')^2+(H_2^0\nu z -iE_3^0j_1')^2-(iH_j^0 j_j')^2)  \label{symbol}
\end{align}
is the symbol of the paradifferential operator $\text{Op}(\Delta)$, and the source term 
$$\mathfrak{A}^m =\sum\limits_{j'\in\mathbb{Z}^2\setminus\{0\}} \textbf{c}_{j'}(\mathfrak{A}^m) e^{ij'\cdot y'},$$
is determined by all previous profiles, and where  $\mathfrak{A}^m \in H^{\infty} ( [0,T]\times \mathbb{T}^2 )$ has zero mean with respect to $y'$.

Note that $\Delta(z, j')$ equals the Lopatinskii determinant in \cite[Proposition~8.1]{Morando2020a}. 
It follows from  \cite[Theorem~1]{Trakhinin2020b} and \cite[Proposition~8.1]{Morando2020a} that 
if $\nu\ll1$ and the stability condition \eqref{H1} holds, then $\Delta(z, j')=0$ has no roots with $\gamma=\mathfrak{R}z>0$ and only simple roots $(z_0,j'_0)\in\Sigma:=\{(z,j')\in \mathbb{C}\times \R^2; \mathfrak{R}z\geq0, |z|^2+|j'|^2=1\}$ with $\gamma_0=\mathfrak{R}z_0=0$. 
Therefore, by the similar arguments as in \cite{Coulombel2004a,Morando2008} on the vanishing of Lopatinskii determinant, 
it holds that for any roots $(z_0,j'_0)\in \Sigma$,
there exists a neighborhood $\mathcal{V}$ of $(z_0,j'_0)$ in $\Sigma$ and a $C^{\infty}$ function $h$ on $\mathcal{V}$ such that 
$$
\Delta(z,j')=(z-z_0) h(z,j'),\quad h(z,j')\neq 0~\quad\forall (z,j')\in \Sigma.
$$
Since $\Delta(z,j')$ is of degree 3 and thus $h$ is of degree 2, we have
\begin{equation*}
	\norm{\Phi^{m+1}}{2,\gamma}^2 \leq \frac{C}{\gamma^2} \norm{\mathfrak{A}^m}{0,\gamma}^2, 
\end{equation*}
where we have set
\begin{equation*}
\Phi^{m+1}:=\sum\limits_{j'\in\mathbb{Z}^2\setminus\{0\}} \textbf{c}_{j'}(\widehat{\varphi}^{m+1}(0)) e^{ij'\cdot y'},
\end{equation*}
and where the weighted norm $\norm{\cdot}{m,\gamma}$ is defined as 
\begin{equation*}
	\norm{v}{m,\gamma}^2:=\int_{0}^{\infty}\int_{\R^2}^{}(\gamma^2 + \delta^2 +|j'|^2 )^m| \fL \textbf{c}_{j'} (v) |^2\mathrm{d}\delta \mathrm{d}j'.
\end{equation*}
%Actually, when $\nu\ll1<\gamma$, $\Delta(z, j')$ can be controlled as 
%\begin{align*}
%		\Delta(z, j') \geq& |j'| [(z+iu_j^0 j_j')^2 +(B_j^0 j_j')^2+(H_j^0 j_j')^2-(E_3^0)^2|j'|^2] \\
%		& - C\nu |z| [|z|^2+ (|u^0|^2+|B^0|^2+|H^0|^2+|E^0|^2) |j'|^2]\\
%		&-C\nu^2 |z|^2 |H^0|^2 |j'|
%\end{align*}
From \cref{symbol} it follows that the symbol of the adjoint operator of $\text{Op}(\Delta)$ is exactly:
\begin{align*}
	\overline{\Delta(z, j')}= &\sqrt{\nu^2 \bar{z}^2 +|j'|^2}	((\bar{z}-iu_j^0 j_j')^2 -(iB_j^0 j_j')^2)  \notag\\
	&\quad + |j'| ((H_1^0\nu \bar{z} -iE_3^0j_2')^2+(H_2^0\nu \bar{z} +iE_3^0j_1')^2-(iH_j^0 j_j')^2).  
\end{align*}
Obviously, the properties of the roots of $\overline{\Delta(z, j')}=0$ are the same as for $\Delta(z, j')=0$. Thus we can easily show that 
\begin{equation}\label{dual}
	\norm{v}{0,\gamma}^2 \leq \frac{C}{\gamma^2} \norm{ \text{Op}(\Delta)^* v   }{-2,\gamma}^2 \qquad  \forall v\in L^2_{\gamma}=e^{\gamma t}L^2.
\end{equation}
Using \cref{dual} the well-posedness of \cref{s-phi*} can be derived by a duality argument and we obtain the following result.
\begin{Lem}\label{lem-s-phi}
	If $\nu\ll1$ and the stability condition \eqref{H1} holds, there exists a unique solution $\Phi^{m+1} \in H^{\infty} ( [0,T]\times \mathbb{T}^2 )$  to \cref{s-phi*} with zero mean with respect to $y'$ and 
	\begin{equation}\label{s-phi-initial}
		\Phi^{m+1}|_{t=0}=0,~\partial_t\Phi^{m+1}|_{t=0}=0.
	\end{equation}
\end{Lem}
Here we just impose the simplest initial data for $\Phi^{m+1} $, which match \cref{WKB-front}, i.e. at the initial time  $\varphi^m$ has zero mean with respect to $\theta$ for all $m\geq2$. We could also consider more general initial data.

\vspace{0.3cm}

$\blacklozenge$ \underline{The analysis of zero Fourier mode}. \\

To complete the construction of 
	$\widehat{\varphi}^{m+1}(0)$ we notice that the compatibility between the source term for the Poisson equation and the Neumann condition in \cref{s-q} gives
	\begin{align*}
		\int_{\Omega_0^+} \Big( (\dt&+2u^0_j\partial_{y_j}) \bF^{+}_7-2B^0_j \partial_{y_j} \bF^{+}_8 -\partial_{y_{\alpha}} \bF^{+}_{\alpha}\Big) dy\\
		=
		&\int_{\Gamma_0} \Big( \bF^{+}_{3}  |_{y_3=0} -(\dt+u^0_j\partial_{y_j}) \bG_1+B^0_j \partial_{y_j} \bG_2 \Big) dy',
	\end{align*}
	which can be reduced to
	\begin{equation}
		\label{s-comp-B}
		\dfrac{d\,}{dt} \Big( \int_{\Omega_0^+} \bF_7^{m,+} dy + \int_{\Gamma_0} \bG_1^{m} dy' \Big) =0.
	\end{equation}
	As shown in \cite{Pierre2021}, see also \cref{comp-source}, \cref{s-comp-B} holds if and only if
	\begin{equation}
		\label{s-comp-phi}
	\dfrac{d^2}{dt^2}	\mathbf{c_0}(\widehat{\varphi}^{m+1}(0))=\dfrac{d^2}{dt^2} \int_{\Gamma_0}  \widehat{\varphi}^{m+1}(0) dy' =0.
	\end{equation}
Here $\mathbf{c_0}$ represents the zero Fourier coefficient with respect to $y'$.
We also recall from \cref{WKB-front}, that at the initial time $\varphi^m$ has zero mean with respect to $\theta$ for any $m\geq2$.
In the subsequent contents, to ensure compatibility for the divergence of the initial data \cref{s-initial-u}, the time derivative of the slow mean of $\varphi^m$ also has to be of zero mean with respect to $y'$ initially, see \cref{s-comp-phit}.
Therefore, we have to impose that 
\begin{equation}\label{s-phi0}
	\mathbf{c_0}(\widehat{\varphi}^{m+1}(0)):=\int_{\Gamma_0\times \mathbb{T} } \varphi^{m+1} (t,y',\theta) dy'd\theta =0, \quad \text{for all~} t\geq0.
\end{equation}
\vspace{0.3cm}

Combining the analysis of the nonzero Fourier modes and zero Fourier mode with respect to $y'$, the construction of the front $\varphi^{m+1}$ is completed as:
\begin{equation*}\label{s-phi-final}
	\varphi^{m+1}(t,y',\theta)=\varphi^{m+1}_\# (t,y',\theta)+\Phi^{m+1}(t,y')
\end{equation*}
where $\Phi^{m+1}$ is given by \cref{lem-s-phi}.

\subsubsection{Solving the equations in the vacuum region}
We consider the Maxwell equations in the vacuum region:
\begin{align}\label{s*vac}
	&\left\{
	\begin{aligned}
		&(\nu\dt H + \nabla \times E)_\alpha=\bF^{-}_\alpha,&&\quad\alpha=1,2,3,\quad y\in \Omega_0^-, \\
		&(\nu\dt E - \nabla \times H)_\alpha=\bF^{-}_{3+\alpha},&&\quad\alpha=1,2,3,\quad y\in \Omega_0^-, \\
		& \nabla\cdot H = \bF_7^{-},    
		\quad \nabla\cdot E = \bF_8^{-}, && \quad y\in \Omega_0^-,
		\\
	&H_3|_{y_3=0}=  \bG_3,\quad  E_2|_{y_3=0}=  \bG_4,  \quad
	-E_1|_{y_3=0}=  \bG_5, 
	\\	
	&(H,E) |_{t=0}=(H_0,E_0). 
\end{aligned}
\right.
\end{align}
To state the necessary compatibility conditions between the boundary data and the initial values we recursively define $(\dt^k H_0, \dt^k E_0)$, $k\geq 1$, by formally taking $k-1$ time derivatives of the evolution equations in \cref{s*vac}, solving for $(\dt^k H, \dt^k E)$ and evaluating it at time $t=0$; for $k=0$, we set $(\dt^0 H_0, \dt^0 E_0)=(H_0, E_0)$.
\begin{Prop}[Solvability of the Maxwell equations]\label{prop-s-vac}
Let $\bF^{-},\bF_7^{-},\bF_8^{-}  \in H^{\infty}([0,T]\times\Omega_0^-)$, $\bG_3,\bG_4, \bG_5 \in H^{\infty}([0,T]\times\mathbb{T}^2)$ and $H_0,E_0 \in H^{\infty}(\Omega_0^-)$. Then \cref{s*vac} has a unique solution $(H, E) \in  H^{\infty} ([0,T]\times\Omega_0^{-})$  if and only if:
\begin{itemize}
	\item 
	$\bF_3^- |_{y_3=0}= \nu \partial_t \bG_3+  \partial_{y_j}\bG_{3+j}$ (compatibility on $\Gamma_0)$;
	\item 
	$\nu\partial_t \bF_7^-=\partial_{y_\alpha} \bF_{\alpha}^-$, $\nu\partial_t \bF_8^-=\partial_{y_\alpha} \bF_{3+\alpha}^-$ (compatibility for the divergence of the magnetic and electric fields);
	\item  $\nabla\cdot H_0 = \bF_7^- |_{t=0}$, $\nabla\cdot E_0 = \bF_8^- |_{t=0}$ in $ \Omega_0^-$(compatibility for the divergence of the initial data);
	\item 	$H_{0,3}|_{y_3=0}= \partial_t^k\bG_3|_{t=0}$, $\partial_t^k (  E_{0,2}, -E_{0,1} )|_{y_3=0}=  \partial_t^k(  \bG_4, \bG_5) |_{t=0}$ for all integers $k\geq0$ %$u_{0,3}|_{\Gamma_0}=  \bG_1 |_{t=0}$,  $B_{0,3}|_{\Gamma_0}=  \bG_2 |_{t=0}$, $H_{0,3}|_{\Gamma_0}=  \bG_3 |_{t=0}$, $E_{0,2}|_{\Gamma_0}=  \bG_4 |_{t=0}$, $-E_{0,1}|_{\Gamma_0}=  \bG_5 |_{t=0}$  
	(compatibility  of the initial data on the boundary). 
	
\end{itemize}
\end{Prop}
\begin{proof}
The necessity of conditions in \cref{prop-s-vac} follows from classical manipulations. For instance, the first one comes from the evolution equation for $H_3$ restricted on the boundary and the substitution of boundary values. The second and third ones are obtained by taking the divergence of the evolution equation of $H$ and $E$, respectively. 

Let us now assume that all conditions in \cref{prop-s-vac} are satisfied.
The boundary matrix $A^-_3$ has eigenvalues $\lambda_{1,2}=1,\lambda_{3,4}=-1$ and $\lambda_{5,6}=0$. Thus there are two incoming characteristics in $\Omega_0^-$ ($\lambda_{3,4}<0)$ and it is enough to prescribe two boundary conditions (for $E_1$ and $E_2$). Since the initial data satisfy the necessary compatibility condition at the boundary for any order, there exists a unique solution $(H, E) \in  H^{\infty} ([0,T]\times\Omega_0^{\pm})$ to the Maxwell equations, see \cite{Rauch1974}. From the equation for $H_3$ restricted to the boundary and the compatibility of the source terms on $\Gamma_0$ we get
\begin{equation*}
	\nu\dt H_3|_{y_3=0}=(\bF_3^-- \partial_{y_1}E_2+\partial_{y_2}E_1)|_{y_3=0}=\bF_3^-|_{y_3=0}-\partial_{y_j}\bG_{3+j}=\nu \partial_t \bG_3.
\end{equation*}
Then $\nu\dt (H_3|_{y_3=0}-\bG_3)=0$ and the compatibility condition on the initial value gives $H_3(t)|_{y_3=0}=\bG_3(t)$ for all $t>0$.
Applying the divergence operator to the equation for $H$ yields
\begin{equation*}
	\nu\dt \nabla\cdot H=\partial_{y_\alpha} \bF_{\alpha}^-=\nu\partial_t \bF_7^-.
\end{equation*}
Then $\nu\dt (\nabla\cdot H-\bF_7^-)=0$ and the compatibility condition on the initial divergence gives $\nabla\cdot H(t)=\bF_7^-(t)$ for all $t>0$. The proof of the constraint on the divergence of $E$ is similar. Thus $(H, E)$ is the solution to \cref{s*vac}.
\end{proof}
We postpone to \cref{lem-comp-ams} and \cref{lem-comp-dms} the proof that the source terms and boundary values of \cref{s}, \cref{s-boundary} satisfy the compatibility conditions of \cref{prop-s-vac}. We will also construct  suitable initial values $\mathfrak{H}_0^{m+1}, \mathfrak{E}_0^{m+1}$ satisfying the initial constraints in \cref{s-initial-values}. By applying \cref{prop-s-vac} it follows that there exists a unique solution $H=\underline{\widehat{H}}^{m+1}(0), E=\underline{\widehat{E}}^{m+1}(0)$ to \cref{s*vac} with source terms $\bF^{m,-},\bF_7^{m,-},\bF_8^{m,-}$, $\bG_3^{m},\bG_{3+j}^{m}$, and initial values $\mathfrak{H}_0^{m+1}, \mathfrak{E}_0^{m+1}$. For the sake of shortness, we continue to denote this solution by $H,E$.

\subsubsection{Solving the equations in the plasma region}
Now we study the equations in the plasma region:
\begin{align}
	&\label{s*plasma}
	\left\{
	\begin{aligned}
		&(\dt+u^0_j\partial_{y_j})u_\alpha+u^0_\alpha\nabla\cdot u-B^0_j\partial_{y_j}B_\alpha -B^0_\alpha\nabla\cdot B +\partial_{x_\alpha}q =\bF^{+}_\alpha,&&\quad\alpha=1,2,3, \\ &(\dt+u^0_j\partial_{y_j})B_\alpha-u^0_\alpha\nabla\cdot B-B^0_j\partial_{y_j}u_\alpha +B^0_\alpha\nabla\cdot u  =\bF^{+}_{3+\alpha},&&\quad\alpha=1,2,3, \\ &\nabla\cdot u = \bF_7^{+}, \quad \nabla\cdot B = \bF_8^{+},  && \quad y\in \Omega_0^+,
		\\
		&
		u_3|_{y_3=0}=  \bG_1,\quad  	B_3|_{y_3=0}=  \bG_2,\\
		&
		(u,B) |_{t=0}=(u_0,B_0),
	\end{aligned}
	\right.
\end{align} 
The evolution equations for $u,B$ in \cref{s*plasma} can be regarded as a hyperbolic system for $u,B$, if $\nabla q$ is a given source term.

\begin{Prop}[Solvability of the plasma equations]\label{prop-s-plasma}
	Let $\bF^{+},\bF_8^{+}  \in H^{\infty}([0,T]\times\Omega_0^+)$, $\bG_1,\bG_2 \in H^{\infty}([0,T]\times\mathbb{T}^2)$ and $u_0,B_0 \in H^{\infty}(\Omega_0^+)$. %Let ${q}\in  H^{\infty} ([0,T]\times\Omega_0^{+})$ be given by \cref{prop-q} and solution to \cref{s*total}. 
	Then \cref{s*plasma} has a unique solution $(u, B, \nabla q) \in  H^{\infty} ([0,T]\times\Omega_0^{+})$  if and only if:
	\begin{itemize}
		\item $\bF_6^+ |_{y_3=0}= (\partial_t +u_j^0 \partial_{y_j})\bG_2-B_j^0 \partial_{y_j}\bG_1$ (compatibility on $\Gamma_0)$;
		\item $\partial_t \bF_8^+=\partial_{y_\alpha} \bF_{3+\alpha}^+$ (compatibility for the divergence of the magnetic field);
		\item $\nabla\cdot u_0 = \bF_7^+ |_{t=0}$, $\nabla\cdot B_0 = \bF_8^+ |_{t=0}$ in $\Omega_0^+$ (compatibility for the divergence of the initial data); 
		\item 	$(u_{0,3}, B_{0,3})|_{y_3=0}=(\bG_1, \bG_2)|_{t=0}$
		(compatibility  of the initial data on the boundary);
	\end{itemize}
\end{Prop} 
\begin{proof}
	The necessity of conditions in \cref{prop-s-plasma} follows from manipulations as in the previous \cref{prop-s-vac}. For the sufficiency, let us assume that all conditions in the statement are satisfied. 
	
	Given $q$ we determine the velocity and magnetic field as in \cite{Pierre2021, Secchi1993}. We consider the symmetric hyperbolic system
	\begin{align}
		&\label{s*plasma-uB}
		\left\{
		\begin{aligned}
			&(\dt+u^0_j\partial_{y_j})u_\alpha-B^0_j\partial_{y_j}B_\alpha   =\bF^{+}_\alpha -u^0_\alpha\bF^{+}_7
			+B^0_\alpha\bF^{+}_8-\partial_{x_\alpha}q,&&\quad\alpha=1,2,3, 
			\\ &(\dt+u^0_j\partial_{y_j})B_\alpha-B^0_j\partial_{y_j}u_\alpha   =\bF^{+}_{3+\alpha}
			+u^0_\alpha\bF^{+}_8
			-B^0_\alpha\bF^{+}_7,&&\quad\alpha=1,2,3,
			\\
			&
			(u,B) |_{t=0}=(u_0,B_0),
		\end{aligned}
		\right.
	\end{align} 
	This system does not need any boundary condition since it only involves tangential operators (recall that $u^0_3=B^0_3=0$). By the classical theory (\cite{Friedrichs1954}), problem \cref{s*plasma-uB} admits a unique solution $(u, B) \in  H^{\infty} ([0,T]\times\Omega_0^{+})$. To verify the remaining equations in \cref{s*plasma} we proceed as in \cite{Pierre2021}, so we omit the details.
	%apply the divergence operator to \cref{s*plasma-uB} {and use \cref{s*total}} to obtain
	%\begin{align*}
	%	&
	%	\left\{
	%	\begin{aligned}
	%		&(\dt+u^0_j\partial_{y_j})(\nabla\cdot u-\bF^{+}_7)-B^0_j\partial_{y_j}(\nabla\cdot B-\bF^{+}_8)   =0, 
	%		\\ &(\dt+u^0_j\partial_{y_j})(\nabla\cdot B-\bF^{+}_8)  -B^0_j\partial_{y_j}(\nabla\cdot u-\bF^{+}_7)   =0,
	%	\end{aligned}
	%	\right.
	%\end{align*} 
	%The compatibility for the divergence of the initial data gives the constraints on the divergences in \cref{s*plasma}. Restricting to $\Gamma_0$ the equations \cref{s*plasma-uB} for $\alpha=3$ and using the compatibility condition on $\bF_6^+$ and the Neumann condition \cref{partial3-q} for $q$ we get
	%\begin{align*}
	%	&
	%	\left\{
	%	\begin{aligned}
	%		&(\dt+u^0_j\partial_{y_j})(u_3-\bG_1)-B^0_j\partial_{y_j}(B_3-\bG_2)   =0, \\ &(\dt+u^0_j\partial_{y_j})(B_3-\bG_2)-B^0_j\partial_{y_j}(u_3-\bG_1)   =0. 
	%		\end{aligned}
	%	\right.
	%\end{align*} 
	%The compatibility of the initial data on the boundary gives the boundary values in \cref{s*plasma}. 
\end{proof}
In \cref{lem-comp-ams} and \cref{lem-comp-dms} we will show that the source terms and boundary values of \cref{s}, \cref{s-boundary} satisfy the compatibility conditions of \cref{prop-s-plasma}. We will also construct  suitable initial values $\mathfrak{u}_0^{m+1}, \mathfrak{B}_0^{m+1}$ satisfying the initial constraints in \cref{s-initial-values}. By applying \cref{prop-s-plasma} it follows that there exists a unique solution $u=\underline{\widehat{u}}^{m+1}(0), B=\underline{\widehat{B}}^{m+1}(0)$ to \cref{s*plasma} with source terms $\bF^{m,+},\bF_8^{m,+}$, $\bG_1^{m},\bG_{2}^{m}$, and initial values $\mathfrak{u}_0^{m+1}, \mathfrak{B}_0^{m+1}$. 
This complete the construction of $\underline{\widehat{U}}^{m+1}(0)$, $\underline{\widehat{V}}^{m+1}(0)$.

\vline

\vspace{0.3cm}

\subsubsection{Compatibility of the initial data}\label{s-initial-values}

For later use we give the following lemma.
\begin{Lem}
	\label{comp-source}	The source terms of \cref{s}, \cref{s-boundary} satisfy
	\begin{align*}
		\forall\, t\in[0,T]\qquad &\int_{\Omega_0^+} \bF_7^{m,+} dy + \int_{\Gamma_0} \bG_1^{m} dy'=\int_{\Gamma_0} \dt \widehat{\varphi}^{m+1}(0) dy',\\
		&\int_{\Omega_0^+} \bF_8^{m,+} dy + \int_{\Gamma_0} \bG_2^{m} dy'=0,\\
		&\int_{\Omega_0^+} \bF_7^{m,-} dy - \int_{\Gamma_0} \bG_3^{m} dy'=0\,.
	\end{align*}
\end{Lem}
\begin{proof}
	The proof for the equalities in the plasma part is given in \cite{Pierre2021}. The proof of the last equality is similar.
\end{proof}

We first consider the plasma part. The initial velocity field  $\mathfrak{u}_0^{m+1} $ in \cref{s-initial} should satisfy
\begin{equation}
	\label{s-initial-u}
	\left\{
	\begin{aligned}
		&\nabla \cdot \mathfrak{u}_0^{m+1} =\bF_7^{m,+}|_{t=0} \quad \text{in~} \Omega_0^+,\\
		& \mathfrak{u}_{0,3}^{m+1}|_{\Gamma_0} =\bG_1^{m}|_{t=0}.
	\end{aligned}
	\right.
\end{equation}
The solution to \cref{s-initial-u} exists if and only if 
\begin{equation*}
	\label{s-comp-u}
	\int_{\Omega_0^+} \bF_7^{m,+} |_{t=0}\,dy + \int_{\Gamma_0} \bG_1^{m} |_{t=0}\,dy' =0,
\end{equation*} 
and, from \cref{comp-source}, the above equation holds if and only if 
\begin{equation*}
	\label{s-comp-phit}
	\int_{\Gamma_0} \dt \widehat{\varphi}^{m+1}(0)|_{t=0}\, dy' =0,
\end{equation*}
which holds as a consequence of the choice \cref{s-phi0}. Since the source terms in \cref{s-initial-u} have $H^\infty$ regularity, we can find the solution $\mathfrak{u}_0^{m+1} \in  H^{\infty} ([0,T]\times\Omega_0^{+})$.

The initial magnetic field in the plasma part $\mathfrak{B}_0^{m+1}$ in \cref{s-initial} should satisfy
\begin{equation}
	\label{s-initial-B}
	\left\{
	\begin{aligned}
		&\nabla \cdot \mathfrak{B}_0^{m+1} =\bF_8^{m,+}|_{t=0} \quad \text{in~} \Omega_0^+,\\
		& \mathfrak{B}_{0,3}^{m+1}|_{\Gamma_0} =\bG_2^{m}|_{t=0}.
	\end{aligned}
	\right.
\end{equation}
The solution to \cref{s-initial-B} exists if and only if 
\begin{equation*}
	\int_{\Omega_0^+} \bF_8^{m,+} |_{t=0}\,dy + \int_{\Gamma_0} \bG_2^{m} |_{t=0}\,dy' =0.
\end{equation*}
From \cref{comp-source}, the above equation holds for all time, not only for $t=0$. Since the source terms in \cref{s-initial-B} have $H^\infty$ regularity, we can find the solution $\mathfrak{B}_0^{m+1} \in  H^{\infty} ([0,T]\times\Omega_0^{+})$.

We now focus on the vacuum part. For the initial magnetic field $\mathfrak{H}_0^{m+1}$ in \cref{s-initial} the construction is similar as for $\mathfrak{B}_0^{m+1}$. In fact it should satisfy 
\begin{equation}
	\label{s-initial-H}
	\left\{
	\begin{aligned}
		&\nabla \cdot \mathfrak{H}_0^{m+1} =\bF_7^{m,-}|_{t=0} \quad \text{in~} \Omega_0^-,\\
		& \mathfrak{H}_{0,3}^{m+1}|_{\Gamma_0} =\bG_3^{m}|_{t=0}.
	\end{aligned}
	\right.
\end{equation}
The solution to \cref{s-initial-H} exists if and only if 
\begin{equation*}
	\label{s-comp-H}
	\int_{\Omega_0^-} \bF_7^{m,-} |_{t=0}\,dy - \int_{\Gamma_0} \bG_3^{m} |_{t=0}\,dy' =0.
\end{equation*} 
From \cref{comp-source}, the above equation holds for all time, not only for $t=0$. Since the source terms in \cref{s-initial-H} have $H^\infty$ regularity, we can find the solution $\mathfrak{H}_0^{m+1} \in  H^{\infty} ([0,T]\times\Omega_0^{+})$.

The initial electric field $\mathfrak{E}_0^{m+1}$ should satisfy a different problem 
\begin{equation}
	\label{s-initial-E}
	\left\{
	\begin{aligned}
		&\nabla \cdot \mathfrak{E}_0^{m+1} =\bF_8^{m,-}|_{t=0} \quad \text{in~} \Omega_0^-;\\
		&\partial_t^k \mathfrak{E}_{0,2}^{m+1}|_{\Gamma_0} =\partial_t^k\bG_4^{m}|_{t=0}, ~ 
		-\partial_t^k \mathfrak{E}_{0,1}^{m+1}|_{\Gamma_0} =\partial_t^k\bG_5^{m}|_{t=0}.
	\end{aligned}
	\right.
\end{equation}
This problem has no compatibility condition for the source terms, thus the solution $\mathfrak{E}_0^{m+1}\in  H^{\infty} ([0,T]\times\Omega_0^{-})$ to \cref{s-initial-E} exists for any smooth $\bF_8^{m,-}, \bG_4^{m,-}, \bG_5^{m,-}$. 

\vspace{0.3cm}
\textbf{Summary.~} With the above compatibility conditions on the initial data and source terms (\cref{lem-comp-ams} and \cref{lem-comp-dms}), 
\cref{prop-s-vac}, \cref{prop-s-plasma} and \cref{lem-s-phi},  that there is a unique solution  $\underline{\widehat{U}}^{m+1}(0)$, $\underline{\widehat{V}}^{m+1}(0) \in  H^{\infty} ([0,T]\times\Omega_0^{\pm})$ and $\widehat{\varphi}^{m+1}(0) \in H^{\infty} ( [0,T]\times \mathbb{T}^2 )$ to \cref{s}, \cref{s-boundary}, \cref{s-initial} and \cref{s-phi-initial}.

The construction of the front $\varphi^{m+1}$ is completed as
\begin{equation*}
	\varphi^{m+1}(t,y',\theta)=\varphi^{m+1}_\# (t,y',\theta)+\widehat{\varphi}^{m+1}(t,y',0),
\end{equation*}
where $\varphi^{m+1}_\#$ is assigned in the induction assumption \cref{H(m)-1} and $\widehat{\varphi}^{m+1}(t,y',0)=\Phi^{m+1}(t,y')$ is given in \cref{lem-s-phi}.

In our constructions of the correctors \cref{decomp-fm+1},
it remains to determine  $\widehat{u}_{1,\star}^{m+1}(0)$, $\widehat{u}_{2,\star}^{m+1}(0)$, $\widehat{B}_{1,\star}^{m+1}(0)$, $\widehat{B}_{2,\star}^{m+1}(0)$, which is done by enforcing H((m+1)-4),
and the mean free part (in $\theta$) of the front profile $\varphi_{\#}^{m+2}$ by enforcing H((m+1)-5).

\subsection{The tangential components of the fast mean in plasma}
The components of the fast mean in vacuum have been determined in \cref{decomp-fm+1} when solving \cref{fm+1}.
Only the tangential components of the fast mean in plasma, i.e. $\widehat{u}_{1,\star}^{m+1}(0)$, $\widehat{u}_{2,\star}^{m+1}(0)$, $ \widehat{B}_{1,\star}^{m+1}(0)$, $\widehat{B}_{2,\star}^{m+1}(0)$, are degrees of freedom.
Our construction is the same as in \cite{Pierre2021}, and is obtained by enforcing condition (H(m+1)-4) and choosing the initial conditions
\begin{equation}\label{tangential-initial}
	\Pi\, \widehat{U}^{m+1}_\star(0) |_{t=0}=0 \quad  \text{for all~}m\geq1.
\end{equation}
We refer to \cite{Pierre2021} for details.

\subsection{The linearized nonlocal Hamilton-Jacobi equation for the front}
At this stage the only unknown parts of $U^{m+1}$ and $V^{m+1}$ in \cref{decomp-fm+1} are 
\begin{align*}
	&\mathring{U}^{m+1}:=\sum\limits_{k\in \mathbb{Z} \setminus {0} }|k| \widehat{\varphi}^{m+2}(t,y',k) \chi(y_3) e^{-|k|Y_3+ik\theta} \sR^+(k),\\
	&\mathring{V}^{m+1}:=\sum\limits_{k\in \mathbb{Z} \setminus \{0\}}\widetilde{\gamma}^{m+1}_j |k| \widehat{\varphi}^{m+2}(t,y',k) \chi(y_3) e^{|k|\sqrt{1-\nu^2\tau^2}Y_3+ik\theta} \sR^-_j(k),
\end{align*}
which depend on the oscillating modes in $\theta$ of the front profile $\varphi^{m+2}$. Such oscillating modes are determined by enforcing H((m+1)-5), i.e. the orthogonality condition .
Similar as deriving the equation \cref{eq-HJ} for the non-zero modes of the leading front $\varphi^{2}$, we can rewrite (H(m+1)-5) in terms of $\widehat{\varphi}^{m+2}$, i.e. replacing $\widehat{F}^{m+1,\pm} (k), \widehat{G}^{m+1}(k)$ by $\widehat{U}^{m+1}(k), \widehat{V}^{m+1}(k), \widehat{\varphi}^{m+2}(k)$, and finally by the only unknown $\widehat{\varphi}^{m+2}(k)$.

Note that during this process, the only terms depending on the lifting function $\chi(y_3)$ are $\sL^+(k) \cdot A^+_3 \partial_{y_3} \mathring{U}^{m+1}$ and  $ \sL^-(k) \cdot A^-_3 \partial_{y_3} \mathring{V}^{m+1}$.
These two terms can be cancelled since  $\sL^\pm(k) \cdot A_3^\pm \sR^\pm(k)=0$ for $k\neq0$.

The calculations of the substitution of \cref{decomp-fm+1} into (H(m+1)-5) are similar as those in \cref{sec-leading}, and thus we omit here. Eventually, taking $(\widetilde{\gamma}^{m+1}_1, \widetilde{\gamma}^{m+1}_2)=(\iota_1, \iota_2)$ (see \cref{def-gamma-special}),
(H(m+1)-5) turns out to be equivalent to 
\begin{equation}
	\begin{aligned}
		&(c^+ + d_0) \dt \widehat{\varphi}^{m+2}(k)  +(c^+ u^0_j -b^+ B^0_j+d_j) \partial_{y_j} \widehat{\varphi}^{m+2}(k)\\
		&\qquad+ 2i \Big( (c^+)^2 -(b^+)^2 - (a_1^-)^2-(a_2^-)^2+(b^-)^2 \Big) 
		\sgnk\\
		&\qquad \cdot \sum\limits_{k_1+k_2 = k}  \frac{|k||k_1||k_2|}{|k|+|k_1|+|k_2|}\widehat{\varphi}^{2} (k_1) \widehat{\varphi}^{m+2}(k_2)= \widehat{\mathfrak{g}}^{m+1}(k)
		\quad \text{for all~} k\in \mathbb{Z}\setminus \{0\},
		\label{eq-HJ-m}
	\end{aligned}
\end{equation}
where constants $d_0, d_j$ are same as those in  \cref{sec-leading}, and $\widehat{\mathfrak{g}}^{m+1}(k)$ of \eqref{eq-HJ-m} are the terms from the previously determined quantities.
The above equation is exactly the linearization of \cref{eq-HJ} about the leading front $\varphi^2_\#$.
It follows from the linear analogue of \cref{thm-HJ}, which is obviously easier to prove, that since $c^+\neq0$ and $\nu \ll1$, there is a unique $\varphi^{m+2}_{\#} \in C^{\infty}([0,T];H^\infty_\#)$ to \cref{eq-HJ-m} with zero initial data, 
with the positive time $T$ being fixed in \cref{thm-HJ} as the lifespan of $\varphi^{2}_{\#}$. 
%So, as quoted at the beginning of this section, the time $T$ in the induction assumption H(m) is determined by the solvability of the leading front to \cref{eq-HJ}. 

With the above regularity of $\varphi^{m+2}_{\#}$, it is easy to verify that $(U^{m+1}, V^{m+1})\in S^{\pm}$ in \cref{decomp-fm+1}, i.e. H((m+1)-1) is fulfilled. 
This completes the proof of our induction, i.e.  H(m+1) is satisfied for all $m\ge1$, with the same positive time $T$ fixed by \cref{thm-HJ}.

\subsection{Summary on the construction of the correctors}

Before we finish proving \cref{thm-main}, let us make a short summary on our construction of the correctors.
After defining H(m) at the beginning of this section and verifying H(1) in \cref{sec-initial step}, we construct 
$(U^{m+1}, V^{m+1}, \varphi^{m+2}_{\#})$ and $\widehat{\varphi}^{m+1}(0)$, for given $(U^{1}, V^{1}, \varphi^{2})$, $(U^{2}, V^{2}, \varphi^{3})$, $\cdots$, $(U^{m-1},V^{m-1},\varphi^{m})$, $(U^m,V^m,\varphi^{m+1}_{\#})$
and the induction assumption H(m) in time $[0,T]$, such that H(m+1) is satisfied for the same time $T$.

\begin{enumerate}
	\item Due to H((m)-3), H((m)-4), H((m)-5), \cref{comp-am} and \cref{comp-dm},  the fast problem (H(m+1)-2) (that is \cref{fm+1}) has solutions $(U^{m+1}, V^{m+1}, \varphi^{m+2}_{\#})$ in the form \cref{decomp-fm+1}.
	
	\item By enforcing H((m+1)-3), we determine $\underline{\widehat{U}}^{m+1}(0)$, $\underline{\widehat{V}}^{m+1}(0)$ and
	$\widehat{\varphi}^{m+1}(0)$ by the slow problem \cref{s}, \cref{s-boundary} and \cref{s-initial}. 
	
	\item By enforcing H((m+1)-4), we determine  $\widehat{u}_{1,\star}^{m+1}(0)$, $\widehat{u}_{2,\star}^{m+1}(0)$, $ \widehat{B}_{1,\star}^{m+1}(0)$, $\widehat{B}_{2,\star}^{m+1}(0)$ as in \cite{Pierre2021}. 
	
	\item By enforcing H((m+1)-5), we determine $\varphi^{m+2}_\#$ as the unique solution to \cref{eq-HJ-m} with zero initial data (see \cref{WKB-front} for $m+2\geq3$). \cref{eq-HJ-m} is the linearization of the nonlocal Hamilton-Jacobi equation \cref{eq-HJ}, when taking $(\widetilde{\gamma}^{m+1}_1, \widetilde{\gamma}^{m+1}_2)=(\iota_1, \iota_2)$.	
\end{enumerate}
With all the above steps, the construction of the corrector $(U^{m+1}, V^{m+1}, \varphi^{m+2}_\#)$ $\in S^{\pm} \times H^{\infty} ([0,T]\times\mathbb{T}^2 \times\mathbb{T} )$ (defined as in \cref{decomp-fm+1}) is completed.

\vspace{0.3cm}
Now we compare our constructions with the statement of \cref{thm-main}.
The initial conditions for the front profiles in \cref{thm-main} have been fulfilled by our settings \cref{WKB-front}, \cref{s-phi-initial}.
% and \cref{s-phi0} (the nonzero modes with respect to $y'$ of $\partial_t \widehat{\varphi}^m(0)$ actually could be chosen arbitrarily, while the mean with respect to $y'$ has to vanish due to the compatibility for the divergence of the initial data).
In \cref{sec-leading}, we have set that  $\underline{U}^1=0$, $\underline{V}^1=0$, $\widehat{U}^1_\star(0)=0$, $\widehat{V}^1_\star(0)=0$, while $\underline{U}^2=0$ and $\underline{V}^2=0$ will be examined in the next section.
When solving the tangential components of the fast mean in plasma, we have set $\Pi\widehat{U}^{m}_\star(0) |_{t=0}=0$ in \cref{tangential-initial}, for any $m\ge1$.

To finish the proof of \cref{thm-main}, it remains to verify that the approximate solutions defined in  \cref{thm-main} nearly satisfy the original plasma--vacuum interface problem \cref{MHD}, \cref{Maxwell} and \cref{interface} at any desired order of accuracy, and also the nontriviality of $\underline{\widehat{U}}^3(0)$ or $\underline{\widehat{V}}^3(0)$.

\subsection{Approximate solutions}
Recall that the approximate solutions are defined in \cref{thm-main} as:
%for any integer $M\geq 1$
\begin{align*}
	U_\e^{\textnormal{app},M}(t,x) & := U^0+\sum\limits_{m= 1}^M \e^m U^m(t,x',y_3,Y_3,\theta),\\
	V_\e^{\textnormal{app},M}(t,x) & := V^0+\sum\limits_{m= 1}^M \e^m V^m(t,x',y_3,Y_3,\theta),\\
	\varphi_\e^{\textnormal{app},M}(t,x') & := \sum\limits_{m= 2}^{M+1} \e^m\varphi^m(t,x',\theta)\qquad \text{for any integer~}M\geq 1,
\end{align*}where $ y_3:= x_3-\varphi_\e^{\textnormal{app},M} (t,x'),~ Y_3:=\frac{y_3}{\e},~ \theta:=\frac{\tau t+\xi'\cdot x'}{\e}$, 
and $(U^{\mu},V^{\mu},\varphi^{\mu+1})_{\mu=0,\dots,M}$ satisfy the condition (H(M)) stated at the beginning of \cref{sec-correctors}.
 
Firstly, we start with the error terms for the jump conditions.
\begin{align*}
	R_{b,\e}^{2}=&q_\e^{\textnormal{app},M}|_{\Gamma_\e^{\textnormal{app},M}(t)}-\frac{1}{2} \abs{H_\e^{\textnormal{app},M}}{2}|_{\Gamma_\e^{\textnormal{app},M}(t)}+\frac{1}{2}\abs{E_\e^{\textnormal{app},M}}{2}|_{\Gamma_\e^{\textnormal{app},M}(t)}\\
	=&\Big( q^0 +  \sum\limits_{m= 1}^M \e^m q^m(\cdot)\Big)-\frac{1}{2} \Big( H_\alpha^0 +  \sum\limits_{m= 1}^M \e^m H_\alpha^m(\cdot)\Big)\Big( H_\alpha^0 +  \sum\limits_{m= 1}^M \e^m H_\alpha^m(\cdot)\Big)\\
	&\quad +\frac{1}{2}\Big( E_\alpha^0 +  \sum\limits_{m= 1}^M \e^m E_\alpha^m(\cdot)\Big)\Big( E_\alpha^0 +  \sum\limits_{m= 1}^M \e^m E_\alpha^m(\cdot)\Big)\\
	=&\Big(q^0-\frac{1}{2}(H_1^0)^2-\frac{1}{2}(H_2^0)^2+\frac{1}{2}(E_3^0)^2\Big)\\
	&\quad\quad +\sum\limits_{m=1}^M \e^m \Big( q^m (\cdot)-H_j^0 H_j^m (\cdot) +E_3^0 E_3^m(\cdot)  - \frac{1}{2} \sum\limits_{\substack{l_1+l_2=m\\l_1,l_2\geq 1}}  H_\alpha^{l_1} (\cdot)H_\alpha^{l_2} (\cdot) \\
	&\quad
	+\frac{1}{2} \sum\limits_{\substack{l_1+l_2=m\\l_1,l_2\geq 1}}  E_\alpha^{l_1} (\cdot)E_\alpha^{l_2} (\cdot)  \Big) +\mathcal{O}(\e^{M+1})\\
	=&\sum\limits_{m=1}^M \e^m \Big( q^m|_{y_3=Y_3=0}-H_j^0 H_j^m |_{y_3=Y_3=0} +E_3^0 E_3^m  |_{y_3=Y_3=0}\\
	&\quad 
	- G_6^{m-1}  \Big)(t,x',\frac{\tau t+\xi'\cdot x'}{\e}) +\mathcal{O}(\e^{M+1})\\
	=&\mathcal{O}(\e^{M+1}),
\end{align*}
where $(\cdot)$ stands for 
$$\left(t,x',0,0,\frac{\tau t+\xi'\cdot x'}{\e}\right).$$
The last equality comes from the jump condition of the fast problem (H(M)-2), see \cref{WKB-jump-q}, \eqref{Gm6}. It follows that
\begin{equation*}
\sup\limits_{ t\in [0,T], x\in \Gamma_\e^{\textnormal{app},M}(t)}   \abs{R_{b,\e}^{2} }{}=\mathcal{O} (\e^{M+1}).
\end{equation*}

Let us estimate the error terms for the partial differential equations in the vacuum part. Specifically, let us consider the error term in the equation for the first component of $H_{\e}^{\textnormal{app},M}$. We have
\begin{align*}
	R_{\e,1}^{1,-}=& \nu \dt H_{\e,1}^{\textnormal{app},M} +  \partial_{x_2}E_{\e,3}^{\textnormal{app},M} -\partial_{x_3} E_{\e,2}^{\textnormal{app},M}\\
	=& \sum\limits_{m= 1}^{M} \e^{m-1} \Big( \nu \tau \partial_{\theta} H_1^m + \xi_2 \partial_{\theta} E_3^m  -\partial_{Y_3}E_2^m  + \nu \dt H_1^{m-1}+ \partial_{x_2} E_3^{m-1}   \\
	&\quad -\partial_{x_3}E_2^{m-1}
	 -\sum\limits_{\substack{l_1+l_2=m+1\\l_1\geq 1}} \partial_{\theta} \varphi^{l_2} \partial_{Y_3}(\nu \tau H_1^{l_1} +\xi_2 E_3^{l_1} )  \\
	&\quad -\sum\limits_{\substack{l_1+l_2=m\\l_1\geq 1}}( \partial_{t} \varphi^{l_2} \partial_{Y_3}H_1^{l_1} + \partial_{x_2}\varphi^{l_2} \partial_{Y_3}  E_3^{l_1} )
	- \sum\limits_{\substack{l_1+l_2=m\\l_1\geq 1}} \partial_{\theta} \varphi^{l_2} \partial_{y_3}(\nu \tau H_1^{l_1} +\xi_2 E_3^{l_1} )\\
	&\quad -\sum\limits_{\substack{l_1+l_2=m-1\\l_1\geq 1}}( \partial_{t} \varphi^{l_2} \partial_{y_3}H_1^{l_1} + \partial_{x_2}\varphi^{l_2} \partial_{y_3}  E_3^{l_1} ) \Big)(\cdot) +\mathcal{O}(\e^{M})\\
	=& \sum\limits_{m= 1}^{M} \e^{m-1} \Big( \nu \tau \partial_{\theta} H_1^m + \xi_2 \partial_{\theta} E_3^m  -\partial_{Y_3}E_2^m  -F_1^{m-1,-} \Big)(\cdot) +\mathcal{O}(\e^{M})\\
	=&\mathcal{O}(\e^{M}),
\end{align*}
where here $(\cdot)$ stands for 
$$\left(t,x',x_3-\varphi_{\e}^{\textnormal{app},M}(t,x') ,\frac{x_3-\varphi_{\e}^{\textnormal{app},M}(t,x')}{\e},\frac{\tau t+\xi'\cdot x'}{\e}\right).$$
The last equality comes from the first equation for the vacuum part of the fast problem (H(M)-2).
It follows that
\begin{align*}
\sup\limits_{ t\in [0,T], x\in \Omega_\e^{\textnormal{app},M,\pm}(t)} \abs{R_{\e,1}^{1,-} }{} =\mathcal{O} (\e^M).
\end{align*}
The remaining estimates of \cref{thm-main} are obtained in a similar way; in particular, in the plasma part the estimates are as those in \cite{Pierre2021}. Thus we feel free to omit the details.
Therefore, the desired order of accuracy of the approximate solutions is verified, and it remains to consider the rectification phenomenon to complete \cref{thm-main}.

%%%%%%%%%%%%%%%%%%%%%%%%%%%%%%%%%%%%%%%%%%%%%%%%%%%%%%%%%%%%%%%%
%%%%%The rectification phenomenon
%%%%%%%%%%%%%%%%%%%%%%%%%%%%%%%%%%%%%%%%%%%%%%%%%%%%%%%%%%%%%%%%

\section{The rectification phenomenon}
In this section, we will study whether there is a nontrivial residual component of one corrector $(U^m,V^m)$.
It will be proved that without any more requirement on $\varphi^2_0$ other than \cref{thm-main}, the rectification phenomenon will not appear at the first corrector, i.e., $\underline{U}^2=0$ and $\underline{V}^2=0$, which is the last remaining property to prove of our construction of the WKB expansion.
It is the same as \cite{Pierre2021}, but different from elastodynamics \cite{Marcou2011}.

In contrast to uncertain nontriviality of the residual component of the second corrector in \cite{Pierre2021}, in our case it is confirmed that $\underline{U}^3=\underline{\widehat{U}}^3(0)\neq0$ or $\underline{V}^3=\underline{\widehat{V}}^3(0)\neq0$ for general $\varphi^2_0$. 

\subsection{The first corrector}
In \cref{sec-leading}, we have settled that  $\underline{U}^1=\widehat{U}_\star^1(0)=0, \underline{V}^1=\widehat{V}_\star^1(0)=0$ and 
\begin{equation}\label{sol-leading-simp}
	\begin{aligned}
		&\widehat{U}^1(t,y',0,Y_3,k)=\widehat{U}^1_\star(t,y',0,Y_3,k) %\notag\\	&\quad 
		= |k| \widehat{\varphi}^2(t,y',k)  e^{-|k|Y_3} \sR^+(k), \notag\\
		&\widehat{V}^1(t,y',0,Y_3,k)=\widehat{V}^1_\star(t,y',0,Y_3,k) %\notag\\ &\quad 
		= |k| \widehat{\varphi}^2(t,y',k)  e^{|k|\sqrt{1-\nu^2\tau^2}Y_3} \sR^-(k), \notag
	\end{aligned}
\end{equation}
where $\sR^-(k)=\iota_j\sR^-_j(k) (j=1,2)$ is defined in \cref{vector-R3}.
Now we solve $\underline{U}^2$ and $\underline{V}^2$ as \cref{sec-fast,sec-correctors} in following several steps.

\subsubsection{The nonzero modes of $\underline{U}^2$ and $\underline{V}^2$} Since $\underline{U}^2,\underline{V}^2$ vanish, so do $\underline{F}^{1,\pm}.$ Then as solving the inhomogeneous fast problem on nonzero Fourier modes (Section 3.2.3), we have 
$$
\underline{\widehat{U}}^2(k)=\frac{1}{ik} (\sA^+)^{-1} \underline{\widehat{F}}^{1,+}(k)=0,\quad 
\underline{\widehat{V}}^2(k)=\frac{1}{ik} (\sA^-)^{-1} \underline{\widehat{F}}^{1,-}(k)=0.
$$
Since $\sA^\pm$ are invertible, $\underline{\widehat{U}}^2(k), 
\underline{\widehat{V}}^2(k)$ are trivial. 

\subsubsection{Deriving the equations for $\underline{\widehat{U}}^2(0)$ and $\underline{\widehat{V}}^2(0)$}
It is easy to obtain $\bF^{1,\pm}=0$, and thus as \cref{s} for $m=1$, $\underline{\widehat{U}}^2$ and $\underline{\widehat{V}}^2$ satisfy 
\begin{equation}
	\label{s2}
	\left\{
	\begin{aligned}
		&L_s^+ \underline{\widehat{U}}^{2}(0)=0,  \nabla\cdot \underline{\widehat{B}}^{2}(0) = 0,  & y\in \Omega_0^+, \\
		&L_s^- \underline{\widehat{V}}^{2}(0)=0,  \nabla\cdot \underline{\widehat{H}}^{2}(0) = 0,
		\nabla\cdot \underline{\widehat{E}}^{2}(0) = 0, & y\in \Omega_0^-,
	\end{aligned}
	\right.
\end{equation}
The boundary conditions \cref{s-boundary} for $m=1$ read explicitly:
\begin{equation}
	\label{s-boundary2}
	\left\{
	\begin{aligned}
		&\,
		\underline{\widehat{u}}_3^{2}(0) |_{y_3=0}= (\dt +u^0_j \partial_{y_j})  \widehat{\varphi}^{2} (0)
		+\mathbf{c_0} \{ \partial_{\theta}\varphi^{2} \xi_j   u^{1}_j |_{y_3=Y_3=0} \} -\widehat{u}_{3,\star}^{2}(0) |_{y_3=Y_3=0} =: \bG^{1}_1,  \\
		&
		\underline{\widehat{B}}_3^{2}(0) |_{y_3=0}=  B^0_j \partial_{y_j}   \widehat{\varphi}^{2} (0)
		+\mathbf{c_0} \{ \partial_{\theta}\varphi^{2} \xi_j   B^{1}_j |_{y_3=Y_3=0} \} -\widehat{B}_{3,\star}^{2}(0) |_{y_3=Y_3=0} =: \bG^{1}_2, \\
		&
		\underline{\widehat{H}}_3^{2}(0) |_{y_3=0}=  H^0_j \partial_{y_j}  \widehat{\varphi}^{2} (0)
		+\mathbf{c_0} \{ \partial_{\theta}\varphi^{2} \xi_j   H^{1}_j |_{y_3=Y_3=0} \}-\widehat{H}_{3,\star}^{2}(0) |_{y_3=Y_3=0} =: \bG^{1}_3, \\
		&
		\underline{\widehat{E}}_2^{2}(0) |_{y_3=0}= (-E^0_3\partial_{y_2} -\nu H^0_1 \dt) \widehat{\varphi}^{2} (0)
		-\mathbf{c_0} \{ \partial_{\theta}\varphi^{2} (\xi_2   E^{1}_3 +\nu \tau H^{1}_1 ) |_{y_3=Y_3=0}\} \\
		&\qquad\qquad\qquad\qquad
		-\widehat{E}_{2,\star}^{2}(0) |_{y_3=Y_3=0}  =: \bG^{1}_4, \\
		&
		-\underline{\widehat{E}}_1^{2}(0) |_{y_3=0}= (E^0_3\partial_{y_1} -\nu H^0_2 \dt)  \widehat{\varphi}^{2} (0)
		+\mathbf{c_0} \{ \partial_{\theta}\varphi^{2} (\xi_1   E^{1}_3-\nu \tau H^{1}_2  ) |_{y_3=Y_3=0}\} \\
		&\qquad\qquad\qquad\qquad
		+\widehat{E}_{1,\star}^{2}(0) |_{y_3=Y_3=0} =: \bG^{1}_5, \\
		&
		\{\underline{\widehat{q}}^{2}(0)-H_j^0\underline{\widehat{H}}_j^{2}(0) + E_3^0 \underline{\widehat{E}}_3^{2}(0) \}  |_{y_3=0}=\mathbf{c_0} \{ \frac{1}{2} ( H^{1}_{\alpha}H^{1}_{\alpha}-  E^{1}_{\alpha}E^{1}_{\alpha}) |_{y_3=Y_3=0} \} \\
		&\qquad\qquad\qquad\qquad
		-\{\widehat{q}_{\star}^{2}(0)-H_j^0\widehat{H}_{j,\star}^{2}(0) + E_3^0 \widehat{E}_{3,\star}^{2}(0) \}  |_{y_3=Y_3=0}
		=: \bG^{1}_6.
	\end{aligned}
	\right.
\end{equation}

Before we solve \cref{s2} and \cref{s-boundary2} directly, we first simplify the source terms of the boundary conditions.

\subsubsection{Simplifying the source terms of \cref{s-boundary2}}
Recalling the expression of $F^{1,-}_1$ in \cref{F01-} and vanishing $\widehat{U}_\star^1(0), \widehat{V}_\star^1(0)$, we have 
\begin{align*}
	-\partial_{Y_3} \widehat{E}_{2,\star}^{2}(0) |_{y_3=0} = &  -(\nu \dt \widehat{H}_{1,\star}^{1} + \partial_{y_2} \widehat{E}_{3,\star}^{1} - \partial_{y_3} \widehat{E}_{2,\star}^{1})(0)  |_{y_3=0}\\
	& + \mathbf{c_0} \{ \partial_{\theta}\varphi^{2} (\nu \tau \partial_{Y_3} H^{1}_{1,\star} +\xi_2 	\partial_{Y_3}  E^{1}_{3,\star} ) |_{y_3=0}\}\\
	=&\partial_{Y_3}\mathbf{c_0} \{ \partial_{\theta}\varphi^{2} (\nu \tau  H^{1}_{1,\star} +\xi_2  E^{1}_{3,\star} ) |_{y_3=0}\}.
\end{align*}
Integrating the above equality with respect to $Y_3$ from $0$ to $-\infty$, it yields that
\begin{align*}
	\widehat{E}_{2,\star}^{2}(0) |_{y_3=0} = -\mathbf{c_0} \{ \partial_{\theta}\varphi^{2} (\nu \tau  H^{1}_{1,\star} +\xi_2  E^{1}_{3,\star} ) |_{y_3=0}\},
\end{align*}
and thus substituting it into \cref{s-boundary2} with $\underline{V}^1=0$ yields that $\underline{\widehat{E}}_2^{2}(0) |_{y_3=0}= (-E^0_3\partial_{y_2} -\nu H^0_1 \dt)  \widehat{\varphi}^{2} (0)$. 
Note that actually it holds that $\widehat{E}_{2,\star}^{2}(0) |_{y_3=0}=0$ since 
\begin{align*}
	&\mathbf{c_0} \{ \partial_{\theta}\varphi^{2} (\nu \tau  H^{1}_{1,\star} +\xi_2  E^{1}_{3,\star} ) |_{y_3=0}\}\\
	=&\Big(\nu\tau \frac{\nu\tau a^-_1-\xi_1 b^-}{\sqrt{1-\nu^2\tau^2 }}+\xi_2  \frac{\xi_1 a^-_2-\xi_2 a^-_1}{\sqrt{1-\nu^2\tau^2 }}  \Big)\\
	& \quad \cdot
	\sum\limits_{k\neq 0} i(-k)|k| \widehat{\varphi}^2(-k) \widehat{\varphi}^2(k) e^{|k|\sqrt{1-\nu^2\tau^2}Y_3}=0,
\end{align*}
where the final equality follows from changing $k$ to $-k$.
The similar arguments also work for $	(\underline{\widehat{u}}_3^{2}, 	\underline{\widehat{B}}_3^{2}, 	\underline{\widehat{H}}_3^{2}, 	\underline{\widehat{E}}_1^{2})(0) |_{y_3=0}$.

Now we consider $\{\underline{\widehat{q}}^{2}(0)-H_j^0\underline{\widehat{H}}_j^{2}(0) + E_3^0 \underline{\widehat{E}}_3^{2}(0) \}  |_{y_3=0}$, and by applying the similar calculations of the above, we have
\begin{align*}
	&\{\underline{\widehat{q}}^{2}(0)-H_j^0\underline{\widehat{H}}_j^{2}(0) + E_3^0 \underline{\widehat{E}}_3^{2}(0) \}  |_{y_3=0}\\
	&=\mathbf{c_0} \{ \frac{1}{2} ( H^{1}_{\alpha}H^{1}_{\alpha}-  E^{1}_{\alpha}E^{1}_{\alpha}) 
	+ \partial_{\theta}\varphi^{2}  (a^-_2  E^{1}_{1,\star} -a^-_1 E^{1}_{2,\star}-b^-  H^{1}_{3,\star}  )   \}|_{y_3=Y_3=0} ,
\end{align*} 
where $\widehat{q}_{\star}^{2}(0)|_{y_3=Y_3=0} =0$ as \cite{Pierre2021}. Furthermore, with $\xi_j a^-_j=\nu\tau b^-$, it holds that
\begin{align*}
	&\mathbf{c_0} \{H^{1}_{\alpha}H^{1}_{\alpha}-  E^{1}_{\alpha}E^{1}_{\alpha}\} |_{y_3=0} \\
	&\quad= \Big(\frac{(\nu\tau a^-_1-\xi_1 b^-)^2+(\nu\tau a^-_2-\xi_2 b^-)^2}{1-\nu^2\tau^2}+(b^-)^2
	-(a^-_1)^2-(a^-_2)^2\\
	&\qquad\qquad-\frac{(\xi_1 a^-_2-\xi_2a^-_1)^2}{1-\nu^2\tau^2   }\Big)
	\sum\limits_{k\neq 0} k^2 \widehat{\varphi}^2(-k) \widehat{\varphi}^2(k) e^{2|k|\sqrt{1-\nu^2\tau^2}Y_3} \\
	&\quad=	\Big(\frac{\nu^2\tau^2(a^-_1)^2+\nu^2\tau^2(a^-_2)^2-2\nu^2\tau^2 (b^-)^2 + (b^-)^2 }{1-\nu^2\tau^2}+(b^-)^2
	-(a^-_1)^2-(a^-_2)^2\\
	&\qquad\qquad- \frac{(a^-_1)^2+(a^-_2)^2-\nu^2\tau^2 (b^-)^2}{1-\nu^2\tau^2}\Big) \sum\limits_{k\neq 0} k^2 \widehat{\varphi}^2(-k) \widehat{\varphi}^2(k) e^{2|k|\sqrt{1-\nu^2\tau^2}Y_3}\\
	&\quad= 
	\Big(-2(a^-_1)^2-2(a^-_2)^2+2(b^-)^2\Big)
	\sum\limits_{k\neq 0} k^2 \widehat{\varphi}^2(-k) \widehat{\varphi}^2(k) e^{2|k|\sqrt{1-\nu^2\tau^2}Y_3},\\
	%%%
	&\mathbf{c_0} \{ \partial_{\theta}\varphi^{2}  (a^-_2  E^{1}_{1,\star} -a^-_1 E^{1}_{2,\star}-b^-  H^{1}_{3,\star}  )  \}|_{y_3=0}\\
	&\quad= \Big((a^-_1)^2+(a^-_2)^2-(b^-)^2\Big)
	\sum\limits_{k\neq 0} k^2 \widehat{\varphi}^2(-k) \widehat{\varphi}^2(k) e^{|k|\sqrt{1-\nu^2\tau^2}Y_3}.
\end{align*}
Therefore, \cref{s-boundary2} can be simplified as 
\begin{equation}
	\label{s-boundary2-simp}
	\left\{
	\begin{aligned}
		&\,
		\underline{\widehat{u}}_3^{2}(0) |_{y_3=0}= (\dt +u^0_j \partial_{y_j})   \widehat{\varphi}^{2} (0)  ,  \quad		\underline{\widehat{B}}_3^{2}(0) |_{y_3=0}=  B^0_j \partial_{y_j}   \widehat{\varphi}^{2} (0) , \\
		&
		\underline{\widehat{H}}_3^{2}(0) |_{y_3=0}=  H^0_j \partial_{y_j}  \widehat{\varphi}^{2} (0) , \quad
		\underline{\widehat{E}}_2^{2}(0) |_{y_3=0}= (-E^0_3\partial_{y_2} -\nu H^0_1 \dt)   \widehat{\varphi}^{2} (0) ,\\
		&-\underline{\widehat{E}}_1^{2}(0) |_{y_3=0}= (E^0_3\partial_{y_1} -\nu H^0_2 \dt)   \widehat{\varphi}^{2} (0) , \\
		&
		\{\underline{\widehat{q}}^{2}(0)-H_j^0\underline{\widehat{H}}_j^{2}(0) + E_3^0 \underline{\widehat{E}}_3^{2}(0) \}  |_{y_3=0}=0.
	\end{aligned}
	\right.
\end{equation}

\subsubsection{Solving the slow problem \cref{s2} and \cref{s-boundary2-simp}} As in \cref{sec-slow means}, we first construct $\widehat{\varphi}(0)$ by considering the following elliptic-hyperbolic problem:
\begin{equation}
	\label{s-front2}
	\left\{ 
	\begin{aligned}
		-\Delta \underline{\widehat{q}}^{2}(0)&=  0,   \\
		\partial_{y_3}\underline{\widehat{q}}^{2}(0)  |_{y_3=0} &=  -(\dt+u^0_i\partial_{y_i})(\dt+u^0_j\partial_{y_j})  \widehat{\varphi}^{2} (0)+B^0_iB^0_j \partial_{y_i}  \partial_{y_j} \widehat{\varphi}^{2} (0),   \\
		(\nu^2 \dt^2-\Delta) \widetilde{\underline{q}}^{2}(0) &= 0,   \\
		\partial_{y_3}\widetilde{\underline{q}}^{2}(0)  |_{y_3=0}& = H^0_iH^0_j \partial_{y_i}  \partial_{y_j} \widehat{\varphi}^{2} (0)  - (H^0_1 \nu \dt +E^0_3 \partial_{y_2} )(H^0_1 \nu \dt +E^0_3 \partial_{y_2} ) \widehat{\varphi}^{2} (0) \\
		&\quad
		-(H^0_2 \nu \dt -E^0_3 \partial_{y_1} )(H^0_2 \nu \dt -E^0_3 \partial_{y_1} )  \widehat{\varphi}^{2} (0) ,\\
		\{ \underline{\widehat{q}}^{2}(0)  -\widetilde{\underline{q}}^{2}(0)& \} |_{y_3=0} =0,
	\end{aligned}
	\right.
\end{equation}
where $\widetilde{\underline{q}}^{2}(0)  =  H_j^0\underline{\widehat{H}}_j^{2}(0) - E_3^0 \underline{\widehat{E}}_3^{2}(0)$.
As in \cref{sec-slow means}, it turns out that the nonzero Fourier modes $\textbf{c}_{j'}$ w.r.t. $y'$ and Laplace transform w.r.t $t$ of $\widehat{\varphi}^2(0)$ must satisfy
\begin{align*}
	&\sqrt{\nu^2 z^2 +|j'|^2}	((z+iu_j^0 j_j')^2 -(iB_j^0 j_j')^2)  + |j'| ((H_1^0\nu z +iE_3^0j_2')^2 \notag\\
	&\quad+(H_2^0\nu z -iE_3^0j_1')^2-(iH_j^0 j_j')^2)  \fL\textbf{c}_{j'} (\widehat{\varphi}^{2}(0))=0,
\end{align*}
and the zero Fourier modes $\textbf{c}_{j'}$ w.r.t. $y'$ of $\widehat{\varphi}^2(0)$ has to be $0$ due to the compatible conditions of \cref{s2} and \cref{s-boundary2-simp}. 
Therefore, it is easy to obtain that $\widehat{\varphi}^2(0)=0$ and thus $\underline{\widehat{q}}^{2}(0)=\widetilde{\underline{q}}^{2}(0)=0$. 
It is obvious that the solution to \cref{s2} and \cref{s-boundary2} is $\underline{\widehat{U}}^2(0)=0, \underline{\widehat{V}}^2(0)=0$, which finishes the proof of \cref{thm-main}.

\subsection{The second corrector}
Since $\underline{U}^1=\underline{U}^2=0, \underline{V}^1=\underline{V}^2=0,$ we have $\underline{F}^{2,\pm}=0$ and thus 
$$\underline{U}^3=\underline{\widehat{U}}^3(0), \quad \underline{V}^3=\underline{\widehat{V}}^3(0).$$ 
In the following, we will show  $\underline{\widehat{U}}^3(0)\neq0$ or $\underline{\widehat{V}}^3(0)\neq0$ for general $\varphi^2_0$ by contradiction arguments.

\subsubsection{Deriving the equations for $\underline{\widehat{U}}^3(0)$ and $\underline{\widehat{V}}^3(0)$}
By the compatible conditions H(m-3) and a direct computation of $\bF^{2,\pm}=0$, the slow means $\underline{\widehat{U}}^3(0), \underline{\widehat{V}}^3(0)$ satisfy 
\begin{equation}
	\label{s3}
	\left\{
	\begin{aligned}
		&L_s^+ \underline{\widehat{U}}^{3}(0)=0,  \nabla\cdot \underline{\widehat{B}}^{3}(0) = 0,  & y\in \Omega_0^+, \\
		&L_s^- \underline{\widehat{V}}^{3}(0)=0,  \nabla\cdot \underline{\widehat{H}}^{3}(0) = 0,
		\nabla\cdot \underline{\widehat{E}}^{3}(0) = 0, & y\in \Omega_0^-,
	\end{aligned}
	\right.
\end{equation}
\begin{equation}
	\label{s-boundary3}
	\left\{
	\begin{aligned}
		&
		\underline{\widehat{u}}_3^{3}(0) |_{y_3=0}= \widehat{G}_1^{2}(0)-\widehat{u}_{3,\star}^{3}(0) |_{y_3=Y_3=0} =: \bG^{2}_1,  \\
		&
		\underline{\widehat{B}}_3^{3}(0) |_{y_3=0}= \widehat{G}_2^{2}(0)-\widehat{B}_{3,\star}^{3}(0) |_{y_3=Y_3=0} =: \bG^{2}_2, \\
		&
		\underline{\widehat{H}}_3^{3}(0) |_{y_3=0}= \widehat{G}_3^{2}(0)-\widehat{H}_{3,\star}^{3}(0) |_{y_3=Y_3=0} =: \bG^{2}_3, \\
		&
		\underline{\widehat{E}}_2^{3}(0) |_{y_3=0}= \widehat{G}_4^{2}(0)-\widehat{E}_{2,\star}^{3}(0) |_{y_3=Y_3=0}  =: \bG^{2}_4, \\
		&
		-\underline{\widehat{E}}_1^{3}(0) |_{y_3=0}= \widehat{G}_5^{2}(0)+\widehat{E}_{1,\star}^{3}(0) |_{y_3=Y_3=0} =: \bG^{2}_5, \\
		&
		\{\underline{\widehat{q}}^{3}(0)-H_j^0\underline{\widehat{H}}_j^{3}(0) + E_3^0 \underline{\widehat{E}}_3^{3}(0) \}  |_{y_3=0} \\
		&
		\qquad\qquad= \widehat{G}_6^{2}(0) -\{\widehat{q}_{\star}^{3}(0)-H_j^0\widehat{H}_{j,\star}^{3}(0) + E_3^0 \widehat{E}_{3,\star}^{3}(0) \}  |_{y_3=Y_3=0}
		=: \bG^{2}_6.
	\end{aligned}
	\right.
\end{equation}
%As the above arguments of the first corrector, the rectification phenomenon occurs at the second corrector only if $\bG^2$ does not reduce to the linear terms of $\widehat{\varphi}^3(0)$. \todo{to change}

Now we simplify the source terms of the boundary conditions.
Recalling the expression of $F^{2,-}_1$ in \cref{Fm-} and %vanishing $\widehat{U}_\star^1(0), \widehat{V}_\star^1(0)$
, we have 
\begin{align*}
	-\partial_{Y_3} \widehat{E}_{2,\star}^{3}(0) |_{y_3=0} 
	= &  -(\nu \dt \widehat{H}_{1,\star}^{2} + \partial_{y_2} \widehat{E}_{3,\star}^{2} - \partial_{y_3} \widehat{E}_{2,\star}^{2})(0) |_{y_3=0} \\
	& + \mathbf{c_0} \Big\{ \partial_{\theta}\varphi^{3} (\nu \tau \partial_{Y_3} H^{1}_{1,\star} +\xi_2 	\partial_{Y_3}  E^{1}_{3,\star} ) + \partial_{\theta}\varphi^{2} (\nu \tau \partial_{Y_3} H^{2}_{1,\star} +\xi_2 	\partial_{Y_3}  E^{2}_{3,\star} )\\
	& +\partial_{\theta}\varphi^{2} (\nu \tau \partial_{y_3} H^{1}_{1,\star} +\xi_2 \partial_{y_3}  E^{1}_{3,\star} ) 
	+  (\nu \partial_{t} \varphi^{2} \partial_{Y_3} H^{1}_{1,\star} +\partial_{y_2}  \varphi^{2}	\partial_{Y_3}  E^{1}_{3,\star} )
	\Big\} |_{y_3=0}.
\end{align*}
Note that $\partial_{y_3} H^{1}_{1,\star} , \partial_{y_3}  E^{1}_{3,\star} $ vanishes at $y_3=0$ since $ H^{1}_{1,\star} ,  E^{1}_{3,\star} $ w.r.t. $y_3$ are totally dependent on the lifting function $\chi(y_3)$, which equals to $1$ near $y_3$.
This property also holds for $F^{1,\pm}_{\star}$ and thus $\widehat{V}_{\star}^{2}(0)$.
Therefore, the above equality can be reduced to 
\begin{align*}
	\widehat{E}_{3,\star}^{2}(0) |_{y_3=Y_3=0} 
	= &  \int_{-\infty}^{0} (\nu \dt \widehat{H}_{1,\star}^{2} + \partial_{y_2} \widehat{E}_{3,\star}^{2} )(0) |_{y_3=0} \, dY_3 \\
	& - \mathbf{c_0} \{ \partial_{\theta}\varphi^{3} (\nu \tau  H^{1}_{1,\star} +\xi_2 	  E^{1}_{3,\star} ) + \partial_{\theta}\varphi^{2} (\nu \tau  H^{2}_{1,\star} +\xi_2 	  E^{2}_{3,\star} )\\
	& 
	+  (\nu \partial_{t} \varphi^{2}  H^{1}_{1,\star} +\partial_{y_2}  \varphi^{2}	 E^{1}_{3,\star} )
	\} |_{y_3=Y_3=0}.
\end{align*}
At the same time, 
\begin{align*}
	\widehat{G}_4^{2}(0)&=(-E^0_3\partial_{y_2} -\nu H^0_1 \dt) \widehat{\varphi}^{3} (0) -  \mathbf{c_0} \{ \partial_{\theta}\varphi^{3} (\nu \tau  H^{1}_{1} +\xi_2 	  E^{1}_{3} ) \\
	& + \partial_{\theta}\varphi^{2} (\nu \tau  H^{2}_{1} +\xi_2 	  E^{2}_{3} )
	+  (\nu \partial_{t} \varphi^{2}  H^{1}_{1} +\partial_{y_2}  \varphi^{2}	 E^{1}_{3} )
	\} |_{y_3=Y_3=0}.
\end{align*}
Since $\underline{V}^1$ vanishes, it holds that
\begin{align}
	\underline{\widehat{E}}_2^{3}(0) |_{y_3=Y_3=0}=&(-E^0_3\partial_{y_2} -\nu H^0_1 \dt) \widehat{\varphi}^{3} (0)
	-\int_{-\infty}^{0} (\nu \dt \widehat{H}_{1,\star}^{2} + \partial_{y_2} \widehat{E}_{3,\star}^{2} )(0) |_{y_3=0} \, dY_3\notag\\
	=&(-E^0_3\partial_{y_2} -\nu H^0_1 \dt) \widehat{\varphi}^{3} (0)
	+\int_{-\infty}^{0} \Big(\nu \dt \mathbf{c_0} \{ \partial_{\theta}\varphi^{2} (\xi_1  H^{1}_{3,\star} +\nu\tau  E^{1}_{2,\star} ) |_{y_3=0}\}\notag\\
	&\quad 
	- \partial_{y_2} \mathbf{c_0} \{ \partial_{\theta}\varphi^{2} (\xi_j  E^{1}_{j,\star} ) |_{y_3=0}\} \Big) \, dY_3  \notag\\
	=&(-E^0_3\partial_{y_2} -\nu H^0_1 \dt) \widehat{\varphi}^{3} (0)+ \Big( -\frac{\nu\tau a^-_1-\xi_1b^-}{\sqrt{1-\nu^2\tau^2}} \nu \partial_t \notag	\\
	&\quad
	-  \frac{\xi_1a^-_2-\xi_2a^-_1}{\sqrt{1-\nu^2\tau^2}} \partial_{y_2} \Big) \sum\limits_{k\neq 0} |k| \widehat{\varphi}^2(-k) \widehat{\varphi}^2(k). 	\label{s-boundary3-simpE2}
\end{align}
Similarly, it holds that
\begin{align}
%	\underline{\widehat{H}}_3^{3}(0) |_{y_3=Y_3=0}=&H^0_j\partial_{y_j} \widehat{\varphi}^{3} (0)
%	+\int_{-\infty}^{0} ( \partial_{y_j} \widehat{H}_{j,\star}^{2} )(0) |_{y_3=0} \, dY_3  \notag \\
%	=&H^0_j\partial_{y_j} \widehat{\varphi}^{3} (0)
%	+ \Big( \frac{\nu\tau a^-_j-\xi_jb^-}{\sqrt{1-\nu^2\tau^2}} \partial_{y_j} 	\Big) \sum\limits_{k\neq 0} |k| \widehat{\varphi}^2(-k) \widehat{\varphi}^2(k), \label{s-boundary3-simpH3}\\
	%%%
	-\underline{\widehat{E}}_1^{3}(0) |_{y_3=Y_3=0}=&(E^0_3\partial_{y_1} -\nu H^0_2 \dt) \widehat{\varphi}^{3} (0)
	-\int_{-\infty}^{0} (\nu \dt \widehat{H}_{2,\star}^{2} - \partial_{y_1} \widehat{E}_{3,\star}^{2} )(0) |_{y_3=0} \, dY_3 \notag \\
	=&(E^0_3\partial_{y_1} -\nu H^0_2 \dt) \widehat{\varphi}^{3} (0)
	+ \Big(- \frac{\nu\tau a^-_2-\xi_2b^-}{\sqrt{1-\nu^2\tau^2}} \nu\partial_t \notag	\\
	&\quad
	+   \frac{\xi_1a^-_2-\xi_2a^-_1}{\sqrt{1-\nu^2\tau^2}} \partial_{y_1} \Big) \sum\limits_{k\neq 0} |k| \widehat{\varphi}^2(-k) \widehat{\varphi}^2(k). \label{s-boundary3-simpE1}
\end{align}

\subsubsection{Contradiction arguments}
Suppose that $\underline{V}^3$ vanishes, and then \textbf{the right-hand side of \cref{s-boundary3-simpE2} and \cref{s-boundary3-simpE1} both have to vanish}.
As a result, Applying $(E^0_3\partial_{y_1} -\nu H^0_2 \dt)$ to \cref{s-boundary3-simpE2}, $(E^0_3\partial_{y_2}+\nu H^0_1 \dt)$ to \cref{s-boundary3-simpE1}, and summing them up, we have
\begin{equation} \label{eq-phi2-rectification}
	\begin{aligned}
	&	\Big\{(\xi_2 H^0_1-\xi_1 H^0_2+\nu\tau E^0_3)b^- \nu^2 \partial_{t} \partial_{t} + 
	    \big(-\nu\tau \xi_2 (H^0_1)^2+\nu\tau \xi_1 H^0_1H^0_2 \\
	&\qquad -(\nu^2\tau^2+\xi_2^2)  H^0_1E^0_3 +\xi_1\xi_2H^0_2E^0_3 -\xi_2 (E^0_3)^2 \big) \nu\partial_{t}\partial_{y_1}\\
	&\quad + (-\nu\tau \xi_2 H^0_1H^0_2 +\nu\tau \xi_1 (H^0_2)^2 +\xi_1\xi_2 H^0_1E^0_3-(\nu^2\tau^2+\xi_1^2) H^0_2E^0_3 \\
	&\qquad+\xi_1 (E^0_3)^2 )\nu\partial_{t}\partial_{y_2} \Big\} \sum\limits_{k\neq 0} |k| \widehat{\varphi}^2(-k) \widehat{\varphi}^2(k) =0.
	\end{aligned}
\end{equation}
Obviously, this is quite different from the equation of $\varphi^2$, \cref{eq-HJ}, and thus it does not hold for general $\varphi^2_0$. This gives a contradiction to the assumption of this paragraph.

\vline

In summary, without any more requirement on $\varphi^2_0$ other than \cref{thm-main}, $\varphi^2$ constructed by the Hamilton-Jacobi equation \cref{eq-HJ} will not satisfy \cref{eq-phi2-rectification}, and thus the rectification phenomenon occurs at the second corrector, i.e., $\underline{U}^3=\underline{\widehat{U}}^3(0)\neq0$ or $\underline{V}^3=\underline{\widehat{V}}^3(0)\neq0$. The proof of \cref{thm-main} is complete.

%%%%%%%%%%%%%%%%%%%%%%%%%%%%%%%%%%%%%%%%%%%%%%%%%%%%%%%%%%%%%%%%
%%%%%Appendix 
%%%%%%%%%%%%%%%%%%%%%%%%%%%%%%%%%%%%%%%%%%%%%%%%%%%%%%%%%%%%%%%%
\newpage
\appendix
\section{Basic algebraic calculations}\label{appen-a}
In this appendix we give the explicit expressions of some quantities that are involved in the analysis of the WKB cascade.

\subsection{Matrices}
The fluxes $f_\alpha^+$ are defined by
\begin{equation*}
\label{f+}
f_1^+:=\begin{pmatrix}
	u_1^2 - B_1^2 +q\\
	u_1u_2-B_1B_2\\
	u_1u_3-B_1B_3\\
	0\\
	u_1B_2-B_1u_2\\
	u_1B_3-B_1u_3\\
	u_1
\end{pmatrix},\quad 
f_2^+:=\begin{pmatrix}
	u_2u_1-B_2B_1\\
	u_2^2-B_2^2+q\\
	u_2u_3-B_2B_3\\
	u_2B_1-B_2u_1\\
	0\\
	u_2B_3-B_2u_
	3\\
	u_2
\end{pmatrix},\quad
f_3^+:=\begin{pmatrix}
	u_3u_1-B_3B_1\\
	u_3u_2-B_3B_2\\
	u_3^2-B_3^2+q\\
	u_3B_1-B_3u_1\\
	u_3B_2-B_3u_2\\
	0\\
	u_3
\end{pmatrix}.
\end{equation*}
Therefore, the Jacobian matrices $A^+_\alpha := \text{d}  f^+_{\alpha} (U^0)$, and $\sA^\pm:=\tau A^\pm_0 +\xi_j A^\pm_j$ are 
\begin{align*}
&A^+_1=\begin{pmatrix}
2u^0_1 &   0   &   0   & -2B^0_1 &   0   &   0   &   1   \\
 u^0_2 & u^0_1 &   0   & - B^0_2 & -B^0_1&   0   &   0   \\
   0   &   0   & u^0_1 &   0     &   0   & -B^0_1&   0   \\
   0   &   0   &   0   &   0     &   0   &   0   &   0   \\
 B^0_2 & -B^0_1&   0   & -u^0_2  & u^0_1 &   0   &   0   \\
   0   &   0   & -B^0_1&   0     &   0   & u^0_1 &   0   \\
   1   &   0   &   0   &   0     &   0   &   0   &   0   
\end{pmatrix},\\
&A^+_2=\begin{pmatrix}
	u^0_2 &  u^0_1  &  0  & -B^0_2 & -B^0_1 &   0   &   0   \\
	0   & 2u^0_2 &   0   &   0     & -2B^0_2&   0   &   1   \\
	0   &   0    & u^0_2 &   0     &   0    & -B^0_2&   0   \\
 -B^0_2 &  B^0_1 &   0   &  u^0_2  & -u^0_1 &   0   &   0   \\
	0   &   0    &   0   &   0     &   0    &   0   &   0   \\
	0   &   0    & -B^0_1&   0     &   0    & u^0_2 &   0   \\
	0   &   1    &   0   &   0     &   0    &   0   &   0   
\end{pmatrix},\\
&A^+_3=\begin{pmatrix}
0   &	0   & u^0_1 &   0   &   0   & -B^0_1 &   0   \\
0   &	0   & u^0_2 &   0   &   0   & -B^0_2 &   0   \\
0   &	0   &   0   &   0   &   0   &   0    &   1   \\
0   &	0   & B^0_1 &   0   &   0   & -u^0_1 &   0   \\
0   &	0   & B^0_2 &   0   &   0   & -u^0_2 &   0   \\
0   &	0   &   0   &   0   &   0   &   0    &   0   \\
0   &	0   &   1   &   0   &   0   &   0    &   0   
\end{pmatrix},\\
&\sA^+=\begin{pmatrix}
c^++\xi_1u^0_1&    \xi_2u^0_1&   0   &-(b^++\xi_1B^0_1) &    -\xi_2B^0_1   &   0   &\xi_1 \\
\xi_1u^0_2&c^++\xi_2u^0_2&   0   &    -\xi_1B^0_2   &-(b^++\xi_2B^0_1) &   0   &\xi_2 \\
0     &    	0  	     &  c^+  &   	0  	        &   	0  	       & -b^+  &   0  \\
-\xi_2B^0_2&    \xi_2B^0_1&   0   &  c^+-\xi_1u^0_1  &  -\xi_2u^0_1     &   0   &   0  \\
\xi_1B^0_2&   -\xi_1B^0_1&   0   &  -\xi_1u^0_2     &  c^+-\xi_2u^0_2  &   0   &   0  \\
0  	  &    	0        & -b^+  &   	0  	        &   	0  	       &  c^+  &   0  \\
\xi_1	  &    	\xi_2    &   0   &   	0  	        &   	0  	       &   0   &   0   
\end{pmatrix}, \\
&A^-_0=\nu \text{I}_6   \,,
\\ %change a line here 
&A^-_1=\begin{pmatrix}
0
&
\begin{matrix}
0 &  0 & 0\\
0 &  0 & -1\\
0 &  1 & 0
\end{matrix}       \\
\begin{matrix}
0 &  0 & 0\\
0 &  0 & 1\\
0 & -1 & 0
\end{matrix} 
&
0
\end{pmatrix},
\quad %another
A^-_2=\begin{pmatrix}
0
&
\begin{matrix}
0  &  0 & 1\\
0  &  0 & 0\\
-1 &  0 & 0
\end{matrix}       \\
\begin{matrix}
0  &  0 & -1\\
0  &  0 & 0\\
1  &  0 & 0
\end{matrix} 
&
0
\end{pmatrix},
\\ %change a line here
&A^-_3=\begin{pmatrix}
0
&
\begin{matrix}
0 & -1 & 0\\
1 &  0 & 0\\
0 &  0 & 0\\
\end{matrix}       \\
\begin{matrix}
0  &  1 & 0\\
-1 &  0 & 0\\
0  &  0 & 0\\
\end{matrix}  
&
0
\end{pmatrix},\quad %another
\sA^-=\begin{pmatrix}
\nu\tau &  0 & 0 &0      &      0	& \xi_2\\
0  & \nu\tau & 0 &0      &      0   &-\xi_1\\
0  &  0 & \nu\tau& -\xi_2&  \xi_1   & 0  \\
0      &      0	 &-\xi_2 &\nu\tau &  0 & 0\\
0      &      0  & \xi_1 &0  & \nu\tau & 0\\
\xi_2  & -\xi_1  & 0     &0  &  0 & \nu\tau
\end{pmatrix} .
\end{align*}

\subsection{Eigenvectors}
By $\tau\neq0$, $(c^+)^2\neq (b^+)^2$ and 
$\nu\ll1$ (see \cref{sec-frequencies}), one can verify that the vector spaces $\text{Ker} (-A_3^+ +i\sA^+)$ and  $\text{Ker} (\sqrt{1-\nu^2\tau^2}A_3^- +i\sA^-)$ are of dimension 1 and 2, respectively, and they are spanned by the vectors
\begin{align*}
	&\textbf{R}^+:= \Big(\xi_1 c^+, \xi_2 c^+, ic^+, \xi_1 b, \xi_2 b, ib^+, (b^+)^2-(c^+)^2\Big)^\top \in\mathbb{C}^7,\\
	&\textbf{R}^-_1:= \Big( -i\nu\tau\sqrt{1-\nu^2\tau^2}, 0, -\nu\tau \xi_1, -\xi_1\xi_2, \nu^2 \tau^2 -\xi_2^2 , i\sqrt{1-\nu^2\tau^2}\xi_2 \Big)^\top \in\mathbb{C}^6,\\
	&\textbf{R}^-_2:= \Big( 0, -i\nu\tau\sqrt{1-\nu^2\tau^2},  -\nu\tau \xi_2,  -(\nu^2 \tau^2 -\xi_1^2) ,\xi_1\xi_2, -i\sqrt{1-\nu^2\tau^2}\xi_1 \Big)^\top \in\mathbb{C}^6,
\end{align*}
Note that the vectors $\textbf{R}^+, \textbf{R}^-_1, \textbf{R}^-_2$ satisfy the additional conditions:
\begin{align*}
	&	\begin{pmatrix}
		0&0&0&i\xi_1&i\xi_2& -1&0
	\end{pmatrix}~\textbf{R}^+ =0,\\
	&   \begin{pmatrix}
		i\xi_1&i\xi_2& -\sqrt{1-\nu^2\tau^2}&0&0&0\\
		0&0&0&i\xi_1&i\xi_2& -\sqrt{1-\nu^2\tau^2}
	\end{pmatrix} ~ \begin{pmatrix}
	\textbf{R}^-_1 & \textbf{R}^-_2
\end{pmatrix}=\mathbf{0},
\end{align*}
which fulfills the divergent-free constraints of the magnectic and electric fields.

Therefore, $\sR^+(k),\sR^-_j(k)$ satisfying
$$
(-A_3^+ +i \sgnk \sA^+ ) \sR^+=0,\quad
(\sqrt{1-\nu^2\tau^2}A_3^- +i \sgnk \sA^- ) \sR^-_j=0,~~j=1,2,
$$
can be defined as 
$$
\forall k\in \mathbb{Z}\setminus \{ 0\},~ \sR^+(k):=\begin{cases}
	\textbf{R}^+, &\text{if~}k>0,\\
	\overline{\textbf{R}^+}, &\text{if~}k<0,
\end{cases} 
~ \sR^-_j(k):=\begin{cases}
	\textbf{R}^-_j, &\text{if~}k>0,\\
	\overline{\textbf{R}^-_j}, &\text{if~}k<0,
\end{cases} ~j=1,2.
$$
More explicitly, 
\begin{align}
	\label{vector-R}
	\begin{aligned}
		&\sR^+(k)=(\xi_1 c^+, \xi_2 c^+,i\sgnk c^+,\xi_1 b^+,\xi_2 b^+,i\sgnk b^+,(b^+)^2-(c^+)^2)^\top,\\
		&\sR^-_1(k)=( \sqrt{1-\nu^2\tau^2 }\nu\tau, 0,-i\sgnk\nu\tau\xi_1,-i\sgnk\xi_1\xi_2,i\sgnk(\nu^2\tau^2-\xi_2^2),-\sqrt{1-\nu^2\tau^2 } \xi_2)^\top,\\
		&\sR^-_2(k)=(0,\sqrt{1-\nu^2\tau^2 }\nu\tau, -i\sgnk\nu\tau\xi_2,-i\sgnk(\nu^2\tau^2-\xi_1^2),i\sgnk\xi_1\xi_2, \sqrt{1-\nu^2\tau^2 } \xi_1)^\top.
	\end{aligned}
\end{align}
For $\iota_j$ defined in \cref{def-gamma-special}, we define 
	\begin{align}
		&\sR^-(k):=\iota_j\sR^-_1(k) 
		\notag \\
		\,=&( \frac{\nu\tau a_1^--\xi_1 b^-}{\sqrt{1-\nu^2\tau^2 }},
		\frac{\nu\tau a_2^--\xi_2 b^-}{\sqrt{1-\nu^2\tau^2 }},  
		i\sgnk b^-, i\sgnk a_2^-, -i\sgnk a_1^-,
		\frac{\xi_1 a_2^--\xi_2 a_1^-}{\sqrt{1-\nu^2\tau^2 }})^\top, \label{vector-R3}
	\end{align} 
	where the relation $\xi_j a^-_j=\nu\tau b^-$ are frequently used.
Similarly, we can obtain $\sL^+(k),\sL^-_j(k)$ satisfying
$$
(-A_3^+ +i \sgnk \sA^+ )^\top \sL^+=0,\quad
(\sqrt{1-\nu^2\tau^2}A_3^- +i \sgnk \sA^- )^\top \sL^-_j=0,~~j=1,2,
$$
and we define $\sL^-(k)$ as the linear combination of $\sL^-_j(k)$, $j=1,2$, satisfying  
$$  \nu\tau \sL^-(k)_5=a^-_1\sL^+(k)_3,\quad\text{and}~
-\nu\tau \sL^-(k)_4 =a^-_2 \sL^+(k)_3.$$
More explicitly, 
\begin{align}
	\label{vector-L}
	\begin{aligned}
		\sL^+(k)=&\left(\xi_1 \tau,\xi_2 \tau, i\sgnk\tau, 2\xi_1 b^+,2\xi_2 b^+, 2i\sgnk b^+,-\tau(a^+ + c^+)   \right)^\top \\
		\sL^-(k)=&\frac{1}{\nu^2\tau(1-\nu^2\tau^2)}\Big( (-\nu\tau a_1^- + \xi_1 b^-) \sR^-_1(k)+(-\nu\tau a_2^- + \xi_2 b^-)\sR^-_2(k) \Big) \\
		=&-\frac{1}{\nu}( \frac{\nu\tau a_1^--\xi_1 b^-}{\sqrt{1-\nu^2\tau^2 }},
			\frac{\nu\tau a_2^--\xi_2 b^-}{\sqrt{1-\nu^2\tau^2 }},  
			i\sgnk b^-,\\
		&\qquad\qquad i\sgnk a_2^-, -i\sgnk a_1^-,
			\frac{\xi_1 a_2^--\xi_2 a_1^-}{\sqrt{1-\nu^2\tau^2 }})^\top.
		%&\frac{1}{\nu\sqrt{1-\nu^2\tau^2} } ( -\nu\tau a_1^- + \xi_1 b^-, -\nu\tau a_2^- + \xi_2 b^-, -i\sgnk \sqrt{1-\nu^2\tau^2}b^-,\\
		%& -i\sgnk\sqrt{1-\nu^2\tau^2}a_2^-, i\sgnk \sqrt{1-\nu^2\tau^2}a_1^-, \xi_2 a^-_1 -\xi_1 a^-_2
		%)^{\top}.
	\end{aligned}
\end{align}

\vspace{0.5cm}

\subsection{Bilinear algebra}
\begin{align*}
\sL^-(k_1) \cdot A_0^- \sR^-_1(k_2) & =-2\nu^2 \tau^2 a^-_1 + (1+\sgnkone\sgnktwo) \big(
(\nu^2 \tau^2 -\xi_2^2) a^-_1 +\xi_1\xi_2 a^-_2 +\nu\tau \xi_1 b^-
\big),
\\
\sL^-(k_1) \cdot A_1^- \sR^-_1(k_2) & =\frac{1}{\nu}\{\nu\tau \xi_1 a^-_1-\nu\tau \xi_2a^-_2+\nu^2 \tau^2 b^-  \\
&\quad - (1+\sgnkone\sgnktwo)\big(
\nu\tau \xi_1 a^-_1  +(\nu^2 \tau^2 -\xi_2^2)b^-
\big)\},
\\
\sL^-(k_1) \cdot A_2^- \sR^-_1(k_2) & =\frac{1}{\nu}\{ 2\nu \tau \xi_2 a^-_1 - (1+\sgnkone\sgnktwo) \big(
\nu \tau \xi_1 a^-_2 +\xi_1\xi_2 b^-
\big) \}
\\
\sL^-(k_1) \cdot \sA^- \sR^-_1(k_2) & =\{1 - \sgnkone\sgnktwo \} \tau(1-\nu^2\tau^2) a^-_1,\\
\sL^-(k_1) \cdot A_0^- \sR^-_2(k_2) & =-2\nu^2 \tau^2 a^-_2 + (1+\sgnkone\sgnktwo) \big(
\xi_1\xi_2  a^-_1 +(\nu^2 \tau^2 -\xi_1^2)a^-_2 +\nu\tau \xi_2 b^-
\big),
\\
\sL^-(k_1) \cdot A_1^- \sR^-_2(k_2) & =\frac{1}{\nu}\{ 2\nu\tau \xi_1 a^-_2 - (1+\sgnkone\sgnktwo)  \big(
\nu \tau \xi_2 a^-_1 +\xi_1\xi_2 b^-
\big)\},
\\
\sL^-(k_1) \cdot A_2^- \sR^-_2(k_2) & =\frac{1}{\nu}\{-\nu \tau \xi_1 a^-_1 +\nu\tau \xi_2 a^-_2+\nu^2 \tau^2 b^- \\
&\quad- (1+\sgnkone\sgnktwo) \big(
\nu\tau \xi_2 a^-_2  +(\nu^2 \tau^2 -\xi_1^2)b^-
\big)\}
\\
\sL^-(k_1) \cdot \sA^- \sR^-_2(k_2) & =\{1 - \sgnkone\sgnktwo \} \tau(1-\nu^2\tau^2) a^-_2,
\end{align*}

%%%%%%%%%%%%%%%%%%%%%%%%%%%%%%%%%%%%%%%%%%%%%%%%%%%%%%%%%%%%%%%%
%%%%%Appendix - compatibility conditions
%%%%%%%%%%%%%%%%%%%%%%%%%%%%%%%%%%%%%%%%%%%%%%%%%%%%%%%%%%%%%%%%
\section{Proofs of \cref{lem-comp-am} and \cref{lem-comp-dm}}\label{appen-proofs}

\begin{proof}[Proof of \cref{lem-comp-am}]
	The proof of the equation for the plasma is the same as in \cite{Pierre2021}. We just consider the equation for the vacuum. By \cref{Fm-} one has
	\begin{align}
		F^{m,-}_3 =& - ( \nu \dt H^m_3 +\partial_{y_1}E^m_2 - \partial_{y_2} E^m_1 ) \notag \\
		&\quad+ \sum\limits_{\substack{l_1+l_2=m+2\\l_1\geq 1}} \uline{\partial_{\theta}\varphi^{l_2} (\nu\tau \partial_{Y_3} H^{l_1}_3 +\xi_1\partial_{Y_3} E^{l_1}_2-\xi_2\partial_{Y_3} E^{l_1}_1)} \notag \\
		&\quad+ \sum\limits_{\substack{l_1+l_2=m+1\\l_1\geq 1}} \uwave{( \dt \varphi^{l_2} \nu \partial_{Y_3} H^{l_1}_3 +\partial_{y_1}\varphi^{l_2}\partial_{Y_3} E^{l_1}_2  -\partial_{y_2}\varphi^{l_2}\partial_{Y_3} E^{l_1}_1 )} \notag \\
		&\quad+ \sum\limits_{\substack{l_1+l_2=m+1\\l_1\geq 1}}\uline{ \partial_{\theta}\varphi^{l_2}(\nu\tau \partial_{y_3} H^{l_1}_3 +\xi_1\partial_{y_3} E^{l_1}_2 -\xi_2\partial_{y_3} E^{l_1}_1)} \notag \\
		&\quad+ \sum\limits_{\substack{l_1+l_2=m\\l_1\geq 1}}\uwave{ ( \dt \varphi^{l_2} \nu \partial_{y_3} H^{l_1}_3 +\partial_{y_1}\varphi^{l_2}\partial_{y_3} E^{l_1}_2  -\partial_{y_2}\varphi^{l_2}\partial_{y_3} E^{l_1}_1 )}  \label{Fm-3}.
	\end{align}
	On the boundary $y_3=Y_3=0$, substituting $H^m_3, E^m_1, E^m_2$ by the boundary conditions in \cref{fm} gives
	\begin{align*}
		&  - ( \nu \dt H^m_3 +\partial_{y_1}E^m_2 - \partial_{y_2} E^m_1 )\\
		=& 
		- \nu b^- \partial_{\theta} \dt \varphi^{m+1}  +a^-_j \partial_{\theta} \partial_{y_j} \varphi^{m+1}
		\\
		&
		-\sum\limits_{\substack{l_1+l_2=m+1\\l_1\geq 1}} \dt \{ \partial_{\theta}  \varphi^{l_2} \nu \xi_j H^{l_1}_j \}-\sum\limits_{\substack{l_1+l_2=m\\l_1\geq 1}} \dt \{ \partial_{y_j}  \varphi^{l_2} \nu  H^{l_1}_j \}
		\\
		&
		+\sum\limits_{\substack{l_1+l_2=m+1\\l_1\geq 1}} \partial_{y_1} \{\partial_{\theta}  \varphi^{l_2} (\nu \tau  H^{l_1}_1 +\xi_2 E^{l_1}_3  )\} +\sum\limits_{\substack{l_1+l_2=m\\l_1\geq 1}} \partial_{y_1} \{ \dt  \varphi^{l_2} \nu   H^{l_1}_1 + \partial_{y_2}  \varphi^{l_2}  E^{l_1}_3 \}
		\\
		&
		+\sum\limits_{\substack{l_1+l_2=m+1\\l_1\geq 1}} \partial_{y_2} \{\partial_{\theta}  \varphi^{l_2} (\nu \tau  H^{l_1}_2 -\xi_1 E^{l_1}_3  )\} +\sum\limits_{\substack{l_1+l_2=m\\l_1\geq 1}} \partial_{y_2} \{ \dt  \varphi^{l_2} \nu   H^{l_1}_2 - \partial_{y_1}  \varphi^{l_2}  E^{l_1}_3 \}
		\\
		=&
		- \nu b^- \partial_{\theta} \dt \varphi^{m+1}  +a^-_j \partial_{\theta} \partial_{y_j} \varphi^{m+1}
		\\
		&
		+\sum\limits_{\substack{l_1+l_2=m+1\\l_1\geq 1}}   \{ \partial_{\theta} \dt \varphi^{l_2} (-\nu \xi_j H^{l_1}_j) + \partial_{\theta} \partial_{y_1}  \varphi^{l_2} (\nu \tau  H^{l_1}_1 +\xi_2 E^{l_1}_3  ) + \partial_{\theta} \partial_{y_2} \varphi^{l_2} (\nu \tau  H^{l_1}_2 -\xi_1 E^{l_1}_3  ) \}
		\\
		&
		+\sum\limits_{\substack{l_1+l_2=m+1\\l_1\geq 1}}\uline{  \partial_{\theta}  \varphi^{l_2} \{  (-\nu \xi_j \dt H^{l_1}_j) +   (\nu \tau \partial_{y_1} H^{l_1}_1 +\xi_2 \partial_{y_1} E^{l_1}_3  ) +  (\nu \tau \partial_{y_2} H^{l_1}_2 -\xi_1 \partial_{y_2} E^{l_1}_3  ) \}}
		\\
		&
		+\sum\limits_{\substack{l_1+l_2=m\\l_1\geq 1}}\uwave{  \{   \dt \varphi^{l_2} (\nu \partial_{y_j}    H^{l_1}_j) +  \partial_{y_1}  \varphi^{l_2} ( -\nu \dt H^{l_1}_1 -\partial_{y_2} E^{l_1}_3  ) +  \partial_{y_2}  \varphi^{l_2} ( -\nu \dt H^{l_1}_2 +\partial_{y_1} E^{l_1}_3  )
			\}	}.
	\end{align*}
	On the other hand, the definition \eqref{Gm} of $G^m$ implies that 
	\begin{align*}
		&\nu\tau \partial_{\theta}G_3^{m}+\xi_j\partial_{\theta}G_{3+j}^{m}\\
		=&
		- \nu b^- \partial_{\theta} \dt \varphi^{m+1}  +a^-_j \partial_{\theta} \partial_{y_j} \varphi^{m+1}
		\\
		&
		+\sum\limits_{\substack{l_1+l_2=m+1\\l_1\geq 1}}  \{\partial_{\theta} \dt \varphi^{l_2} (-\nu \xi_j H^{l_1}_j)\uuline{ +   \dt \varphi^{l_2} (-\nu \xi_j \partial_{\theta} H^{l_1}_j) }  \}
		\\
		&
		+\sum\limits_{\substack{l_1+l_2=m+1\\l_1\geq 1}} \{ \partial_{\theta} \partial_{y_1}  \varphi^{l_2} (\nu \tau  H^{l_1}_1 +\xi_2 E^{l_1}_3  ) \uuline{+  \partial_{y_1}  \varphi^{l_2} (\nu \tau \partial_{\theta} H^{l_1}_1 +\xi_2 \partial_{\theta} E^{l_1}_3  ) } \}
		\\
		&
		+\sum\limits_{\substack{l_1+l_2=m+1\\l_1\geq 1}}  \{\partial_{\theta} \partial_{y_2} \varphi^{l_2}(\nu\tau  H^{l_1}_2 -\xi_1  E^{l_1}_3  ) \uuline{ + \partial_{y_2}\varphi^{l_2} (\nu \tau \partial_{\theta} H^{l_1}_2 -\xi_1 \partial_{\theta} E^{l_1}_3  ) } \}.
	\end{align*}
	Here the sums in $G^m_3,G^m_4,G^m_5$ with ``$l_1+l_2=m+2$'' are cancelled. So in order to prove $\eqref{comp-am}_2$, it suffices to prove that 
	\begin{equation}
		\begin{aligned}
			& \text{~``the underlining part''~}+\text{~ ``the underwaving part ''~}\\
			=& \text{~``the double underlining part''}
			%=& \sum\limits_{\substack{l_1+l_2=m+1\\l_1\geq 1}} \{   \dt \varphi^{l_2} (-\nu \xi_j \partial_{\theta} H^{l_1}_j) 
			%+  \partial_{y_1}  \varphi^{l_2} (\nu \tau \partial_{\theta} H^{l_1}_1 +\xi_2 \partial_{\theta} E^{l_1}_3  )
			%+ \partial_{y_2}\varphi^{l_2} (\nu \tau \partial_{\theta} H^{l_1}_2 -\xi_1 \partial_{\theta} E^{l_1}_3  ) \}.
		\end{aligned}
		\label{comp-am-simp}
	\end{equation}

	It follows from the 1st, 2nd and 7th (divergence-free condition for $H^{\mu}$) equations for the vacuum in \cref{fm} that
	\begin{align*}
		&\text{~ ``the underwaving part ''~}-\text{~``the double underlining part''}\\
		%the 1st equality
		=&\sum\limits_{\substack{l_1+l_2=m\\l_1\geq 1}}  \dt \varphi^{l_2} (\nu F^{l_1,-}_7 +\nu \partial_{y_\alpha} H^{l_1}_\alpha )\\
		&+\sum\limits_{\substack{l_1+l_2=m\\l_1\geq 1}}  \partial_{y_1} \varphi^{l_2} ( F^{l_1,-}_1 +\nu \dt H^{l_1}_1+\partial_{y_2} E^{l_1}_3-\partial_{y_3} E^{l_1}_2 )\\
		&+\sum\limits_{\substack{l_1+l_2=m\\l_1\geq 1}}  \partial_{y_2} \varphi^{l_2} ( F^{l_1,-}_7 +\nu \dt H^{l_1}_2 -\partial_{y_1} E^{l_1}_3-\partial_{y_3} E^{l_1}_1  )\\
		%the 2nd equality
		=&\sum\limits_{\substack{l_1+l_2+l_3=m+2\\l_1\geq 1}}  \dt \varphi^{l_2} \partial_{\theta} \varphi^{l_3} \nu \xi_j \partial_{Y_3} H^{l_1}_j 
		+
		\sum\limits_{\substack{l_1+l_2+l_3=m+1\\l_1\geq 1}} \underbrace{ \dt \varphi^{l_2} \partial_{y_j} \varphi^{l_3} \nu  \partial_{Y_3} H^{l_1}_j   }\\
		&+\sum\limits_{\substack{l_1+l_2+l_3=m+1\\l_1\geq 1}}  \dt \varphi^{l_2} \partial_{\theta} \varphi^{l_3} \nu \xi_j \partial_{y_3} H^{l_1}_j 
		+
		\sum\limits_{\substack{l_1+l_2+l_3=m\\l_1\geq 1}}  \underbrace{\dt \varphi^{l_2} \partial_{y_j} \varphi^{l_3} \nu  \partial_{y_3} H^{l_1}_j  }\\
		&-\sum\limits_{\substack{l_1+l_2+l_3=m+2\\l_1\geq 1}}  \partial_{y_1} \varphi^{l_2} \partial_{\theta} \varphi^{l_3} (\nu \tau \partial_{Y_3} H^{l_1}_1+\xi_2 \partial_{Y_3} E^{l_1}_3) 
		\\
		&\quad\quad -
		\sum\limits_{\substack{l_1+l_2+l_3=m+1\\l_1\geq 1}}  \underbrace{\partial_{y_1} \varphi^{l_2}  (\dt \varphi^{l_3} \nu  \partial_{Y_3} H^{l_1}_1+ \partial_{y_2} \varphi^{l_3} \partial_{Y_3} E^{l_1}_3)  } 
		\\
		&-\sum\limits_{\substack{l_1+l_2+l_3=m+1\\l_1\geq 1}}  \partial_{y_1} \varphi^{l_2} \partial_{\theta} \varphi^{l_3} (\nu \tau \partial_{y_3} H^{l_1}_1+\xi_2 \partial_{y_3} E^{l_1}_3) 
		\\
		&\quad\quad -
		\sum\limits_{\substack{l_1+l_2+l_3=m\\l_1\geq 1}}  \underbrace{\partial_{y_1} \varphi^{l_2}  (\dt \varphi^{l_3} \nu  \partial_{y_3} H^{l_1}_1+ \partial_{y_2} \varphi^{l_3} \partial_{y_3}E^{l_1}_3 )  }
		\\
		&-\sum\limits_{\substack{l_1+l_2+l_3=m+2\\l_1\geq 1}}  \partial_{y_2} \varphi^{l_2} \partial_{\theta} \varphi^{l_3} (\nu \tau \partial_{Y_3} H^{l_1}_2-\xi_1 \partial_{Y_3} E^{l_1}_3) 
		\\
		&\quad\quad -
		\sum\limits_{\substack{l_1+l_2+l_3=m+1\\l_1\geq 1}}  \underbrace{\partial_{y_2} \varphi^{l_2}  (\dt \varphi^{l_3} \nu  \partial_{Y_3} H^{l_1}_2- \partial_{y_1} \varphi^{l_3} \partial_{Y_3}E^{l_1}_3  ) }
		\\
		&-\sum\limits_{\substack{l_1+l_2+l_3=m+1\\l_1\geq 1}}  \partial_{y_2} \varphi^{l_2} \partial_{\theta} \varphi^{l_3} (\nu \tau \partial_{y_3} H^{l_1}_2-\xi_1 \partial_{y_3} E^{l_1}_3) 
		\\
		&\quad\quad -
		\sum\limits_{\substack{l_1+l_2+l_3=m\\l_1\geq 1}}  \underbrace{\partial_{y_2} \varphi^{l_2}  (\dt \varphi^{l_3} \nu  \partial_{y_3} H^{l_1}_2- \partial_{y_1} \varphi^{l_3} \partial_{y_3}E^{l_1}_3  ) } \\
		%the 3rd equality
		=& \sum\limits_{\substack{l_1+l_2+l_3=m+2\\l_1\geq 1}}  \partial_{\theta} \varphi^{l_3}   \{  \dt \varphi^{l_2} ( \nu \xi_j \partial_{Y_3} H^{l_1}_j  ) -  \partial_{y_1} \varphi^{l_2} (\nu \tau \partial_{Y_3} H^{l_1}_1+\xi_2 \partial_{Y_3} E^{l_1}_3)   \\
		&\quad\quad  -  \partial_{y_2} \varphi^{l_2} (\nu \tau \partial_{Y_3} H^{l_1}_2-\xi_1 \partial_{Y_3} E^{l_1}_3)  \}\\
		& + \sum\limits_{\substack{l_1+l_2+l_3=m+1\\l_1\geq 1}}  \partial_{\theta} \varphi^{l_3}   \{  \dt \varphi^{l_2} ( \nu \xi_j \partial_{y_3} H^{l_1}_j  ) -  \partial_{y_1} \varphi^{l_2} (\nu \tau \partial_{y_3} H^{l_1}_1+\xi_2 \partial_{y_3} E^{l_1}_3)   \\
		&\quad\quad  -  \partial_{y_2} \varphi^{l_2} (\nu \tau \partial_{y_3} H^{l_1}_2-\xi_1 \partial_{y_3} E^{l_1}_3)  \},
	\end{align*}
	where the underbrace part are cancelled during the calculation. Moreover, similar substitutions as above give 
	\begin{align*}
		&\text{~ ``the underlining part ''~}\\
		%the 1st equality
		=& \sum\limits_{\substack{l_1+l_2=m+2\\l_1\geq 1}} \partial_{\theta} \varphi^{l_2}  \{ \nu \tau (-\xi_j \partial_{\theta} H^{l_1}_j) +\xi_1 (\nu\tau \partial_{\theta} H^{l_1}_1 +\xi_2 \partial_{\theta} E^{l_1}_3) + \xi_2 (\nu\tau \partial_{\theta} H^{l_1}_2-\xi_1 \partial_{\theta} E^{l_1}_3)  \} \\
		& + \sum\limits_{\substack{l_1+l_2+l_3=m+3\\l_1\geq 1}} \partial_{\theta} \varphi^{l_2}  \partial_{\theta} \varphi^{l_3} \{ \nu \tau (\xi_j \partial_{Y_3} H^{l_1}_j) -\xi_1 (\nu\tau \partial_{Y_3} H^{l_1}_1 +\xi_2 \partial_{Y_3} E^{l_1}_3) \\
		&\quad\quad - \xi_2 (\nu\tau \partial_{Y_3} H^{l_1}_2-\xi_1 \partial_{Y_3} E^{l_1}_3)  \} \\
		& + \sum\limits_{\substack{l_1+l_2+l_3=m+2\\l_1\geq 1}} \partial_{\theta} \varphi^{l_2} \{  \partial_{y_j} \varphi^{l_3}  \nu \tau  \partial_{Y_3} H^{l_1}_j -\xi_1 ( \dt \varphi^{l_3} \nu \partial_{Y_3} H^{l_1}_1 +  \partial_{y_2} \varphi^{l_3} \partial_{Y_3} E^{l_1}_3)\\
		&\quad \quad - \xi_2 (\dt \varphi^{l_3} \nu \partial_{Y_3} H^{l_1}_2 -  \partial_{y_1} \varphi^{l_3} \partial_{Y_3} E^{l_1}_3)  \} \\
		& + \sum\limits_{\substack{l_1+l_2+l_3=m+2\\l_1\geq 1}} \partial_{\theta} \varphi^{l_2}  \partial_{\theta} \varphi^{l_3} \{ \nu \tau (\xi_j \partial_{y_3} H^{l_1}_j) -\xi_1 (\nu\tau \partial_{y_3} H^{l_1}_1 +\xi_2 \partial_{y_3} E^{l_1}_3) \\
		&\quad\quad - \xi_2 (\nu\tau \partial_{y_3} H^{l_1}_2-\xi_1 \partial_{y_3} E^{l_1}_3)  \} \\
		& + \sum\limits_{\substack{l_1+l_2+l_3=m+1\\l_1\geq 1}} \partial_{\theta} \varphi^{l_2} \{  \partial_{y_j} \varphi^{l_3}  \nu \tau  \partial_{y_3} H^{l_1}_j -\xi_1 ( \dt \varphi^{l_3} \nu \partial_{y_3} H^{l_1}_1 +  \partial_{y_2} \varphi^{l_3} \partial_{y_3} E^{l_1}_3)\\
		&\quad \quad - \xi_2 (\dt \varphi^{l_3} \nu \partial_{y_3} H^{l_1}_2 -  \partial_{y_1} \varphi^{l_3} \partial_{y_3} E^{l_1}_3)  \} \\
		%the 2nd equality
		=&  \sum\limits_{\substack{l_1+l_2+l_3=m+2\\l_1\geq 1}} \partial_{\theta} \varphi^{l_2} \{  \partial_{y_j} \varphi^{l_3}  \nu \tau  \partial_{Y_3} H^{l_1}_j -\xi_1 ( \dt \varphi^{l_3} \nu \partial_{Y_3} H^{l_1}_1 +  \partial_{y_2} \varphi^{l_3} \partial_{Y_3} E^{l_1}_3)\\
		&\quad \quad - \xi_2 (\dt \varphi^{l_3} \nu \partial_{Y_3} H^{l_1}_2 -  \partial_{y_1} \varphi^{l_3} \partial_{Y_3} E^{l_1}_3)  \} \\
		& + \sum\limits_{\substack{l_1+l_2+l_3=m+1\\l_1\geq 1}} \partial_{\theta} \varphi^{l_2} \{  \partial_{y_j} \varphi^{l_3}  \nu \tau  \partial_{y_3} H^{l_1}_j -\xi_1 ( \dt \varphi^{l_3} \nu \partial_{y_3} H^{l_1}_1 +  \partial_{y_2} \varphi^{l_3} \partial_{y_3} E^{l_1}_3)\\
		&\quad \quad - \xi_2 (\dt \varphi^{l_3} \nu \partial_{y_3} H^{l_1}_2 -  \partial_{y_1} \varphi^{l_3} \partial_{y_3} E^{l_1}_3)  \},
	\end{align*}
	which is exactly the opposite of the simplification of ``the underwaving part ''$-$``the double underlining part''. 
	Therefore, \eqref{comp-am-simp} is proved, and so is \eqref{comp-am}, that is the compatibility condition \eqref{comp-a} of problem \eqref{fm+1}.
\end{proof}

\begin{proof}[Proof of \cref{lem-comp-dm}]
	The proof of the equation for the plasma is the same as in \cite{Pierre2021}. We just consider one of the equations for the vacuum, $\partial_{Y_3} F_3^{m,-} + \xi_j \partial_{\theta} F_{j}^{m,-}-\nu\tau \partial_{\theta} F_7^{m,-}=0,$ the proof of the other one is similar.
	
	We first express $F_1^{m,-},F_2^{m,-},F_7^{m,-}$ as $F_3^{m,-}$ in \cref{Fm-3}:
	\begin{align}
		%1st eq
		F^{m,-}_1 =& - ( \nu \dt H^m_1 +\partial_{y_2}E^m_3 - \partial_{y_3} E^m_2 )+ \sum\limits_{\substack{l_1+l_2=m+2\\l_1\geq 1}} \partial_{\theta}\varphi^{l_2} (\nu\tau \partial_{Y_3} H^{l_1}_1 +\xi_2\partial_{Y_3} E^{l_1}_3) 
		\notag \\
		&\quad
		+ \sum\limits_{\substack{l_1+l_2=m+1\\l_1\geq 1}} ( \dt \varphi^{l_2} \nu \partial_{Y_3} H^{l_1}_1 +\partial_{y_2}\varphi^{l_2}\partial_{Y_3} E^{l_1}_3 ) 
		\notag \\
		&\quad
		+ \sum\limits_{\substack{l_1+l_2=m+1\\l_1\geq 1}} \partial_{\theta}\varphi^{l_2} (\nu\tau \partial_{y_3} H^{l_1}_1 +\xi_2\partial_{y_3} E^{l_1}_3) 
		\notag \\
		&\quad
		+ \sum\limits_{\substack{l_1+l_2=m\\l_1\geq 1}} ( \dt \varphi^{l_2} \nu \partial_{y_3} H^{l_1}_1 +\partial_{y_2}\varphi^{l_2}\partial_{y_3} E^{l_1}_3 )  \label{Fm-1},\\  
		%2nd eq
		F^{m,-}_2 =& - ( \nu \dt H^m_2 -\partial_{y_1}E^m_3 + \partial_{y_3} E^m_1 )+ \sum\limits_{\substack{l_1+l_2=m+2\\l_1\geq 1}} \partial_{\theta}\varphi^{l_2} (\nu\tau \partial_{Y_3} H^{l_1}_2 -\xi_1\partial_{Y_3} E^{l_1}_3) 
		\notag \\
		&\quad
		+ \sum\limits_{\substack{l_1+l_2=m+1\\l_1\geq 1}} ( \dt \varphi^{l_2} \nu \partial_{Y_3} H^{l_1}_2 -\partial_{y_1}\varphi^{l_2}\partial_{Y_3} E^{l_1}_3 ) \notag \\
		&\quad
		+ \sum\limits_{\substack{l_1+l_2=m+1\\l_1\geq 1}} \partial_{\theta}\varphi^{l_2} (\nu\tau \partial_{y_3} H^{l_1}_2 -\xi_1\partial_{y_3} E^{l_1}_3) 
		\notag \\
		&\quad
		+ \sum\limits_{\substack{l_1+l_2=m\\l_1\geq 1}} ( \dt \varphi^{l_2} \nu \partial_{y_3} H^{l_1}_2 -\partial_{y_1}\varphi^{l_2}\partial_{y_3} E^{l_1}_3 )  \label{Fm-2},\\
		%3rd eq
		F^{m,-}_7 =&  -( \partial_{y_\alpha} H^m_\alpha )+ \sum\limits_{\substack{l_1+l_2=m+2\\l_1\geq 1}} \partial_{\theta}\varphi^{l_2} ( \xi_j \partial_{Y_3} H^{l_1}_j ) 
		+ \sum\limits_{\substack{l_1+l_2=m+1\\l_1\geq 1}} ( \partial_{y_j} \varphi^{l_2}  \partial_{Y_3} H^{l_1}_j ) 
		\notag \\
		&\quad
		+ \sum\limits_{\substack{l_1+l_2=m+1\\l_1\geq 1}} \partial_{\theta}\varphi^{l_2} ( \xi_j \partial_{y_3} H^{l_1}_j ) 
		+ \sum\limits_{\substack{l_1+l_2=m\\l_1\geq 1}} ( \partial_{y_j} \varphi^{l_2}  \partial_{y_3} H^{l_1}_j )  \label{Fm-7}.
	\end{align}
	
	With these formulas, it follows from \cref{fm} (substituting proper terms by $F^{\nu,-}$) that 
	\begin{align}
		&\partial_{Y_3} F_3^{m,-} + \xi_j \partial_{\theta} F_{j}^{m,-}-\nu\tau \partial_{\theta} F_7^{m,-}  
		\notag\\
		%=& -\nu \dt ( \partial_{Y_3} H^{m}_3 + \xi_j \partial_{\theta} H^m_j ) + \partial_{y_1} (-\partial_{Y_3} E^m_2 + \nu\tau \partial_{\theta} H^m_1 +\xi_2 \partial_{\theta} E^m_3 )   \notag\\
		%& + \partial_{y_2} (\partial_{Y_3} E^m_1 + \nu\tau \partial_{\theta} H^m_2 -\xi_1 \partial_{\theta} E^m_3 ) 
		%\partial_{y_3} ( \nu\tau \partial_{\theta} H^m_3 +\xi_1 \partial_{\theta} E^m_2-\xi_2 \partial_{\theta} E^m_1 )  \notag\\
		%&+ \sum\limits_{\substack{l_1+l_2=m+2\\l_1\geq 1}}  \partial_{\theta}\varphi^{l_2} \partial_{Y_3} (\nu\tau \partial_{Y_3} H^{l_1}_3 +\xi_1\partial_{Y_3} E^{l_1}_2-\xi_2\partial_{Y_3} E^{l_1}_1) \notag \\
		%+ \sum\limits_{\substack{l_1+l_2=m+1\\l_1\geq 1}} \Big\{  \dt \varphi^{l_2}( \nu \partial_{Y_3}^2  H^{l_1}_3 + \nu \xi_j \partial_{Y_3}\partial_{\theta}  H^{l_1}_j  ) -\partial_{y_1}\varphi^{l_2}(-\partial_{Y_3}^2 E^m_2 + \nu\tau \partial_{Y_3}\partial_{\theta} H^m_1 +\xi_2 \partial_{Y_3}\partial_{\theta} E^m_3)  \notag\\
		%&\quad -\partial_{y_2}\varphi^{l_2} (\partial_{Y_3}^2 E^m_1 + \nu\tau \partial_{Y_3}\partial_{\theta} H^m_2 -\xi_1 \partial_{Y_3}\partial_{\theta} E^m_3)
		% \Big\} \notag \\
		&\quad= -\nu \dt F^{m-1,-}_7   + \partial_{y_\alpha} F^{m-1,-}_\alpha  
		\notag \\
		&\quad\quad~ 
		+\sum\limits_{\substack{l_1+l_2=m+1\\l_1\geq 1}}  \partial_{\theta}\varphi^{l_2} \partial_{Y_3} (\nu\tau F^{l_1,-}_7 -\xi_j F^{l_1,-}_j )
		+ \sum\limits_{\substack{l_1+l_2=m\\l_1\geq 1}} \Big\{  \dt \varphi^{l_2}\nu \partial_{Y_3} F^{l_1,-}_7  -\partial_{y_j}\varphi^{l_2} \partial_{Y_3} F^{l_1,-}_j   
		\Big\} 
		\notag \\
		&\quad\quad~
		+ \sum\limits_{\substack{l_1+l_2=m+1\\l_1\geq 1}} \uline{\Big\{  \partial_{\theta}\dt \varphi^{l_2} (\nu \xi_j \partial_{Y_3} H^{l_1}_j) -\partial_{\theta}\partial_{y_1} \varphi^{l_2} (\nu \tau \partial_{Y_3} H^{l_1}_1 +\xi_2 \partial_{Y_3} E^{l_1}_3  )}
		\notag\\
		&\quad\quad~
		\quad \quad 
		\uline{-\partial_{\theta}\partial_{y_2} \varphi^{l_2} (\nu \tau \partial_{Y_3} H^{l_1}_2 -\xi_1 \partial_{Y_3} E^{l_1}_3  )
			\Big\}} 
		\notag \\
		&\quad\quad~
		+\sum\limits_{\substack{l_1+l_2=m\\l_1\geq 1}}  \partial_{\theta}\varphi^{l_2} \partial_{y_3} (\nu\tau F^{l_1,-}_7 -\xi_j F^{l_1,-}_j )
		+ \sum\limits_{\substack{l_1+l_2=m-1\\l_1\geq 1}} \Big\{  \dt \varphi^{l_2}\nu \partial_{y_3} F^{l_1,-}_7  -\partial_{y_j}\varphi^{l_2} \partial_{y_3} F^{l_1,-}_j   
		\Big\} 
		\notag \\
		&\quad\quad~
		+ \sum\limits_{\substack{l_1+l_2=m\\l_1\geq 1}} \uline{\Big\{  \partial_{\theta}\dt \varphi^{l_2} (\nu \xi_j \partial_{y_3} H^{l_1}_j) -\partial_{\theta}\partial_{y_1} \varphi^{l_2} (\nu \tau \partial_{y_3} H^{l_1}_1 +\xi_2 \partial_{y_3} E^{l_1}_3  )}
		\notag\\
		&\quad\quad~
		\quad \quad 
		\uline{-\partial_{\theta}\partial_{y_2} \varphi^{l_2} (\nu \tau \partial_{y_3} H^{l_1}_2 -\xi_1 \partial_{y_3} E^{l_1}_3  )
			\Big\} }
		\label{comp-dm-1}
	\end{align}
	Note that here $\varphi$ is independent of $y_3, Y_3$. By the definition of $F^{m,-}$, i.e. \cref{Fm-}, we can extend the above formula one by one: 
	\begin{align}
		&-\nu \dt F^{m-1,-}_7   + \partial_{y_\alpha} F^{m-1,-}_\alpha  
		\notag \\
		&\quad= \sum\limits_{\substack{l_1+l_2=m+1\\l_1\geq 1}} \uwave{\partial_{\theta} \varphi^{l_2} \Big\{ -\nu \dt  (\xi_j \partial_{Y_3} H^{l_1}_j ) + \partial_{y_1} ( \nu\tau \partial_{Y_3} H^{l_1}_1 +\xi_2 \partial_{Y_3} E^{l_1}_3  ) }
		\notag \\
		&\quad\quad~\quad \quad 
		\uwave{+ \partial_{y_2} ( \nu\tau \partial_{Y_3} H^{l_1}_2 - \xi_1 \partial_{Y_3} E^{l_1}_3  ) 
			+ \partial_{y_3} ( \nu\tau \partial_{Y_3} H^{l_1}_3 +\xi_1 \partial_{Y_3} E^{l_1}_2 - \xi_2\partial_{Y_3} E^{l_1}_1 ) 
			\Big\} } 
		\notag \\
		&\quad\quad~
		+\sum\limits_{\substack{l_1+l_2=m\\l_1\geq 1}} \uuline{\Big\{  
			\dt \varphi^{l_2}(\nu \partial_{Y_3}\partial_{y_{\alpha}} H^{l_1}_{\alpha}  )  - \partial_{y_1} \varphi^{l_2}(\nu \partial_{Y_3} \dt H^{l_1}_{1} +  \partial_{Y_3} \partial_{y_2} E^{l_1}_{3} -\partial_{Y_3} \partial_{y_3} E^{l_1}_{2} ) }
		\notag \\
		&\quad\quad~ \quad \quad 
		\uuline{- \partial_{y_2} \varphi^{l_2}(\nu \partial_{Y_3} \dt H^{l_1}_{1} -  \partial_{Y_3} \partial_{y_1} E^{l_1}_{3} +\partial_{Y_3} \partial_{y_3} E^{l_1}_{1} ) 
			\Big\} }
		\notag \\
		&\quad\quad~
		+ \sum\limits_{\substack{l_1+l_2=m\\l_1\geq 1}} \partial_{\theta} \varphi^{l_2} \Big\{ -\nu \dt  (\xi_j \partial_{y_3} H^{l_1}_j ) + \partial_{y_1} ( \nu\tau \partial_{y_3} H^{l_1}_1 +\xi_2 \partial_{y_3} E^{l_1}_3  ) 
		\notag \\
		&\quad\quad~\quad \quad 
		+ \partial_{y_2} ( \nu\tau \partial_{y_3} H^{l_1}_2 - \xi_1 \partial_{y_3} E^{l_1}_3  ) 
		+ \partial_{y_3} ( \nu\tau \partial_{y_3} H^{l_1}_3 +\xi_1 \partial_{y_3} E^{l_1}_2 - \xi_2\partial_{y_3} E^{l_1}_1 ) 
		\Big\}  
		\notag \\
		& \quad\quad~
		+\sum\limits_{\substack{l_1+l_2=m-1\\l_1\geq 1}} \Big\{  
		\dt \varphi^{l_2}(\nu \partial_{y_3}\partial_{y_{\alpha}} H^{l_1}_{\alpha}  )  - \partial_{y_1} \varphi^{l_2}(\nu \partial_{y_3} \dt H^{l_1}_{1} +  \partial_{y_3} \partial_{y_2} E^{l_1}_{3} -\partial_{y_3} \partial_{y_3} E^{l_1}_{2} ) 
		\notag \\
		&\quad\quad~ \quad \quad 
		- \partial_{y_2} \varphi^{l_2}(\nu \partial_{y_3} \dt H^{l_1}_{1} -  \partial_{y_3} \partial_{y_1} E^{l_1}_{3} +\partial_{y_3} \partial_{y_3} E^{l_1}_{1} ) 
		\Big\} 
		\notag \\
		&\quad\quad~
		- ``\text{the underlining part}'' ,
		\label{comp-dm-2}
	\end{align}
	\begin{align}
		&\sum\limits_{\substack{l_1+l_2=m+1\\l_1\geq 1}}  \partial_{\theta}\varphi^{l_2} \partial_{Y_3} (\nu\tau F^{l_1,-}_7 -\xi_j F^{l_1,-}_j ) 
		\notag\\
		&\quad=
		-``\text{the underwaving part}''
		+\sum\limits_{\substack{l_1+l_2+l_3=m+2\\l_1\geq 1}}\underbrace{
			\Big\{  
			\partial_{\theta}\varphi^{l_2}\dt\varphi^{l_3}  (-\nu\xi_j \partial_{Y_3}^2 H^{l_1}_j)}
		\notag\\
		&\quad\quad~\quad\quad
		\underbrace{
			+\partial_{\theta}\varphi^{l_2}\partial_{y_1}\varphi^{l_3} (\nu\tau \partial_{Y_3}^2 H^{l_1}_1+\xi_2 \partial_{Y_3}^2 E^{l_1}_3)
			+\partial_{\theta}\varphi^{l_2}\partial_{y_2}\varphi^{l_3} (\nu\tau \partial_{Y_3}^2 H^{l_1}_2-\xi_1\partial_{Y_3}^2 E^{l_1}_3)
			\Big\}}
		\notag \\
		&\quad
		+\sum\limits_{\substack{l_1+l_2+l_3=m+1\\l_1\geq 1}}  \underbrace{
			\Big\{  
			\partial_{\theta}\varphi^{l_2}\dt\varphi^{l_3}  (-\nu\xi_j \partial_{Y_3}\partial_{y_3} H^{l_1}_j)
			+\partial_{\theta}\varphi^{l_2}\partial_{y_1}\varphi^{l_3} (\nu\tau \partial_{Y_3}\partial_{y_3} H^{l_1}_1+\xi_2 \partial_{Y_3}\partial_{y_3} E^{l_1}_3)   }
		\notag\\
		&\quad\quad~\quad\quad
		\underbrace{
			+\partial_{\theta}\varphi^{l_2}\partial_{y_2}\varphi^{l_3} (\nu\tau \partial_{Y_3}\partial_{y_3} H^{l_1}_2-\xi_1\partial_{Y_3}\partial_{y_3} E^{l_1}_3)
			\Big\} },
		\label{comp-dm-3}
	\end{align}
	\begin{align}
		&\sum\limits_{\substack{l_1+l_2=m\\l_1\geq 1}} \Big\{  \dt \varphi^{l_2}\nu \partial_{Y_3} F^{l_1,-}_7  -\partial_{y_j}\varphi^{l_2} \partial_{Y_3} F^{l_1,-}_j   
		\Big\} 
		\notag\\
		&\quad=-``\text{the double underlining part}''
		+\sum\limits_{\substack{l_1+l_2+l_3=m+2\\l_1\geq 1}} \underbrace{
			\Big\{  
			\dt\varphi^{l_2}\partial_{\theta}\varphi^{l_3}  (\nu\xi_j \partial_{Y_3}^2 H^{l_1}_j)}
		\notag\\
		&\quad\quad~\quad\quad
		\underbrace{
			-\partial_{y_1}\varphi^{l_2}\partial_{\theta}\varphi^{l_3} (\nu\tau \partial_{Y_3}^2 H^{l_1}_1+\xi_2 \partial_{Y_3}^2 E^{l_1}_3)
			-\partial_{y_2}\varphi^{l_2}\partial_{\theta}\varphi^{l_3} (\nu\tau \partial_{Y_3}^2 H^{l_1}_2-\xi_1\partial_{Y_3}^2 E^{l_1}_3)
			\Big\}}
		\notag \\
		&\quad
		+\sum\limits_{\substack{l_1+l_2+l_3=m+1\\l_1\geq 1}} \underbrace{
			\Big\{  
			\dt\varphi^{l_2}\partial_{\theta}\varphi^{l_3}  (\nu\xi_j \partial_{Y_3}\partial_{y_3} H^{l_1}_j)
			-\partial_{y_1}\varphi^{l_2}\partial_{\theta}\varphi^{l_3} (\nu\tau \partial_{Y_3}\partial_{y_3} H^{l_1}_1+\xi_2 \partial_{Y_3}\partial_{y_3} E^{l_1}_3)}
		\notag\\
		&\quad\quad~\quad\quad
		\underbrace{
			-\partial_{y_2}\varphi^{l_2}\partial_{\theta}\varphi^{l_3} (\nu\tau \partial_{Y_3}\partial_{y_3} H^{l_1}_2-\xi_1\partial_{Y_3}\partial_{y_3} E^{l_1}_3)
			\Big\}}.
		\label{comp-dm-4}
	\end{align}
	The two underbraced parts are canceled. The remaining calculations are similar:
	\begin{align}
		&\sum\limits_{\substack{l_1+l_2=m\\l_1\geq 1}}  \partial_{\theta}\varphi^{l_2} \partial_{y_3} (\nu\tau F^{l_1,-}_7 -\xi_j F^{l_1,-}_j )
		+ \sum\limits_{\substack{l_1+l_2=m-1\\l_1\geq 1}} \Big\{  \dt \varphi^{l_2}\nu \partial_{y_3} F^{l_1,-}_7  -\partial_{y_j}\varphi^{l_2} \partial_{y_3} F^{l_1,-}_j   
		\Big\} 
		\notag \\
		&\quad\quad~
		=- \sum\limits_{\substack{l_1+l_2=m\\l_1\geq 1}} \partial_{\theta} \varphi^{l_2} \Big\{ -\nu \dt  (\xi_j \partial_{y_3} H^{l_1}_j ) + \partial_{y_1} ( \nu\tau \partial_{y_3} H^{l_1}_1 +\xi_2 \partial_{y_3} E^{l_1}_3  ) 
		\notag \\
		&\quad\quad~\quad \quad 
		+ \partial_{y_2} ( \nu\tau \partial_{y_3} H^{l_1}_2 - \xi_1 \partial_{y_3} E^{l_1}_3  ) 
		+ \partial_{y_3} ( \nu\tau \partial_{y_3} H^{l_1}_3 +\xi_1 \partial_{y_3} E^{l_1}_2 - \xi_2\partial_{y_3} E^{l_1}_1 ) 
		\Big\}  
		\notag \\
		& \quad\quad~\quad
		-\sum\limits_{\substack{l_1+l_2=m-1\\l_1\geq 1}} \Big\{  
		\dt \varphi^{l_2}(\nu \partial_{y_3}\partial_{y_{\alpha}} H^{l_1}_{\alpha}  )  - \partial_{y_1} \varphi^{l_2}(\nu \partial_{y_3} \dt H^{l_1}_{1} +  \partial_{y_3} \partial_{y_2} E^{l_1}_{3} -\partial_{y_3} \partial_{y_3} E^{l_1}_{2} ) 
		\notag \\
		&\quad\quad~ \quad \quad 
		- \partial_{y_2} \varphi^{l_2}(\nu \partial_{y_3} \dt H^{l_1}_{1} -  \partial_{y_3} \partial_{y_1} E^{l_1}_{3} +\partial_{y_3} \partial_{y_3} E^{l_1}_{1} ) 
		\Big\} .
		\label{comp-dm-5}
	\end{align}
	Summing up \eqref{comp-dm-1}, \eqref{comp-dm-2}, \eqref{comp-dm-3}, \eqref{comp-dm-4} and \eqref{comp-dm-5} gives
	$$\partial_{Y_3} F_3^{m,-} + \xi_j \partial_{\theta} F_{j}^{m,-}-\nu\tau \partial_{\theta} F_7^{m,-}=0.$$
	This completes the proof of \eqref{comp-dm}, that is the compatibility condition \eqref{comp-d} of problem \eqref{fm+1}.    
\end{proof}

\section{Compatibilities of the slow problem}\label{appen-comp-s}
In this part, we will prove the remaining two compatibilities of the slow problem \cref{s}, \cref{s-boundary} and \cref{s-initial} for the slow means of the correctors in \cref{sec-slow means}.
The proofs are similar as \cref{lem-comp-am} and \cref{lem-comp-dm}, where substitutions coming from the residue components of H((m+1)-2) are frequently used.

\begin{Lem}\label{lem-comp-ams}
	The source terms and boundary terms in \cref{s-source} and \cref{s-boundary} satisfy:
	\begin{equation}
		\label{comp-ams}
		\begin{aligned}
			&\bF_6^{m,+} |_{\Gamma_0}= (\partial_t +u_j^0 \partial_{y_j})\bG^m_2-H_j^0 \partial_{y_j}\bG^m_1, \\ 
			&\bF_3^{m,-} |_{\Gamma_0}= \nu \partial_t \bG^m_3+  \partial_{y_j}\bG^m_{3+j},
		\end{aligned}
	\end{equation}
independently of the choice of the slow mean $\widehat{\varphi}^{m+1}(0)$.
\end{Lem}
\begin{proof}[Proof of \cref{lem-comp-ams}]
	The proof of the equation for the plasma is the same as \cite{Pierre2021}. We just consider the equation for the vacuum.
	Recalling the definition of the boundary terms \cref{s-boundary} and \eqref{Gm}, we need to show that 
	\begin{align}
	    &\bF_3^{m,-} |_{\Gamma_0} +  \Big( \nu \dt \widehat{H}^{m+1}_{3,\star}(0) +\partial_{y_1} \widehat{E}^{m+1}_{2,\star} (0) -\partial_{y_2} \widehat{E}^{m+1}_{1,\star}(0)\Big) |_{y_3=Y_3=0}   \notag\\
\qquad=	&\nu \partial_t \widehat{G}^{m}_3(0)+  \partial_{y_j}\widehat{G}^{m}_{3+j}(0)  \notag\\
\qquad= &\mathbf{c_0} \Big\{
           \sum\limits_{\substack{l_1+l_2=m+2\\l_1\geq 1}} \partial_{\theta} \varphi^{l_2} \Big( 
	\nu \dt (\xi_j H^{l_1}_j) -\partial_{y_1} (\nu\tau H^{l_1}_1 +\xi_2 E^{l_1}_3)-\partial_{y_2} (\nu\tau H^{l_1}_2 -\xi_1 E^{l_1}_3)
	\Big) \notag\\
	    &  +\sum\limits_{\substack{l_1+l_2=m+2\\l_1\geq 1}} \Big(\nu \dt  \partial_{\theta} \varphi^{l_2} 
	      (\xi_j H^{l_1}_j) -\partial_{y_1} \partial_{\theta} \varphi^{l_2} (\nu\tau H^{l_1}_1 +\xi_2 E^{l_1}_3)-\partial_{y_2}\partial_{\theta} \varphi^{l_2}  (\nu\tau H^{l_1}_2 -\xi_1 E^{l_1}_3)
	      \Big) \notag\\
	    &  + \sum\limits_{\substack{l_1+l_2=m+1\\l_1\geq 1}}\Big( \partial_{y_j} \varphi^{l_2}  
	      \nu \dt H^{l_1}_j - (\dt \varphi^{l_2} \nu \partial_{y_1} H^{l_1}_1 +\partial_{y_2} \varphi^{l_2}  \partial_{y_1} E^{l_1}_3) \notag\\
	    &  \qquad-(\dt \varphi^{l_2} \nu \partial_{y_2} H^{l_1}_2 -\partial_{y_1} \varphi^{l_2}  \partial_{y_2} E^{l_1}_3)
	      \Big)
         	\Big\}. \label{comp-ams-1}
	\end{align}
    Here $\mathbf{c_0}$ represents the constant for the mean with respect to $\theta$ on $\mathbb{T}$, namely the integration with respect to $\theta$ on $\mathbb{T}$.
    We will prove that in \cref{comp-ams-1}, the two parts in the first line are the residue and surface wave components of the last four lines, respectively.
    
     \underline{\textbf{Step 1: Residue components.}}
    Recalling the definitions of $F^{m,-}_7, F^{m,-}_j$ (\cref{Fm-7}, \cref{Fm-1} and \cref{Fm-2}), substituting $\partial_{y_3}  \underline{H}^{l_1}_3$, $-\partial_{y_3}  \underline{E}^{l_2}_3$, $\partial_{y_3}  \underline{E}^{l_1}_1$ by $\underline{F}^{l_1}_7 $, $\underline{F}^{l_1}_1 $, $\underline{F}^{l_1}_2 $ and other terms from H((m+1)-2), and integrating by parts with respect to $\theta$, we have 
    \begin{align*}
    	\bF^{m,-}_3=
    	&\mathbf{c_0} \Big\{ \sum\limits_{\substack{l_1+l_2=m+2\\l_1\geq 1}} \partial_{\theta} \varphi^{l_2}  \partial_{y_3}  \Big(\nu\tau \underline{H}^{l_1}_3 +\xi_1 \underline{E}^{l_1}_2 -\xi_2\underline{E}^{l_1}_1 \Big)\\
    	&\qquad
    	+\sum\limits_{\substack{l_1+l_2=m+1\\l_1\geq 1}} \Big(\dt \varphi^{l_2} \nu\partial_{y_3}  \underline{H}^{l_1}_3  + \partial_{y_1} \varphi^{l_2} \partial_{y_3} \underline{E}^{l_1}_2 -\partial_{y_2} \varphi^{l_2} \partial_{y_3} \underline{E}^{l_1}_1   \Big)
    	\Big\}=:I_1+I_2, 
    \end{align*}
    	%2nd
    \begin{align*}
    	I_1=& \mathbf{c_0} \Big\{ 
    	\sum\limits_{\substack{l_1+l_2=m+2\\l_1\geq 1}} \uwave{\partial_{\theta} \varphi^{l_2}   \Big(-\nu\tau (\xi_j \partial_{\theta}  \underline{H}^{l_1+1}_j  )   +\xi_1  (\nu \tau \partial_{\theta}  \underline{H}^{l_1+1}_1 +\xi_2 \partial_{\theta}  \underline{E}^{l_1+1}_3)} \\
    	&\qquad\uwave{+ \xi_2  (\nu \tau \partial_{\theta}  \underline{H}^{l_1+1}_2 -\xi_1 \partial_{\theta}  \underline{E}^{l_1+1}_3) \Big)}\\
        &+\sum\limits_{\substack{l_1+l_2=m+2\\l_1\geq 1}}\partial_{\theta} \varphi^{l_2}   \Big(-\nu\tau ( \partial_{y_j}  \underline{H}^{l_1}_j  )   +\xi_1  (\nu  \dt \underline{H}^{l_1}_1 + \partial_{y_2}  \underline{E}^{l_1}_3) \\
        &\qquad +\xi_2  (\nu  \dt \underline{H}^{l_1}_2 - \partial_{y_1}  \underline{E}^{l_1}_3) \Big) \\
        &+\sum\limits_{\substack{l_1+l_2+l_3=m+2\\l_1\geq 1}}  \partial_{\theta} \varphi^{l_2} \partial_{\theta} \varphi^{l_3}  \uwave{  \Big(\nu\tau (\xi_j \partial_{y_3}  \underline{H}^{l_1}_j  )   -\xi_1  (\nu \tau \partial_{y_3}  \underline{H}^{l_1}_1 +\xi_2 \partial_{y_3}  \underline{E}^{l_1}_3)} \\
        &\qquad\uwave{- \xi_2  (\nu \tau \partial_{y_3}  \underline{H}^{l_1}_2 -\xi_1 \partial_{y_3}  \underline{E}^{l_1}_3) \Big) }\\
        &+\sum\limits_{\substack{l_1+l_2+l_3=m+1\\l_1\geq 1}}\partial_{\theta} \varphi^{l_2}   \Big( \uuline{\nu\tau (\partial_{y_j} \varphi^{l_3}  \partial_{y_3}  \underline{H}^{l_1}_j  )   -\xi_1  (\dt \varphi^{l_3} \nu \partial_{y_3}  \underline{H}^{l_1}_1 +\partial_{y_2} \varphi^{l_3} \partial_{y_3}  \underline{E}^{l_1}_3)} \\
        &\qquad \uwave{ -\xi_2  (\dt \varphi^{l_3} \nu \partial_{y_3}  \underline{H}^{l_1}_2 -\partial_{y_1} \varphi^{l_3} \partial_{y_3}  \underline{E}^{l_1}_3)} \Big)
    	\Big\}\\
    =& \mathbf{c_0} \Big\{ 
    	 \sum\limits_{\substack{l_1+l_2=m+2\\l_1\geq 1}} \partial_{\theta} \varphi^{l_2} \Big( 
    	\nu \dt (\xi_j \underline{H}^{l_1}_j) -\partial_{y_1} (\nu\tau \underline{H}^{l_1}_1 +\xi_2 \underline{E}^{l_1}_3)-\partial_{y_2} (\nu\tau \underline{H}^{l_1}_2 -\xi_1 \underline{E}^{l_1}_3)
    	\Big)\\
    	&+\sum\limits_{\substack{l_1+l_2+l_3=m+1\\l_1\geq 1}}\uuline{ \partial_{\theta} \varphi^{l_2}   \Big( \nu\tau (\partial_{y_j} \varphi^{l_3}  \partial_{y_3}  \underline{H}^{l_1}_j  )   -\xi_1  (\dt \varphi^{l_3} \nu \partial_{y_3}  \underline{H}^{l_1}_1 +\partial_{y_2} \varphi^{l_3} \partial_{y_3}  \underline{E}^{l_1}_3)} \\
    	&\qquad \uuline{ -\xi_2  (\dt \varphi^{l_3} \nu \partial_{y_3}  \underline{H}^{l_1}_2 -\partial_{y_1} \varphi^{l_3} \partial_{y_3}  \underline{E}^{l_1}_3)} \Big)
        \Big\},  
    \end{align*}
      %3nd
    \begin{align*}
    I_2=& \mathbf{c_0} \Big\{ 
      \sum\limits_{\substack{l_1+l_2=m+1\\l_1\geq 1}}   \Big(-\dt \varphi^{l_2} \nu (\xi_j \partial_{\theta}  \underline{H}^{l_1+1}_j  )   +  \partial_{y_1} \varphi^{l_2}   (\nu \tau \partial_{\theta}  \underline{H}^{l_1+1}_1 +\xi_2 \partial_{\theta}  \underline{E}^{l_1+1}_3) \\
        &\qquad+ \partial_{y_2} \varphi^{l_2}   (\nu \tau \partial_{\theta}  \underline{H}^{l_1+1}_2 -\xi_1 \partial_{\theta}  \underline{E}^{l_1+1}_3) \Big)\\
        &+\sum\limits_{\substack{l_1+l_2=m+1\\l_1\geq 1}}   \Big(- \dt \varphi^{l_2}\nu ( \partial_{y_j}  \underline{H}^{l_1}_j  )   +\partial_{y_1} \varphi^{l_2}   (\nu  \dt \underline{H}^{l_1}_1 + \partial_{y_2}  \underline{E}^{l_1}_3) \\
        &\qquad +\partial_{y_2} \varphi^{l_2}   (\nu  \dt \underline{H}^{l_1}_2 - \partial_{y_1}  \underline{E}^{l_1}_3) \Big) \\
        &+\sum\limits_{\substack{l_1+l_2+l_3=m+1\\l_1\geq 1}}   \uuline{  \Big( \dt \varphi^{l_2} \partial_{\theta} \varphi^{l_3}  \nu (\xi_j \partial_{y_3}  \underline{H}^{l_1}_j  )   -\partial_{y_1} \varphi^{l_2} \partial_{\theta} \varphi^{l_3}   (\nu \tau \partial_{y_3}  \underline{H}^{l_1}_1 +\xi_2 \partial_{y_3}  \underline{E}^{l_1}_3)} \\
        &\qquad\uuline{- \partial_{y_2} \varphi^{l_2} \partial_{\theta} \varphi^{l_3}   (\nu \tau \partial_{y_3}  \underline{H}^{l_1}_2 -\xi_1 \partial_{y_3}  \underline{E}^{l_1}_3) \Big) }\\
        &+\sum\limits_{\substack{l_1+l_2+l_3=m\\l_1\geq 1}}   \Big( \uwave{ \dt \varphi^{l_2} \nu (  \partial_{y_j} \varphi^{l_3} \partial_{y_3}  \underline{H}^{l_1}_j  )   -\partial_{y_1} \varphi^{l_2} ( \dt \varphi^{l_3} \nu \partial_{y_3}  \underline{H}^{l_1}_1 +\partial_{y_2} \varphi^{l_3} \partial_{y_3}  \underline{E}^{l_1}_3)} \\
       &\qquad \uwave{ -\partial_{y_2} \varphi^{l_2}  (\dt \varphi^{l_3} \nu \partial_{y_3}  \underline{H}^{l_1}_2 -\partial_{y_1} \varphi^{l_3} \partial_{y_3}  \underline{E}^{l_1}_3)} \Big)
    \Big\}\\
    =& \mathbf{c_0} \Big\{ 
    \sum\limits_{\substack{l_1+l_2=m+2\\l_1\geq 1}}   \Big(\nu \dt\partial_{\theta} \varphi^{l_2}  (\xi_j   \underline{H}^{l_1}_j  )   -  \partial_{y_1}\partial_{\theta}  \varphi^{l_2}   (\nu \tau  \underline{H}^{l_1}_1 +\xi_2   \underline{E}^{l_1}_3) \\
    &\qquad- \partial_{y_2}\partial_{\theta}  \varphi^{l_2}   (\nu \tau  \underline{H}^{l_1}_2 -\xi_1  \underline{E}^{l_1}_3) \Big)\\
    &+\sum\limits_{\substack{l_1+l_2=m+1\\l_1\geq 1}}   \Big(- \dt \varphi^{l_2}\nu ( \partial_{y_j}  \underline{H}^{l_1}_j  )   +\partial_{y_1} \varphi^{l_2}   (\nu  \dt \underline{H}^{l_1}_1 + \partial_{y_2}  \underline{E}^{l_1}_3) \\
    &\qquad +\partial_{y_2} \varphi^{l_2}   (\nu  \dt \underline{H}^{l_1}_2 - \partial_{y_1}  \underline{E}^{l_1}_3) \Big) \\
    &+\sum\limits_{\substack{l_1+l_2+l_3=m+1\\l_1\geq 1}}   \uuline{  \Big( \dt \varphi^{l_2} \partial_{\theta} \varphi^{l_3}  \nu (\xi_j \partial_{y_3}  \underline{H}^{l_1}_j  )   -\partial_{y_1} \varphi^{l_2} \partial_{\theta} \varphi^{l_3}   (\nu \tau \partial_{y_3}  \underline{H}^{l_1}_1 +\xi_2 \partial_{y_3}  \underline{E}^{l_1}_3)} \\
    &\qquad\uuline{- \partial_{y_2} \varphi^{l_2} \partial_{\theta} \varphi^{l_3}   (\nu \tau \partial_{y_3}  \underline{H}^{l_1}_2 -\xi_1 \partial_{y_3}  \underline{E}^{l_1}_3) \Big) }
    \Big\},
    \end{align*}
    where the underwave terms are cancelled themselves in the summations, while the double underline terms are cancelled after summing up $I_1$ and $I_2$. 
    With the above calculations, it is found that $\bF^{m,-}_3$ is the residue component of the last four lines of \cref{comp-ams-1}.
    
    \underline{\textbf{Step 2: Surface wave components.}}
    Since $\partial_{Y_3} \widehat{H}^{m+1}_{3,\star}(0)=\widehat{F}^{m,-}_{7,\star}(0)$,   $-\partial_{Y_3} \widehat{E}^{m+1}_{2,\star}(0)=\widehat{F}^-_{1,\star}(0)$, $\partial_{Y_3} \widehat{E}^{m+1}_{1,\star}(0)=\widehat{F}^-_{2,\star}(0)$, and also the induction assumption (H(m)-3), we have that for any $y_3, Y_3<0$, 
    \begin{align}
    	&\partial_{Y_3} \Big( \nu \dt \widehat{H}^{m+1}_{3,\star}(0) +\partial_{y_1} \widehat{E}^{m+1}_{2,\star} (0) -\partial_{y_2} \widehat{E}^{m+1}_{1,\star}(0)\Big) \notag\\
       =&\nu \dt \widehat{F}^{m,-}_{7,\star}(0) - \partial_{y_j}\widehat{F}^{m,-}_{7,\star}(0)\notag\\
       =&\mathbf{c_0} \Big\{ 
       \nu \dt F^{m,-}_7-\partial_{y_j} F^{m,-}_j
       \Big\}\notag\\
       =&\mathbf{c_0} \Big\{ 
       -\Big(\nu \dt \partial_{y_3} H^m_3 + \partial_{y_1}\partial_{y_3} E^m_3 -\partial_{y_2}\partial_{y_3} E^m_1 \Big)\notag\\
        % Y_3 part
        &+ \sum\limits_{\substack{l_1+l_2=m+2\\l_1\geq 1}} \partial_{\theta} \varphi^{l_2} \partial_{Y_3} \Big( 
        \nu \dt (\xi_j H^{l_1}_j) -\partial_{y_1} (\nu\tau H^{l_1}_1 +\xi_2 E^{l_1}_3)-\partial_{y_2} (\nu\tau H^{l_1}_2 -\xi_1 E^{l_1}_3)
        \Big) \notag\\
        &+\sum\limits_{\substack{l_1+l_2=m+2\\l_1\geq 1}} \Big(\nu \dt  \partial_{\theta} \varphi^{l_2}  \partial_{Y_3}
        (\xi_j H^{l_1}_j) -\partial_{y_1} \partial_{\theta} \varphi^{l_2}  \partial_{Y_3}(\nu\tau H^{l_1}_1 +\xi_2 E^{l_1}_3)\notag\\
        &\qquad-\partial_{y_2}\partial_{\theta} \varphi^{l_2}  \partial_{Y_3} (\nu\tau H^{l_1}_2 -\xi_1 E^{l_1}_3)
        \Big)\notag\\
        &+\sum\limits_{\substack{l_1+l_2=m+1\\l_1\geq 1}} \Big(  (\partial_{y_j} \varphi^{l_2}  \nu \dt\partial_{Y_3}
         H^{l_1}_j) -(\dt \varphi^{l_2} \nu \partial_{y_1} \partial_{Y_3} H^{l_1}_1 +\partial_{y_2} \varphi^{l_2}  \partial_{y_1} \partial_{Y_3} E^{l_1}_3) \notag\\
        &  \qquad-(\dt \varphi^{l_2} \nu \partial_{y_2} \partial_{Y_3}H^{l_1}_2 -\partial_{y_1} \varphi^{l_2}  \partial_{y_2}\partial_{Y_3} E^{l_1}_3)
        \Big)\notag\\
        % y_3 part
        &+ \sum\limits_{\substack{l_1+l_2=m+1\\l_1\geq 1}} \partial_{\theta} \varphi^{l_2} \partial_{y_3} \Big( 
        \nu \dt (\xi_j H^{l_1}_j) -\partial_{y_1} (\nu\tau H^{l_1}_1 +\xi_2 E^{l_1}_3)-\partial_{y_2} (\nu\tau H^{l_1}_2 -\xi_1 E^{l_1}_3)\Big)\notag \\
        &+\sum\limits_{\substack{l_1+l_2=m+1\\l_1\geq 1}} \Big(\nu \dt  \partial_{\theta} \varphi^{l_2}  \partial_{y_3}
        (\xi_j H^{l_1}_j) -\partial_{y_1} \partial_{\theta} \varphi^{l_2}  \partial_{y_3}(\nu\tau H^{l_1}_1 +\xi_2 E^{l_1}_3)\notag\\
        &\qquad-\partial_{y_2}\partial_{\theta} \varphi^{l_2}  \partial_{y_3} (\nu\tau H^{l_1}_2 -\xi_1 E^{l_1}_3) \Big) \notag\\
        &+\sum\limits_{\substack{l_1+l_2=m\\l_1\geq 1}} \Big(  (\partial_{y_j} \varphi^{l_2}  \nu \dt\partial_{y_3}
        H^{l_1}_j) -(\dt \varphi^{l_2} \nu \partial_{y_1} \partial_{y_3} H^{l_1}_1 +\partial_{y_2} \varphi^{l_2}  \partial_{y_1} \partial_{y_3} E^{l_1}_3) \notag\\
        &  \qquad-(\dt \varphi^{l_2} \nu \partial_{y_2} \partial_{y_3}H^{l_1}_2 -\partial_{y_1} \varphi^{l_2}  \partial_{y_2}\partial_{y_3} E^{l_1}_3) \Big)
       \Big\}\notag\\
       =:& \mathbf{c_0} \Big\{ 
       -\Big(\nu \dt \partial_{y_3} H^m_3 + \partial_{y_1}\partial_{y_3} E^m_3 -\partial_{y_2}\partial_{y_3} E^m_1 \Big) \Big\} + \sum\limits_{n=1}^6 I_n. \label{comp-ams-2}
    \end{align}
     In the end of the above equality, we use $I_n$ ($n=1,2, \cdots, 6$) to denote the six summations, respectively.
     
     With (H(m)-3) and \cref{lem-comp-dm}, it holds that $\widehat{F}^{m,-}_3 (0) =0$,  which gives that 
     \begin{align*}
     	&\mathbf{c_0} \Big\{ 
     	\nu \dt  H^m_3 + \partial_{y_1}E^m_3 -\partial_{y_2} E^m_1 \Big) \\
     	=& \mathbf{c_0} \Big\{ 
     	\sum\limits_{\substack{l_1+l_2=m+2\\l_1\geq 1}} \partial_{\theta}\varphi^{l_2} (\nu\tau \partial_{Y_3} H^{l_1}_3 +\xi_1\partial_{Y_3} E^{l_1}_2-\xi_2\partial_{Y_3} E^{l_1}_1) \notag \\
     	&\quad+ \sum\limits_{\substack{l_1+l_2=m+1\\l_1\geq 1}} ( \dt \varphi^{l_2} \nu \partial_{Y_3} H^{l_1}_3 +\partial_{y_1}\varphi^{l_2}\partial_{Y_3} E^{l_1}_2  -\partial_{y_2}\varphi^{l_2}\partial_{Y_3} E^{l_1}_1 ) \notag \\
     	&\quad+ \sum\limits_{\substack{l_1+l_2=m+1\\l_1\geq 1}} \partial_{\theta}\varphi^{l_2}(\nu\tau \partial_{y_3} H^{l_1}_3 +\xi_1\partial_{y_3} E^{l_1}_2 -\xi_2\partial_{y_3} E^{l_1}_1) \notag \\
     	&\quad+ \sum\limits_{\substack{l_1+l_2=m\\l_1\geq 1}} ( \dt \varphi^{l_2} \nu \partial_{y_3} H^{l_1}_3 +\partial_{y_1}\varphi^{l_2}\partial_{y_3} E^{l_1}_2  -\partial_{y_2}\varphi^{l_2}\partial_{y_3} E^{l_1}_1 ) 
     	\Big\} \\
     	=:&J_1 +J_2+J_3+J_4. 
     \end{align*}
     In the end of the above equality, we use $J_1, J_2, J_3, J_4$ to denote the four summations, respectively. Now we take the parital derivative with repsect to $y_3$ of the above terms, and substitute it back to \cref{comp-ams-2}, subtitute proper terms by  $\underline{F}^{l_1}_7 $, $\underline{F}^{l_1}_1 $, $\underline{F}^{l_1}_2 $ and other terms from H((m+1)-2), and integrate by parts with respect to $\theta$ (for $I_5$). Then it yields that 
     \begin{align}
      &	-\partial_{y_3} J_2 + I_5 + I_6-\partial_{y_3} J_4   \notag\\
      % 1st =
     =&  \mathbf{c_0} \Big\{ 
     \sum\limits_{\substack{l_1+l_2=m+1\\l_1\geq 1}} \Big(  -\nu \dt \varphi^{l_2}  \partial_{y_3} (\partial_{Y_3} H^{l_1}_3  + \xi_j\partial_{\theta} H^{l_1}_j 
     +\partial_{y_{\alpha}} H^{l_1-1}_{\alpha}  ) \notag\\
     &+  \partial_{y_1} \varphi^{l_2}  \partial_{y_3} (-\partial_{Y_3} E^{l_1}_2 +\nu\tau \partial_{\theta} H^{l_1}_1 +\xi_2 \partial_{\theta} E^{l_1}_3
      +\nu \dt H^{l_1-1}_1 +\partial_{y_2} E^{l_1-1}_3 -\partial_{y_3} E^{l_1-1}_2  ) \notag\\
     &+  \partial_{y_2} \varphi^{l_2}  \partial_{y_3} (\partial_{Y_3} E^{l_1}_2 +\nu\tau \partial_{\theta} H^{l_1}_1 -\xi_1 \partial_{\theta} E^{l_1}_3
     +\nu \dt H^{l_1-1}_2 -\partial_{y_1} E^{l_1-1}_3 +\partial_{y_3} E^{l_1-1}_1  ) 
      \Big) \Big\} \notag\\
      % 2st =
      =&  \mathbf{c_0} \Big\{ 
      \sum\limits_{\substack{l_1+l_2=m+1\\l_1\geq 1}} \Big(  -\nu \dt \varphi^{l_2}  \partial_{y_3} (F^{l_1-1}_7 
      +\partial_{y_{\alpha}} H^{l_1-1}_{\alpha}  ) \notag\\
      &+  \partial_{y_1} \varphi^{l_2}  \partial_{y_3} ( F^{l_1-1}_1
      +\nu \dt H^{l_1-1}_1 +\partial_{y_2} E^{l_1-1}_3 -\partial_{y_3} E^{l_1-1}_2  ) \notag\\
      &+  \partial_{y_2} \varphi^{l_2}  \partial_{y_3} (F^{l_1-1}_2
      +\nu \dt H^{l_1-1}_2 -\partial_{y_1} E^{l_1-1}_3 +\partial_{y_3} E^{l_1-1}_1  ) 
      \Big) \Big\} \notag\\
      %3rd =
    =& \mathbf{c_0} \Big\{ 
    \sum\limits_{\substack{l_1+l_2+l_3=m+2\\l_1\geq 1}}     \partial_{\theta} \varphi^{l_3}\Big( -\nu\dt \varphi^{l_2}  \partial_{y_3}   (\xi_j \partial_{Y_3}  H^{l_1}_j  )  +\partial_{y_1} \varphi^{l_2}  \partial_{y_3} (\nu \tau \partial_{Y_3}  H^{l_1}_1 +\xi_2 \partial_{Y_3} E^{l_1}_3) \notag\\
    &\qquad+ \partial_{y_2} \varphi^{l_2}   \partial_{y_3} (\nu \tau \partial_{Y_3}  H^{l_1}_2 -\xi_1 \partial_{Y_3}  E^{l_1}_3) \Big)  \notag\\
    &+\sum\limits_{\substack{l_1+l_2+l_3=m+1\\l_1\geq 1}}  \partial_{\theta} \varphi^{l_3}  \Big( -\nu\dt \varphi^{l_2}   \partial_{y_3}   (\xi_j \partial_{y_3}  H^{l_1}_j  )  +\partial_{y_1} \varphi^{l_2} \partial_{y_3} (\nu \tau \partial_{y_3}  H^{l_1}_1 +\xi_2 \partial_{y_3} E^{l_1}_3) \notag \\
    &\qquad+ \partial_{y_2} \varphi^{l_2} \partial_{y_3} (\nu \tau \partial_{y_3}  H^{l_1}_2 -\xi_1 \partial_{y_3}  E^{l_1}_3) \Big)
    \Big\} \notag,
     \end{align}
  \begin{align}
 	&	-\partial_{y_3} J_1 + I_4 -\partial_{y_3} J_3   \notag\\
 	%1st =
   =&  \mathbf{c_0} \Big\{ 
   \sum\limits_{\substack{l_1+l_2=m+2\\l_1\geq 1}}  \partial_{\theta} \varphi^{l_2}  \partial_{y_3}\Big(  -\nu \tau (\partial_{Y_3} H^{l_1}_3  
   +\partial_{y_{\alpha}} H^{l_1-1}_{\alpha}  ) \notag\\
   &+  \xi_1  (-\partial_{Y_3} E^{l_1}_2 
   +\nu \dt H^{l_1-1}_1 +\partial_{y_2} E^{l_1-1}_3 -\partial_{y_3} E^{l_1-1}_2  ) \notag\\
   &+  \xi_2 (\partial_{Y_3} E^{l_1}_2 
   +\nu \dt H^{l_1-1}_2 -\partial_{y_1} E^{l_1-1}_3 +\partial_{y_3} E^{l_1-1}_1  ) 
   \Big) \Big\} \notag\\
   % 2st =
   =&  \mathbf{c_0} \Big\{ 
   \sum\limits_{\substack{l_1+l_2=m+2\\l_1\geq 1}}  \partial_{\theta} \varphi^{l_2}  \partial_{y_3} \Big(  -\nu \tau (F^{l_1-1}_7  + \uwave{\xi_j\partial_{\theta} H^{l_1}_j }
   +\partial_{y_{\alpha}} H^{l_1-1}_{\alpha}  ) \notag\\
   &+  \xi_1 ( F^{l_1-1}_1
   +\uwave{\nu\tau \partial_{\theta} H^{l_1}_1 +\xi_2 \partial_{\theta} E^{l_1}_3}
   +\nu \dt H^{l_1-1}_1 +\partial_{y_2} E^{l_1-1}_3 -\partial_{y_3} E^{l_1-1}_2  ) \notag\\
   &+  \xi_2 (F^{l_1-1}_2
   +\uwave{\nu\tau \partial_{\theta} H^{l_1}_1 -\xi_1 \partial_{\theta} E^{l_1}_3}
   +\nu \dt H^{l_1-1}_2 -\partial_{y_1} E^{l_1-1}_3 +\partial_{y_3} E^{l_1-1}_1  ) 
   \Big) \Big\} \notag\\
   %3rd =
   =& \mathbf{c_0} \Big\{ 
   \sum\limits_{\substack{l_1+l_2+l_3=m+2\\l_1\geq 1}}     \partial_{\theta} \varphi^{l_2}\Big( -\nu\tau (\partial_{y_j} \varphi^{l_3}  \partial_{y_3}   \partial_{Y_3}  H^{l_1}_j  )  +\xi_1  ( \dt \varphi^{l_3}\nu \partial_{y_3}\partial_{Y_3}  H^{l_1}_1 +\partial_{y_2} \varphi^{l_3} \partial_{y_3} \partial_{Y_3} E^{l_1}_3) \notag\\
   &\qquad+ \xi_2 ( \dt \varphi^{l_3}\nu \partial_{y_3} \partial_{Y_3}  H^{l_1}_2 -\partial_{y_1} \varphi^{l_3} \partial_{y_3} \partial_{Y_3}  E^{l_1}_3) \Big)  \notag\\
   &+\sum\limits_{\substack{l_1+l_2+l_3=m+1\\l_1\geq 1}}     \partial_{\theta} \varphi^{l_2}\Big( -\nu\tau (\partial_{y_j} \varphi^{l_3}  \partial_{y_3}^2  H^{l_1}_j  )  +\xi_1  (\nu \dt \varphi^{l_3} \partial_{y_3}^2 H^{l_1}_1 +\partial_{y_2} \varphi^{l_3} \partial_{y_3}^2 E^{l_1}_3) \notag\\
   &\qquad+ \xi_2 (\nu \dt \varphi^{l_3} \partial_{y_3}^2  H^{l_1}_2 -\partial_{y_1} \varphi^{l_3} \partial_{y_3}^2 E^{l_1}_3) \Big) 
   \Big\} \notag.
   \end{align}
   By exchanging the summation orders of $l_2$ and $l_3$, we have
   \begin{align*}
   	&\mathbf{c_0} \Big\{ 
   	-\Big(\nu \dt \partial_{y_3} H^m_3 + \partial_{y_1}\partial_{y_3} E^m_3 -\partial_{y_2}\partial_{y_3} E^m_1 \Big) \Big\} +I_4 + I_5 + I_6 \\
   	=& (-\partial_{y_3} J_2 + I_5 + I_6-\partial_{y_3} J_4) +(-\partial_{y_3} J_1 + I_4 -\partial_{y_3} J_3)=0.
   \end{align*}
   Therefore, rearranging the terms in $I_3$ gives that
   \begin{align}
   	&\partial_{Y_3} \Big( \nu \dt \widehat{H}^{m+1}_{3,\star}(0) +\partial_{y_1} \widehat{E}^{m+1}_{2,\star} (0) -\partial_{y_2} \widehat{E}^{m+1}_{1,\star}(0)\Big) = I_1+I_2+I_3 \notag \\
   	=&\partial_{Y_3} \mathbf{c_0} \Big\{  \sum\limits_{\substack{l_1+l_2=m+2\\l_1\geq 1}} \partial_{\theta} \varphi^{l_2}  \Big( 
   	\nu \dt (\xi_j H^{l_1}_j) -\partial_{y_1} (\nu\tau H^{l_1}_1 +\xi_2 E^{l_1}_3)\notag\\
   	 &\qquad-\partial_{y_2} (\nu\tau H^{l_1}_2 -\xi_1 E^{l_1}_3)
   	\Big) \notag\\
   	&+\sum\limits_{\substack{l_1+l_2=m+2\\l_1\geq 1}} \Big(\nu \dt  \partial_{\theta} \varphi^{l_2}  
   	(\xi_j H^{l_1}_j) -\partial_{y_1} \partial_{\theta} \varphi^{l_2} (\nu\tau H^{l_1}_1 +\xi_2 E^{l_1}_3)\notag\\
   	&\qquad-\partial_{y_2}\partial_{\theta} \varphi^{l_2}   (\nu\tau H^{l_1}_2 -\xi_1 E^{l_1}_3)
   	\Big)\notag\\
   	&+\sum\limits_{\substack{l_1+l_2=m+1\\l_1\geq 1}} \Big(  (\partial_{y_j} \varphi^{l_2}  \nu \dt
   	H^{l_1}_j) -(\dt \varphi^{l_2} \nu \partial_{y_1}  H^{l_1}_1 +\partial_{y_2} \varphi^{l_2}  \partial_{y_1}  E^{l_1}_3) \notag\\
   	&  \qquad-(\dt \varphi^{l_2} \nu \partial_{y_2}H^{l_1}_2 -\partial_{y_1} \varphi^{l_2}  \partial_{y_2} E^{l_1}_3)
   	\Big) \Big\}\notag
   \end{align}
   It  follows directly from integrating the above equality with respect to $Y_3$ from $-\infty$ to $0$ that 
    $\Big( \nu \dt \widehat{H}^{m+1}_{3,\star}(0) +\partial_{y_1} \widehat{E}^{m+1}_{2,\star} (0) -\partial_{y_2} \widehat{E}^{m+1}_{1,\star}(0)\Big) |_{y_3=Y_3=0}$ is the surface wave component of the last four lines of \cref{comp-ams-1}.
    
    With the two steps above, the proof is completed.
\end{proof}

\begin{Lem}\label{lem-comp-dms}
	The source terms in \cref{s-source} satisfy:
	\begin{equation}
		\label{comp-dms}
		\partial_t \bF_8^{m,+}=\partial_{y_\alpha} \bF_{3+\alpha}^{m,+},\quad 
		\nu\partial_t \bF_7^{m,-}=\partial_{y_\alpha} \bF_{\alpha}^{m,-},\quad 
		\nu\partial_t \bF_8^{m,-}=\partial_{y_\alpha} \bF_{3+\alpha}^{m,-}.
	\end{equation}
\end{Lem}
\begin{proof}[Proof of \cref{lem-comp-dms}]
The proof of the equation for the plasma is the same as in \cite{Pierre2021}. We just consider one of the equality for the vacuum, $\nu\partial_t \bF_7^{m,-}=\partial_{y_\alpha} \bF_{\alpha}^{m,-}$, the proof of the other one is similar.
Recalling \cref{s-source} and using integration by parts with respect to $\theta$, it yields that
\begin{align*}
	&\partial_{y_\alpha} \bF_{\alpha}^{m,-}-\nu\partial_t \bF_7^{m,-} \notag\\
	%1st =
   =& \mathbf{c_0} \Big\{  \sum\limits_{\substack{l_1+l_2=m+2\\l_1\geq 1}} \Big( 
        \partial_{y_1} \partial_{\theta} \varphi^{l_2} \partial_{y_3} (\nu\tau \underline{H}^{l_1}_1 +\xi_2 \underline{E}^{l_1}_3 )
       \notag \\
    &\qquad +\partial_{y_2} \partial_{\theta} \varphi^{l_2} \partial_{y_3} (\nu\tau \underline{H}^{l_1}_2 -\xi_1 \underline{E}^{l_1}_3 )
            -\nu \dt \partial_{\theta} \varphi^{l_2} \partial_{y_3} (\xi_j \underline{H}^{l_1}_j )
   \Big) \notag\\
    &+ \sum\limits_{\substack{l_1+l_2=m+2\\l_1\geq 1}}\partial_{\theta} \varphi^{l_2} \partial_{y_3}  \Big( \nu\tau (\partial_{y_{\alpha}} \underline{H}^{l_1}_{\alpha}  ) -\xi_1 (\nu\dt \underline{H}^{l_1}_1 + \partial_{y_2} \underline{E}^{l_1}_3 - \partial_{y_3} \underline{E}^{l_1}_2 )\\
    &\qquad -\xi_2 (\nu\dt \underline{H}^{l_1}_2 - \partial_{y_1} \underline{E}^{l_1}_3 + \partial_{y_3} \underline{E}^{l_1}_1 )
   \Big) \notag\\
    &+ \sum\limits_{\substack{l_1+l_2=m+1\\l_1\geq 1}}  \Big( \nu \dt \varphi^{l_2} \partial_{y_3} (\partial_{y_{\alpha}} \underline{H}^{l_1}_{\alpha}  ) -\partial_{y_1} \varphi^{l_2} \partial_{y_3} (\nu\dt \underline{H}^{l_1}_1 + \partial_{y_2} \underline{E}^{l_1}_3 - \partial_{y_3} \underline{E}^{l_1}_2 )\\
    &\qquad -\partial_{y_2} \varphi^{l_2} \partial_{y_3} (\nu\dt \underline{H}^{l_1}_2 - \partial_{y_1} \underline{E}^{l_1}_3 + \partial_{y_3} \underline{E}^{l_1}_1 )
    \Big) 
   \Big\}\notag\\
   =:&I_1+I_2+I_3,\notag
\end{align*}
\begin{align*}
   I_1+I_3
   %1st =
   =& \mathbf{c_0} \Big\{  \sum\limits_{\substack{l_1+l_2=m+1\\l_1\geq 1}}  \Big( \nu \dt \varphi^{l_2} \partial_{y_3} ( \underline{F}^{l_1,-}_7+ \partial_{y_{\alpha}} \underline{H}^{l_1}_{\alpha}  ) \notag\\
   &\qquad -\partial_{y_1} \varphi^{l_2} \partial_{y_3} (\underline{F}^{l_1,-}_1+\nu\dt \underline{H}^{l_1}_1 + \partial_{y_2} \underline{E}^{l_1}_3 - \partial_{y_3} \underline{E}^{l_1}_2 )\\
   &\qquad -\partial_{y_2} \varphi^{l_2} \partial_{y_3} (\underline{F}^{l_1,-}_2+\nu\dt \underline{H}^{l_1}_2 - \partial_{y_1} \underline{E}^{l_1}_3 + \partial_{y_3} \underline{E}^{l_1}_1 ) \Big)  
   \Big\}\notag\\
   %2nd =
   =& \mathbf{c_0} \Big\{  \sum\limits_{\substack{l_1+l_2+l_3=m+2\\l_1\geq 1}}  \Big( \nu \dt \varphi^{l_2} \partial_{\theta} \varphi^{l_3} \partial_{y_3} (  \xi_j\underline{H}^{l_1}_j  )  -\partial_{y_1} \varphi^{l_2} \partial_{\theta} \varphi^{l_3} \partial_{y_3} (\nu\tau \underline{H}^{l_1}_1 + \xi_2 \underline{E}^{l_1}_3 )\\
   &\qquad -\partial_{y_2} \varphi^{l_2} \partial_{\theta} \varphi^{l_3} \partial_{y_3} (\nu\tau \underline{H}^{l_1}_1 + \xi_2 \underline{E}^{l_1}_3 ) \Big)    
   \Big\},\notag
\end{align*}
\begin{align*}
   I_2=& \mathbf{c_0} \Big\{  \sum\limits_{\substack{l_1+l_2=m+2\\l_1\geq 1}} \partial_{\theta} \varphi^{l_2} \partial_{y_3}  \Big( \uwave{\nu\tau ( - \xi_j \partial_{\theta} \underline{H}^{l_1+1}_{j} } + \underline{F}^{l_1,-}_{7}   ) \notag\\
   &\qquad \uwave{-\xi_1 (- \nu\tau  \partial_{\theta} \underline{H}^{l_1+1}_1 - \xi_2 \partial_{\theta}\underline{E}^{l_1+1}_3 }+ \underline{F}^{l_1,-}_{1}   )\\
   &\qquad \uwave{-\xi_2 (- \nu\tau  \partial_{\theta} \underline{H}^{l_1+1}_2 + \xi_1 \partial_{\theta}\underline{E}^{l_1+1}_3 +} \underline{F}^{l_1,-}_{2}   ) \Big)
   \Big\}\notag\\
%   =& \mathbf{c_0} \Big\{  \sum\limits_{\substack{l_1+l_2+l_3=m+2\\l_1\geq 1}} \partial_{\theta} \varphi^{l_2}  \partial_{\theta} \varphi^{l_3}\partial_{y_3}  \Big( \nu\tau (  \xi_j  \underline{H}^{l_1}_{j} ) -\xi_1 ( \nu\tau   \underline{H}^{l_1}_1 + \xi_2 \underline{E}^{l_1}_3  )\\
%   &\qquad -\xi_2 ( \nu\tau   \underline{H}^{l_1}_2 -\xi_1\underline{E}^{l_1}_3  ) \Big)
%   \Big\}\notag
     =& \mathbf{c_0} \Big\{  \sum\limits_{\substack{l_1+l_2+l_3=m+2\\l_1\geq 1}} \partial_{\theta} \varphi^{l_2}  \Big( \nu\tau (   \partial_{y_j} \varphi^{l_3}  \partial_{y_3} \underline{H}^{l_1}_{j} ) -\xi_1 ( \dt \varphi^{l_3} \nu  \partial_{y_3}    \underline{H}^{l_1}_1 + \partial_{y_2} \varphi^{l_3}  \partial_{y_3}  \underline{E}^{l_1}_3  )\\
   &\qquad -\xi_2 ( \dt \varphi^{l_3}\nu  \partial_{y_3}    \underline{H}^{l_1}_2 -\partial_{y_1} \varphi^{l_3}  \partial_{y_3} \underline{E}^{l_1}_3  ) \Big)
   \Big\} = -(I_1+I_3),\notag
\end{align*}
where the underwave terms are cancelled themselves during the summations.
Therefore, the proof of $\nu\partial_t \bF_7^{m,-}=\partial_{y_\alpha} \bF_{\alpha}^{m,-}$ is completed.
\end{proof}

\vline

\section*{Acknowledgements}
Yuan Yuan thanks the Department of Mathematics, University of Brescia for its hospitality during his postdoctoral stay, in which the research of this work was initiated.

The research of Paolo Secchi was supported in part by the Italian MIUR Project PRIN 2015YCJY3A-004.
The research of Yuan Yuan is supported by the National Natural Science Foundation of China (Grants No. 11901208 and 11971009), the Natural Science Foundation of Guangdong Province, China (Grant No. 2021A1515010247) and by the Italian MIUR Project PRIN 2015YCJY3A-004.

%%%%%%%%%%%%%%%%%%%%%%%%%%%%%%%%%%%%%%%%%%%%%%%%%%%%%%%%%%%%%%%%
%%%%% Bib
%%%%%%%%%%%%%%%%%%%%%%%%%%%%%%%%%%%%%%%%%%%%%%%%%%%%%%%%%%%%%%%%
%\bibliographystyle{amsplain}
%\bibliography{ref}

\end{document}